\let\chooseClass1   %tac journal
\let\chooseClass3   %article, 12pt, for arXiv

\ifx\chooseClass1
\documentclass{tac}
\newtheoremrm{conjecture}{Conjecture}
\newtheoremrm{exercise}{Exercise}
\fi

\ifx\chooseClass2
\IfFileExists{extarticle.cls}{
\documentclass[14pt]{extarticle}
\emergencystretch 7 pt
\setlength{\oddsidemargin}{-17 mm} %-17     %-19    %-19     %-19
\setlength{\textwidth}{192 mm}     %192     %196    %196     %196
\setlength{\textheight}{246 mm}    %246     %246    %270     %266
\setlength{\topmargin}{-31 mm}     %-31     %-31    %-31     %-31
}{
\documentclass[12pt]{article}
\textheight = 8.5in
\textwidth 6.3in
\setlength{\oddsidemargin}{0mm}
\setlength{\topmargin}{-15 mm}
}{}
\fi

\ifx\chooseClass3
\documentclass[12pt]{article}
\fi

\ifx\chooseClass1
	\else
\makeatletter
\def\@seccntformat#1{\csname the#1\endcsname.\quad}
\renewcommand\section{\@startsection {section}{1}{\z@}%
                                   {-3.5ex \@plus -1ex \@minus -.2ex}%
                                   {2.3ex \@plus.2ex}%
                                   {\normalfont\large\bfseries}}
\renewcommand\subsection{\@startsection{subsection}{2}{\z@}%
                        {3.25ex plus 1ex minus .2ex}{-.5em}%
                        {\normalfont\normalsize\bfseries}}
\renewcommand\subsubsection{\@startsection{subsubsection}{3}{\z@}%
                        {3.25ex plus 1ex minus .2ex}{-.5em}%
                        {\normalfont\normalsize\it}}
\@addtoreset{equation}{section}
\makeatother

\usepackage{amsthm}
\newtheoremstyle{boldhead}%     name
{\topsep}%                      abovespace
{\topsep}%                      belowspace
{\slshape}%                     bodyfont
{}%                             indentation=noindent
{\bfseries}%                    headfont
{.}%                            headpunctuation
{ }%                            headspace=interword space
{\thmname{#1}\thmnumber{ #2}\thmnote{ (#3)}}%   custom head specification

\newtheoremstyle{boldremark}%   name
{\topsep}%                      abovespace
{\topsep}%                      belowspace
{\upshape}%                     bodyfont
{}%                             indentation=noindent
{\bfseries}%                    headfont
{.}%                            headpunctuation
{ }%                            headspace=interword space
{\thmname{#1}\thmnumber{ #2}\thmnote{ (#3)}}%   custom head specification

\swapnumbers

\theoremstyle{boldhead}
\newtheorem{theorem}[subsection]{Theorem}
\newtheorem{corollary}[subsection]{Corollary}
\newtheorem{lemma}[subsection]{Lemma}
\newtheorem{proposition}[subsection]{Proposition}

\theoremstyle{boldremark}
\newtheorem{conjecture}[subsection]{Conjecture}
\newtheorem{definition}[subsection]{Definition}
\newtheorem{example}[subsection]{Example}
\newtheorem{examples}[subsection]{Examples}
\newtheorem{exercise}[subsection]{Exercise}
\newtheorem{remark}[subsection]{Remark}

\fi%\chooseClass2,3

\usepackage{    %abbrevs,
				amsmath,
                amsfonts,
                amssymb,
                amsxtra,
                fancybox,
                stmaryrd,
                varioref,
}

\numberwithin{equation}{section}

\IfFileExists{url.sty}{\usepackage{url}}{}
\providecommand{\url}[1]{{\tt #1}}

\usepackage{ifpdf}
\ifpdf
    \IfFileExists{hyperref.sty}{\usepackage[pdftex]{hyperref}}{}
\else
    \IfFileExists{hyperref.sty}{\usepackage[hypertex]{hyperref}}{}
\fi

    \ifpdf
        \usepackage[pdftex]{graphicx}
    \else
        \usepackage[dvips]{graphicx}
    \fi

\IfFileExists{euscript.sty}{\usepackage[mathcal]{euscript}}{}
\IfFileExists{mathrsfs.sty}{\usepackage{mathrsfs}}{\let\mathscr\mathfrak}

\usepackage{diagrams}
\diagramstyle[height=2em,balance,righteqno,PostScript=dvips,nohug]

\IfFileExists{dsfont.sty}{\usepackage[sans]{dsfont}%
\newcommand\1{{\mathds 1}}}{\newcommand\1{{1\mkern-5mu {\mathrm I}}}}

\newlength{\mylabelwidths}
\newenvironment{myitemize}{\begin{list}{}{%
\setlength{\labelwidth}{\mylabelwidths}%
\setlength{\leftmargin}{\mylabelwidths}\addtolength{\leftmargin}{0.5em}%
\setlength{\itemsep}{-0.2\baselineskip}}}%
{\end{list}}

\newarrowhead{pto}{\mkern-6mu\raise.162em\hbox{$\sss\succ$}\mkern6mu}%
{\mkern6mu\raise.162em\hbox{$\sss\prec$}\mkern-6mu}%
{\sss\curlyvee}{\sss\curlywedge}%vertical heads are unusable yet

\newarrowtail{mfrom}{\mkern6mu\lower.018em\hbox{$\sss\prec$}\mkern-6mu}%
{\mkern-6mu\lower.018em\hbox{$\sss\succ$}\mkern6mu}%
{\sss\curlywedge}{\sss\curlyvee}%vertical tails are unusable yet

\newarrowtail{sscurlyvee}{\raise.18ex\hbox{$\sss\succ$}}%
{\mkern-.5mu\raise.18ex\hbox{$\sss\prec$}\mkern.4mu}%
{\sss\curlyvee}{\sss\curlywedge}

\newarrow{Cof}{sscurlyvee}---{->}
\newarrow{Epi}----{triangle}
\newarrow{Fib}----{->>}
\newarrow{FromTo}{mfrom}==={pto}
\newarrow{Gyrlyga}{boldhook}----
\newarrow{MapsTo}{mapsto}---{->}
\newarrow{Mono}{boldhook}---{->}
\newarrow{TTo}----{->}
\newarrow{Twoar}===={=>}

\newcommand\NN{{\mathbb N}}
\newcommand\ZZ{{\mathbb Z}}

\newcommand{\ca}{{\mathcal A}}
\newcommand{\cb}{{\mathcal B}}
\newcommand{\cc}{{\mathcal C}}
\newcommand{\cd}{{\mathcal D}}
\newcommand{\ce}{{\mathcal E}}
\newcommand{\cg}{{\mathcal G}}
\newcommand{\CG}{{\mathscr G}}
\newcommand{\ci}{{\mathcal I}}
\newcommand{\ck}{{\mathcal K}}
\newcommand{\cl}{{\mathcal L}}
\newcommand{\cm}{{\mathcal M}}
\newcommand{\cn}{{\mathcal N}}
\newcommand{\co}{{\mathcal O}}
\newcommand{\cp}{{\mathcal P}}
\newcommand{\cq}{{\mathcal Q}}
\newcommand{\cs}{{\mathcal S}}
\newcommand{\CS}{{\mathscr S}}
\newcommand{\cu}{{\mathcal U}}
\newcommand{\cv}{{\mathcal V}}
\newcommand{\cw}{{\mathcal W}}
\newcommand{\cx}{{\mathcal X}}
\newcommand{\Rho}{{\mathrm P}}
\newcommand{\rmF}{{\mathrm F}}

\newcommand{\yi}{\ddot{\imath}}
\newcommand{\one}{{\mathsf1}}
\newcommand{\sfi}{{\mathsf i}}
\newcommand{\sfj}{{\mathsf j}}
\newcommand{\sfv}{{\mathsf v}}
\newcommand{\fu}{{\mathscr U}}

\newcommand{\0}{\phantom0}
\newcommand{\bi}{{\mathbf i}}
\newcommand{\bj}{{\mathbf j}}
\newcommand{\bn}{{\mathbf n}}
\newcommand{\bone}{{\mathbf1}}
\newcommand{\bv}{{\mathbf v}}
\newcommand{\bull}{{\scriptscriptstyle\bullet}}
\newcommand{\colim}{\qopname\relax m{colim}}
\newcommand{\Com}{{{\mathsf C}_\kk}}
\newcommand{\cOm}{{\mathsf C}_\kk}
\newcommand{\uCom}{{\underline{\mathsf C}_\kk}}
\newcommand{\pmQuiver}{{\mathcal{P\!M\!Q}}}
\newcommand{\Sim}{{\mkern-4.5mu\sim\mkern1.5mu}}
\newcommand{\sk}{{\mathsf{sk}}}
\newcommand{\SMCcat}{\mathcal{S\!M\!C}_{\Cat}}
\newcommand{\smQuiver}{\mathcal{S\!M\!Q}_{\Cat}}
\newcommand{\su}{{\mathsf{su}}}
\newcommand{\tdt}{\otimes\dots\otimes}
\newcommand{\tR}[4]{\sS{_{#1,#2}^{\phantom{#1,}#3}}#4}

\newcommand{\sou}{\sS{_{\mathsf s}}}
\newcommand{\soup}{\sS{^{\mathsf s}}}
\newcommand{\tar}{\sS{_{\mathsf t}}}
\newcommand{\sS}[2]{\vphantom{#2}#1 #2}
\newcommand{\n}[1]{\nobreakdash-\hspace{0pt}}
\newcommand{\ainf}[1]{$A_\infty$\nobreakdash-\hspace{0pt}}
\newcommand{\ainfm}[1]{$\mathrm{A}_\infty$\nobreakdash-\hspace{0pt}}
\newcommand{\alinf}{\mathsf a_\infty}
\newcommand{\mainf}{\mathrm{A}_\infty}
\newcommand{\Cat}{{\mathcal C}at}
\newcommand{\ucatspan}{\underline{\textit{cat}\textup{-span}}}
\newcommand{\Catspan}{{\mathcal C}at\textup{-span}}
\DeclareMathOperator\END{{\mathcal E}{\it nd}}
\newcommand{\nOp}[1]{{}_{#1}\mkern-4.5mu\Op}
\newcommand{\StrictMonCat}{\mathrm{StrictMon}{\mathcal C}at}

\let\boxt\boxtimes
\let\con\triangleright
\let\emptyset\varnothing
\let\eps\varepsilon
\let\ge\geqslant
\let\kk\Bbbk
\let\le\leqslant
\let\lto\xleftarrow
\let\mb\mathbf
\let\rto\xrightarrow
\let\sss\scriptstyle
\let\tens\otimes
\let\ttt\textstyle
\let\und\underline
\let\wh\widehat
\let\wt\widetilde

\newcommand{\aA}{\mathbb{A}}
\newcommand\BB{\mathbb{B}}
\newcommand\KK{\mathbb{K}}

\newcommand\mcDG{{\mathsf{DG}}}

\newcommand\cA{{\mathcal A}}
\newcommand\cB{{\mathcal B}}
\newcommand\cC{{\mathcal C}}
\newcommand\cD{{\mathcal D}}
\newcommand\cF{{\mathcal F}}
\newcommand\cH{{\mathcal H}}
\newcommand\cI{{\mathcal I}}
\newcommand\cM{{\mathcal M}}
\newcommand\cO{{\mathcal O}}
\newcommand\cP{{\mathcal P}}
\newcommand\cR{{\mathcal R}}
\newcommand\cS{{\mathcal S}}

\newcommand\mcA{{\mathsf A}}
\newcommand\mcB{{\mathsf B}}
\newcommand\mcC{{\mathsf C}}
\newcommand\mcD{{\mathsf D}}
\newcommand\mcF{{\mathsf F}}
\newcommand\mcG{{\mathsf G}}
\newcommand\mcH{{\mathsf H}}
\newcommand\mcL{{\mathsf L}}
\newcommand\mcM{{\mathsf M}}

\newcommand\sff{{\mathsf f}}
\newcommand\sfg{{\mathsf g}}

\newcommand\dd[1]{{}^<#1}

\newcommand{\imodi}{\textup{-mod-}}
\newcommand{\modul}{\textup{-mod}}

\DeclareMathOperator\gr{\mathbf{gr}}
\DeclareMathOperator\dg{\mathbf{dg}}
\DeclareMathOperator\udg{\underline{\dg}}

\DeclareMathOperator\Set{\mathcal Set}

\newcommand{\dgac}{\mathbf{dgac}}

\DeclareMathOperator\Mat{Mat}
\DeclareMathOperator\Span{\mathcal Span}

\DeclareMathOperator\Ab{Ab}
\DeclareMathOperator\as{{\it as}}
\DeclareMathOperator\As{{\it As}}
\DeclareMathOperator\AS{\mathsf{As}}
\DeclareMathOperator\ass{{\it as1}}
\DeclareMathOperator\Ass{{\it As1}}
\DeclareMathOperator{\can}{can}
\DeclareMathOperator\COM{{\it Com}}
\DeclareMathOperator\COMM{{\it Com1}}
\DeclareMathOperator\comp{comp}
\DeclareMathOperator\Cone{Cone}
\DeclareMathOperator\edges{e}

\DeclareMathOperator\ev{ev}
\DeclareMathOperator\FAs{{\it FAs}}
\DeclareMathOperator\FAss{{\it FAs1}}

\DeclareMathOperator\hoM{{\it hom}}
\DeclareMathOperator\HOM{{\mathcal H}{\it om}}
\DeclareMathOperator\id{id}
\DeclareMathOperator\Id{Id}
\DeclareMathOperator\im{Im}
\DeclareMathOperator{\In}{In}
\DeclareMathOperator\inj{in}
\DeclareMathOperator\Inp{Inp}
\DeclareMathOperator\Iso{Iso}
\DeclareMathOperator\iV{\overline{v}}
\DeclareMathOperator\IV{v}
\DeclareMathOperator\Ker{Ker}
\DeclareMathOperator\Mor{Mor}
\DeclareMathOperator\Ob{Ob}
\newcommand{\op}{{\operatorname{op}}}
\DeclareMathOperator\Op{Op}
\DeclareMathOperator\pr{pr}
\DeclareMathOperator\troot{root}
\DeclareMathOperator\sign{sign}
\DeclareMathOperator\src{src}
\DeclareMathOperator\supp{supp}
\DeclareMathOperator\tgt{tgt}
\DeclareMathOperator\Tor{Tor}

\newcommand{\corref}[1]{Corollary~\ref{#1}}
\newcommand{\defref}[1]{Definition~\ref{#1}}
\newcommand{\exaref}[1]{Example~\ref{#1}}
\newcommand{\figref}[1]{Fig.~\ref{#1}}
\newcommand{\lemref}[1]{Lemma~\ref{#1}}
\newcommand{\propref}[1]{Proposition~\ref{#1}}
\newcommand{\remref}[1]{Remark~\ref{#1}}
\newcommand{\secref}[1]{Section~\ref{#1}}
\newcommand{\thmref}[1]{Theorem~\ref{#1}}

\newlength{\texthigh}
\begin{document}
\bibliographystyle{amsalpha}
\title{$A_\infty$-morphisms with several entries}
\author{Volodymyr Lyubashenko}
\ifx\chooseClass1
\address{Institute of Mathematics,
National Academy of Sciences of Ukraine, \\
3 Tereshchenkivska st.,
Kyiv-4, 01601 MSP, Ukraine
}
\copyrightyear{2011}
\eaddress{lub@imath.kiev.ua}
\fi
\maketitle

\allowdisplaybreaks[1]

\begin{abstract}
We show that morphisms from $n$ \ainf-algebras to a single one are maps
over an operad module with $n+1$ commuting actions of the operad
$A_\infty$, whose algebras are conventional \ainf-algebras.
Similar statement holds for homotopy unital \ainf-algebras.
The operad $A_\infty$ and modules over it have two useful gradings
related by isomorphisms which change the degree.
The composition of \ainf-morphisms with several entries is presented as
a convolution of a coalgebra-like and an algebra-like structures.
For this sake notions of lax $\Cat$-span multicategories and
multifunctors are introduced.
They are lax versions of strict multicategories and multifunctors
associated with the monad of free strict monoidal category.
\end{abstract}

It is well-known that operads play a prominent part in the study of
\ainfm-algebras.
In particular, \ainfm-algebras in the conventional sense 
\cite{Stasheff:HomAssoc} are algebras over the $\dg$\n-operad $\mainf$,
a resolution (a cofibrant replacement) of the $\dg$\n-operad $\As$ of 
associative non-unital $\dg$\n-algebras.
What about morphisms $A\to B$ of \ainfm-algebras?
It is shown in \cite{Lyu-Ainf-Operad} that they are maps over certain
bimodule over the $\dg$\n-operad $\mainf$.
This bimodule is a resolution (a cofibrant replacement) of the
corresponding $\As$\n-bimodule.

The current article addresses morphisms with several arguments
\(f:A_1,\dots,A_n\to B\) of \ainfm-algebras.
We explain that they are maps over certain $n\wedge1$-operad
$\mainf$-module $\rmF_n$.
The latter means an $\NN^n$\n-graded complex with $n$ left and one
right pairwise commuting actions of $\mainf$.
Furthermore, it is a resolution (a cofibrant replacement) of the
corresponding notion for associative $\dg$\n-algebras without unit.

The unital case is quite similar to the non-unital one.
There is an operad \(\mainf^{hu}\) governing homotopy unital
\ainfm-algebras.
Homotopy unital morphisms $A\to B$ are controlled by an operad 
\(\mainf^{hu}\)\n-bimodule \cite{Lyu-Ainf-Operad}.
In the current article we describe the $n\wedge1$-operad
\(\mainf^{hu}\)\n-module $\rmF_n^{hu}$ responsible for homotopy unital
\ainfm-morphisms \(f:A_1,\dots,A_n\to B\).
We see that it is a resolution (a cofibrant replacement) of the
corresponding $n\wedge1$-operad module over the operad of associative
unital $\dg$\n-algebras.

The $\dg$-operad of \ainf-algebras has two useful forms.
The first, $\mainf$, is already presented as a resolution of the
operad $\As$.
The second, $A_\infty$, is easy to remember, because all generators
have the same degree 1 and the expression for the differential
contains no oscillating signs.
These two are related by an isomorphism of operads that changes the
degrees in a prescribed way.
Structure equations for this isomorphism use certain signs.
These signs reappear in the formula for the differential in the first
operad, $\mainf$.
There are similar duplicates of other considered operads and modules
over them: $F_n$, \(A_\infty^{hu}\), \(F_n^{hu}\), etc.

The definition and main properties of $n\wedge1$-operad modules pop out
in the study of lax $\Cat$-span multicategories -- one more direction
treated in the article as a category theory base of the whole subject.
These lax multicategories generalize strict multicategories associated
with the monad of free strict monoidal category.
Composition of \ainf- and \(A_\infty^{hu}\)\n-morphisms of several
arguments is presented as convolution of a certain colax $\Cat$-span 
multifunctor viewed as a coalgebra and the lax $\Cat$-span multifunctor
$\HOM$ viewed as an algebra.
This gives the multicategory of \ainf- or \(A_\infty^{hu}\)\n-algebras.

\subsection{Notations and conventions.}\label{sec-Notations-conventions}
We denote by $\NN$ the set of non-negative integers $\ZZ_{\ge0}$. 
By norm on $\NN^n$ we mean the function \(|\cdot|:\NN^n\to\NN\), 
\(j\mapsto|j|=\sum_{i=1}^nj^i\).
Let $\kk$ denote the ground commutative ring. 
Tensor product \(\tens_\kk\) will be denoted simply $\tens$. 
When a $\kk$\n-linear map $f$ is applied to an element $x$, the result
is typically written as $x.f=xf$. 
The tensor product of two maps of graded $\kk$\n-modules $f$, $g$ of
certain degree is defined so that for elements $x$, $y$ of arbitrary 
degree
\[ (x\tens y).(f\tens g) = (-1)^{\deg y\cdot\deg f}x.f\tens y.g.
\]
In other words, we strictly follow the Koszul rule. 
Composition of $\kk$\n-linear maps $X\rto{f}Y\rto{g}Z$ is usually
denoted $f\cdot g=fg:X\to Z$.
For other types of maps composition is often written as \(g\circ f=gf\).

We assume that each set is an element of some universe.
This universe is not fixed through the whole article.
For instance, the category of categories $\Cat$ means the category of
$\fu'$\n-small, locally $\fu$\n-small categories for some universes 
\(\fu\in\fu'\) (thus, a category $\cc$ is in $\Cat$ iff
\(\Ob\cc\in\fu'\) and \(\cc(X,Y)\in\fu\) for all $X,Y\in\Ob\cc$).
These universes are used tacitly, without being explicitly mentioned.

We consider the category of totally ordered finite sets and their 
non-decreasing maps.
An arbitrary totally ordered finite set is isomorphic to a unique set 
\(\bn=\{1<2<\dots<n\}\) via a unique isomorphism, $n\ge0$.
Functions of totally ordered finite set that we use in this article
are assumed to \emph{depend only on the isomorphism class of the set}.
Thus, it suffices to define them only for skeletal totally ordered
finite sets $\bn$.
The full subcategory of such sets and their non-decreasing maps is
denoted $\co_\sk$.
The full subcategory of $\Set$ formed by $\bn$, $n\ge0$, is denoted
$\cS_\sk$.

Whenever $I\in\Ob\co_\sk$, there is another totally ordered set 
\([I]=\{0\}\sqcup I\) containing $I$, where element $0$ is the smallest 
one.
Thus, \([n]=[\bn]=\{0<1<2<\dots<n\}\).

The list $A$, \dots, $A$ consisting of $n$ copies of the same object
$A$ is denoted $\sS{^n}A$.

For any graded \(\kk\)-module $M$ denote by $sM=M[1]$ the same module
with the grading shifted by 1: \(M[1]^k=M^{k+1}\).
Denote by \(\sigma:M\to M[1]\), \(M^k\ni x\mapsto x\in M[1]^{k-1}\)
the ``identity map'' of degree \(\deg\sigma=-1\).

\subsection{Motivation.}\label{sec-Motivation}
\ainf-algebras and \ainf-categories arise in symplectic topology in the
studies of Floer cohomology of Lagrangian submanifolds of symplectic
manifolds, see the monograph of K.~Fukaya, Y.-G.~Oh, H.~Ohta and K.~Ono 
\cite{FukayaOhOhtaOno:Anomaly}.
This article is devoted to \ainf-algebras, which in particular are
cochain complexes of $\kk$\n-modules (differential graded $\kk$-modules)
\[ \dots \rTTo^d X^{-1} \rTTo^d X^0 \rTTo^d X^1 \rTTo^d \dots\;,
\quad d^2 =0.
\]
Denote by $(\dg,\tens_\kk)$ the monoidal category of cochain complexes
of $\kk$-modules with chain maps as morphisms.
Tensor product of complexes is denoted $\tens=\tens_\kk$ as usual.
Sometimes we denote the same product by $\boxt$ instead.
The reason lies in expected generalization from \ainf-algebras to 
\ainf-categories.
The latter have the underlying structure of a $\dg$\n-quiver.
Unlike $\dg$\n-modules $\dg$\n-quivers $\ca$ and $\cb$ admit two
products: the external product $\ca\boxt\cb$ with the set of objects 
\(\Ob\ca\boxt\cb=\Ob\ca\times\Ob\cb\) and the tensor product
$\ca\tens\cb$ defined if and only if \(\Ob\ca=\Ob\cb\).
In the latter case \(\Ob\ca\tens\cb=\Ob\ca\).
For $\dg$\n-modules viewed as a particular case of $\dg$\n-quivers 
(\(\Ob\ca=\Ob\cb=\{1\}\)) both products coincide.

\ainfm-algebras are complexes $A\in\dg$ with the differential $m_1$
and operations \(m_n:A^{\tens n}\to A\), \(\deg m_n=2-n\), for $n\ge2$
such that
\begin{equation}
\sum_{j+p+q=n}(-1)^{jp+q}(1^{\tens j}\tens m_p\tens1^{\tens q})\cdot m_{j+1+q}=0
\label{eq-A8-equation}
\end{equation}
for all $n\ge2$. 
For instance, binary multiplication $m_2$ is a chain map, it is
associative up to the boundary of the homotopy $m_3$:
\[ (m_2\tens1)m_2 -(1\tens m_2)m_2 =m_3m_1
+(1\tens1\tens m_1 +1\tens m_1\tens1 +m_1\tens1\tens1)m_3,
\]
and so on.

\subsection{\texorpdfstring{$\dg$}{dg}-operads.}
Informally, non-symmetric operads are collections of operations, which
can be performed without permuting the arguments, in algebras of a
certain type. 
In this article an operad will mean a \emph{non-symmetric} operad. 

A \emph{(non-symmetric) operad} $\co$ is a collection of sets 
\((\co(n))_{n\in\NN}\) -- operations of arity $n$, an associative
family of substitution compositions
\[
\mu:\co(n_1)\times\dots\times\co(n_k)\times\co(k)\to\co(n_1+\dots+n_k)
\]
(one for each $\kk$\n-tuple \((n_1,\dots,n_k)\in\NN^k\), $k\in\NN$),
which has a two-sided unit \(\eta\in\co(1)\) -- the identity operation.

\begin{example}
The operad $\as$ of semigroups without unit has precisely one operation 
\(m^{(n)}:(x_1,\dots,x_n)\mapsto x_1\dots x_n\) for each $n>0$.
Thus, $\as(n)=\{m^{(n)}\}$ for positive $n$ and \(\as(0)=\emptyset\).

Similarly, there is the operad $\ass$ in $\Set$ with
$\ass(n)=\{m^{(n)}\}$ for all $n\ge0$.
$\ass$\n-algebras are monoids -- semigroups with a unit,which is
implemented by the nullary operation $m^{(0)}$.
\end{example}

Operations from a $\dg$\n-operad have in addition a degree and a 
boundary.
Category $\cm=\dg^\NN$ of collections of complexes
\((\cw(n))_{n\in\NN}\) is equipped with the tensor product $\odot$:
\[ (\cu\odot\cw)(n) 
=\bigoplus_{n_1+\dots+n_k=n}^{k\ge0} \cu(n_1)\tdt\cu(n_k)\tens\cw(k).
\]
The unit object $\1$ has \(\1(1)=\kk\), $\1(n)=0$ for $n\ne1$.

A \emph{(non-symmetric) $\dg$-operad} $\co$ is a monoid in
$\cm=(\dg^\NN,\odot)$, say
\((\co,\mu:\co\odot\co\to\co,\eta:\1\to\co)\).
Multiplication consists of maps
\begin{equation}
\mu:\co(n_1)\tdt\co(n_k)\tens\co(k) \to \co(n_1+\dots+n_k).
\label{eq-mu-OOOOOOOOO-O}
\end{equation}

\begin{example}
For any complex $X\in\dg$ there is the $\dg$-operad $\END X$ of its 
endomorphisms.
It has \((\END X)(n)=\und\dg(X^{\tens n},X)\), the complex of
$\kk$\n-linear maps \(X^{\tens n}\to X\) of certain degree.
Here $\und\dg$ is the category enriched in $\dg$ due to closedness of 
$(\dg,\tens)$.
\end{example}

\begin{definition}
An algebra $X$ over a $\dg$-operad $\co$ is a complex $X$ together with
a morphism of operads \(\co\to\END X\) (morphism of monoids in $\cm$).
\end{definition}

\begin{example}
The $\dg$\n-operad $\As$ is the $\kk$\n-linear envelope of the operad
$\as$ in $\Set$. 
They have \(\As(0)=0\) and $\As(n)=\kk m^{(n)}=\kk$ for $n>0$. 
The identity operation $m^{(1)}$ is the unit of the operad, and
$m^{(2)}=m$ is the binary multiplication.
$\As$\n-algebras are associative differential graded $\kk$\n-algebras
without unit.

Similarly, the operad $\ass$ in $\Set$ has the $\kk$\n-linear
envelope -- the $\dg$\n-operad $\Ass$ with $\Ass(n)=\kk m^{(n)}=\kk$
for all $n\ge0$. 
Clearly, $\Ass$\n-algebras are associative differential graded
$\kk$\n-algebras with the multiplication $m^{(2)}$ and the unit
$m^{(0)}=1^\su$.
\end{example}

\subsection{Model category structures.}
	\label{sec-Model-category-structures}
The following theorem is proved by Hinich in 
\cite[Section~2.2]{Hinich:q-alg/9702015}, except that he relates a
category with the category of complexes $\dg$, not with its power 
$\dg^S$. 
A generalization is given in \cite[Theorem~1.2]{Lyu-Hinich-Thm}.
It has the same formulation as below, however, $\dg$ means there the category of differential graded modules over a graded commutative ring.

\begin{theorem}[{\cite{Hinich:q-alg/9702015,Lyu-Hinich-Thm}}]\label{thm-Hinich-model}
Suppose that $S$ is a set, a category $\cc$ is complete and cocomplete
and \(F:\dg^S\rightleftarrows\cc:U\) is an adjunction. 
Assume that $U$ preserves filtering colimits. 
For any $x\in S$, $p\in\ZZ$ consider the object $\KK[-p]_x$ of $\dg^S$, 
\(\KK[-p]_x(x)=\bigl(0\to\kk\rTTo^1 \kk\to0\bigr)\) (concentrated in
degrees $p$ and $p+1$), \(\KK[-p]_x(y)=0\) for $y\ne x$. 
Assume that the chain map \(U(\inj_2):UA\to U(F(\KK[-p]_x)\sqcup A)\)
is a quasi-isomorphism for all objects $A$ of $\cc$ and all $x\in S$, 
$p\in\ZZ$. 
Equip $\cc$ with the classes of weak equivalences (resp. fibrations)
consisting of morphisms $f$ of $\cc$ such that $Uf$ is a
quasi-isomorphism (resp. an epimorphism). 
Then the category $\cc$ is a model category.
\end{theorem}

We shall recall also several constructions used in the proof of this
theorem.
They describe cofibrations and trivial cofibrations in $\cc$. 
Assume that \(M\in\Ob\dg^S\), \(A\in\Ob\cc\),
\(\alpha:M\to UA\in\dg^S\).
Denote by \(C=\Cone\alpha=(M[1]\oplus UA,d_{\Cone})\in\Ob\dg^S\) the
cone taken pointwise, that is, for any $x\in S$ the complex 
\(C(x)=\Cone\bigl(\alpha(x):M(x)\to(UA)(x)\bigr)\) is the usual cone.
Denote by \(\bar\imath:UA\to C\) the obvious embedding.
Let \(\eps:FU(A)\to A\) be the adjunction counit.
Following Hinich \cite[Section~2.2.2]{Hinich:q-alg/9702015} define an
object \(A\langle M,\alpha\rangle\in\Ob\cc\) as the pushout
\begin{diagram}[w=4em,LaTeXeqno]
FU(A) & \rTTo^\eps &A
\\
\dTTo<{F\bar\imath} &&\dTTo>{\bar\jmath}
\\
FC &\rTTo^g &\NWpbk A\langle M,\alpha\rangle
\label{dia-pushout-A<M-alpha>}
\end{diagram}
If $\alpha=0$, then \(A\langle M,0\rangle\simeq F(M[1])\sqcup A\) and 
$\bar\jmath=\inj_2$ is the canonical embedding.
We say that $M$ consists of free $\kk$\n-modules if for any $x\in S$,
$p\in\ZZ$ the $\kk$\n-module $M(x)^p$ is free.

The proof contains the following important statements.
If $M$ consists of free $\kk$\n-modules and $d_M=0$, then 
\(\bar\jmath:A\to A\langle M,\alpha\rangle\) is a cofibration.
It might be called an \emph{elementary standard cofibration}. If
\[ A \to A_1 \to A_2 \to \cdots
\]
is a sequence of elementary standard cofibrations, $B$ is a colimit
of this diagram, then the ``infinite composition'' map $A\to B$ is a
cofibration called a \emph{standard cofibration}
\cite[Section~2.2.3]{Hinich:q-alg/9702015}.

Assume that \(N\in\Ob\dg^S\) consists of free $\kk$\n-modules, $d_N=0$
and \(M=\Cone\bigl(1_{N[-1]}\bigr)=(N\oplus N[-1],d_{\Cone})\).
Then for any morphism \(\alpha:M\to UA\in\dg^S\) the morphism 
\(\bar\jmath:A\to A\langle M,\alpha\rangle\) is a trivial cofibration
in $\cc$ and a standard cofibration, composition of two elementary
standard cofibrations.
It is called a \emph{standard trivial cofibration}.
Any (trivial) cofibration is a retract of a standard (trivial)
cofibration \cite[Remark~2.2.5]{Hinich:q-alg/9702015}.

When \(F:\dg^S\to\cc\) is the functor of constructing a free
$\dg$\n-algebra of some kind, the maps $\bar\jmath$ are interpreted
as ``adding variables to kill cycles''.

The category $\Op$ of operads admits an adjunction 
\(F:\dg^\NN\rightleftarrows\Op:U\).
Applying \thmref{thm-Hinich-model} to this category one gets \cite[Proposition~1.8]{Lyu-Ainf-Operad}

\begin{proposition}%\label{pro-Operads-model-category}
Define weak equivalences (resp. fibrations) in $\Op$ as morphisms $f$
of $\Op$ such that $Uf$ is a quasi-isomorphism (resp. an epimorphism). 
These classes make $\Op$ into a model category.
\end{proposition}

This statement was proven previously in \cite{Hinich:q-alg/9702015}, \cite[Remark~2]{math/0101102} and follows from \cite[Theorem~1.1]{MR2821434}.

\begin{example}
Using Stasheff's associahedra one proves that there is a cofibrant
replacement \(\mainf\to\As\) where the graded operad $\mainf$ is freely 
generated by $n$\n-ary operations $m_n$ of degree $2-n$ for $n\ge2$. 
The differential is found as
\begin{equation}
m_n\partial =-\sum_{j+p+q=n}^{1<p<n}
(-)^{jp+q}(1^{\tens j}\tens m_p\tens1^{\tens q})\cdot m_{j+1+q}.
\label{eq-md-(-)(1m1)m}
\end{equation}
Basis \((m(t))\) of \(\mainf=T\bigl(\kk\{m_n\mid n\ge2\}\bigr)\) over
$\kk$ is indexed by isomorphism classes of ordered rooted trees $t$
without unary vertices (those with one incoming edge).
The tree $t^1$ which has just one vertex (the root and the leaf)
corresponds to the unit from $\mainf(1)$.

Algebras over the $\dg$\n-operad $\mainf$ are precisely \ainf-algebras
in the sense of \eqref{eq-A8-equation}.

Furthermore, the chain map \(\mainf(n)\to\As(n)\) is homotopy
invertible for each \(n\ge1\). 
One way to prove it is implied by a remark of Markl 
\cite[Example~4.8]{Markl:ModOp}. 
Another proof uses the operad of Stasheff associahedra
\cite{Stasheff:HomAssoc} and the configuration space of
$(n+1)$\n-tuples of points on a circle considered by Seidel in his
book \cite{SeidelP-book-Fukaya}. 
Details can be found in \cite[Proposition~1.19]{BesLyuMan-book}.
\end{example}

\subsection{Morphisms of operads.}
Besides usual homomorphisms of operads, which are chain maps of
degree 0, we consider also maps that change the degree.

\begin{definition}
A $\dg$\n-operad homomorphism \(t:\co\to\cp\) of degree
\(\bar{t}=r\in\ZZ\) is a collection of homogeneous $\kk$\n-linear maps
\(t(n):\co(n)\to\cp(n)\), $n\ge0$, of degree \(g(n)=(1-n)r\) such that
\begin{itemize}
\item \(1_\co.t(1)=1_\cp\);

\item for all \(k,n_1,\dots,n_k\in\NN\) the following square commutes
up to a certain sign:
\begin{diagram}[LaTeXeqno]
\co(n_1)\tdt\co(n_k)\tens\co(k) &\rTTo^\mu &\co(n_1+\dots+n_k)
\\
\dTTo<{t(n_1)\tdt t(n_k)\tens t(k)}&\sss(-1)^c &\dTTo>{t(n_1+\dots+n_k)}
\\
\cp(n_1)\tdt\cp(n_k)\tens\cp(k) &\rTTo^\mu &\cp(n_1+\dots+n_k)
\label{dia-ttttttttttttttt-t}
\end{diagram}
where the tensor product of homogeneous right maps $t(\_)$ is that of
the $\dg$\n-category $\und\dg$ and the sign is determined by
\end{itemize}
\begin{equation}
c= r\sum_{i=1}^k (i-1)(1-n_i)
+\frac{r(r-1)}2\sum_{1\le i<j\le k}(1-n_i)(1-n_j) 
+\frac{r(r-1)}2(1-k)\sum_{i=1}^k (1-n_i);
\label{eq-c-r-rr/2}
\end{equation}
\begin{itemize}
\item for all $n\in\ZZ$ 
\[ d\cdot t(n) =(-1)^{g(n)}t(n)\cdot d: \co(n) \to \cp(n). 
\]
\end{itemize}
\end{definition}

Notice that the only functions \(g:\NN\to\ZZ\) that satisfy the 
equations
\[ g(1)=0, \qquad g(n_1)+ \dots +g(n_k) +g(k) =g(n_1+\dots+n_k)
\]
are functions \(g(n)=(1-n)r\) for some $r\in\ZZ$.

\begin{example}
Let \((X,d_X)\) be a complex of $\kk$\n-modules and let 
\((X[1],d_{X[1]}=-\sigma^{-1}\cdot d_X\cdot\sigma)\) be its suspension.
There is a $\dg$\n-operad morphism
\[ \varSigma =\HOM(\sigma;\sigma^{-1}) 
=\HOM(\sigma;1)\cdot\HOM(1;\sigma^{-1}): \END(X[1]) \to \END X
\]
of degree 1. 
Thus, the mapping
\(f\mapsto(-1)^{nf}\sigma^{\tens n}\cdot f\cdot\sigma^{-1}\),
\[ \varSigma(n)=\und\dg(\sigma^{\tens n};1)\cdot\und\dg(1;\sigma^{-1}):
\und\dg(X[1]^{\tens n},X[1]) \to \und\dg(X^{\tens n},X),
\]
has degree $1-n$.
The sign $(-1)^c$, \(c=\sum_{i=1}^k(i-1)(1-n_i)\), pops out in the
following procedure.
Write down the tensor product corresponding to the left vertical arrow
of \eqref{dia-ttttttttttttttt-t} for $t=\varSigma$:
\[ (\sigma^{\tens n_1}\tens\sigma^{-1}) 
\tens (\sigma^{\tens n_2}\tens\sigma^{-1})
\tdt (\sigma^{\tens n_k}\tens\sigma^{-1})
\tens (\sigma^{\tens k}\tens\sigma^{-1});
\]
move factors $\sigma^{-1}$ using the Koszul rule to their respective
opponents, factors $\sigma$ of \(\sigma^{\tens k}\), in order to 
cancel them and to obtain finally 
\(\sigma^{\tens(n_1+\dots+n_k)}\tens\sigma^{-1}\).

Maps \(\varSigma(n)\) commute with the differential in the graded
sense because their factors $\sigma^{\pm1}$ do.
\end{example}

\begin{remark}\label{rem-homogeneous-operad-homomorphisms}
Summands
	\(r(r-1)/2\sum_{1\le i<j\le k}(1-n_i)(1-n_j)
	+(1-k)r(r-1)/2\sum_{i=1}^k(1-n_i)\)
of $c$ make sure that the composition of two morphisms of operads
\(t:\co\to\cp\) and \(u:\cp\to\cq\) of degree $\bar t$ and $\bar u$ 
respectively be an operad morphism of degree $\bar t+\bar u$.
Furthermore, if all homogeneous maps \(t(n):\co(n)\to\cp(n)\) for 
\(t:\co\to\cp\) are invertible, than there is an inverse morphism of 
operads \(t^{-1}:\cp\to\co\) of degree $-\bar t$ with
\(t^{-1}(n)=t(n)^{-1}\).

Let $\co$ be a $\dg$\n-operad, $\cp$ be a graded operad and
\(t:\co\to\cp\) be an invertible graded operad homomorphism of degree 
$r$ (\(1_\co.t(1)=1_\cp\) and \eqref{dia-ttttttttttttttt-t} holds).
Then $\cp$ has a unique differential $d$ which turns it into a
$\dg$\n-operad and \(t:\co\to\cp\) into a $\dg$\n-operad homomorphism
of degree $r$.
\end{remark}

\begin{example}
There is a version of the $\dg$\n-operad $\mainf$ denoted $A_\infty$.
This is a $\dg$\n-operad freely generated as a graded operad by
$n$\n-ary operations $b_n$ of degree $1$ for $n\ge2$.
The differential is defined as
\[  b_n\partial =-\sum_{j+p+q=n}^{1<p<n}
(1^{\tens j}\tens b_p\tens1^{\tens q})\cdot b_{j+1+q}.
\]
Comparing the differentials we find that these two operads are 
isomorphic via an isomorphism of degree 1
\[ \varSigma: A_\infty \to \mainf, \qquad b_i \mapsto m_i.
\]
In fact, due to \eqref{dia-ttttttttttttttt-t}
\begin{align*}
&[(1^{\tens j}\tens b_p\tens1^{\tens q})b_{j+1+q}].\varSigma(j+p+q)
\\
&= (-1)^{j(1-p)+1-p}
[(1^{\tens j}\tens b_p\tens1^{\tens q}).
(\varSigma(1)^{\tens j}\tens\varSigma(p)\tens\varSigma(1)^{\tens q})] 
\cdot [b_{j+1+q}.\varSigma(j+1+q)]
\\
&=(-1)^{(j+1)(1-p)} [(1^{\tens j}\tens m_p\tens1^{\tens q})m_{j+1+q}].
\end{align*}
Therefore,
\begin{align*}
m_n\partial
&=(b_n.\varSigma(n))\partial =(-1)^{1-n}(b_n\partial).\varSigma(n)
\\
&=(-1)^n\sum_{j+p+q=n}^{1<p<n}
[(1^{\tens j}\tens b_p\tens1^{\tens q})b_{j+1+q}].\varSigma(n)
\\
&=(-1)^n\sum_{j+p+q=n}^{1<p<n} (-1)^{(j+1)(1-p)}
(1^{\tens j}\tens m_p\tens1^{\tens q})m_{j+1+q}
\\
&=-\sum_{j+p+q=n}^{1<p<n}
(-1)^{jp+q}(1^{\tens j}\tens m_p\tens1^{\tens q})\cdot m_{j+1+q}
\end{align*}
which coincides with \eqref{eq-md-(-)(1m1)m}.
This fixes the differential on $\mainf$ since $m_n$ generate the
graded operad.
An easy lemma shows that it suffices to verify graded commutation of
the differential and any operad homomorphism on generators.
In particular, \(\varSigma:A_\infty\to\mainf\) commutes with
$\partial$ in the graded sense.

Knowing that $\mainf$ is homotopy isomorphic to its cohomology $\As$,
we conclude that $A_\infty$ is homotopy isomorphic to its cohomology as 
well.
There is an isomorphism of degree~1 between graded operads 
\(\varSigma:H^\bull(A_\infty)\to\As\).
Hence, \(H^\bull(A_\infty(n))=\kk[1-n]\) for $n\ge1$.

For any algebra \(A\in\dg\) over the $\dg$\n-operad $\mainf$ the 
$\dg$\n-module $A[1]$ becomes an algebra over the $\dg$\n-operad 
$A_\infty$ so that the square of operad homomorphisms commutes:
\begin{diagram}[nobalance,bottom,LaTeXeqno]
A_\infty &\rTTo &\END A[1]
\\
\dTTo<\varSigma &&\dTTo>{\HOM(\sigma;\sigma^{-1})} &&\qquad, \qquad 
(-1)^n\sigma^{\tens n}\cdot b_n\cdot\sigma^{-1} =m_n: A^{\tens n} \to A, 
\quad n\ge1.
\\
\mainf &\rTTo &\END A
\label{dia-A8-EndA1-A8-EndA}
\end{diagram}
Verification is straightforward.
\end{example}

\begin{example}
Approaching homotopy unital \ainfm-algebras we start with strictly 
unital ones. 
They are governed by the operad $\mainf^\su$ generated over $\mainf$ by 
a nullary degree 0 cycle $\one^\su$ subject to the following relations:
\[ (1\tens\one^\su)m_2 =1, \quad (\one^\su\tens1)m_2 =1, \quad
(1^{\tens a}\tens\one^\su\tens1^{\tens c})m_{a+1+c} =0
\text{ \ if \ } a+c>1.
\]
There is a standard trivial cofibration and a homotopy isomorphism
	$\mainf^\su\rCof~\Sim \mainf^\su\langle\one^\su-\sfi,
	\sfj\rangle=\mainf^\su\langle \sfi,\sfj\rangle$,
where $\sfi$, $\sfj$ are two nullary operations, $\deg \sfi=0$,
$\deg \sfj=-1$, with \(\sfi\partial=0\), \(\sfj\partial=\one^\su-\sfi\). 

A cofibrant replacement \(\mainf^{hu}\to\Ass\) is constructed as a 
$\dg$\n-suboperad of \(\mainf^\su\langle\sfi,\sfj\rangle\) generated as 
a graded operad by $\sfi$ and $n$\n-ary operations of degree $4-n-2k$
\[ m_{n_1;n_2;\dots;n_k} =(1^{\tens n_1}\tens \sfj\tens1^{\tens n_2} 
\tens\sfj\tdt1^{\tens n_{k-1}}\tens 
\sfj\tens1^{\tens n_k})m_{n+k-1},
\]
where \(n=\sum_{q=1}^kn_q\), $k\ge1$, $n_q\ge0$, \(n+k\ge3\).
Notice that the graded operad \(\mainf^{hu}\) is free.
See \cite[Section~1.11]{Lyu-Ainf-Operad} for the proofs.

One can perform all the above steps also for the operad $A_\infty$:

1) Adding to $A_\infty$ a nullary degree $-1$ cycle $\bone^\su$ subject 
to the relations:
\begin{equation}
(1\tens\bone^\su)b_2 =1, \quad (\bone^\su\tens1)b_2 =-1, \quad
(1^{\tens a}\tens\bone^\su\tens1^{\tens c})b_{a+1+c} =0
\text{ \ if \ } a+c>1.
\label{eq-(1bone)b2-1-(bone1)b2-1}
\end{equation}
The resulting operad is denoted $A_\infty^\su$.

2) Adding to $A_\infty^\su$ two nullary operations $\bi$, $\bj$, 
$\deg\bi=-1$, $\deg\bj=-2$, with \(\bi\partial=0\), 
\(\bj\partial=\bi-\bone^\su\).
The standard trivial cofibration
$A_\infty^\su\rCof~\Sim A_\infty^\su\langle\bi,\bj\rangle$ is a homotopy 
isomorphism.

3) $A_\infty^{hu}$ is a $\dg$\n-suboperad of
$A_\infty^\su\langle\bi,\bj\rangle$ generated as a graded operad by 
$\bi$ and $n$\n-ary operations of degree $3-2k$
\[ b_{n_1;n_2;\dots;n_k} 
=(1^{\tens n_1}\tens\bj\tens1^{\tens n_2}\tens\bj\tdt
1^{\tens n_{k-1}}\tens\bj\tens1^{\tens n_k})b_{n+k-1},
\]
where \(n=\sum_{q=1}^kn_q\), $k\ge1$, $n_q\ge0$, \(n+k\ge3\).

The obtained operads are related to the previous ones by invertible 
homomorphisms of degree 1, extending
\(\varSigma:b_n\mapsto m_n\),
\[ \varSigma: A_\infty^\su \to \mainf^\su, \; \bone^\su\mapsto \one^\su; 
\quad \varSigma: A_\infty^\su\langle\bi,\bj\rangle \to 
\mainf^\su\langle\sfi,\sfj\rangle,\; \bi\mapsto \sfi,\; \bj\mapsto \sfj; 
\quad \varSigma: A_\infty^{hu} \to \mainf^{hu}.
\]
The latter is a restriction of the previous map.
For algebras $A$ over operads $\mainf^\su$,
\(\mainf^\su\langle\sfi,\sfj\rangle\), $\mainf^{hu}$ the complex $A[1]$ 
obtains a structure of an algebra over the operad $A_\infty^\su$, 
\(A_\infty^\su\langle\bi,\bj\rangle\) or $A_\infty^{hu}$ due to a 
property similar to \eqref{dia-A8-EndA1-A8-EndA}, in particular,
\[ \bone^\su\sigma^{-1} =\one^\su, \qquad \bi\sigma^{-1} =\sfi, \qquad 
\bj\sigma^{-1} =\sfj: \kk \to A.
\]
\end{example}

\subsection{Morphisms of \texorpdfstring{$A_\infty$}{A8}-algebras.}
Objects of the category \(\cOm^{\NN\sqcup\NN\sqcup\NN}\) are written as 
triples of collections
\((\ca;\cp;\cb)=(\ca(n);\cp(n);\cb(n))_{n\in\NN}\) of complexes.
An \emph{operad bimodule} is defined as a triple \((\ca;\cp;\cb)\), 
consisting of operads $\ca$, $\cb$ and an $\ca$-$\cb$-bimodule $\cp$
in the monoidal category \((\dg^\NN,\odot)\).
The actions \(\lambda:\ca\odot\cp\to\cp\) and \(\rho:\cp\odot\cb\to\cp\) 
consist of chain maps
\begin{align*}
\lambda_{n_1,\dots,n_k}:\ca(n_1)\tdt\ca(n_k)\tens\cp(k) &\to 
\cp(n_1+\dots+n_k),
\\
\rho_{n_1,\dots,n_k}:\cp(n_1)\tdt\cp(n_k)\tens\cb(k) &\to 
\cp(n_1+\dots+n_k).
\end{align*}

The category of operad bimodules $\nOp1_1$ has morphisms 
\((f;h;g):(\ca;\cp;\cb)\to(\cc;\cq;\cd)\), where \(f:\ca\to\cc\), 
\(g:\cb\to\cd\) are morphisms of $\dg$-operads and
\(h:\cp\to\sS{_f}\cq_g\) is an $\ca$-$\cb$-bimodule morphism, where 
\(\sS{_f}\cq_g=\cq\) obtains its $\ca$-$\cb$-bimodule structure via
$f$, $g$.

\begin{example}
Let $X$, $Y$ be objects of $\dg$ (complexes of $\kk$\n-modules).
Define a collection \(\HOM(X,Y)\) as
\(\HOM(X,Y)(n)=\und\dg(X^{\tens n},Y)\).
Substitution composition \(\HOM(X,Y)\odot\HOM(Y,Z)\to\HOM(X,Z)\) and 
obvious units \(\1\to\HOM(X,X)\) make the category of complexes enriched 
in the monoidal category $(\dg^\NN,\odot)$. 
In particular, \(\END X=\HOM(X,X)\) are algebras in $\dg^\NN$
($\dg$-operads). 
The collection \(\HOM(X,Y)\) is an \(\END X\)-\(\END Y\)-bimodule. 
The multiplication and the actions are induced by substitution 
composition. 
\end{example}

In the nearest sections we use the shorthand \((\co,\cp)\) for an operad 
$\co$\n-bimodule \((\co;\cp;\co)\).

\subsubsection{Morphisms come from bimodules.}
Consider the operad bimodule \((\As,\As)\), where the first term is an 
operad and the second term is a regular bimodule.
One easily checks that a morphism of operad bimodules
\[ (\As;\As;\As) \to (\END X;\HOM(X,Y);\END Y)
\]
amounts to a morphism \(f:X\to Y\) of associative differential graded
$\kk$-algebras without units.

There is a pair of adjoint functors 
\(T:\dg^{\NN\sqcup\NN\sqcup\NN}\rightleftarrows\nOp1_1:U\), 
\(T(\cu;\cx;\cw)=(T\cu;T\cu\odot\cx\odot T\cw;T\cw)\).
Applying Hinich's \thmref{thm-Hinich-model} we get

\begin{proposition}[Proposition~2.2 \cite{Lyu-Ainf-Operad}]
Define weak equivalences (resp. fibrations) in $\nOp1_1$ as morphisms 
$f$ of $\nOp1_1$ such that $Uf$ is a quasi-isomorphism (resp. an 
epimorphism). 
These classes make $\nOp1_1$ into a model category.
\end{proposition}

Cofibrant replacement of a $\dg$-operad bimodule $(\co,\cp)$ is a 
trivial fibration \((\ca,\cF)\to(\co,\cp)\) (surjective mapping inducing 
isomorphism in cohomology) such that the only map \((\1,0)\to(\ca,\cF)\) 
is a cofibration in $\nOp1_1$, for instance, the graded operad bimodule 
\((\ca,\cF)\) is freely generated.

\begin{definition}
A $\dg$-operad bimodule homomorphism
\((u,v,w):(\ca;\cp;\cb)\to(\cc;\cq;\cd)\) of degree \(r\in\ZZ\) is a 
pair of $\dg$\n-operad homomorphisms \(u:\ca\to\cc\), \(w:\cb\to\cd\) of 
degree $r$ and a collection of homogeneous $\kk$\n-linear maps
\(v(n):\cp(n)\to\cq(n)\), $n\ge0$, of degree \(r(1-n)\) such that
\begin{itemize}
\item for all \(k,n_1,\dots,n_k\in\NN\) the following squares commute up 
to the sign given in \eqref{eq-c-r-rr/2}:
\begin{diagram}
\ca(n_1)\tdt\ca(n_k)\tens\cp(k) &\rTTo^\lambda &\cp(n_1+\dots+n_k)
\\
\dTTo<{u(n_1)\tdt u(n_k)\tens v(k)} &\sss(-1)^c&\dTTo>{v(n_1+\dots+n_k)}
\\
\cc(n_1)\tdt\cc(n_k)\tens\cq(k) &\rTTo^\lambda &\cq(n_1+\dots+n_k)
\end{diagram}
\begin{diagram}
\cp(n_1)\tdt\cp(n_k)\tens\cb(k) &\rTTo^\rho &\cp(n_1+\dots+n_k)
\\
\dTTo<{v(n_1)\tdt v(n_k)\tens w(k)} &\sss(-1)^c&\dTTo>{v(n_1+\dots+n_k)}
\\
\cq(n_1)\tdt\cq(n_k)\tens\cd(k) &\rTTo^\rho &\cq(n_1+\dots+n_k)
\end{diagram}

\item for all $n\in\NN$
\[ d\cdot v(n) =(-1)^{r(1-n)}v(n)\cdot d: \cp(n) \to \cq(n). 
\]
\end{itemize}
\end{definition}

Properties of $\dg$-operad bimodule homomorphisms are quite similar to 
those of $\dg$\n-operad homomorphisms, described in
\remref{rem-homogeneous-operad-homomorphisms}.

\begin{proposition}[cf. Proposition~2.7 \cite{Lyu-Ainf-Operad}]
There is an operad bimodule \((A_\infty,F_1)\) freely generated by
$n$-ary elements $f_n$ of degree $0$ over the graded operad $A_\infty$.
The differential for it is given by 
\begin{equation}
f_k\partial
=\sum_{r+n+t=k}^{n>1} (1^{\tens r}\tens b_n\tens1^{\tens t}) f_{r+1+t}
-\sum^{l>1}_{i_1+\dots+i_l=k} (f_{i_1}\tens f_{i_2}\tdt f_{i_l}) b_l.
\label{eq-fk-partial-1b1f-fffb}
\end{equation}
$F_1$-maps are \ainf-algebra morphisms (for algebras written with 
operations $b_n$).
There is an isomorphic form of this bimodule -- the operad bimodule 
\((\mainf,\rmF_1)\) freely generated by $n$-ary elements $\sff_n$ of 
degree $1-n$ over the graded operad $\mainf$.
The differential for it is given by 
\begin{gather*}
\sff_k\partial =\sum_{r+n+t=k}^{n>1} 
(-1)^{(r+1)n+t-1}(1^{\tens r}\tens m_n\tens1^{\tens t}) \sff_{r+1+t} 
-\!\! \sum^{l>1}_{i_1+\dots+i_l=k} \!\! (-1)^\sigma 
(\sff_{i_1}\tens\sff_{i_2}\tdt\sff_{i_l}) m_l, 
\\
\sigma =k-1 +\sum_{j=1}^k j(1-i_j).
\end{gather*}
The isomorphism between these bimodules
\[ (\varSigma,\varSigma): (A_\infty,F_1) \to (\mainf,\rmF_1),
\qquad b_i \mapsto m_i, \quad f_k \mapsto \sff_k
\]
has degree 1.
$\rmF_1$-maps are \ainf-algebra morphisms $A\to B$ (for algebras written 
with operations $m_n$).
The two notions of \ainf-morphisms agree in the sense that the square of 
operad bimodule maps
\begin{diagram}[nobalance]
(A_\infty;F_1;A_\infty) &\rTTo &(\END A[1];\HOM(A[1],B[1]);\END B[1])
\\
\dTTo<{(\varSigma;\varSigma;\varSigma)}
&&\dTTo>{(\HOM(\sigma;\sigma^{-1});
\HOM(\sigma;\sigma^{-1});\HOM(\sigma;\sigma^{-1}))}
\\
(\mainf;\rmF_1;\mainf) &\rTTo &(\END A;\HOM(A,B);\END B)
\end{diagram}
commutes.
The bimodule \((\mainf,\rmF_1)\) is a cofibrant replacement of
$(\As,\As)$.
\((\mainf,\rmF_1)\to(\As,\As)\) is a homotopy isomorphism in 
\(\dg^{\NN\sqcup\NN}\).
\end{proposition}

\begin{proof}
There is a degree preserving $\dg$-operad bimodule isomorphism 
\((\mainf,\rmF_1)\to(\mainf,\rmF'_1)\),
\(\sff_k\mapsto(-1)^{1-k}\sff'_k\), to the bimodule presented in 
\cite{Lyu-Ainf-Operad}, whose generators are denoted here by 
\(\sff'_k\).
Thus all previously proven properties are inherited by the
\ainfm-bimodule described above. 
\end{proof}

\subsection{Tensor coalgebra.}
The \emph{tensor $\kk$-module} of $A[1]$ is 
\(T(A[1])=\oplus_{n\ge0}T^n(A[1])=\oplus_{n\ge0}A[1]^{\tens n}\).
Multiplication in an \ainf-algebra $A$ is given by the operations of 
degree $+1$
\[ b_n: T^n(A[1]) =A[1]^{\tens n} \to A[1], \qquad n\ge1.
\]

Recall that $\kk$\n-linear maps are composed \emph{from left to right}. 
Operations $b_n$ have to satisfy the \ainf-equations, $n\ge1$:
\begin{equation*}
\sum_{r+k+t=n} (1^{\tens r}\tens b_k\tens1^{\tens t})b_{r+1+t} =0:
T^n(A[1]) \to A[1].
%\label{eq-b-b-0}
\end{equation*}

Tensor $\kk$-module \(T(A[1])\) has a coalgebra structure: the cut 
coproduct
\[ \Delta(x_1\tens x_2\tdt x_n)
=\sum_{k=0}^n x_1\tdt x_k\bigotimes x_{k+1}\tdt x_n.
\]

An \ainf-structure on a graded $\kk$-module $A$  is equivalent to 
\(b^2=0\), where \(b:T(A[1])\to T(A[1])\) is a coderivation of degree 
$+1$ given by the formula
\[ b =\sum_{r+k+t=n} 1^{\tens r}\tens b_k\tens1^{\tens t}:
T^n(A[1]) \to T(A[1]), \qquad b_0 =0.
\]
In particular, \(b\Delta=\Delta(1\tens b+b\tens1)\).

\subsection{Comultiplication and composition.}
In order to have an associative composition of \((\co,\cF)\)-morphisms, 
we postulate an associative counital comultiplication 
\(\Delta:\cF\to\cF\odot_\co\cF\) in the category of $\co$-bimodules.

Suppose that $A$, $B$, $C$ are $\co$-algebras and $g:\cF\to\HOM(A,B)$, 
$h:\cF\to\HOM(B,C)$ are \((\co,\cF)\)-morphisms.
Then their composition is defined as the convolution
\begin{equation*}
g\cdot h =\bigl[ \cF \rTTo^\Delta \cF\odot_\co\cF \rTTo^{g\odot h}
\HOM(A,B)\odot_{\END B}\HOM(B,C) \to \HOM(A,C) \bigr].
\end{equation*}
For \((\As,\As)\) the comultiplication is the identity map.
For $A_\infty$-morphisms the comultiplication is chosen as
\[ \Delta:F_1\to F_1\odot_{A_\infty}F_1, \quad f_n\Delta 
=\sum_{i_1+\dots+i_k=n} (f_{i_1}\tens f_{i_2}\tdt f_{i_k})\tens f_k,
\]
see \cite[Section~2.16]{Lyu-Ainf-Operad}, which results in the 
composition
\[ (g\cdot h)_n
=\sum_{i_1+\dots+i_k=n}(g_{i_1}\tens g_{i_2}\tdt g_{i_k})h_k.
\]

\subsection{Homomorphisms with \texorpdfstring{$n$}n arguments.}
\ainf-morphisms with several arguments \(f:A_1,\dots,A_n\to B\) are 
defined as augmented $\dg$\n-coalgebra morphisms
\[ \hat f: T(A_1[1])\tdt T(A_n[1]) \to T(B[1]).
\]
Here both augmented graded coalgebras \((C,\Delta,\eps,\eta)\) are of 
the form
	\((\kk\oplus\bar C,\Delta(x)=1\tens x+x\tens1+\bar{\Delta}(x)\;
	\forall\,x\in\bar{C},\pr_1,\inj_1)\),
where the non-counital coassociative coalgebra \((\bar C,\bar\Delta)\) 
is conilpotent (cocomplete \cite[Section~1.1.2]{Lefevre-Ainfty-these}, 
\cite[Section~4.3]{math.RT/0510508}).
Thus \((\bar C,\bar\Delta)\) is identified with a $T^{\ge1}$\n-coalgebra 
\cite[Proposition~6.8]{BesLyuMan-book}, see also
\propref{pro-Tge1-coalgebra} of the current article.
In the category of such augmented graded coalgebras the target 
\(\kk\oplus T^{\ge1}(B[1])\) is cofree, see
\corref{cor-Tge1-coalgebra}, hence, augmented graded coalgebra morphisms 
$\hat f$ are in bijection with the degree 0 $\kk$\n-linear maps
\[ f: T(A_1[1])\tdt T(A_n[1]) \to B[1],
\]
whose restriction to \(T^0(A_1[1])\tdt T^0(A_n[1])\simeq\kk\) vanishes.
The morphism $\hat f$ will be a chain map, $\hat fb=b\hat f$,
if and only if
\begin{multline}
\sum_{q=1}^n \sum_{r+c+t=\ell^q}^{c>0}
\bigl[ \boxt^{i\in\bn}T^{\ell^i}sA_i
\rTTo^{1^{\boxt(q-1)}\boxt(1^{\tens r}\tens b_c
\tens1^{\tens t})\boxt1^{\boxt(n-q)}}
\\
T^{\ell^1}sA_1\boxt\dots\boxt T^{\ell^{q-1}}sA_{q-1}\boxt T^{r+1+t}sA_q
\boxt T^{\ell^{q+1}}sA_{q+1}\boxt\dots\boxt T^{\ell^n}sA_n
\rTTo^{f_{\ell-(c-1)e_q}} sB \bigr]
\\
\hskip\multlinegap
=\sum_{\substack{j_1,\dots,j_k\in\NN^n-0\\j_1+\dots+j_k=\ell}}^{k>0}
\bigl[ \boxt^{i\in\bn}T^{\ell^i}sA_i \rTTo^\sim
\boxt^{i\in\bn}\tens^{p\in\mb k}T^{j_p^i}sA_i \hfill
\\
\rTTo^\sim \tens^{p\in\mb k}\boxt^{i\in\bn}T^{j_p^i}sA_i 
\rTTo^{\tens^{p\in\mb k} f_{j_p}}\tens^{p\in\mb k}sB \rTTo^{b_k} 
sB\bigr].
\label{eq-sff-b-commute}
\end{multline}

\subsection{Composition of \texorpdfstring{$A_\infty$}{A8}-morphisms.}
A \emph{tree} $t$ is a composable sequence of non-\hspace{0pt}decreasing 
maps of totally ordered finite sets \(\mb m=\{1<2<\dots<m\}\)
(objects of $\co_\sk$)
\begin{equation}
t =\bigl(t(0) \rTTo^{t_1} t(1)\rTTo^{t_2} \dots\;t(l-1) \rTTo^{t_l} t(l) 
=\mb1\bigr).
\label{eq-t-t(0)-t(l)-1}
\end{equation}
Composition of a family of $A_\infty$-morphisms
$(g_h^b)_{h>0}^{b\in t(h)}$ indexed by vertices of a tree $t$ is
\begin{multline*}
g_j =\comp(t)(g_h^b)_j =\bigl[ \tens^{a\in t(0)}T^{j^a}sA_0^a
\rTTo^{\tens^{b_1\in t(1)}\wh{g_1^{b_1}}} 
\tens^{b_1\in t(1)}T^{j_1^{b_1}}sA_1^{b_1}
\\
\rTTo^{\tens^{b_2\in t(2)}\wh{g_2^{b_2}}}
\tens^{b_2\in t(2)}T^{j_2^{b_2}}sA_2^{b_2} \to\dots \rTTo^{g_l^1} 
sA_l^1\bigr],
\end{multline*}
and $g_h^b$ is given via its components
\(g_{h,j}^b:\tens^{a\in t_h^{-1}b}T^{j^a}sA_{h-1}^a\to sA_h^b\).
Here $j$ belongs to \(\NN^{t_h^{-1}b}\).

Explicit formula for the composition is
\begin{align*}
g_j =\sum^{t-\text{tree }\tau}_{\forall a\in t(0)\,|\tau(0,a)|=j^a}
\tens^{h\in I} \tens^{b\in t(h)} \tens^{p\in\tau(h,b)}
g^b_{h,|\tau_{(h-1,a)}^{-1}(p)|_{a\in t_h^{-1}b}}.
\end{align*}
A $t$\n-tree is a functor \(\tau:t\to\co_\sk\), \(\tau(\troot)=\mb1\), 
where the poset $t$ is the free category (of paths) built on the quiver 
$t$ oriented towards the root.
It has the set of objects $\iV(t)$, the set of vertices of $t$, and the 
root is the terminal object of $t$.
A \emph{symmetric tree} is a sequence \eqref{eq-t-t(0)-t(l)-1}, where 
maps $t_h$ are not supposed to be non-decreasing.
A \emph{braided tree} is a symmetric tree (a functor)
\(t:[l]\to\cS_\sk\) for which the maps $t_h$ satisfy an extra condition: 
for all \(0\le p<q<r\le l\), \(a,b\in t(p)\) inequalities $a<b$ and 
$t_{p\to q}(a)>t_{p\to q}(b)$ imply $t_{p\to r}(a)\ge t_{p\to r}(b)$ 
(see \cite[Definition~2.3]{BesLyuMan-book}).
Labelling of a (symmetric) tree in a set $S$ consists of functions 
\(\ell:t(h)\to S\), \(0\le h\le l\).

\subsection{Graded complexes.}
Now let us explain what is new in the main part of the article.
In \secref{sec-Lax-Cat-multicategories} we define such notions as lax 
$\Cat$-span multicategories and their particular cases: lax $\Cat$-span 
operads, lax $\Cat$-multicategories and lax $\Cat$-operads.
These are accompanied with definitions of lax $\Cat$-span multifunctors 
and $\Cat$-span multinatural transformations.
Together they form a 2-category, whose objects are lax $\Cat$-span 
multicategories.

\subsubsection{\texorpdfstring{$\Cat$}{Cat}-operad of graded 
\texorpdfstring{$\kk$}k-modules.}
 \label{sec-Cat-operad-graded-k-modules}
Let us describe examples $\mcG$, $\mcDG$ of a weak $\Cat$-operad: that 
of (differential) graded $\kk$-modules. 
We define $\mcG(n)=\gr^{\NN^n}$, $\mcDG(n)=\dg^{\NN^n}$.
The structure of one is obtained from the structure of the other by 
forgetting or introducing the differential.
So we describe only one of them.

A functor is given for a tree $t$
\[ \circledast(t): \prod_{(h,b)\in\IV(t)}\mcG(t_h^{-1}b) \to \mcG(t(0)), 
\qquad (\cp_h^b)_{(h,b)\in\IV(t)} \mapsto \circledast(t)
(\cp_h^b)_{(h,b)\in\IV(t)}.
\]

Namely, for a tree $t:[n]\to\co_\sk$ with \([n]=\{0,1,2,\dots,n\}\),
\begin{equation*}
\circledast(t)(\cp_h^b)_{(h,b)\in\IV(t)}(z) 
=\bigoplus^{t-\text{tree }\tau}_{\forall a\in t(0)\,|\tau(0,a)|=z^a}
\bigotimes^{h\in I} \bigotimes^{b\in t(h)} \bigotimes^{p\in \tau(h,b)} 
\cp_h^b\Bigl(
\bigl(|\tau_{(h-1,a)\to(h,b)}^{-1}(p)|\bigr)_{a\in t_h^{-1}b} \Bigr).
\end{equation*}

Let $[I]=[n]$, $[J]=[m]$.
Correspondingly we denote \(I=[I]-\{0\}=\bn\), \(J=[J]-\{0\}=\mb m\).
Let \(f:I\to J\) be an isotonic map.
Let isotonic map \(\psi=[f]:[J]\to[I]\) viewed as a functor be right 
adjoint to the functor \([f]^*=(0\mapsto0)\sqcup f:[I]\to[J]\).
This means that for any \(x\in[I]\), \(y\in[J]\) the following 
inequalities are equivalent:
\begin{equation}
x\le[f](y) \Longleftrightarrow [f]^*(x)\le y.
\label{eq-x<fy-fx<y}
\end{equation}
Formula reads
\[ [f]: [J] \to [I], \qquad y \mapsto 
[f](y) \overset{\text{def}}= \max([f]^*)^{-1}([0,y]).
\]
Here \([0,y]=\{z\in[J]\mid z\le y\}\subset[J]\).
If \(t:[I]\to\cS_\sk\) is a (plain, symmetric or braided) tree, then the 
composite functor
\(t_\psi=\bigl([J]\rTTo^\psi \relax[I]\rTTo^t \cS_\sk\bigr)\) is also a 
tree.
If $a,b\in[I]$, $a\le b$, $c\in t(b)$, then the tree
\(t^{|c}_{[a,b]}:[a,b]\to\cS_\sk\) is the subtree of $t$ consisting of  
vertices in the preimage of $c$, whose level $k$ is above $a$:
\(t^{|c}_{[a,b]}(k)=t_{k\to b}^{-1}(c)\) for \(k\in[a,b]\).

For $f$ and \(\psi=[f]\) as above a natural bijection is constructed:
\[ \lambda^f: \circledast(t)(\cp_h^b)_{(h,b)\in\IV(t)} \to
\circledast(t_\psi)\bigl( \circledast(t^{|c}_{[\psi(g-1),\psi(g)]})
(\cp_h^b)_{(h,b)\in\IV(t^{|c}_{[\psi(g-1),\psi(g)]})} 
\bigr)_{(g,c)\in\IV(t_\psi)}.
\]

In \secref{sec-Morphisms-several-entries} we explain that morphisms with 
$n$ entries of algebras over operads form an $n\wedge1$-operad module.
In particular, we find this module for \ainf-algebras.

\subsection{\texorpdfstring{$n\wedge1$}{n1}-operad modules.}
	\label{sec-n1-operad-modules}
An $n\wedge1$-operad module is a sequence
\((\ca_1,\dots,\ca_n;\cp;\cb)\), where $\cb$, $\ca_i$ are operads for 
\(i\in\bn\), and \(\cp\in\Ob\dg^{\NN^n}\) is equipped with unital, 
associative and commuting actions
\begin{gather}
\rho=\rho_{(k_r)}: \circledast(\bn\to\mb1\to\mb1)(\cp;\cb)(\tau_{(k_r)}) 
=\Bigl(\bigotimes_{r=1}^m\cp(k_r)\Bigr)\tens\cB(m)
\to \cP\biggl(\sum_{r=1}^mk_r\biggr),
\notag
\\
\tau_\rho =\tau_{(k_r)} =
\begin{diagram}[inline,h=0.6em,w=3em,nobalance]
\mb k_1^n+\dots+\mb k_m^n & \\
&&\rdTTo(2,3) \\
\cdots & \\
&&\rdTTo(2,1) &\mb m &\rTTo &\mb1 &, \\
\mb k_1^2+\dots+\mb k_m^2 &&\ruTTo(2,1) \ruTTo(2,3) \\
\\
\mb k_1^1+\dots+\mb k_m^1 &
\end{diagram}
\label{eq-tau-kkk-m-1}
\\
\begin{split}
\lambda =\lambda_{k,(j_p^i)} &:
\circledast(\bn\rto1\bn\to\mb1)((\ca_i)_{i\in\bn};\cp)(\tau_{k,(j_p^i)})
\\
&=\biggl[\bigotimes_{i=1}^n \bigotimes_{p=1}^{k^i} \cA_i(j_p^i)\biggr]
\tens\cP\bigl((k^i)_{i=1}^n\bigr) \to
\cP\biggl(\Bigl(\sum_{p=1}^{k^i}j_p^i\Bigr)_{i=1}^n\biggr), 
\end{split}
\notag
\\
\tau_\lambda =\tau_{k,(j_p^i)} =
\begin{diagram}[inline,h=0.8em,w=3em,nobalance]
\ttt\sum_{p=1}^{k^n}\mb j_p^n &\rTTo & \mb k^n \\
&&&\rdTTo(2,3) \\
\ \cdots\ &\rTTo &\cdots \\
&&&\rdTTo(2,1) &\mb1 &. \\
\ttt\sum_{p=1}^{k^2}\mb j_p^2 &\rTTo &\mb k^2&\ruTTo(2,1) \ruTTo(2,3) \\
\\
\ttt\sum_{p=1}^{k^1}\mb j_p^1 &\rTTo & \mb k^1
\end{diagram}
\label{eq-tau-jjj-kkk-1}
\end{gather}
Commutativity and associativity of the above actions can be incorporated 
into a single requirement: existence. unitality and associativity of the 
actions
\begin{gather}
\begin{split}
\alpha &:
\circledast(\bn\rto1 \bn\to\mb1\to\mb1)((\ca_i)_{i\in\bn};\cp;\cb)(\tau)
\\
&=\Bigl(\bigotimes_{i=1}^n \bigotimes_{p=1}^{k_1^i+\dots+k_m^i} 
\cA_i(j_p^i)\Bigr)\tens\Bigl(\bigotimes_{r=1}^m\cp(k_r)\Bigr)\tens\cb(m)
\to 
\cp\biggl(\Bigl(\sum_{p=1}^{k_1^i+\dots+k_m^i}j_p^i\Bigr)_{i=1}^n\biggr),
\end{split}
\notag
\\
\tau_\alpha =
\begin{diagram}[inline,h=0.8em,w=3em,nobalance]
\ttt\sum_{p=1}^{k_1^n+\dots+k_m^n}\mb j_p^n &\rTTo
&\mb k_1^n+\dots+\mb k_m^n & \\
&&&&\rdTTo(2,3) \\
\ \ \cdots\ \ &\rTTo &\ \ \cdots\ \ & \\
&&&&\rdTTo(2,1) &\mb m &\rTTo &\mb1 \\
\ttt\sum_{p=1}^{k_1^2+\dots+k_m^2}\mb j_p^2 &\rTTo
&\mb k_1^2+\dots+\mb k_m^2 &&\ruTTo(2,1) \ruTTo(2,3) \\
\\
\ttt\sum_{p=1}^{k_1^1+\dots+k_m^1}\mb j_p^1 &\rTTo
&\mb k_1^1+\dots+\mb k_m^1 &
\end{diagram}
\label{eq-j-kk-m-1}
\end{gather}
for each $m\in\NN$, each family \(k_1,\dots,k_m\in\NN^n\) and each 
family of non-negative integers 
\(\bigl((j_p^i)_{p=1}^{k_1^i+\dots+k_m^i}\bigr)_{i=1}^n\).
Associativity of $\alpha$ is formulated via contraction of trees.

An example of an $n\wedge1$-operad module is given by
\[ (\END A_1,\dots,\END A_n;\HOM(A_1,\dots,A_n;B);\END B),
\]
where $B$, $A_i$, $i\in\bn$, are complexes of $\kk$\n-modules.
The object \(\HOM((A_i)_{i\in I};B)\) of \(\dg^{\NN^I}\), 
\(I=\mb m\in\co_\sk\), is specified by
\[ \HOM((A_i)_{i\in I};B)((n^i)_{i\in I})
=\uCom\bigl((\sS{^{n^i}}A_i)_{i\in I};B\bigr),
\]
where $\uCom$ is the symmetric $\dg$\n-multicategory associated with the 
symmetric closed monoidal category $\dg$.
The actions $\rho$, $\lambda$ are particular cases of the composition 
for $\HOM$'s, which is defined as the multiplication \(\mu_{\uCom}\) in 
the symmetric $\dg$\n-multicategory $\uCom$.
Recall the latter.
For a labelled symmetric tree
\[ T =\bigl(J \rTTo^\phi P \rTTo^\con \mb1; 
(X_j)_{j\in J}, (Y_p)_{p\in P}, Z \mid X_j, Y_p, Z \in \Ob\dg \bigr),
\]
$J,P\in\co_\sk$, \(\phi\in\Set\), the composition map
\[ \mu_{\uCom}^T: 
\Bigl(\tens^{p\in P} \uCom\bigl((X_j)_{j\in\phi^{-1}p};Y_p\bigr)\Bigr) 
\tens \uCom\bigl((Y_p)_{p\in P};Z\bigr)
\to \uCom\bigl((X_j)_{j\in J};Z\bigr)
\]
takes into account the Koszul rule (see
\secref{sec-Notations-conventions}) and the symmetry in $\udg$.
The composition for $\HOM$'s given for a tree
\(t=\bigl(A\rTTo^\theta B\rTTo^\con \mb1\bigr)\), \(\theta\in\co_\sk\), 
labelled with \(A\to\Ob\dg\), \(a\mapsto X_a\), \(B\to\Ob\dg\),
\(b\mapsto Y_b\), \(Z\in\Ob\dg\) and a $t$\n-tree \(\tau:t\to\co_\sk\), 
$a\mapsto J_a$, $b\mapsto P_b$,
\(a\mapsto(\tau_a:J_a\to P_{\theta(a)})\) is
\begin{multline*}
\comp_\tau: \Bigl[\bigotimes^{b\in B} \bigotimes^{p\in P_b} 
\HOM\bigl((X^a)_{a\in\theta^{-1}b};Y^b\bigr)
\bigl((|\tau_a^{-1}(p)|)_{a\in\theta^{-1}b}\bigr) \Bigl]
\tens \HOM\bigl((Y^b)_{b\in B};Z\bigr)\bigl((|P_b|)_{b\in B}\bigr)
\\
\hskip\multlinegap
=\Bigl[\bigotimes^{b\in B} \bigotimes^{p\in P_b} 
\uCom\bigl((\sS{^{|\tau_a^{-1}(p)|}}X^a)_{a\in\theta^{-1}b};Y^b\bigr) 
\Bigl] \tens \uCom\bigl((\sS{^{|P_b|}}Y^b)_{b\in B};Z\bigr) \hfill
\\
\to \uCom\bigl((\sS{^{|J_a|}}X^a)_{a\in A};Z\bigr)
=\HOM\bigl((X^a)_{a\in A};Z\bigr)\bigl((|J_a|)_{a\in A}\bigr).
\end{multline*}
It is defined as the multiplication \(\mu_{\uCom}^{\tilde\tau}\), where 
the labelled symmetric tree
\[ \tilde\tau =\Bigl(\bigsqcup_{a\in A}J_a \rTTo^{\tilde\phi}
\bigsqcup_{b\in B}P_b \rTTo^\con \mb1\Bigr)
\]
is formed by
	\(\tilde\phi|_{J_a}=\bigl(J_a\rTTo^{\tau_a} P_{\theta(a)}\rMono 
	\sqcup_{b\in B}P_b\bigr)\).
The label associated to any \(j\in J_a\) (resp. \(p\in P_b\), 
\(1\in\mb1\)) is $X^a$ (resp. $Y^b$, $Z$).
The trees \(t_\rho=(\bn\to\mb1\to\mb1)\), $\tau_\rho$ from
\eqref{eq-tau-kkk-m-1} determine \(\rho=\mu_{\uCom}^{\tilde\tau_\rho}\) 
and the trees \(t_\lambda=(\bn\rto1 \bn\to\mb1)\), $\tau_\lambda$ from
\eqref{eq-tau-jjj-kkk-1} determine 
\(\lambda=\mu_{\uCom}^{\tilde\tau_\lambda}\).

Given an operad $\co$ and an $n\wedge1$-operad $\co$\n-module $\cF_n$ 
for each $n\ge0$ we define a \emph{morphism of $\co$\n-algebras with $n$ 
arguments} \(X_1,\dots,X_n\to Y\) as a morphism of $\nOp n_1$
\begin{equation*}
(\co,\dots,\co;\cF_n;\co) 
\to (\END X_1,\dots,\END X_n;\HOM(X_1,\dots,X_n;Y);\END Y).
\end{equation*}

\subsection{Resolution \texorpdfstring{$F_n\to\FAs_n$}{Fn->FAsn}.}
Consider a $n\wedge1$-operad $\As$\n-module $\FAs_n$ having 
\(\FAs_n(j^1,\dots,j^n)=\kk\) for all non--vanishing 
\((j^1,\dots,j^n)\in\NN^n\), while \(\FAs_n(0,\dots,0)=0\).
The actions for $\FAs_n$ are given by multiplication in $\kk$.
A morphism of $n\wedge1$-operad modules
\begin{gather*}
(\As,\dots, \As;\FAs_n;\As)
\to (\END A_1,\dots,\END A_n;\HOM(A_1,\dots,A_n;B);\END B),
\\
\HOM(A_1,\dots,A_n;B)(j^1,\dots,j^n) 
=\und\dg(A_1^{\tens j^1}\tdt A_n^{\tens j^n},B),
\end{gather*}
amounts to a family of morphisms \(f_i:A_i\to B\) of associative 
differential graded $\kk$\n-algebras without units, \(i\in\bn\), such 
that the following diagrams commute for all \(1\le i<j\le n\):
\begin{diagram}
A_i\tens A_j &\rTTo^c_\sim &A_j\tens A_i &\rTTo^{f_j\tens f_i} &B\tens B
\\
\dTTo<{f_i\tens f_j} &&&&\dTTo>{m_B}
\\
B\tens B &&\rTTo^{m_B} &&B
\end{diagram}
see \exaref{exa-Com-FAsn}.

\begin{theorem}[Propositions \ref{pro-Fn-FAsn-Fn},
\ref{pro-Fn-FAsn-rmFn} and \thmref{thm-Fn-FAsn-cofibrant-replacement}]
The $n\wedge1$-operad module $(\As,\FAs_n)$ admits a cofibrant 
replacement \((\mainf,\rmF_n)\), where the $n\wedge1$-operad
\ainfm-module
	$\rmF_n
	=\boxdot_{\ge0}(\sS{^n}\mainf;\kk\{\sff_j\mid j\in\NN^n-0\};\mainf)$ 
is freely generated as a graded module by elements
$\sff_{j^1,\dots,j^n}\in\rmF_n(j^1,\dots,j^n)$,
\((j^1,\dots,j^n)\in\NN^n-0\), of degree $1-j^1-\dots-j^n=1-|j|$.
When these generators are taken into \(\HOM(A_1,\dots,A_n;B)\) they 
become linear maps \(\boxt^{i\in\bn}T^{\ell^i}A_i\to B\).
The definition of the differential is
\begin{multline*}
\sff_\ell\partial =\sum_{q=1}^n \sum_{r+x+t=\ell^q}^{x>1}
(-1)^{(1-x)(\ell^1+\dots+\ell^{q-1}+r)+1-|\ell|}
\lambda^q_{(\sS{^r}1,x,\sS{^t}1)}(\sS{^r}1,m_x,\sS{^t}1;
\sff_{\ell-(x-1)e_q})
\\
+\sum_{\substack{j_1,\dots,j_k\in\NN^n-0\\j_1+\dots+j_k=\ell}}^{k>1}
(-1)^{k+\sum_{1\le b<a\le k}^{1\le c<d\le n}j_a^cj_b^d
+\sum_{p=1}^k(p-1)(|j_p|-1)}
\rho_{(j_p^i)}((\sff_{j_p})_{p=1}^k;m_k).
\end{multline*}
The first arguments of $\lambda$ are all $1=\id$ except $m_x$ on the 
only place $p=r+1$. 
Moreover, \((\mainf,\rmF_n)\to(\As,\FAs_n)\) is a homotopy isomorphism.
\end{theorem}

The differential interpreted on \ainf-algebras $A_1$, \dots, $A_n$, $B$ 
means
\begin{multline*}
\sff_\ell\partial =\sum_{q=1}^n \sum_{r+x+t=\ell^q}^{x>1}
(-1)^{(1-x)(\ell^1+\dots+\ell^{q-1}+r)+1-|\ell|}
\bigl[ \tens^{i\in\bn}T^{\ell^i}A_i
\rTTo^{1^{\tens(q-1)}\tens(1^{\tens r}\tens m_x
\tens1^{\tens t})\tens1^{\tens(n-q)}}
\\
\hfill T^{\ell^1}A_1\tdt T^{\ell^{q-1}}A_{q-1}\tens 
T^{r+1+t}A_q\tens T^{\ell^{q+1}}A_{q+1}\tdt T^{\ell^n}A_n
\rTTo^{\sff_{\ell-(x-1)e_q}} B \bigr] \quad
\\
\hskip\multlinegap
+\sum_{\substack{j_1,\dots,j_k\in\NN^n-0\\j_1+\dots+j_k=\ell}}^{k>1}
(-1)^{k+\sum_{1\le b<a\le k}^{1\le c<d\le n}j_a^cj_b^d
+\sum_{p=1}^k(p-1)(|j_p|-1)}
\bigl[ \tens^{i\in\bn}T^{\ell^i}A_i \rTTo^\sim
\tens^{i\in\bn}\tens^{p\in\mb k}T^{j_p^i}A_i \rTTo^\sim \hfill
\\
\tens^{p\in\mb k}\tens^{i\in\bn}T^{j_p^i}A_i
\rTTo^{\tens^{p\in\mb k}\sff_{j_p}} \tens^{p\in\mb k}B \rTTo^{m_k} 
B\bigr].
\end{multline*}

\subsection{Morphisms of
	\texorpdfstring{$\mathrm{A}_\infty$}{A8}-algebras with 
	\texorpdfstring{$n$}n arguments.}
$\rmF_n$\n-algebra maps consist of \ainfm-algebras $A_1$, \dots, $A_n$, 
$B$, and an \ainfm-morphism \((\sff_j)_{j\in\NN^n-0}\). 
The latter means that the equation holds for all $\ell\in\NN^n-0$:
\begin{equation*}
\sff_\ell m_1 +(-1)^{|\ell|} \biggl[ \sum_{q=1}^n \sum_{r+1+t=\ell^q}
1^{\tens(q-1)}\tens(1^{\tens r}\tens m_1\tens1^{\tens t})
\tens1^{\tens(n-q)} \biggr]\sff_\ell
=\sff_\ell\partial.
\end{equation*}
This is equivalent to equation \eqref{eq-sff-b-commute}.

Composition of \ainfm-morphisms can be obtained from comultiplication in 
the system $\rmF_n$.
We view it as a coalgebra and $\HOM$ as an algebra.
Then homomorphisms between them form an algebra as well.
This way multiquiver $\alinf$ whose objects are \ainfm-algebras
\((B,\alpha_B:\mainf\to\END B)\) with the set of morphisms
\begin{multline*}
\alinf((A_i,\alpha_{A_i})_{i\in I};(B,\alpha_B)) 
\\
=\{ ((\alpha_{A_i})_{i\in I};\phi;\alpha_B):(\sS{^I}\mainf;\rmF_I;mainf) 
\to ((\END A_i)_{i\in I};\HOM((A_i)_{i\in I};B);\END B) \}
\end{multline*}
becomes a multicategory.

For any tree $t$ and any collection of \ainfm-algebras
\(\alpha_h^b:\mainf\to\END A_h^b\), \((h,b)\in\iV(t)\), assume given 
\(t_h^{-1}b\wedge1\)\n-operad module morphisms for \((h,b)\in\IV(t)\):
\begin{multline*}
\bigl((\alpha_{h-1}^a)_{a\in t_h^{-1}b};g_h^b;\alpha_h^b\bigr): 
\bigl(\sS{^{t_h^{-1}b}}\mainf;\rmF_{t_h^{-1}b};\mainf\bigr) 
\\
\to \bigl((\END A_{h-1}^a)_{a\in t_h^{-1}b};
\HOM((A_{h-1}^a)_{a\in t_h^{-1}b};A_h^b); \END A_h^b\bigr).
\end{multline*}
Their composition can be defined as
\(\bigl((\alpha_0^a)_{a\in t(0)};\comp(t)(g_h^b);\alpha_l^1\bigr)\), 
where
\begin{multline*}
\comp(t)(g_h^b) =\bigl[\rmF_{t(0)} \rTTo^{\Delta(t)}
\circledast(t)(\rmF_{t_h^{-1}b})_{(h,b)\in\IV(t)}
\rTTo^{\circledast(t)(g_h^b)}
\\
\circledast(t)
\bigl(\HOM((A_{h-1}^a)_{a\in t_h^{-1}b};A_h^b)\bigr)_{(h,b)\in\IV(t)}
\rTTo^{\comp(t)} \HOM((A_0^a)_{a\in t(0)};A_l^1) \bigr].
\end{multline*}
Explicit form of the comultiplication is
\begin{align*}
\Delta(t)(\sff_j)
=\sum^{t-\text{tree }\tau}_{\forall a\in t(0)\,|\tau(0,a)|=j^a}
(-1)^{c(\tilde\tau)} \tens^{h\in I} \tens^{b\in t(h)} 
\tens^{p\in\tau(h,b)} \sff_{|\tau_{(h-1,a)}^{-1}(p)|_{a\in t_h^{-1}b}}.
\end{align*}
The sign $c(\tilde\tau)$ is computed via
recipe~\eqref{eq-c(tilde-tau)SS-SS} through the Koszul rule.

\section{Lax \texorpdfstring{$\Cat$}{Cat}-span multicategories}
 \label{sec-Lax-Cat-multicategories}
In this section we describe the categorical background to the main 
subject of morphisms with several entries.

\subsection{\texorpdfstring{$\Cat$}{Cat}-spans.}
A \emph{$\Cat$-span} $\mcC$ is a pair of functors with a common source, 
which we denote \(\soup\mcC\lTTo^\src \mcC\rTTo^\tgt \tar\mcC\). 
We say that \(\soup\mcC\) is the source category and \(\tar\mcC\) is the 
target category.
For example, anafunctors (see Makkai \cite{MakkaiM:avoacc}) are 
particular instances of $\Cat$-spans.
% Avoiding the axiom of choice in general category theory, Michael Makkai \url{http://www.math.mcgill.ca/makkai/anafun/}, \url{http://ncatlab.org/nlab/show/anafunctor}.
A \emph{morphism of $\Cat$-spans} \(F:\mcA\to\mcB\) is a triple of 
functors
\begin{equation*}
\soup F: \soup\mcA \to \soup\mcB, \qquad F: \mcA \to \mcB, 
\qquad \tar F: \tar\mcA \to \tar\mcB, 
\end{equation*}
strictly commuting with the source functor $\src$ and the target functor 
$\tgt$:
\begin{diagram}[LaTeXeqno]
\soup\mcA &\lTTo^\src &\mcA &\rTTo^\tgt &\tar\mcA
\\
\dTTo<{\soup F} &= &\dTTo<F &= &\dTTo>{\tar F}
\\
\soup\mcB &\lTTo^\src &\mcB &\rTTo^\tgt &\tar\mcB
\label{dia-sF-F-tF}
\end{diagram}

So described category $\Catspan$ has arbitrary products.
In fact,
\[ \prod_{i\in I}(\soup\mcC_i \lTTo^\src \mcC_i \rTTo^\tgt \tar\mcC_i)
=\bigl( \prod_{i\in I}\soup\mcC_i \lTTo^{\hspace*{0.2em}\prod\src} 
\prod_{i\in I}\mcC_i \rTTo^{\prod\tgt} \prod_{i\in I}\tar\mcC_i \bigr).
\]

First proof of the following statement was given by Sergiy Slobodianiuk 
(unpublished).
The proof presented here is devised by the author.

\begin{proposition}
The category $\Catspan$ is Cartesian closed.
\end{proposition}

\begin{proof}
The inner homomorphisms object \(\und\Catspan(\mcA,\mcB)\) for 
\((\Catspan,\times)\) is given by the pair of functors 
 \(\und\Cat(\soup\mcA,\soup\mcB)\lTTo_{\pr_1}^\src \ucatspan(\mcA,\mcB)
 \rTTo_{\pr_3}^\tgt \und\Cat(\tar\mcA,\tar\mcB)\),
where objects of the category \(\ucatspan(\mcA,\mcB)\) are triples of 
functors \((\soup F,F,\tar F)\) such that diagram~\eqref{dia-sF-F-tF} 
commutes.
Morphisms \((\soup F,F,\tar F)\to(\soup G,G,\tar G)\) of the category 
$\ucatspan(\mcA,\mcB)$ are triples of natural transformations
\[ \bigl( (\soup\phi:\soup F\to\soup G:\soup\mcA\to\soup\mcB),
(\phi:F\to G:\mcA\to\mcB),
(\tar\phi:\tar F\to\tar G:\tar\mcA\to\tar\mcB)\bigr)
\]
such that
\begin{equation}
\begin{split}
\begin{diagram}[inline]
\mcA &\rTTo^\src &\soup\mcA
\\
\dTTo<F \overset\phi\Rightarrow \dTTo>G &= &\dTTo>{\soup G}
\\
\mcB &\rTTo^\src &\soup\mcB
\end{diagram}
\quad &=\quad
\begin{diagram}[inline]
\mcA &\rTTo^\src &\soup\mcA
\\
\dTTo<F &=&\dTTo<{\soup F} \overset{\soup\phi}\Rightarrow\dTTo>{\soup G}
\\
\mcB &\rTTo^\src &\soup\mcB
\end{diagram}
\quad,
\\
\begin{diagram}[inline]
\mcA &\rTTo^\tgt &\tar\mcA
\\
\dTTo<F \overset\phi\Rightarrow \dTTo>G &= &\dTTo>{\tar G}
\\
\mcB &\rTTo^\tgt &\tar\mcB
\end{diagram}
\quad &=\quad
\begin{diagram}[inline]
\mcA &\rTTo^\tgt &\tar\mcA
\\
\dTTo<F &= &\dTTo<{\tar F} \overset{\tar\phi}\Rightarrow \dTTo>{\tar G}
\\
\mcB &\rTTo^\tgt &\tar\mcB
\end{diagram}
\quad.
\end{split}
\label{eq-4-squares-4-phi}
\end{equation}
The source and target functors are projections 
\(\src=\pr_1:(\soup F,F,\tar F)\mapsto\soup F\), 
\((\soup\phi,\phi,\tar\phi)\mapsto\soup\phi\), and
\(\tgt=\pr_3:(\soup F,F,\tar F)\mapsto\tar F\), 
\((\soup\phi,\phi,\tar\phi)\mapsto\tar\phi\).
The evaluation morphism \(\ev:\mcA\times\und\Catspan(\mcA,\mcB)\to\mcB\) 
consists of three functors
\[ \bigl( \soup\mcA\times\und\Cat(\soup\mcA,\soup\mcB) \rTTo^\ev 
\soup\mcB, \mcA\times\mcD\rTTo^{1\times\pr_2} 
\mcA\times\und\Cat(\mcA,\mcB)\rTTo^\ev \mcB,
\tar\mcA\times\und\Cat(\tar\mcA,\tar\mcB) \rTTo^\ev \tar\mcB \bigr),
\]
where \(\mcD=\ucatspan(\mcA,\mcB)\).

Notice that there is a 2\n-category $\ucatspan$, whose objects are 
$\Cat$-spans and whose categories of morphisms are
\(\ucatspan(\mcA,\mcB)\).
All morphisms in this 2\n-category, including left and right whiskering, 
are compositions in $\Cat$, performed simultaneously in three places: 
the source, the main body and the target.
Such compositions preserve commutation relations \eqref{dia-sF-F-tF}, \eqref{eq-4-squares-4-phi}, hence, they give well-defined compositions in $\ucatspan$.
Standard equations involving them hold in $\ucatspan$, since they hold in the 2\n-category $\und{\Cat}$.

The underlying category of objects and 1\n-morphisms in $\ucatspan$ is 
precisely $\Catspan$.
Furthermore, finite products in $\Catspan$ extend to a (symmetric) 
monoidal 2\n-category structure of $\ucatspan$, just as they do for 
$\Cat$ and $\und{\Cat}$.

We have to prove that the mapping
\begin{multline}
\varphi =\bigl[ \Catspan(\mcC,\und\Catspan(\mcA,\mcB)) 
\rTTo^{(1_\mcA\times-)\times\ev}
\\
\Catspan(\mcA\times\mcC,\mcA\times\und\Catspan(\mcA,\mcB))
	\times\Catspan(\mcA\times\und\Catspan(\mcA,\mcB),\mcB)
\\
\rTTo^\comp \Catspan(\mcA\times\mcC,\mcB) \bigr]
\label{eq-varphi-ev}
\end{multline}
is a bijection.
Let us construct a map inverse to $\varphi$.

The unit object of $(\Catspan,\times)$ and $(\ucatspan,\times)$ is the 
terminal $\Cat$-span \(\1^3=(\1\leftarrow\1\to\1)\), where $\1$ is the 
terminal category.
A morphism \(\1^3\to\mcC\) of $\Catspan$ has the form 
\((*\mapsto\src X,*\mapsto X,*\mapsto\tgt X)\) for some object 
\(X\in\Ob\mcC\).
Denote it by \(\dot{X}:\1^3\to\mcC\).
A 2\n-morphism \(\dot{X}\to\dot{Y}:\1^3\to\mcC\) of $\ucatspan$ 
identifies with a morphism \(f:X\to Y\in\mcC\).
Denote it by \(\dot{f}:\dot{X}\to\dot{Y}:\1^3\to\mcC\), 
\((*\mapsto\src f,*\mapsto f,*\mapsto\tgt f)\).

Given a $\Cat$-span morphism \(F:\mcA\times\mcC\to\mcB\), let us 
construct a $\Cat$-span morphism
\(\Psi(F):\mcC\to\und\Catspan(\mcA,\mcB)\).
First of all, for any object $X$ of the category $\mcC$ there is a 
morphism of $\Cat$-spans
\[ \mcA \rTTo^\sim \mcA\times\1^3 \rTTo^{1\times\dot X} \mcA\times\mcC 
\rTTo^F \mcB,
\]
given by three functors
\[ (\soup F(-,\src X),F(-,X),\tar F(-,\tgt X)).
\]
This is an object \(\psi(F)(X)\) of \(\ucatspan(\mcA,\mcB)\).
Secondly, for any morphism \(f:X\to Y\) of the category $\mcC$ there is 
a 2\n-morphism of $\ucatspan$
\[ \biggl( \mcA \rTTo^\sim \mcA\times\1^3 
\pile{\rTTo^{1\times\dot X} \\ \sss 1\times\dot f 
\Downarrow \\ \rTTo_{1\times\dot Y} }
\mcA\times\mcC \rTTo^F \mcB \biggr)
=\biggl( \psi(F)(f): \psi(F)(X) \to \psi(F)(Y) \biggr)
\]
given by three natural transformations
\[ (\soup F(-,\src f),F(-,f),\tar F(-,\tgt f)).
\]
Clearly, this determines a functor
\[ \psi(F): \mcC \to \ucatspan(\mcA,\mcB).
\]
On the other hand, given functors 
\(\soup F:\soup\mcA\times\soup\mcC\to\soup\mcB\) and 
\(\tar F:\tar\mcA\times\tar\mcC\to\tar\mcB\) induce by closedness of 
$(\Cat,\times)$ the functors
\[
\begin{aligned}
\overline{\soup F}: \soup\mcC &\longrightarrow
\und\Cat(\soup\mcA,\soup\mcB)
\\
U &\longmapsto \soup F(-,U)
\end{aligned}
\quad,\qquad
\begin{aligned}
\overline{\tar F}: \tar\mcC &\longrightarrow \und\Cat(\tar\mcA,\tar\mcB)
\\
V &\longmapsto \tar F(-,V)
\end{aligned}
\quad.
\]
Explicit expressions show that both squares of the diagram
\begin{diagram}
\soup\mcC &\lTTo^\src &\mcC &\rTTo^\tgt &\tar\mcC
\\
\dTTo<{\overline{\soup F}} &&\dTTo>{\psi(F)} &&\dTTo>{\overline{\tar F}}
\\
\und\Cat(\soup\mcA,\soup\mcB) &\lTTo^\src_{\pr_1} &\ucatspan(\mcA,\mcB) 
&\rTTo^\tgt_{\pr_3} &\und\Cat(\tar\mcA,\tar\mcB)
\end{diagram}
commute.
Therefore, there is a $\Cat$-span morphism
\[ \Psi(F) =(\overline{\soup F},\psi(F),\overline{\tar F}): 
\mcC \to \und\Catspan(\mcA,\mcB).
\]
This gives a map
\[ \Psi: \Catspan(\mcA\times\mcC,\mcB) \to
\Catspan(\mcC,\und\Catspan(\mcA,\mcB)).
\]
Without much efforts one checks that it is inverse to mapping $\varphi$ 
given by \eqref{eq-varphi-ev}.
\end{proof}

Thus there is a category $\und\Catspan$ enriched in $\Catspan$.
The underlying functor \(\pr_2:\Catspan\to\Cat\), 
\((\soup\mcC\lTTo^\src \mcC\rTTo^\tgt \tar\mcC)\mapsto\mcC\) turns 
$\Catspan$ into a $\Cat$-\hspace{0pt}category.
This 2\n-category structure of $\Catspan$ coincides with that of 
$\ucatspan$ discussed in the above proof.
In particular, a 2\n-morphism of $\Catspan$ is a triple 
\((\soup\phi,\phi,\tar\phi)\) which satisfies
\eqref{eq-4-squares-4-phi}.

\subsection{\texorpdfstring{$\Cat$}{Cat}-span multiquivers.}
Consider the category $\StrictMonCat$ of strict monoidal categories and 
strict monoidal functors.
The underlying functor $U$ and the `free strict monoidal category' 
functor $F$ form an adjunction
\(F:\Cat\rightleftarrows\StrictMonCat:U\).
The monad \(\text-^*=U\circ F:\Cat\to\Cat\) associated with this 
adjunction takes a category $\cc$ to \(\cc^*=\sqcup_{n\in\NN}\,\cc^n\).
This `free strict monoidal category' monad is cartesian, see 
Definition~4.1.1 and Example~4.1.15 of \cite{math.CT/0305049}.
Thus, the monad \(\text-^*\) is suitable for introducing a kind of 
multicategories.

A \emph{$\Cat$-span multiquiver} is a $\Cat$-span $\mcC$ together with 
$\sou\mcC$ such that
\(\soup\mcC=(\sou\mcC)^*=\sqcup_{n\in\NN}\,\sou\mcC^n\), thus,
\[ \soup\mcC((X_i)_{i\in\bn},(Y_j)_{j\in\mb m}) =
\begin{cases}
\prod_{i\in\bn}\sou\mcC(X_i,Y_i), &\text{ if } n=m,
\\
\emptyset, &\text{ otherwise.}
\end{cases}
\]
Since $\soup\mcC$ is a disjoint union of categories \(\sou\mcC^n\), the 
mere existence of a functor \(\src:\mcC\to\soup\mcC\) implies that 
\(\mcC=\sqcup_{n\in\NN}\,\mcC(n)\) is a disjoint union as well, where 
\(\mcC(n)=\src^{-1}(\sou\mcC^n)\).
Thus we can view a $\Cat$-span multiquiver as a sequence of pairs of 
functors
\[
\bigl(\sou\mcC^n \lTTo^\src \mcC(n)\rTTo^\tgt \tar\mcC\mid n\ge0 \bigr).
\]
A \emph{morphism of $\Cat$-span multiquivers} is a morphism of
$\Cat$-spans \(F:\mcC\to\mcD\) together with
\(\sou F:\sou\mcC\to\sou\mcD\) such that \(\soup F=(\sou F)^*\).
Necessarily \(F=\bigl(F(n):\mcC(n)\to\mcD(n)\mid n\ge0\bigr)\).
This defines the category \(\smQuiver\) of $\Cat$-span multiquivers.

Moreover, \(\smQuiver\) is a 2\n-category and has a strict 2\n-functor 
\(\smQuiver\to\Catspan\), injective on objects, 1\n-morphisms and 
2\n-morphisms.
Namely, a \emph{2-morphism of $\Cat$-span multiquiver morphisms} 
\(\phi:F\to G:\mcC\to\mcD\) is a triple of natural transformations 
\(\sou\phi:\sou F\to\sou G:\sou\mcC\to\sou\mcD\), $\phi$ and 
\(\tar\phi:\tar F\to\tar G:\tar\mcC\to\tar\mcD\) such that the triple 
\(((\sou\phi)^*,\phi,\tar\phi)\) satisfies \eqref{eq-4-squares-4-phi}:
\begin{equation}
\begin{split}
\begin{diagram}[inline]
\mcC &\rTTo^\tgt &\tar\mcC
\\
\dTTo<F \overset\phi\Rightarrow \dTTo>G &= &\dTTo>{\tar G}
\\
\mcD &\rTTo^\tgt &\tar\mcD
\end{diagram}
\quad &=\quad
\begin{diagram}[inline]
\mcC &\rTTo^\tgt &\tar\mcC
\\
\dTTo<F &= &\dTTo<{\tar F} \overset{\tar\phi}\Rightarrow \dTTo>{\tar G}
\\
\mcD &\rTTo^\tgt &\tar\mcD
\end{diagram}
\quad,
\\
\begin{diagram}[inline]
\mcC &\rTTo^\src &(\sou\mcC)^*
\\
\dTTo<F \overset\phi\Rightarrow \dTTo>G &= &\dTTo>{(\sou G)^*}
\\
\mcD &\rTTo^\src &(\sou\mcD)^*
\end{diagram}
\quad &=
\begin{diagram}[inline]
\mcC &\rTTo^\src &(\sou\mcC)^*
\\
\dTTo<F &= &\dTTo<{(\sou F)^*} 
\overset{(\sou\phi)^*}\Longrightarrow \dTTo>{(\sou G)^*}
\\
\mcD &\rTTo^\src &(\sou\mcD)^*
\end{diagram}
\quad.
\end{split}
\label{eq-4dia-phi-tphi}
\end{equation}

\subsection{Spans and related structures.}
Consider the category strict-2-$\Cat=\Cat\text-\Cat$ of strict
2\n-categories and strict 2\n-functors.
It is complete, so there exists a bicategory $\Span\Cat\text-\Cat$ of 
spans in it.
We shall need the category \(\CS=\Span\Cat\text-\Cat(\Cat,\Cat)\), whose 
objects are pairs of strict 2\n-functors
\(\cs=\bigl(\Cat\lTTo^\src \cs\rTTo^\tgt \Cat\bigr)\), and morphisms
$F:\cs\to\cp$ are commutative diagrams of strict 2\n-functors
\begin{diagram}
\Cat &\lTTo^\src &\cs &\rTTo^\tgt &\Cat
\\
\dEq &= &\dTTo>F &= &\dEq
\\
\Cat &\lTTo^\src &\cp &\rTTo^\tgt &\Cat
\end{diagram}
By the general theory of bicategories the category $\CS$ is Monoidal, 
the term is suggested in \cite[Definition~2.5]{BesLyuMan-book}.
Its monoidal product is
\begin{equation}
\boxbox^{h\in\bn}\cs_h =\lim( \Cat \lTTo^\src \cs_1 \rTTo^\tgt \Cat 
\lTTo^\src \cs_2 \rTTo^\tgt \dots\Cat \lTTo^\src \cs_n \rTTo^\tgt \Cat).
\label{eq-power-SMQCat-lim}
\end{equation}
In this diagram the limit is taken in strict-2-$\Cat$.
The canonical strict 2\n-functors to the first and the last 
$\Cat$-\hspace{0pt}vertices are the horizontal source and target of 
\(\boxbox^{h\in\bn}\cs_h\).
Of course,
\(\boxbox^\emptyset=\bigl(\Cat\lTTo^\Id \Cat\rTTo^\Id \Cat\bigr)\).
For any non-decreasing map \(\phi:\bn\to\mb k\) there are isomorphisms 
in strict-2-$\Cat$
	\(\boxbox^{h\in\bn}\cs_h\rTTo^\sim
	\boxbox^{m\in\mb k}\boxbox^{h\in\phi^{-1}m}\cs_h\)
due to presenting repeated limits as a single limit.
The isomorphism \(\boxbox^1\cs=\lim(\Cat\leftarrow\cs\to\Cat)=\cs\) is 
chosen to be the identity.

Now we shall equip $\smQuiver$ with some structures, which will turn it 
later into a weak triple $\Cat$-category.
The latter notion will be defined below as close to the strict one as 
possible.
Such a version is just what we need for the purposes of this article.
Thus, \(\smQuiver\) is equipped with the following data:

1. Horizontal source and target, a pair of strict 2\n-functors
\begin{diagram}
\Cat &\lTTo^{\pr_1} &\smQuiver &\rTTo^{\pr_3} &\Cat
\\
\sou\mcB &\lMapsTo &((\sou\mcB)^* \leftarrow \mcB \to \tar\mcB) 
&\rMapsTo &\tar\mcB
%\label{dia-Cat-SMQCat-Cat}
\end{diagram}
Thus, \(\sou\smQuiver=\tar\smQuiver=\Cat\) and \(\smQuiver\) is an 
object of $\CS$.
Explicitly we may define objects, 1\n-morphisms and 2\n-morphisms of 
\(\boxbox^{\bn}\smQuiver\) as \(\Ob\boxbox^\emptyset=\Cat\),
\begin{gather*}
\Ob\boxbox^{\bn}\smQuiver =\{(\cc_{h-1}^* \lTTo^\src \mcC_h
\rTTo^\tgt \cc_h)_{h\in\bn} \in \Ob\smQuiver^n\}, \qquad n>0,
%\label{eq-Ob-Mor-smQuiver}
\\
\Mor\boxbox^\emptyset =\{ F:\cC\to\cD \text{ -- functor } \},
\\
\begin{split}
\Mor &\boxbox^{\bn}\smQuiver 
\\
&=\Biggl\{ \bigl( (F_h:\mcC_h\to\mcD_h)_{h\in\bn},
(\tar F_h:\cC_h\to\cD_h)_{h\in[n]} \bigr) \;\Bigg|\; \forall\, h\in\bn
\begin{diagram}[inline,w=2em]
\cc_{h-1}^* &\lTTo^\src &\mcC_h &\rTTo^\tgt &\cc_h
\\
\dTTo~{\tar F_{h-1}^*} &= &\dTTo>{F_h} &= &\dTTo<{\tar F_h}
\\
\cd_{h-1}^* &\lTTo^\src &\mcD_h &\rTTo^\tgt &\cd_h
\end{diagram}
\Biggr\},
\end{split}
\\
2\text-\Mor \boxbox^\emptyset
=\{\phi: F\to G: \cC\to\cD \text{ -- natural transformation } \},
\\
\begin{split}
2\text-\Mor \boxbox^{\bn}\smQuiver 
=\bigl\{ \bigl( (\phi_h:F_h\to G_h:\mcC_h\to\mcD_h)_{h\in\bn}, 
(\tar\phi_h:\tar F_h\to\tar G_h:\cC_h\to\cD_h)_{h\in[n]} \bigr) \mid 
\\
\forall\, h\in\bn \; (\tar\phi_{h-1},\phi_h,\tar\phi_h) \text{ satisfy } 
\eqref{eq-4dia-phi-tphi} \bigr\}.
\end{split}
\end{gather*}

2. A morphism $\boxdot^I:\boxbox^I\smQuiver\to\smQuiver$ of $\CS$ for 
every set $I\in\Ob\co_\sk$.
For \(I=\emptyset\) the morphism
\(\boxdot^\emptyset:\boxbox^\emptyset=\Cat\to\smQuiver\), takes a 
category $\cc$ to the $\Cat$-span multiquiver 
\(\mcC=\boxdot^\emptyset\cc\) with \(\sou\mcC=\tar\mcC=\cc\), 
\(\mcC(n)=\emptyset\) for $n\ne1$, and
\(\cC\lTTo^{\src=\id} \mcC(1)=\cc\rTTo^{\tgt=\id} \cC\).
For
	\((\mcC_h)_{h\in I}
	=(\cc_{h-1}^*\lTTo^\src \mcC_h\rTTo^\tgt \cc_h)_{h\in I}\)
with non-empty $I$ there is a multiquiver defined as
\begin{equation*}
\boxdot^{h\in I}\mcC_h =\bigsqcup_{\text{tree }t:[I]\to\co_\sk} 
\lim(D_t:t^* \to \Cat),
%\label{eq-boxICh-Ulim}
\end{equation*}
where the diagram shape (category) $t^*$ has
\(\Ob t^*=\IV(t)\sqcup\iV(t)=\IV(\bar t)\sqcup\edges(\bar t)\).
The tree $\bar t$ is obtained from the tree $t$ by adding an extra edge 
starting from the root.
The set $\IV(\bar t)$ of internal vertices of $\bar t$ in the usual 
sense (without the new added root) coincides with the set $\IV(t)$ of 
internal vertices of $t$.
The set $\edges(\bar t)$ of all edges of $\bar t$ is identified with the 
set $\iV(t)$ of all vertices of $t$.
Arrows of the diagram start in an internal vertex of $\bar t$ and end up 
in (the middle of) an adjacent edge.
Thus,
	\(\Mor t^*=\{v\to e\mid(v,e)\in\IV(\bar t)\times\edges(\bar t),
	\ e\text{ is adjacent to }v\}\)
is the set of flags of $\bar t$.
Examples are given below:
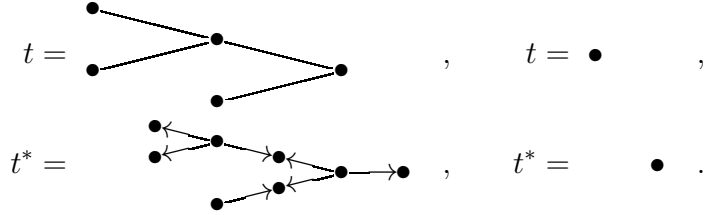
\begin{alignat*}2
t &=\; 
\begin{diagram}[inline,tight,abut,w=2em,h=0.5em]
\bullet \\
&\rdLine \\
&&\bullet \\
&\ruLine &&\rdLine \\
\bullet &&&&\bullet \\
&&&\ruLine \\
&&\bullet
\end{diagram}
\hspace*{2em} \quad, &\qquad t &=\; \bullet \hspace*{2em} \quad,
\\
t^* &=\; \hspace*{2em}
\begin{diagram}[inline,tight,abut,w=1em,h=0.25em]
\bullet \\
&\luTTo \\
&&\bullet \\
&\ldTTo &&\rdTTo \\
\bullet &&&&\bullet \\
&&&&&\luTTo \\
&&&&&&\bullet &\rTTo &\bullet \\
&&&&&\ldTTo \\
&&&&\bullet \\
&&&\ruTTo \\
&&\bullet
\end{diagram}
\quad, &\qquad t^* &=\; \hspace*{2em} \bullet \quad.
\end{alignat*}
The category assigned by $D_t$ to \(v=(h,b)\in\IV(\bar t)\) is 
\(\mcC_h(t_h^{-1}b)\), the category assigned to
\(e_{(h-1,a)}:(h-1,a)\to(h,b)\in\edges(\bar t)\), \(1\le h\le I+1\), is 
$\cc_{h-1}$.
The functor assigned by $D_t$ to the arrow
\((h,b)\to e_{(h,b)}\in\Mor t^*\) is 
\(\mcC_h(t_h^{-1}b)\rTTo^\tgt \cc_h\), \(h\in I\).
The functor assigned to the arrow 
\(e_{(h-1,a)}\leftarrow(h,b)\in\Mor t^*\) for \(a\in t_h^{-1}b\),
$h\in I$, is 
	\(\cc_{h-1}\lTTo^{\pr_a} \cc_{h-1}^{t_h^{-1}b}\lTTo^\src 
	\mcC_h(t_h^{-1}b)\).
Explicitly we may describe the disjoint union of limits as
\begin{multline}
\boxdot^{h\in I}\mcC_h =\bigsqcup_{\text{tree }t:[I]\to\co_\sk} 
\bigl\{(\ell,p) \mid 
\ell=(\ell_h^b)_{h\in[I]}^{b\in t(h)} \in\prod_{h\in[I]}\cc_h^{t(h)},\;
p=(p_h^b)_{h\in I}^{b\in t(h)} 
\in\prod_{h\in I} \prod_{b\in t(h)} \mcC_h(t_h^{-1}b),
\\
\forall h\in I\; \forall b\in t(h)\; \tgt p_h^b=\ell_h^b,\; 
\src p_h^b =(\ell_{h-1}^a)_{a\in t_h^{-1}b} \bigr\},
\label{eq-boxICh-U(lp)}
\end{multline}
where elements of a category mean its morphisms and, in particular, its 
objects.
The source and the target functors are given by
\(\src:\boxdot^{h\in I}\mcC_h\to\cc_0^*\),
\(\src(\ell,p)=(\ell_0^a)_{a\in t(0)}\) and
\(\tgt:\boxdot^{h\in I}\mcC_h\to\cc_{\max[I]}\), 
\(\tgt(\ell,p)=\ell_{\max[I]}^1\).

3. For any non-decreasing map \(\phi:\bn\to\mb k\) and the induced 
\(\psi=[\phi]:[k]\to[n]\) (see \eqref{eq-x<fy-fx<y}) there are
2\n-isomorphisms
\begin{align*}
\Lambda^\phi_{\smQuiver}:\boxdot^{l\in\bn}\mcC_l &\longrightarrow
\boxdot^{m\in\mb k}\boxdot^{l\in\phi^{-1}m}\mcC_l, 
%\label{eq-Lambda-phi-smQuiver}
\\
(t,\ell,(p_l^j)_{l\in\bn}^{j\in t(l)}) &\longmapsto
(t_\psi,\ell_\psi,(t^{|j}_{[\psi(m-1),\psi(m)]},\ell^{|j},
(p_l^i)_{l\in\phi^{-1}m}^{i\in t^{|j}_{[\psi(m-1),\psi(m)]}(l)} 
)_{m\in\mb k}^{j\in t(\psi(m))}). 
\end{align*}

4. The 2\n-isomorphism 
	\(\Rho:\boxdot^{\mb1}\to\Iso:\boxbox^{\mb1}\smQuiver\rTTo^\sim 
	\smQuiver\)
is the obvious one.

Fix a category $\cc$. 
Consider the $\Cat$\n-subcategory \(\sS{^\cc}\smQuiver\) of the 
$\Cat$\n-category \(\smQuiver\) whose objects are $\Cat$-span 
multiquivers \(\mcC\) with \(\sou\mcC=\tar\mcC=\cc\), whose
1\n-morphisms satisfy \(\sou F=\tar F=\id_\cc\) and 2\n-morphisms 
satisfy \(\sou\phi=\tar\phi=\id_{\Id}\).

\begin{definition}\label{def-lax-Cat-span-multicategory}
A \emph{lax $\Cat$-span multicategory} $\mcC$ is the collection of
\begin{enumerate}
\item a $\Cat$-span multiquiver 
\(\mcC=\bigl(\cc^*\lTTo^\src \mcC\rTTo^\tgt \cc\bigr)\);

\item a 1\n-morphism $\circledast^I:\boxdot^I\mcC\to\mcC$ in 
\(\sS{^\cc}\smQuiver\), for every set $I\in\Ob\co_\sk$.

For a map $f:I\to J$ in $\Mor\co_\sk$ introduce a 1\n-morphism
\[ \circledast^f_\mcC = \bigl(\boxdot^I\mcC \rTTo^{\Lambda^f}
\boxdot^{j\in J}\boxdot^{f^{-1}j}\mcC
\rTTo^{\boxdot^{j\in J}\circledast^{f^{-1}j}} \boxdot^J\mcC\bigr);
\]

\item a 2\n-morphism $\lambda^f$ in \(\sS{^\cc}\smQuiver\) for every map 
$f:I\to J$ in $\Mor\co_\sk$:
\[
\begin{diagram}[width=4.7em,inline]
\boxdot^I\mcC & \rTTo^{\circledast^f_\mcC} & \boxdot^J\mcC \\
&\rdTTo_{\circledast^I} \ruTwoar(1,1)_{\lambda^f} &\dTTo>{\circledast^J} 
\\
&& \mcC
\end{diagram}
\quad=\quad
\begin{diagram}[width=4.7em,inline]
\boxdot^{j\in J}\boxdot^{f^{-1}j}\mcC
& \rTTo^{\boxdot^{j\in J}\circledast^{f^{-1}j}} & \boxdot^J\mcC \\
\uTTo<{\Lambda^f} & \uTwoar<{\lambda^f} & \dTTo>{\circledast^J} \\
\boxdot^I\mcC &\rTTo^{\circledast^I} & \mcC
\end{diagram}
\]

\item a 2-morphism 
\(\rho:\circledast^{\mb1}\to\Rho:\boxdot^{\mb1}\mcC\to\mcC\) in 
\(\sS{^\cc}\smQuiver\),
\end{enumerate}
such that
\begin{enumerate}
\renewcommand{\labelenumi}{(\roman{enumi})}
\item for all sets \(I\in\Ob\co_\sk\)
\begin{align*}
\begin{diagram}[width=4.4em,h=3em,inline,nobalance]
\boxdot^{i\in I}\boxdot^{\{i\}}\mcC 
&\pile{\rTTo^{\boxdot^{i\in I}\Rho^{\{i\}}} 
\\ \Uparrow\sss\boxdot^{i\in I}\rho^{\{i\}} 
\\ \rTTo_{\boxdot^{i\in I}\circledast^{\{i\}}}}
&\boxdot^I\mcC 
\\
\uTTo<{\Lambda^{\id_I}} &\uTwoar>{\lambda^{\id_I}} 
&\dTTo>{\circledast^I} 
\\
\boxdot^I\mcC &\rTTo^{\circledast^I} &\mcC
\end{diagram}
&=
\begin{diagram}[width=2.2em,h=3em,inline,nobalance]
\boxdot^{i\in I}\boxdot^{\{i\}}\mcC
&&\rTTo^{\boxdot^{i\in I}\Rho^{\{i\}}} &&\boxdot^I\mcC 
\\
\uTTo<{\Lambda^{\id_I}} &= &&\ruTTo(4,2)<1 = &\dTTo>{\circledast^I} 
\\
\boxdot^I\mcC &&\rTTo^{\circledast^I} &&\mcC
\end{diagram}
\equiv\id: \circledast^I \to \circledast^I,
\\
\begin{diagram}[width=4em,inline,nobalance]
\boxdot^{\mb1}\boxdot^I\mcC &\rTTo^{\boxdot^{\mb1}\circledast^I} 
&\boxdot^{\mb1}\mcC
\\
\uTTo<{\Lambda^{I\to {\mb1}}} &\uTwoar<{\lambda^{I\to {\mb1}}} 
&\dTTo<{\circledast^{\mb1}} \, \overset{\rho}{\sss\Rightarrow} \, 
\dTTo>{\Rho} 
\\
\boxdot^I\mcC &\rTTo^{\circledast^I} &\mcC
\end{diagram}
&=\,
\begin{diagram}[width=2.2em,inline,nobalance]
\boxdot^{\mb1}\boxdot^I\mcC &&\rTTo^{\boxdot^{\mb1}\circledast^I} 
&&\boxdot^{\mb1}\mcC 
\\
\uTTo<{\Lambda^{I\to{\mb1}}} &\rdTTo<{\ttt=}>{\Rho} &&= &\dTTo>{\Rho} 
\\
\boxdot^I\mcC &\rTTo^1 &\boxdot^I\mcC &\rTTo^{\circledast^I} &\mcC
\end{diagram}
\equiv\id: \circledast^I \to \circledast^I;
\end{align*}

\item for any pair of composable maps $I \rTTo^f J \rTTo^g K$ from 
$\co_\sk$ this equation holds:
\end{enumerate}
\begin{equation*}%\label{lax-sym-Mon-cat-lambda-fg}
\begin{diagram}[inline,width=4em,height=4em]
\boxdot^J\mcC & \rTTo^{\circledast^g_\mcC} & \boxdot^K\mcC \\
\uTTo<{\circledast^f_\mcC} &\luTwoar(1,2)<{\lambda^f}
\rdTTo_{\circledast^J} \ruTwoar(1,1)^{\lambda^g} &\dTTo>{\circledast^K} 
\\
\boxdot^I\mcC &\rTTo_{\circledast^I} &\mcC
\end{diagram}
\quad=\quad
\begin{diagram}[inline,width=6em,height=4em]
\boxdot^J\mcC & \rTTo^{\circledast^g_\mcC} & \boxdot^K\mcC \\
\uTTo<{\circledast^f_\mcC} &
\luTwoar(1,1)^{\boxdot^{k\in K}\lambda^{f:f^{-1}g^{-1}k\to g^{-1}k}}
\ruTTo_{\circledast_\mcC^{g\circ f}} \ruTwoar(1,2)>{\lambda^{g\circ f}}
& \dTTo>{\circledast^K} \\
\boxdot^I\mcC & \rTTo_{\circledast^I} & \mcC
\end{diagram}
\end{equation*}
A \emph{weak $\Cat$-span multicategory} $\mcC$ is a lax $\Cat$-span 
multicategory $\mcC$ such that $\lambda^f$, $\rho$ are invertible.
\end{definition}

Here \(\boxdot^{k\in K}\lambda^{f:f^{-1}g^{-1}k\to g^{-1}k}\) means the
2\n-morphism
\begin{diagram}[nobalance,width=7.5em,height=2.2em]
&& \boxdot^J\mcC && \\
& \ruTTo^{\boxdot^{j\in J}\circledast^{f^{-1}j}} & =
& \rdTTo^{\Lambda^g} & \\
\boxdot^{j\in J}\boxdot^{f^{-1}j}\mcC & \rTTo^{\Lambda^g}
& \boxdot^{k\in K}\boxdot^{j\in g^{-1}k}\boxdot^{f^{-1}j}\mcC
& \rTTo^{\boxdot^{k\in K}\boxdot^{j\in g^{-1}k}\circledast^{f^{-1}j}}
& \boxdot^{k\in K}\boxdot^{g^{-1}k}\mcC \\
\uTTo<{\Lambda^f} & =
& \uTTo~{\boxdot^{k\in K}\Lambda^{f:f^{-1}g^{-1}k\to g^{-1}k}}
& \uTwoar~{\boxdot^{k\in K}\lambda^{f:f^{-1}g^{-1}k\to g^{-1}k}}
& \dTTo>{\boxdot^{k\in K}\circledast^{g^{-1}k}} \\
\boxdot^I\mcC & \rTTo_{\Lambda^{g\circ f}}
& \boxdot^{k\in K}\boxdot^{f^{-1}g^{-1}k}\mcC
& \rTTo_{\boxdot^{k\in K}\circledast^{f^{-1}g^{-1}k}} & \boxdot^K\mcC
\end{diagram}

The top quadrilateral in above diagram is the identity 2\n-morphism due
to the 2\n-transformation $\Lambda^g$ being strict. The left square is
the tetrahedron equation for $\Lambda$.

\begin{examples}
1) Assume that \(\mcC=\bigl(\cc\lTTo^\src \mcC\rTTo^\tgt \cc\bigr)\) is 
a lax $\Cat$-span multicategory.
In particular, we assume that \(\mcC(n)=\emptyset\) for $n\ne1$.
Then there is only one summand in \eqref{eq-boxICh-U(lp)} -- the one 
indexed by the linear tree \(t:[I]\to\co_\sk\), \(t(i)=\mb1\).
Thus,
\begin{equation}
\boxdot^{\mb m}\mcC =\lim( \cC \lTTo^\src \mcC \rTTo^\tgt \cC \lTTo^\src 
\mcC \rTTo^\tgt \dots \cC \lTTo^\src \mcC \rTTo^\tgt \cC),
\label{eq-lim(C-C-C-C-C-C)}
\end{equation}
where the number of vertices $\mcC$ is $m$ and the limit is taken in 
$\Cat$.
In these assumptions we describe Examples 1.1--1.3:

1.1) Denote by $\1$ the terminal category.
A lax $\Cat$-span multicategory \(\1\lTTo^\src \mcC\rTTo^\tgt \1\) is 
nothing else but a lax Monoidal category, see
\cite[Definition~2.5]{BesLyuMan-book}.
In fact, here \(\boxdot^I\mcC=\mcC^I\) and
\(\circledast^I:\mcC^I\to\mcC\) becomes the Monoidal product.

1.2) A strict $\Cat$-span multicategory
	\((\cc\lTTo^\src \mcC\rTTo^\tgt \cc,\circledast^I,\lambda^f=\id,
	\rho=\id)\)
is the same thing as a category internal to $\Cat$.
Equivalently, it is a double category and $\cc$ is its category of 
vertical morphisms.

1.3) A lax $\Cat$-span multicategory
\((\cc\lTTo^\src \mcC\rTTo^\tgt \cc,\circledast^I,\lambda^f,\rho)\) can 
be called a lax double category.
Fix an object $A$ of $\cc$ and consider the subcategory $\mcD$ of 
$\mcC$, whose objects $X$ satisfy \(\src X=\tgt X=A\) and morphisms $f$ 
satisfy \(\src f=\tgt f=1_A\).
Then \(\1=\{A\}\leftarrow\mcD\to\{A\}=\1\) is a lax $\Cat$-span 
multicategory of the type considered in Example~1.1. 
Thus, \((\mcD,\circledast^I,\lambda^f,\rho)\) is a lax Monoidal 
category.

2) Assume that $\mcC$ is a lax $\Cat$-span multicategory, $\cc$, $\mcC$ 
are strict 2\n-categories, \(\cc\lto{\src}\mcC\rto{\tgt}\cc\) and 
\(\circledast^I:\boxdot^I\mcC\to\mcC\) are strict 2\n-functors, 
$\lambda^f$ and $\rho$ are strict 2\n-morphisms.
In other words, all data are enriched in $\Cat$.
Notice as above that \(\mcC(n)=\emptyset\) for $n\ne1$.
In these assumptions we have Examples 2.1--2.3:

2.1) If $\cc=\1$ is the terminal 2\n-category (that with a unique
2\n-morphism), then $\mcC$ is nothing else but a lax Monoidal
$\Cat$-category, see \cite[Definition~2.10]{BesLyuMan-book}.

2.2) If $\lambda^f=\id$, $\rho=\id$, then $\mcC$ is the same thing as an 
internal category in strict-2-$\Cat=\Cat\text-\Cat$.
Equivalently, it is a triple category, whose one direction category is 
discrete.
Its 3\n-cells can be visualized as cylinders.

2.3) Dropping the restrictions of 2.2) we may call $\mcC$ a lax triple 
category, whose one direction category is discrete.
For instance, $\smQuiver$ is such a weak triple category.
Notice that its powers are given by \(\boxbox^{\bn}\) in place of 
\(\boxdot^{\bn}\), compare \eqref{eq-power-SMQCat-lim} and 
\eqref{eq-lim(C-C-C-C-C-C)}.
Also \((\smQuiver,\boxdot^I,\Lambda_{\smQuiver}^\phi,\Rho)\) uses 
$\boxdot^I$, $\Lambda_{\smQuiver}^\phi$, $\Rho$ in place of 
$\circledast^I$, $\lambda^\phi$, $\rho$.

Fix an object $A$ of $\cc$ and consider the 2\n-subcategory $\mcD$ of 
$\mcC$, whose objects $X$ (resp. 1\n-morphisms $f$, 2\n-morphisms 
$\alpha$) satisfy \(\src X=\tgt X=A\) (resp. \(\src f=\tgt f=1_A\), 
\(\src\alpha=\tgt\alpha=1_{1_A}\)).
Then \(\1=\{A\}\leftarrow\mcD\to\{A\}=\1\) is a lax $\Cat$-span 
submulticategory of the type considered in Example~2.1.
Thus, \((\mcD,\circledast^I,\lambda^\phi,\rho)\) is a lax Monoidal 
$\Cat$-\hspace{0pt}category.
In particular,
\((\sS{^\cc}\smQuiver,\boxdot^I,\Lambda_{\smQuiver}^\phi,\Rho)\) is a 
(weak) Monoidal $\Cat$-\hspace{0pt}category.
Now we may say that \emph{a lax $\Cat$-span multicategory
\((\mcC,\circledast^I,\lambda^\phi,\rho)\) is a lax-Monoidal-category 
inside the Monoidal $\Cat$-\hspace{0pt}category
\((\sS{^\cc}\smQuiver,\boxdot^I,\Lambda_{\smQuiver}^\phi,\Rho)\)}.

3) A strict $\Cat$-span multicategory $\mcC$ (the one with 
$\lambda^\phi=\id$, $\rho=\id$) is the same thing as a
${}^*$-multicategory in the sense of Leinster
\cite[Definition~4.2.2]{math.CT/0305049} for the Cartesian monad 
\(\text-^*:\Cat\to\Cat\), \(\cc\mapsto\cc^*\) of free strict monoidal 
category.

4) Let us discuss also the particular case of a $\Cat$-span multiquiver 
$\mcC$ for which $\sou\mcC$ and $\tar\mcC$ are discrete categories.
In that case the category $\mcC$ is a disjoint union of full 
subcategories \(\mcC((X_i)_{i\in I};Y)\), \(X_i\in\Ob\sou\mcC\), 
\(Y\in\Ob\tar\mcC\).
Therefore, the notion of a $\Cat$-span multiquiver coincides in the 
mentioned case with the notion of a $\Cat$-multiquiver 
\cite[Definition~3.2]{BesLyuMan-book}.
The latter notion is a particular case of $\cv$\n-multiquivers defined 
for an arbitrary monoidal category $\cv$, not only for $\cv=\Cat$.

On the other hand, the notion a lax $\Cat$-span multicategory $\mcC$ 
with the discrete category of objects $\tar\mcC$, called shortly a lax 
$\Cat$\n-multicategory, comprises more examples than the notion a 
$\Cat$\n-multicategory (a particular case of $\cv$\n-multicategories 
\cite[Definition~3.7]{BesLyuMan-book}).
In fact, natural transformations $\lambda^f$ and $\rho$ have to be 
identity transformations in the latter case.
\end{examples}

\subsection{Lax \texorpdfstring{$\Cat$}{Cat}-span operads.}
A \emph{lax $\Cat$-span operad} is a particular case of a lax
$\Cat$-span multicategory $\mcC$ -- that for which \(\tar\mcC=\1\) is 
the terminal category.
Lax $\Cat$\n-operad is a shorthand for the lax $\Cat$-span operad.
Let us give more comments on its structure.
First of all, the monoidal product of multiquivers $\mcC_h$ in the 
category \(\sS{^\1}\smQuiver\) is
\begin{align*}
\boxdot^{h\in I}\mcC_h &=\bigsqcup_{\text{tree }t:[I]\to\co_\sk} 
\prod_{h\in I}(\mcC_h)^{t(h)},
\\
(\boxdot^{h\in I}\mcC_h)(n)
&=\bigsqcup_{\text{tree }t:[I]\to\co_\sk}^{|t(0)|=n}
\prod_{h\in I} \prod_{b\in t(h)} \mcC_h(t_h^{-1}b).
\end{align*}
The 1\n-morphism $\circledast^I:\boxdot^I\mcC\to\mcC$ takes the form of 
a collection of functors
\begin{equation}
\circledast^I(n): \bigsqcup_{\text{tree }t:[I]\to\co_\sk}^{|t(0)|=n}
\prod_{h\in I} \prod_{b\in t(h)} \mcC_h(t_h^{-1}b) \to \mcC(n).
\label{eq-*I(n)-C(n)}
\end{equation}

\begin{proposition}
Let $\mcC$ be a lax $\Cat$-span operad.
This induces a Monoidal structure on the category $\mcC(1)$.
\end{proposition}

\begin{proof}
Define a \emph{linear tree} corresponding to a set \(I\in\Ob\co_\sk\) as 
the functor \(lt_I:[I]\to\co_\sk\),
\([I]\ni i\mapsto\mb1\).
We may view $lt_I=(\mb1\to\mb1\to\dots\to\mb1)$ as a synonym for $[I]$.
Restricting functor~\eqref{eq-*I(n)-C(n)} to the linear tree $lt_I$ we 
get a functor
\[ \odot^I \overset{\text{def}}= \circledast^I(lt_I)(1)|:
\mcC(1)^I \to\mcC(1).
\]
Notice that for a linear tree $t$ and any \(f:I\to J\) the trees 
$t_\psi$ and \(t^{|1}_{[\psi(m-1),\psi(m)]}\) from
\secref{sec-Cat-operad-graded-k-modules} are also linear.
Therefore, \(\circledast^f(1)\) maps the component indexed by the linear 
tree $lt_I$ to the component indexed by $lt_J$, and this restriction 
\(\circledast^f(1)|\) coincides with $\odot^f$.

Restricting given 2-morphism $\lambda^f$ (for $\boxdot$) to linear trees 
we get a natural transformation $\lambda_\odot^f$ (for $\odot$).
The given 2\n-morphism
\(\rho:\circledast^{\mb1}\to\Rho:\boxdot^{\mb1}\mcC\to\mcC\) restricted 
to the linear tree $lt_\1=(\mb1\to\mb1)$ gives a transformation 
\(\rho_\odot:\odot^{\mb1}\to\Iso:\mcC(1)^{\mb1}\to\mcC(1)\).

Two equations for $\lambda_\odot^f$ and $\rho_\odot$ follow from two 
equations~\ref{def-lax-Cat-span-multicategory}(i) for $\lambda_\mcC^f$ 
and $\rho_\mcC$.
The tetrahedron equation for $\lambda_\odot^f$ follows from tetrahedron
equation~\ref{def-lax-Cat-span-multicategory}(ii) for $\lambda_\mcC^f$.
\end{proof}

\subsection{Lax \texorpdfstring{$\Cat$}{Cat}-operad of graded 
	\texorpdfstring{$\kk$}k-modules.}
Let us construct examples $\mcG$, $\mcDG$ of a weak $\Cat$-operad: that 
of (differential) graded $\kk$-modules.
An operad is a multicategory with one object $*$.
Instead of $\mcG(\sS{^I}{*};*)$, \(I\in\Ob\co_\sk\), the notation 
$\mcG(I)$ is used.
We define $\mcG(I)=\gr^{\NN^I}$, $\mcDG(I)=\dg^{\NN^I}$.
The two cases are similar, differing only by absence or presence of the 
differential.
So we give our formulae only for one of them.
A tree $t$ gives rise to a functor
\[ \circledast(t): \prod_{(h,b)\in\IV(t)}\mcG(t_h^{-1}b) \to \mcG(t(0)), 
\qquad (\cp_h^b)_{(h,b)\in\IV(t)} \mapsto
\circledast(t)(\cp_h^b)_{(h,b)\in\IV(t)}.
\]

A $t$\n-tree is a functor \(\tau:t\to\co_\sk\), \(\tau(\troot)=\mb1\), 
where the poset $t$ is the free category built on the quiver $t$ 
oriented from the root.
It has the set of objects $\iV(t)$, the set of vertices of $t$, 
morphisms are oriented paths, and the root is the terminal object of 
$t$.

A tree can be defined as a successor map \(S:V\to V\), where $V$ is a 
non-empty finite set (of vertices), such that \(\im(S^k)\) contains only 
one element for some $k\in\NN$.
There is only one vertex $v\in V$ such that \(S(v)=v\), it is called the 
root.
An oriented graph without loops $G$ is constructed out of $S$, whose set 
of vertices is $V$ and arrows are \(v\to S(v)\) if vertex $v$ is not the 
root.
Since $G$ is a connected graph, whose number of edges is one less than 
the number of vertices, it is a tree.
The only oriented path connecting a vertex $v$ with the root consists of 
$v$, $S(v)$, $S^2(v)$, \dots, the root.
For any tree $T$ denote by \(S_T:\iV(T)\to\iV(T)\) its successor map.
A morphism of \(S:V\to V\) to \(S':V'\to V'\) is a mapping \(f:V\to V'\) 
such that \(f\cdot S'=S\cdot f\).
A $t$\n-tree $\tau$ comes with a morphism \(S_\tau\to S_t\) of successor 
maps.

Thus, for a tree $t:[I]\to\co_\sk$ with $I=\bn$,
\([I]=[n]=\{0,1,\dots,n\}\),
\begin{equation}
\circledast(t)(\cp_h^b)_{(h,b)\in\IV(t)}(z) 
=\bigoplus^{t-\text{tree }\tau}_{\forall a\in t(0)\,|\tau(0,a)|=z^a}
\bigotimes^{h\in I} \bigotimes^{b\in t(h)} \bigotimes^{p\in \tau(h,b)} 
\cp_h^b\Bigl(
\bigl(|\tau_{(h-1,a)\to(h,b)}^{-1}(p)|\bigr)_{a\in t_h^{-1}b}\Bigr).
\label{eq-O(t)(P)(z)}
\end{equation}

The following 2\n-isomorphism has to be specified for a non-decreasing 
map \(f:I\to J\), a tree $t$ of height $|I|$ and the induced map 
\(\psi=[f]:[J]\to[I]\) and trees $t_\psi$,
\(t^{|c}_{[\psi(g-1),\psi(g)]}\) (see \eqref{eq-x<fy-fx<y}) 
\begin{diagram}[h=2.2em,nobalance]
\prod_{(h,b)\in\IV(t)} \mcDG(t_h^{-1}b) &\rTTo^{\Lambda^f}_\sim 
&\prod_{(g,c)\in\IV(t_\psi)}
\prod_{(h,b)\in\IV(t^{|c}_{[\psi(g-1),\psi(g)]})} \mcDG(t_h^{-1}b)
\\
\dTTo<{\circledast(t)} &\rTwoar^{\lambda^f} 
&\dTTo>{\prod_{(g,c)\in\IV(t_\psi)}
	\circledast(t^{|c}_{[\psi(g-1),\psi(g)]})}
\\
\mcDG(t(0)) &\lTTo^{\circledast(t_\psi)} 
&\prod_{(g,c)\in\IV(t_\psi)} \mcDG((t_\psi)_g^{-1}c) 
\end{diagram}
This is an invertible natural transformation
\begin{diagram}[h=2.3em,nobalance]
\prod_{(h,b)\in\IV(t)} \dg^{\NN^{t_h^{-1}b}} &\rTTo^{\Lambda^f}_\sim 
&\prod_{(g,c)\in\IV(t_\psi)}
\prod_{(h,b)\in\IV(t^{|c}_{[\psi(g-1),\psi(g)]})} \dg^{\NN^{t_h^{-1}b}}
\\
\dTTo<{\circledast(t)} &\rTwoar^{\lambda^f} 
&\dTTo>{\prod_{(g,c)\in\IV(t_\psi)}
	\circledast(t^{|c}_{[\psi(g-1),\psi(g)]})}
\\
\dg^{\NN^{t(0)}} &\lTTo^{\circledast(t_\psi)} 
&\prod_{(g,c)\in\IV(t_\psi)} \dg^{\NN^{(t_\psi)_g^{-1}c}}
\end{diagram}
On collections a natural bijection has to be constructed:
\[ \lambda^f: \circledast(t)(\cp_h^b)_{(h,b)\in\IV(t)} \to
\circledast(t_\psi)\bigl( \circledast(t^{|c}_{[\psi(g-1),\psi(g)]})
(\cp_h^b)_{(h,b)\in\IV(t^{|c}_{[\psi(g-1),\psi(g)]})} 
\bigr)_{(g,c)\in\IV(t_\psi)}.
\]

Denote 
\[ \cq_g^c =\circledast(t^{|c}_{[\psi(g-1),\psi(g)]})
(\cp_h^b)_{(h,b)\in\IV(t^{|c}_{[\psi(g-1),\psi(g)]})}. 
\]
For an arbitrary $t$\n-tree $\tau$ define a $t_\psi$\n-tree 
$\sS{_\psi}\tau$ as the composition \(t_\psi\to t\rTTo^\tau \co_\sk\); 
the first functor being constructed from the maps
\(\id:t_\psi(g)\to t(\psi(g))\).
For arbitrary \(g\in J\), \(c\in t_\psi(g)\), \(q\in \tau(\psi(g),c)\) 
define a tree \(\tR gcq\tau:t^{|c}_{[\psi(g-1),\psi(g)]}\to\co_\sk\) 
with \(\tR gcq\tau(h,b)\simeq \tau_{(h,b)\to(\psi(g),c)}^{-1}(q)\).
Define $\lambda^f$ as the composition
\begin{multline}
\circledast(t)(\cp_h^b)_{(h,b)\in\IV(t)}(z)
=\bigoplus^{t-\text{tree }\tau}_{\forall a\in t(0)\,|\tau(0,a)|=z^a}
\bigotimes^{h\in I} \bigotimes^{b\in t(h)} \bigotimes^{p\in \tau(h,b)} 
\cp_h^b\Bigl(\bigl(|\tau_{a\to b}^{-1}(p)|\bigr)_{a\in t_h^{-1}b}\Bigr) 
\\
\rto\sim 
\bigoplus^{t_\psi-\text{tree }\sS{_\psi}\tau}
 _{\substack{\forall a\in t(0)\\|\sS{_\psi}\tau(0,a)|=z^a}}
\bigoplus^{
 \substack{\forall g\in J,\,c\in t_\psi(g),\,q\in\sS{_\psi}\tau(g,c)
 \\t^{|c}_{[\psi(g-1),\psi(g)]}-\text{tree }\tR gcq\tau}}
 _{\substack{\forall d\in t_{\psi,g}^{-1}(c)
 \\ \tR gcq\tau(\psi(g-1),d)\simeq\sS{_\psi}\tau_{d\to c}^{-1}(q)}}
\bigotimes^{\substack{g\in J\\c\in t_\psi(g)}}_{q\in\sS{_\psi}\tau(g,c)} 
\bigotimes^{\substack{h\in J\\ h\le\psi(g)}}_{h>\psi(g-1)} 
\bigotimes^{b\in t^{-1}_{h\to\psi(g)}(c)}_{p\in\tR gcq\tau(h,b)} 
\cp_h^b\Bigl(
\bigl(|\tR gcq\tau_{a\to b}^{-1}(p)|\bigr)_{a\in t_h^{-1}b}\Bigr)
\\
\rto\sim 
\bigoplus^{t_\psi-\text{tree }\sS{_\psi}\tau}
 _{\substack{\forall a\in t(0)\\|\sS{_\psi}\tau(0,a)|=z^a}}
\bigotimes^{g\in J}_{\substack{c\in t_\psi(g)\\q\in\sS{_\psi}\tau(g,c)}}
\bigoplus^{t^{|c}_{[\psi(g-1),\psi(g)]}-\text{tree }\tR gcq\tau}
 _{\substack{\forall d\in t_{\psi,g}^{-1}(c)
 \\ \tR gcq\tau(\psi(g-1),d)\simeq\sS{_\psi}\tau_{d\to c}^{-1}(q)}}
\bigotimes^{h\in J}_{\substack{h\le\psi(g)\\h>\psi(g-1)}}
\bigotimes^{b\in t^{-1}_{h\to\psi(g)}(c)}_{p\in\tR gcq\tau(h,b)} 
\cp_h^b\Bigl(
\bigl(|\tR gcq\tau_{a\to b}^{-1}(p)|\bigr)_{a\in t_h^{-1}b}\Bigr)
\\
=\bigoplus^{t_\psi-\text{tree }\sS{_\psi}\tau}
 _{\forall a\in t(0)\,|\sS{_\psi}\tau(0,a)|=z^a}
\bigotimes^{g\in J} \bigotimes^{c\in t_\psi(g)}
\bigotimes^{q\in\sS{_\psi}\tau(g,c)}
\cq_g^c\Bigl(\bigl(|\sS{_\psi}\tau_{d\to c}^{-1}(q)|\bigr)
 _{d\in t_{\psi,g}^{-1}c}\Bigr)
\\
=\circledast(t_\psi)\bigl( \circledast(t^{|c}_{[\psi(g-1),\psi(g)]})
(\cp_h^b)_{(h,b)\in\IV(t^{|c}_{[\psi(g-1),\psi(g)]})} 
\bigr)_{(g,c)\in\IV(t_\psi)}.
\label{eq-circledast(tpsi)circledast}
\end{multline}
The second isomorphism is the distributivity isomorphism between the 
tensor product and the direct sum.
The first mapping is well defined, since for any $t$\n-tree $\tau$ the 
trees constructed out of it satisfy
\[ \tR gcq\tau(\psi(g-1),d) \simeq 
\tau_{(\psi(g-1),d)\to(\psi(g),c)}^{-1}(q)
=\sS{_\psi}\tau_{(g-1,d)\to(g,c)}^{-1}(q).
\]
Thus $\tau$ is mapped to the set of collections of trees of the 
following type
\begin{equation}
\bigl( \sS{_\psi}\tau, \tR gcq\tau \mid g\in J,\, c\in t_\psi(g),\, 
q\in\sS{_\psi}\tau(g,c)\; \forall d\in t_{\psi,g}^{-1}(c)\, 
\tR gcq\tau(\psi(g-1),d) \simeq
\sS{_\psi}\tau_{(g-1,d)\to(g,c)}^{-1}(q) \bigr).
\label{eq-(tau-tau)}
\end{equation}
The last condition means that we are given a morphism
\(\tR gcq\tau(\psi(g-1),d)\to\sS{_\psi}\tau(g-1,d)\in\co_\sk\), whose 
image in $\Set$ coincides with the preimage of $q$ under the map
	\(\sS{_\psi}\tau_{(g-1,d)\to(g,c)}:
	\sS{_\psi}\tau(g-1,d)\to\sS{_\psi}\tau(g,c)\).
This gives a bijection from the set of trees $\tau$ to the set of such 
collections.
The inverse mapping is given by the formula
\begin{equation}
\tau(h,b) =\bigsqcup_{q\in\sS{_\psi}\tau(\phi h,t_{h\to\psi\phi h}(b))} 
\tR{\phi h}{t_{h\to\psi\phi h}(b)}q\tau(h,b),
\label{eq-T(hb)=bigsqcup}
\end{equation}
where the disjoint union is taken in the sense of $\co_\sk$.

Disjoint union of a family of sets \((S_k)_{k\in K}\) is 
\begin{equation}
\bigsqcup_{k\in K}S_k 
=\bigl\{(k,s)\in K\times\bigcup_{k\in K}S_k \mid s\in S_k\bigr\}.
\label{eq-Disjoint-union-family-sets}
\end{equation}
If $K$ and all $S_k$ are ordered, then \(\sqcup_{k\in K}S_k\) is 
lexicographically ordered, \(k_1<k_2\) implies \((k_1,s_1)<(k_2,s_2)\), 
and \((k,s)\le(k,s')\) is equivalent to \(s\le s'\).
When $K$ and all $S_k$ are totally ordered, then so is the disjoint 
union.
If $K$ and all $S_k$ are objects of $\co_\sk$, we denote by 
\(\sqcup_{k\in K}S_k\) the only object of $\co_\sk$ in bijection with 
\eqref{eq-Disjoint-union-family-sets}.
It comes equipped with morphisms of $\co_\sk$
\[ S_k \rTTo^{\inj_k} \sqcup_{k\in K}S_k \rTTo K.
\]
We may view \eqref{eq-Disjoint-union-family-sets} as an interpretation 
(a generalization) of a sum \(\sum_{k\in K}s_k\) of non-negative 
integers.

Let us define the successor map 
\(\tau_{(h-1,b)\to(h,t_h(b))}:\tau(h-1,b)\to\tau(h,t_h(b))\) for the 
tree $\tau$ constructed out from the collection
\((\sS{_\psi}\tau,\tR gcq\tau)\) in \eqref{eq-T(hb)=bigsqcup}.
Assume first that \(\phi(h-1)=\phi(h)\). 
Then \(\psi\phi(h-1)=\psi\phi(h)\), the indexing sets coincide, and the 
sought successor map is a disjoint union of successor maps for each 
summand:
\begin{multline*}
\tau_{(h-1,b)\to(h,t_h(b))} =\bigsqcup_{\id}
\tR{\phi h}{t_{h-1\to\psi\phi h}(b)}q\tau_{(h-1,b)\to(h,t_h(b))}:
\\
\tau(h-1,b)
=\bigsqcup_{q\in\sS{_\psi}\tau(\phi h,t_{h-1\to\psi\phi h}(b))} 
\tR{\phi h}{t_{h-1\to\psi\phi h}(b)}q\tau(h-1,b) 
\\
\to \bigsqcup_{q\in\sS{_\psi}\tau(\phi h,t_{h\to\psi\phi h}(t_h(b)))} 
\tR{\phi h}{t_{h\to\psi\phi h}(t_h(b))}q\tau(h,t_h(b)) =\tau(h,t_h(b)).
\end{multline*}
Assume now that \(\phi(h-1)<\phi(h)\). 
Equivalence~\eqref{eq-x<fy-fx<y} written as 
\[ x\le\psi(y) \Longleftrightarrow \phi(x)\le y
\]
for any \(x\in[I]\), \(y\in[J]\), implies the counit
\(\phi\psi y\le y\), the unit \(x\le\psi\phi x\), the identities 
\(\psi\phi\psi=\psi\) and \(\phi\psi\phi=\phi\).
Its equivalent form
\[ x>\psi(y) \Longleftrightarrow \phi(x)>y
\]
implies in our case inequalities 
\(h-1\le\psi\phi(h-1)\le\psi(\phi h-1)<h\le\psi\phi h\).
Hence, \(h-1=\psi\phi(h-1)=\psi(\phi h-1)\).

If \(y<g\in[J]\) and \(\psi y=\psi g\), then for each collection 
\((\sS{_\psi}\tau,\tR gcq\tau)\) and for any
\(c\in t_\psi(y)=t_\psi(g)\) we have
\[ \sS{_\psi}\tau_{(y,c)\to(g,c)}=\id:
\sS{_\psi}\tau(y,c) \to \sS{_\psi}\tau(g,c).
\]
In fact, it suffices to see it for $y=g-1$.
The preimage of any \(q\in\sS{_\psi}\tau(g,c)\) consists of one point, 
since \(\tR gcq\tau(\psi(g-1),c)=\tR gcq\tau(\psi g,c)=\mb1\).
Hence, \(\sS{_\psi}\tau_{(g-1,c)\to(g,c)}\) is an order preserving 
bijection.

For \(b\in t(h-1)\) define 
\(c=t_{h\to\psi\phi h}(t_h(b))=t_{h-1\to\psi\phi h}(b)\).
Let us define the successor map in the second case:
\begin{equation*}
\tau(h-1,b) =\bigsqcup_{q\in\sS{_\psi}\tau(\phi(h-1),b)} 
\tR{\phi(h-1)}bq\tau(h-1,b) 
\to \bigsqcup_{p\in\sS{_\psi}\tau(\phi h,c)} \tR{\phi h}cp\tau(h,t_h(b)) 
=\tau(h,t_h(b)).
\end{equation*}
It will map summand to summand in accordance with the mapping of 
indexing sets chosen as
\[ \sS{_\psi}\tau_{(\phi(h-1),b)\to(\phi h,c)}:
\sS{_\psi}\tau(\phi(h-1),b) \to \sS{_\psi}\tau(\phi h,c), \qquad 
q\mapsto p.
\]
We have chosen an arbitrary $q$ and took its image $p$.
In order to specify a mapping 
\(\gamma:\tR{\phi(h-1)}bq\tau(h-1,b)\to\tR{\phi h}cp\tau(h,t_h(b))\) we 
notice that the source is $\mb1$ embedded into 
\(\sS{_\psi}\tau(\phi(h-1),b)=\sS{_\psi}\tau(\phi h-1,b)\) via 
\(1\mapsto q\).
By definition of $p$ 
\[ q\in \sS{_\psi}\tau_{(\phi h-1,b)\to(\phi h,c)}^{-1}(p) \simeq 
\tR{\phi h}cp\tau(\psi(\phi h-1),b) =\tR{\phi h}cp\tau(h-1,b)
\]
and $q$ is represented by an element $q'$ of the latter set.
Finally, \(\gamma(1)\overset{\text{def}}=\tR{\phi h}cp\tau_h(q')\).

In the particular case of $h=\psi g$ we take into account that 
\(\psi\phi\psi g=\psi g\) and this formula gives
\[ \tau(\psi g,b) 
=\bigsqcup_{q\in\sS{_\psi}\tau(\phi\psi g,b)}
\tR{\phi\psi g}{b}q\tau(\psi g,b)
=\bigsqcup_{q\in\sS{_\psi}\tau(\phi\psi g,b)} \mb1
=\sS{_\psi}\tau(\phi\psi g,b) =\sS{_\psi}\tau(g,b),
\]
as it should be.
Further verifications show that the disassembling and assembling maps 
for $t$\n-trees are indeed inverse to each other, and invertibility of 
$\lambda^f$ follows.

\begin{example}
The Monoidal product in the category $\mcDG(1)=\dg^\NN$ of collections 
$\ca_h$ is isomorphic to \(\circledast(lt_I)(\ca_h)_{h\in I}(z)\):
\begin{equation*}
(\odot^{h\in I}\ca_h)(z) 
=\bigoplus^{\text{tree }\tau:[I]\to\co_\sk}_{|\tau(0)|=z}
\bigotimes^{h\in I} \bigotimes^{p\in \tau(h)} \ca_h(|\tau_h^{-1}p|).
\end{equation*}
This turns $\dg^\NN$ into a quite familiar Monoidal category.
Algebras in this Monoidal category are precisely $\dg$\n-operads.
\end{example}

\begin{definition}\label{def-lax-Cat-span-multifunctor}
A \emph{lax $\Cat$-span multifunctor}
	\((F,\phi^I):(\mcL,\circledast_\mcL^I,\lambda_\mcL^f,\rho_\mcL)
	\to(\mcM,\circledast_\mcM^I,\lambda_\mcM^f,\rho_\mcM)\) 
between lax $\Cat$-span multicategories is
\begin{enumerate}
\renewcommand{\labelenumi}{\roman{enumi})}
\item a 1\n-morphism
	\(F=(\tar F,F,\tar F):\mcL=(\tar\mcL,\mcL,\tar\mcL)
	\to(\tar\mcM,\mcM,\tar\mcM)=\mcM\)
in \(\smQuiver\);

\item a 2\n-morphism for each set $I\in\Ob\co_\sk$
\begin{diagram}[LaTeXeqno]
\boxdot^I\mcL &\rTTo^{\boxdot^IF} &\boxdot^I\mcM
\\
\dTTo<{\circledast_\mcL^I} &\ldTwoar^{\phi^I} &\dTTo>{\circledast_\mcM^I} 
\\
\mcL &\rTTo^F &\mcM
\label{dia-phi-I-2-morphism}
\end{diagram}
\end{enumerate}
such that
\begin{equation*}
\begin{diagram}[inline,w=4em]
\boxdot^{\mb1}\mcL &\rTTo^{\boxdot^{\mb1}F} &\boxdot^{\mb1}\mcM
\\
\dTTo<{\Rho} \overset{\rho_\mcL}\Leftarrow 
\dTTo>{\circledast_\mcL^{\mb1}} 
&\ldTwoar^{\phi^{\mb1}} &\dTTo>{\circledast_\mcM^{\mb1}} 
\\
\mcL &\rTTo^F &\mcM
\end{diagram}
=
\begin{diagram}[inline,w=4em]
\boxdot^{\mb1}\mcL &\rTTo^{\boxdot^{\mb1}F} &\boxdot^{\mb1}\mcM
\\
\dTTo<{\Rho} &= &\dTTo<{\Rho} \overset{\rho_\mcM}\Leftarrow 
\dTTo>{\circledast_\mcM^{\mb1}} 
\\
\mcL &\rTTo^F &\mcM
\end{diagram}
\end{equation*}
and for every map $f:I\to J$ of $\co_\sk$ the following equation holds:
\begin{equation*}
\begin{diagram}[inline,nobalance,tight,width=2.6em]
&&\boxdot^I\mcL &&\rTTo^{\boxdot^IF} &&\boxdot^I\mcM
\\
&\ldTTo^{\circledast_\mcL^f}
&&&&\ldTwoar(6,2)^{\boxdot^{j\in J}\phi^{f^{-1}j}} 
\ldTTo_{\circledast_\mcM^f} &
\\
\boxdot^J\mcL &&\rTTo_{\boxdot^JF} &&\boxdot^J\mcM 
&\lTwoar_{\lambda_\mcM^f} &\dTTo>{\circledast_\mcM^I} 
\\
&\rdTTo_{\circledast_\mcL^J} &&\ldTwoar^{\phi^J} 
&&\rdTTo_{\circledast_\mcM^J} & 
\\
&&\mcL &&\rTTo^F &&\mcM
\end{diagram}
\quad=\quad
\begin{diagram}[inline,nobalance,tight,width=2.6em]
&&\boxdot^I\mcL &\rTTo^{\boxdot^IF} &\boxdot^I\mcM
\\
&\ldTTo^{\circledast_\mcL^f} &&\ldTwoar(2,4)<{\phi^I} &
\\
\boxdot^J\mcL &\lTwoar^{\lambda_\mcL^f} &\dTTo>{\circledast_\mcL^I} 
&&\dTTo>{\circledast_\mcM^I} 
\\
&\rdTTo_{\circledast_\mcL^J}
\\
&&\mcL &\rTTo^F &\mcM
\end{diagram}
\label{lax-sym-Mon-functor-f}
\end{equation*}
\end{definition}

Here 2\n-morphism \(\boxdot^{j\in J}\phi^{f^{-1}j}\) means the pasting
\begin{diagram}
\boxdot^I\mcL &\rTTo^{\boxdot^IF} &\boxdot^I\mcM
\\
\dTTo<{\Lambda^f} &= &\dTTo>{\Lambda^f}
\\
\boxdot^{j\in J}\boxdot^{f^{-1}j}\mcL
&\rTTo^{\boxdot^{j\in J}\boxdot^{f^{-1}j}F}
&\boxdot^{j\in J}\boxdot^{f^{-1}j}\mcM
\\
\dTTo<{\boxdot^{j\in J}\circledast_\mcL^{f^{-1}j}} 
&\ldTwoar^{\boxdot^{j\in J}\phi^{f^{-1}j}} 
&\dTTo>{\boxdot^{j\in J}\circledast_\mcM^{f^{-1}j}}
\\
\boxdot^J\mcL &\rTTo^{\boxdot^JF} &\boxdot^J\mcM
\end{diagram}

\begin{proposition}
	\label{pro-lax-Cat-span-multifunctor-induces-lax-Monoidal-functor}
Any lax $\Cat$-span multifunctor
\[ (F,\phi^I):(\mcL,\circledast_\mcL^I,\lambda_\mcL^f,\rho_\mcL)
\to (\mcM,\circledast_\mcM^I,\lambda_\mcM^f,\rho_\mcM)
\]
between lax $\Cat$-span operads induces a lax Monoidal functor 
\[ (F(1),\phi^I(1)|):
(\mcL(1),\odot_{\mcL(1)}^I,\lambda_{\mcL(1)}^f,\rho_{\mcL(1)})
\to(\mcM(1),\odot_{\mcM(1)}^I,\lambda_{\mcM(1)}^f,\rho_{\mcM(1)}).
\]
\end{proposition}

\begin{proof}
Let us restrict the given transformation
\begin{diagram}
\boxdot^I\mcL &\rTTo^{\boxdot^IF} &\boxdot^I\mcM
\\
\dTTo<{\circledast_\mcL^I} &\ldTwoar^{\phi^I} 
&\dTTo>{\circledast_\mcM^I} 
\\
\mcL &\rTTo^F &\mcM
\end{diagram}
to the component indexed by the linear tree $lt_I$.
This yields a natural transformation
\begin{diagram}
\mcL(1)^I &\rTTo^{F(1)^I} &\mcM(1)^I
\\
\dTTo<{\odot_{\mcL(1)}^I} &\ldTwoar^{\phi^I(1)|} 
&\dTTo>{\odot_{\mcM(1)}^I}
\\
\mcL(1) &\rTTo^{F(1)} &\mcM(1)
\end{diagram}
Two equations for $\phi^I$ from \defref{def-lax-Cat-span-multifunctor} 
imply two equations for $\psi^I$.
\end{proof}

\subsection{\texorpdfstring{$\Cat$}{Cat}-span multinatural 
	transformations.}
\begin{definition}%\label{def-lax-Cat-span-multinatural-transformation}
A \emph{$\Cat$-span multinatural transformation}
	\(\xi:(F,\phi^I)\to(G,\psi^I):
	(\mcL,\circledast_\mcL^I,\lambda_\mcL^f,\rho_\mcL)\to
	(\mcM,\circledast_\mcM^I,\lambda_\mcM^f,\rho_\mcM)\) 
between lax $\Cat$-span multifunctors is a 2\n-morphism
	\((\tar\xi,\xi,\tar\xi):(\tar F,F,\tar F)\to(\tar G,G,\tar G):
	(\tar\mcL,\mcL,\tar\mcL)\to(\tar\mcM,\mcM,\tar\mcM)\)
of $\smQuiver$ such that for all $I\in\Ob\co_\sk$
\begin{equation}
\begin{diagram}[w=4.5em,h=1.7em]
(\boxdot^I\mcL,\tar\mcL) 
&\pile{\rTTo^{(\boxdot^IF,\tar F)} 
\\ \Downarrow\sss(\boxdot^I\xi,\tar\xi) \\ \rTTo_{(\boxdot^IG,\tar G)}}
&(\boxdot^I\mcM,\tar\mcM)
\\
\dTTo<{(\circledast_\mcL^I,\Id)} &&\dTTo>{(\circledast_\mcM^I,\Id)} 
\\
&\ldTwoar^{(\psi^I,\id)} &
\\
(\mcL,\tar\mcL) &\rTTo_{(G,\tar G)} &(\mcM,\tar\mcM)
\end{diagram}
\quad =\quad
\begin{diagram}[w=4.5em,h=1.7em]
(\boxdot^I\mcL,\tar\mcL) &\rTTo^{(\boxdot^IF,\tar F)}
&(\boxdot^I\mcM,\tar\mcM)
\\
\dTTo<{(\circledast_\mcL^I,\Id)} &\ldTwoar_{(\phi^I,\id)}
&\dTTo>{(\circledast_\mcM^I,\Id)} 
\\
&&
\\
(\mcL,\tar\mcL) 
&\pile{\rTTo^{(F,\tar F)} \\ \Downarrow\sss(\xi,\tar\xi) 
\\ \rTTo_{(G,\tar G)}}
&(\mcM,\tar\mcM)
\end{diagram}
\quad.
\label{eq-dia-Cat-span-multinatural-transformation}
\end{equation}
\end{definition}

\begin{proposition}
The collection $\SMCcat$ of lax $\Cat$-span multicategories, lax
$\Cat$-span multifunctors and $\Cat$-span multinatural transformations 
is a 2\n-category.
\end{proposition}

\begin{proof}
Composition of lax $\Cat$-span multifunctors and identity $\Cat$-span 
multifunctors are the obvious ones.
The horizontal and the vertical compositions of 2\n-morphisms are 
implied by the underlying mapping \(\SMCcat\to\smQuiver\), which is a 
functor at the level of objects and 1\n-morphisms and is injective on 
2\n-morphisms.
One only has to verify that the identity 2\n-morphisms are in $\SMCcat$ 
and that if the horizontal or the vertical composition of two
2\n-morphisms in $\SMCcat$ makes sense, then its value in $\smQuiver$ 
belongs actually to $\SMCcat$.
This is clear from the shape of 
equation~\eqref{eq-dia-Cat-span-multinatural-transformation}.
\end{proof}

\begin{proposition}
\label{pro-Cat-span-multinatural-transformation-Monoidal-transformation}
Let $\mcL$, $\mcM$ be lax $\Cat$-span operads.
Any $\Cat$-span multinatural transformation
\[ \xi: (F,\phi^I) \to (G,\psi^I):
(\mcL,\circledast_\mcL^I,\lambda_\mcL^f,\rho_\mcL)
\to (\mcM,\circledast_\mcM^I,\lambda_\mcM^f,\rho_\mcM)
\]
between lax $\Cat$-span multifunctors induces a Monoidal transformation 
\[ \xi(1):(F(1),\phi^I(1))\to(G(1),\psi^I(1)):
(\mcL(1),\odot_{\mcL(1)}^I,\lambda_{\mcL(1)}^f,\rho_{\mcL(1)})
\to(\mcM(1),\odot_{\mcM(1)}^I,\lambda_{\mcM(1)}^f,\rho_{\mcM(1)}).
\]
\end{proposition}

\begin{proof}
Restricting the transformations to linear trees we get
\begin{equation}
\begin{diagram}[w=4.5em,h=1.7em]
\mcL(1)^I
&\pile{\rTTo^{F(1)^I} \\ \Downarrow\sss\xi(1)^I \\ \rTTo_{G(1)^I}}
&\mcM(1)^I
\\
\dTTo<{\odot^I} &&\dTTo>{\odot^I} 
\\
&\ldTwoar^{\psi^I(1)} &
\\
\mcL(1) &\rTTo_{G(1)} &\mcM(1)
\end{diagram}
\quad =\quad
\begin{diagram}[w=4.5em,h=1.7em]
\mcL(1)^I &\rTTo^{F(1)^I} &\mcM(1)^I
\\
\dTTo<{\odot^I} &\ldTwoar_{\phi^I(1)} &\dTTo>{\odot^I} 
\\
&&
\\
\mcL(1)
&\pile{\rTTo^{F(1)} \\ \Downarrow\sss\xi(1) \\ \rTTo_{G(1)}}
&\mcM(1)
\end{diagram}
\quad.
\label{eq-xi(1)-Monoidal-transformation}
\end{equation}
These equations say that the transformation $\xi(1)$ is Monoidal, see 
\cite[Definition~2.20]{BesLyuMan-book}.
\end{proof}

\subsection{\texorpdfstring{$\dg$}{dg}-operads are lax 
	\texorpdfstring{$\Cat$}{Cat}-span multifunctors 
	\texorpdfstring{$\one\to\mcDG$}{1->DG}.}
Let $\mcC$ be an arbitrary lax $\Cat$\n-span operad.
Denote by $\one$ the unit object of the Monoidal category 
\((\sS{^\1}\smQuiver,\boxdot^I,\Lambda^f,\Rho)\) of lax $\Cat$\n-span 
operads.
This is a lax $\Cat$\n-span operad with
\[ \one(n) =
\begin{cases}
\1 =\text{terminal (1-morphism) category}, \quad &\text{ if } n=1,
\\
\emptyset =\text{initial (empty) category}, &\text{ if } n\ne1.
\end{cases}
\]
A multiquiver morphism \(\one\to\mcC\) is a functor \(\1\to\mcC(1)\), so 
it is just an object of $\mcC(1)$.
In particular, a $\Cat$\n-multiquiver map \(\one\to\mcDG\) is the same 
as a functor \(\1\to\mcDG(1)=\dg^\NN\).

\propref{pro-lax-Cat-span-multifunctor-induces-lax-Monoidal-functor} 
implies that a lax $\Cat$\n-span multifunctor \(\one\to\mcDG\) is the 
same as a Monoidal functor \(\1\to\mcDG(1)=(\dg^\NN,\odot^I)\). 
By Definition~2.25 and Proposition~2.28 of \cite{BesLyuMan-book} this is 
the same as an algebra in \((\dg^\NN,\odot^I)\), that is, an operad.

By
\propref{pro-Cat-span-multinatural-transformation-Monoidal-transformation}
a multinatural transformation \(\xi:(F,\phi^I)\to(G,\psi^I):\one\to\mcDG\)
is the same as a Monoidal transformation 
	\(\xi(1):\co=(F(1),\phi^I(1))\to(G(1),\psi^I(1))
	=\cp:\1\to(\dg^\NN,\odot^I)\).
Here operads $\co$ and $\cp$ are identified with the image of 
$1\in\Ob\1$ under corresponding functors.
Equations~\eqref{eq-xi(1)-Monoidal-transformation} for $\xi(1)$ 
translate to
\[ \bigr( \odot^I\co \rTTo^{\odot^I\xi(1)} \odot^I\cp \rTTo^{\psi^I(1)} 
\cp \bigr)
=\bigr( \odot^I\co \rTTo^{\phi^I(1)} \co \rTTo^{\xi(1)} \cp \bigr),
\]
that is, $\xi(1):\co\to\cp$ is a morphism of operads.

\subsection{\texorpdfstring{$n\wedge1$}{n1}-operad modules are lax 
	\texorpdfstring{$\Cat$}{Cat}-span multifunctors.}
	\label{sec-operad-modules-are-lax-Cat-span-multifunctors}
Consider a $\Cat$-\hspace{0pt}multicategory $\mcL_n$ with 
\(\Ob\mcL_n=\{0,1,2,\dots,n\}\) such that
\begin{alignat*}2
\mcL_n(i;i) &= \1 &\quad &\text{for } 0\le i\le n,
\\
\mcL_n(1,2,\dots,n;0) &= \1 &&\text{and}
\\
\mcL_n(k_1,k_2,\dots,k_m;k_0) &= \emptyset
&&\text{for other lists of arguments.}
\end{alignat*}
Thus, $\mcL_n$ is a discrete category with the list of objects $(0;0)$, 
$(1;1)$, \dots, $(n;n)$ and $(1,2,\dots,n;0)$.
Components of the category \(\boxdot^{\mb k}\mcL_n\) either are empty or 
indexed by trees of two kinds:
\begin{itemize}
\item labelled linear trees 
\((lt_{\mb k},\sS{^k}i)=(i\to i\to\dots\to i)\), whose all vertices are 
labelled with the same \(i\in[n]\);

\item labelled trees with $1+l$ vertices labelled with 1 (or 2, or $n$) 
and $k-l$ vertices labelled with 0
\[ (t,\ell) =\Yright_l^k =
\begin{diagram}[inline,h=0.6em,w=2em,nobalance]
n &\rTTo &\cdots &\rTTo &n &\rTTo &n \\
&&&&&&&\rdTTo(2,3) \\
\cdots &\rTTo &\cdots &\rTTo &\cdots &\rTTo &\cdots \\
&&&&&&&\rdTTo(2,1) &0 &\rTTo &0 &\rTTo &\cdots &\rTTo &0 \\
2 &\rTTo &\cdots &\rTTo &2 &\rTTo &2 &\ruTTo(2,1) \ruTTo(2,3) \\
\\
1 &\rTTo &\cdots &\rTTo &1 &\rTTo &1
\end{diagram}
\]
\end{itemize}
The non-vanishing components are terminal categories $\1$.
The 1\n-morphism \(\circledast^{\mb k}:\boxdot^{\mb k}\mcL_n\to\mcL_n\) 
is $\Id_\1$ on any non-vanishing component of the source.
It takes the component indexed by a tree of the first kind to the full 
subcategory $\{(i;i)\}\subset\mcL_n$.
The component indexed by a tree of the second kind goes to the full 
subcategory \(\{(1,2,\dots,n;0)\}\subset\mcL_n\).
The transformation $\lambda^f$ and $\rho$ are \(\id_{\Id_\1}\) on
non-vanishing components.

\begin{definition}\label{def-M-operad-polymodules}
An \emph{$n\wedge1$-operad module} is a lax $\Cat$-span multifunctor 
\(\cp:\mcL_n\to\mcDG\).
A \emph{morphism of $n\wedge1$-operad modules \(r:\cp\to\cq\)} is a 
$\Cat$-span multinatural transformation \(r:\cp\to\cq:\mcL_n\to\mcDG\).
The disjoint union $\mcM$ over $n\ge0$ of so defined categories
$\nOp n_1$ of $n\wedge1$-\hspace{0pt}operad modules is the category of 
operad polymodules.
\end{definition}

Let us describe the structure of an $n\wedge1$-operad module.
A multiquiver 1\n-morphism \(\mcL_n\to\mcDG\) amounts to a sequence 
\((\ca_1,\dots,\ca_n;\cp;\cb)\), where \(\cb,\ca_i\in\Ob\dg^\NN\) for 
\(i\in\bn\), and \(\cp\in\Ob\dg^{\NN^n}\).

\begin{example}
Particular cases of $\circledast$ for $\mcDG$ will obtain a special 
notation.
In addition to the above assume that \(\ca_i^h\in\Ob\dg^\NN\).
We denote
\begin{align*}
\cp\odot_0\cb &=\circledast(\bn \to \mb1 \to \mb1)(\cp; \cb),
\\
\odot_{>0}((\ca_i)_{i\in\bn};\cp)
&=\circledast(\bn \rto1 \bn \to \mb1)((\ca_i)_{i\in\bn}; \cp),
\\
\odot_{\ge0}((\ca_i)_{i\in\bn};\cp;\cb)
&=\circledast(\bn\rto1 \bn\to \mb1\to \mb1)((\ca_i)_{i\in\bn};\cp;\cb).
\end{align*}
All these expressions describe actions of several copies of the category 
$\dg^\NN$ on the category $\dg^{\NN^n}$.
In isomorphic form these actions are given by the graded components, 
$\ell\in\NN^n$,
\begin{gather*}
(\cp\odot_0\cb)(\ell) 
\simeq \bigoplus_{t_1+\dots+t_m=\ell}^{m\ge0} 
\Bigl(\bigotimes_{r=1}^m\cp(t_r)\Bigr) \tens \cb(m),
%\label{eq-Podot0B(l)}
\\
\odot_{>0}(\ca_1,\dots,\ca_n;\cp)(\ell) \simeq
\bigoplus_{k\in\NN^n}
\bigoplus_{j_1^i+\dots+j_{k^i}^i=\ell^i}^{\forall\,i\in\bn}
\biggl[\bigotimes_{i=1}^n\bigotimes_{p=1}^{k^i} 
\cA_i(j_p^i)\biggr]\tens\cP(k),
\\
\begin{split}
\odot_{\ge0}&(\ca_1,\dots,\ca_n;\cp;\cb)(\ell)
\\
&\simeq \bigoplus_{m=0}^\infty \bigoplus_{k_1,\dots,k_m\in\NN^n}
\bigoplus_{\sum_{p=1}^{k_1^i+\dots+k_m^i}j_p^i=\ell^i}^{\forall\,i\in\bn}
\Bigl(\bigotimes_{i=1}^n \bigotimes_{p=1}^{k_1^i+\dots+k_m^i} 
\cA_i(j_p^i)\Bigr)\tens\Bigl(\bigotimes_{r=1}^m\cp(k_r)\Bigr)\tens\cb(m)
\\
&\simeq \bigoplus_{m=0}^\infty \bigoplus_{t_1+\dots+t_m=\ell} 
\bigoplus_{k_1,\dots,k_m\in\NN^n}
\bigoplus_{y_{r,1}^i+\dots+y_{r,k_r^i}^i=t_r^i}^{\forall\,i\in\bn,\,r\in\mb m}
\biggl[\bigotimes_{r=1}^m
\Bigl(\bigotimes_{i=1}^n\bigotimes_{v=1}^{k_r^i} \cA_i(y_{r,v}^i) \Bigr)
\tens\cP(k_r)\biggr] \tens\cb(m).
\end{split}
\end{gather*}
Actually, the action $\odot_{>0}$ can be presented as a combination of 
partial actions $\odot_i$ for \(1\le i\le n\) defined as
\[ \ca\odot_i\cp = \odot_{>0}(\1,\dots,\1,\ca,\1,\dots,\1;\cp), \qquad 
\ca\text{ on $i$-th place}.
\]
Explicit presentation of this action is
\[ (\ca\odot_i\cp)(\ell) 
=\bigoplus_{j_1+\dots+j_q=\ell^i}^{q\ge0} 
\Bigl(\bigotimes_{p=1}^q\ca(j_p)\Bigr) \tens \cp(\ell,\ell^i\mapsto q),
\]
where 
	\((\ell,\ell^i\mapsto q)
	=(\ell^1,\dots,\ell^{i-1},q,\ell^{i+1},\dots,\ell^n)\).

Iterating these actions $I$ times we get the following expressions
\begin{align*}
\odot_0^{[I]}(\cp; (\cb_i)_{i\in I}) &= 
\circledast(\bn\to \underbrace{\mb1\to\mb1\dots\to\mb1}_{\mb1\sqcup I})
(\cp;(\cb_i)_{i\in I}),
\\
\odot_{>0}^{I\sqcup\mb1}(((\ca_i^h)_{h\in\bn})_{i\in I}; \cp)
&= \circledast(
\underbrace{\bn\rto1\bn\rto1\bn\dots\rto1\bn\rto1\bn}_{\mb1\sqcup I}
\to\mb1) (((\ca_i^h)_{h\in\bn})_{i\in I}; \cp),
\\
\odot_{\ge0}^{[I]}(((\ca_i^h)_{h\in\bn})_{i\in I}; \cp;(\cb_i)_{i\in I})
&= \circledast(\underbrace{\bn\dots\rto1 \bn}_{[I]}\to 
\underbrace{\mb1\to \dots \mb1}_{[I]})
(((\ca_i^h)_{h\in\bn})_{i\in I}; \cp; (\cb_i)_{i\in I}).
\end{align*}
Here the last three trees are functors \(\mb1\sqcup[I]\to\co_\sk\), 
\(\mb1\sqcup[I]^\op\to\co_\sk\) and 
\(\mb1\sqcup[I]^\op\cup_{0\sim0}[I]\to\co_\sk\) respectively, where 
\([I]^\op\cup_{0\sim0}[I]\) is obtained by identifying elements 
\(0\in[I]^\op\) and \(0\in[I]\).

One can show that these actions are Monoidal and that actions $\odot_i$ 
for different \(0\le i\le n\) commute up to isomorphisms that satisfy 
coherence conditions.
Furthermore, the action $\odot_{\ge0}$ can be presented as a combination 
of partial actions $\odot_i$ for \(0\le i\le n\).
To be rigorous this approach requires more definitions, and we have 
chosen to avoid it.
\end{example}

Thus, a lax $\Cat$-span multifunctor 
\((\ca_1,\dots,\ca_n;\cp;\cb):\mcL_n\to\mcDG\) consists of a coherent 
system of action maps
\[ \circledast(\underbrace{\bn\rto1 \dots\rto1 \bn}_{[I]}\to 
\underbrace{\mb1\to \dots\to \mb1}_{[J]})
(((\ca_i^h)_{h\in\bn})_{i\in I}; \cp; (\cb_j)_{j\in J}) \to \cp.
\]
Since the actions are Monoidal, this is equivalent to giving a coherent 
system of actions for $J=I$ only, that is, a system of morphisms
	\(\alpha^I:\odot_{\ge0}^{[I]}
	(((\ca_i^h)_{h\in\bn})_{i\in I};\cp;(\cb_i)_{i\in I})\to\cp\).
Equivalently, we consider the monad 
\(\cp\mapsto\odot_{\ge0}(\ca_1,\dots,\ca_n;\cp;\cb)\).

\begin{definition}
An \emph{$n\wedge1$-operad module} can be defined also as a family 
\((\ca_1,\dots,\ca_n;\cp;\cb)\), consisting of $n+1$ operads $\ca_i$, 
$\cb$ and an object $\cp\in\gr^{\NN^n}$ (resp. $\cp\in\dg^{\NN^n}$), 
equipped with an algebra structure
\[ \alpha: \odot_{\ge0}(\ca_1,\dots,\ca_n;\cp;\cb) \to \cp
\]
for the monad \(\cq\mapsto\odot_{\ge0}(\ca_1,\dots,\ca_n;\cq;\cb)\).
\end{definition}

An action is specified by a collection of maps given for each $m\in\NN$, 
each family \(k_1,\dots,k_m\in\NN^n\) and each family of non-negative 
integers \(\bigl((j_p^i)_{p=1}^{k_1^i+\dots+k_m^i}\bigr)_{i=1}^n\)
\begin{equation}
\alpha:
\Bigl(\bigotimes_{i=1}^n \bigotimes_{p=1}^{k_1^i+\dots+k_m^i} 
\cA_i(j_p^i)\Bigr)\tens\Bigl(\bigotimes_{r=1}^m\cp(k_r)\Bigr)\tens\cb(m)
\to 
\cp\biggl(\Bigl(\sum_{p=1}^{k_1^i+\dots+k_m^i}j_p^i\Bigr)_{i=1}^n\biggr).
\label{eq-alpha-full-action}
\end{equation}

Assume that \(f_i:\cc_i\to\ca_i\), \(g:\cd\to\cb\) are morphisms of 
operads. 
They imply a morphism of monads 
	\(\odot_{\ge0}(\cc_1,\dots,\cc_n;\cq;\cd)
	\to\odot_{\ge0}(\ca_1,\dots,\ca_n;\cq;\cb)\).
An algebra $\cp$ over the latter monad becomes an algebra over the 
former monad denoted \(\sS{_{f_1,\dots,f_n}}\cp_g\).

Restricting the action $\alpha$ to submonads 
 \(\odot_{>0}(\ca_1,\dots,\ca_n;\cq)\hookrightarrow
 \odot_{\ge0}(\ca_1,\dots,\ca_n;\cq;\cb)\),
\(\cq\odot_0\cb\hookrightarrow\odot_{\ge0}(\ca_1,\dots,\ca_n;\cq;\cb)\), 
\(\ca_i\odot_i\cq\hookrightarrow\odot_{\ge0}(\ca_1,\dots,\ca_n;\cq;\cb)\),
\(1\le i\le n\), obtained via insertion of operad units $\eta$, we get 
the partial actions
\begin{gather}
\lambda =\lambda_{k,(j_p^i)}:
\biggl[\bigotimes_{i=1}^n \bigotimes_{p=1}^{k^i} \cA_i(j_p^i)\biggr]
\tens\cP\bigl((k^i)_{i=1}^n\bigr) \to
\cP\biggl(\Bigl(\sum_{p=1}^{k^i}j_p^i\Bigr)_{i=1}^n\biggr), \notag
\\
\rho =\rho_{(k_r)}: \Bigl(\bigotimes_{r=1}^m\cp(k_r)\Bigr)\tens\cB(m)
\to \cP\biggl(\sum_{r=1}^mk_r\biggr),
\label{eq-rho-rhokr}
\\
\lambda^i =\lambda^i_{k,(j_p)}:
\biggl[ \bigotimes_{p=1}^{k^i} \cA_i(j_p)\biggr] \tens\cP(k) \to
\cP\biggl(k^1,\dots,k^{i-1},\sum_{p=1}^{k^i}j_p,k^{i+1},
\dots,k^n\biggr). 
\notag
\end{gather}
On the other hand, given $n$ left actions $\lambda^i$ of $\ca_i$ and a 
right action $\rho$ of $\cb$, all pairwise commuting, we can restore the 
total action $\alpha$.

The category of $n\wedge1$-operad modules $\nOp n_1$ has morphisms
\[ (f_1,\dots,f_n;h;f_0): (\ca_1,\dots,\ca_n;\cp;\ca_0) \to
(\cc_1,\dots,\cc_n;\cq;\cc_0),
\]
where \(f_i:\ca_i\to\cc_i\), \(0\le i\le n\), are morphisms of
$\dg$\n-operads and
\(h:\cp\to\sS{_{f_1,\dots,f_n}}\cq_{f_0}\in\dg^{\NN^n}\) is a module 
morphism with respect to actions of all $\ca_i$.
In fact, a morphism of $n\wedge1$-operad modules is by definition a pair 
of natural transformations
	\(\xi:\cp=(\ca_1,\dots,\ca_n;\cp;\ca_0)\to
	(\cc_1,\dots,\cc_n;\cq;\cc_0)=\cq:\mcL_n\to\mcDG\), 
\(\id:\con\to\con:\{0,1,\dots,n\}\to\mb1\) which satisfy certain 
equations.
Since $\mcL_n$ is a discrete category with $n+2$ objects, the 
transformation $\xi$ consists of morphisms
\(f_i:\ca_i\to\cc_i\in\dg^\NN\), \(0\le i\le n\), and
\(h:\cp\to\cq\in\dg^{\NN^n}\).
The category \(\boxdot^{\mb k}\mcL_n\) is also discrete and its objects 
are \((lt_{\mb k},\sS{^k}i)\) and \(\Yright_l^k\), see
\secref{sec-operad-modules-are-lax-Cat-span-multifunctors}.
Equation~\eqref{eq-dia-Cat-span-multinatural-transformation} on 
\(\boxdot^{\mb k}\mcL_n\) reads on objects \((lt_{\mb k},\sS{^k}i)\) as
\[ \bigl( \odot^{\mb k}\ca_i \rTTo^{\odot^{\mb k}f_i} \odot^{\mb k}\cb_i 
\rTTo^{\mu^k} \cb_i \bigr)
=\bigl( \odot^{\mb k}\ca_i \rTTo^{\mu^k} \ca_i \rTTo^{f_i} \cb_i \bigr),
\]
that is, $f_i$ is a morphism of operads.
On the object \(\Yright_l^k\) the equation gives
\begin{multline*}
\bigl[\circledast_\mcDG^{\mb k}(\Yright_l^k)(\ca_1,\dots,\ca_n;\cp;\ca_0)
\rTTo^{\circledast_\mcDG^{\mb k}(\Yright_l^k)(f_1,\dots,f_n;h;f_0)}
\circledast_\mcDG^{\mb k}(\Yright_l^k)(\cc_1,\dots,\cc_n;\cq;\cc_0) 
\rTTo^\alpha \cq \bigr]
\\
=\bigl[\circledast_\mcDG^{\mb k}(\Yright_l^k)(\ca_1,\dots,\ca_n;\cp;\ca_0)
\rTTo^\alpha \cp \rTTo^h \cq \bigr],
\end{multline*}
that is, \(h:\cp\to\sS{_{f_1,\dots,f_n}}\cq_{f_0}\) is a module morphism 
with respect to actions of all $\ca_i$, \(0\le i\le n\).

Objects of the category \(\dg^{n\NN\sqcup\NN^n\sqcup\NN}\) are also 
written as tuples \((\cu_1,\dots,\cu_n;\cx;\cw)\). 
The free algebra functor for the monad 
\(\odot_{\ge0}(\ca_1,\dots,\ca_n;\text-;\cb)\) is the functor 
\(\dg^{\NN^n}\to\ca_1\text-\cdots\text-\ca_n\)-mod-$\cb$,  
\(\cx\mapsto\odot_{\ge0}(\ca_1,\dots,\ca_n;\cx;\cb)\), left adjoint to 
the underlying functor 
\(\ca_1\text-\cdots\text-\ca_n\)-mod-\(\cb\to\dg^{\NN^n}\).
Hence, there is also a pair of adjoint functors
\(F:\dg^{n\NN\sqcup\NN^n\sqcup\NN}\rightleftarrows\nOp n_1:U\),
\[ F(\cu_1,\dots,\cu_n;\cx;\cw)=
(T\cu_1,\dots,T\cu_n;\odot_{\ge0}(T\cu_1,\dots,T\cu_n;\cx;T\cw);T\cw).
\]
The module part is indexed by trees with the top floor describing 
\(\cx(\text-_1)\tdt\cx(\text-_k)\), lower floors indexed by 
$\cw(\text-)$ and $n$ forests indexed by $\cu_i(\text-)$ attached to 
each of $k$ leaves.

In particular,
\[ F(\cu_1,\dots,\cu_n;0;\cw) =(T\cu_1,\dots,T\cu_n;(T\cw)(0);T\cw).
\]

Recall that for any object \(A=(\ca_1,\dots,\ca_n;\cp;\cb)\) and any 
collection of complexes $M$ the object \(A\langle M,0\rangle\) is 
isomorphic to \(F(M[1])\sqcup A\), see
\secref{sec-Model-category-structures}.

\subsection{The monad of free \texorpdfstring{$n\wedge1$}{n1}-operad modules.}
Recall \cite[Section~3.3.6]{BarrWells:TopTT} that a parallel pair of morphisms \(f,g:A\to B\in\cc\) is called \emph{reflexive} if there is a morphism \(r:B\to A\in\cc\) such that \(f\circ r=\id_B=g\circ r\).
Recall that a \emph{contractible coequalizer} \cite[Section~3.3.3]{BarrWells:TopTT} (= a \emph{split fork} \cite[Section~VI.6]{MacLane}) is a diagram in a category $\cd$
\begin{diagram}
A' &\pile{\rTTo^{d^0} \\ \lTTo~t \\ \rTTo_{d^1}} &B' &\pile{\rTTo^d \\ \lTTo_s} &C'
\end{diagram}
such that \(d^0\circ t=\id\), \(d^1\circ t=s\circ d\), \(d\circ s=\id\), and \(d\circ d^0=d\circ d^1\).
Suppose there is a functor \(U:\cc\to\cd\).
Then a pair \(f,g:A\to B\in\cc\) is called \emph{$U$\n-contractible coequalizer pair} if \(d^0=Uf,d^1=Ug:UA\to UB\) extend to a contractible coequalizer in $\cd$.
One says that $U$ \emph{creates} $U$\n-contractible coequalizers \cite[Section~VI.7]{MacLane} if for any pair \(f,g:A\to B\in\cc\) and any contractible coequalizer in $\cd$
\begin{diagram}
UA &\pile{\rTTo^{Uf} \\ \lTTo~t \\ \rTTo_{Ug}} &UB &\pile{\rTTo^d \\ \lTTo_s} &C'
\end{diagram}
\begin{itemize}
\item[\dag] there is a unique morphism \(h:B\to C\in\cc\) such that $C'=UC$, $d=Uh$, and
\label{pg-dag-ddag}
\item[\ddag] $h$ is a coequalizer of $(f,g)$ in $\cc$.
\end{itemize}

The following statement is Exercise~3.3.(PPTT) of \cite{BarrWells:TopTT}.
\begin{theorem}\label{thm-PPTT}
Let \(U:\cc\to\cd\) be a functor which has a left adjoint $F$.
Then the comparison functor \(\Phi:\cc\to\cd^\top\), \(A\mapsto(UA,U\eps:UFUA\to UA)\), for the monad \(\top=U\circ F\) in $\cd$ is an isomorphism of categories if and only if $U$ creates coequalizers of reflexive $U$\n-contractible coequalizer pairs in $\cc$.
\end{theorem}

Here $\cd^\top$ is the category of $\top$\n-algebras.
The condition of the theorem applied to \(f=\id=g\) implies that $U$ reflects isomorphisms.
The proof of this theorem is contained in the proof of (PTT), Beck's Precise Tripleability Theorem \cite[Theorem~3.3.14]{BarrWells:TopTT}.
We shall use the following corollary to \thmref{thm-PPTT}.

One says that (\textit{cf.} \cite[Section~3.5]{BarrWells:TopTT})
\begin{myitemize}
\item[(CTT$'$)] \(U:\cc\to\cd\) creates coequalizers for reflexive pairs $(f,g)$ for which $(Uf,Ug)$ has a coequalizer,
\end{myitemize}
if for any reflexive pair \(f,g:A\to B\in\cc\) and any coequalizer \(d:UB\to C'\) of $(Uf,Ug)$ in $\cd$ conclusions \dag{} and \ddag{} before \thmref{thm-PPTT} hold.

\begin{corollary}[Crude Tripleability Theorem]
Let \(U:\cc\to\cd\) be a functor which satisfies \textup{(CTT$'$)} and has a left adjoint $F$.
Then the comparison functor \(\Phi:\cc\to\cd^\top\), \(\top=U\circ F\) is an isomorphism of categories.
\end{corollary}

\begin{definition}
An \emph{ideal} of an object \(\tilde\cA=(\ca_1,\dots,\ca_n;\cp;\ca_0)\in\nOp n_1\) is a subobject \(\tilde\cI=(\cI_1,\dots,\cI_n;\ck;\cI_0)\) of \(U\tilde\cA\) in $\dg^S$,
\(S=n\NN\sqcup\NN^n\sqcup\NN=\NN\sqcup\dots\sqcup\NN\sqcup\NN^n\sqcup\NN\), 
stable under all multiplications in operads $\cA_j$, $0\le j\le n$, and under the action on $\cp$ from \eqref{eq-alpha-full-action}.
Namely if at least one $\tens$\n-argument of multiplication or action is in $\tilde\cI$, then the result is in $\tilde\cI$ as well.
Equivalently, for all values of indices
\begin{gather*}
\lambda^i_{k,(j_p)}\Biggl(
\biggl[ \bigotimes_{p=1}^{k^i} \cA_i(j_p)\biggr] \tens\ck(k) \Biggr) 
\subset \ck\biggl(k^1,\dots,k^{i-1},\sum_{p=1}^{k^i}j_p,k^{i+1},
\dots,k^n\biggr),
\\
\rho\biggl(\biggl[\Bigl(\bigotimes_{r=1}^{t-1}\cp(k_r)\Bigr) 
\tens\ck(k_t)\tens \Bigl(\bigotimes_{r=t+1}^m\cp(k_r)\Bigr) \biggr] 
\tens\cA_0(m) \biggr)
\subset \ck\biggl(\sum_{r=1}^mk_r\biggr),
\end{gather*}
and $\cI_j$ are ideals of operads $\cA_j$ in a similar sense.
\end{definition}

Ideals are precisely kernels in $\dg^S$ of $Uh$ for morphisms \(h:\tilde\cA\to\tilde\cB\in\nOp n_1\).
If $\tilde\cI$ is an ideal of $\tilde\cA$, then the quotient \(U\tilde\cA/U\tilde\cI\) in the abelian category $\dg^S$ admits a unique structure of an $n\wedge1$-operad module such that the quotient map \(q:\cA\to\tilde\cA/\tilde\cI\) is in \(\nOp n_1\).

For any subcomplex \(\cn\subset\tilde\cA\in\dg^S\) there is the smallest ideal $\tilde\cI$ of $\tilde\cA$ containing $\cn$.
It is spanned as a graded $\kk$\n-submodule of $\tilde\cA$ by results of multiplications or actions containing an element of $\cn$ among its $\tens$\n-arguments.
So obtained $\tilde\cI$ is indeed an ideal due to associativity of the action.
In particular, for a pair of parallel arrows \(f,g:\tilde\cA\to\tilde\cB\in\nOp n_1\) there is the image \(\cn=\im(f-g)\) in the abelian category $\dg^S$.
If $\tilde\cI$ is the smallest ideal of $\tilde\cA$ containing $\cn$, then the quotient $\tilde\cA/\tilde\cI$ is the coequalizer of $f$ and $g$ in $\nOp n_1$.

\begin{proposition}
The comparison functor for the underlying functor \(U:\nOp n_1\to\dg^{n\NN\sqcup\NN^n\sqcup\NN}\) is an isomorphism of categories.
\end{proposition}

Thus \(\nOp n_1\) is isomorphic to the category of $\top$\n-algebras for the monad \(\top=U\circ F\) in \(\dg^{n\NN\sqcup\NN^n\sqcup\NN}\).

\begin{proof}
Let us prove that $U$ satisfies condition (CTT$'$).
First (as a warm-up) we show it for \(U:\Op\to\dg^\NN\).
Let a reflexive pair \(f,g:\ca\leftrightarrows\cb:r\) in $\Op$ be given together with a coequalizer \(d:UB\to\cc'\) in $\dg^\NN$ of \((Uf,Ug)\).
The subobject \(\ck=\im(f-g)=\Ker d\in\dg^\NN\) of $UB$ is an ideal of the operad $\cb$.
In fact, for any multiplication $\mu$ from \eqref{eq-mu-OOOOOOOOO-O} for the operads $\cb$ and $\cA$ we have
\begin{align*}
&\mu^\cb(b_1\tdt b_{i-1}\tens(f-g)a\tens b_{i+1}\tdt b_k\tens b)
\\
&=\mu^\cb(frb_1\tdt frb_{i-1}\tens fa\tens frb_{i+1}\tdt frb_k\tens frb)
\\
&-\mu^\cb(grb_1\tdt grb_{i-1}\tens ga\tens grb_{i+1}\tdt grb_k\tens grb)
\\
&=(f-g)\mu^\cA(rb_1\tdt rb_{i-1}\tens a\tens rb_{i+1}\tdt rb_k\tens rb) \in\ck
\end{align*}
for all \(1\le i\le k\).
Similarly,
\[ \mu^\cb(b_1\tdt b_k\tens(f-g)a) =(f-g)\mu^\cA(rb_1\tdt rb_k\tens a) \in\ck.
\]
Thus the quotient operad $\cb/\ck$ exists together with the quotient map \(q:\cb\to\cb/\ck\in\Op\).
The map $d$ factorises as \(d=\bigl(U\cb\rTTo^{Uq} U(\cb/\ck)\rTTo^\phi_\sim \cc'\bigr)\) for a unique isomorphism \(\phi\in\dg^\NN\).
Transferring the operad structure from $\cb/\ck$ to $\cc'$ along $\phi$ we make the latter into an operad $\cc$, make $d$ into a morphism of operads.
Clearly, properties \dag{} and \ddag{} \vpageref{pg-dag-ddag} hold true.

Consider now a reflexive pair in $\nOp n_1$
\[ \tilde f=((f_i)_0^n;f), \tilde g=((g_i)_0^n;g): \tilde\cA=((\cA_i)_1^n;\cp;\cA_0) \leftrightarrows \tilde\cB=((\cB_i)_1^n;\cq;\cB_0) : \tilde r=((r_i)_0^n;r).
\]
Then \(\tilde\cI=((\cI_i)_1^n;\ck;\cI_0)=\im(U\tilde f-U\tilde g)\) is an ideal in \(\tilde\cB\).
In fact it suffices to take in the source of action map~\eqref{eq-alpha-full-action} one of the $\tens$\n-arguments equal to \((\tilde f-\tilde g)x\) for $x\in\cA_i$ or $\cp$ or $\cA_0$.
Then
\begin{align*}
&\alpha(\dots\tens b_i^p\tdt(\tilde f-\tilde g)x\tdt q_j\tdt b_0)
\\
&=\alpha(\dots\tens f_ir_ib_i^p\tdt\tilde fx\tdt frq_j\tdt f_0r_0b_0)
\\
&-\alpha(\dots\tens g_ir_ib_i^p\tdt\tilde gx\tdt grq_j\tdt g_0r_0b_0)
\\
&=(f-g)\alpha(\dots\tens r_ib_i^p\tdt x\tdt rq_j\tdt r_0b_0) \in\ck.
\end{align*}
Thus \(\tilde\cB/\tilde\cI\) is an $n\wedge1$-operad module and the rest of the proof goes similarly to the case of operads.
\end{proof}

\begin{corollary}\label{cor-category-nOp1-complete-cocomplete}
The category $\nOp n_1$ is complete and cocomplete.
\end{corollary}

\begin{proof}
In fact, \(\cd=\dg^S\) is complete and cocomplete.
Completeness of \(\nOp n_1\simeq\cd^\top\) follows by \cite[Corollary~3.4.3]{BarrWells:TopTT}.
We have seen that the category $\nOp n_1$ has coequalizers.
Therefore $\nOp n_1$ is cocomplete by \cite[Corollary~9.3.3]{BarrWells:TopTT}.
\end{proof}

\subsection{Peculiar features.}\label{sec-Peculiar-features}
The category \((\ca_1,\dots,\ca_n)\imodi\cb\) of operad 
\((\ca_1,\dots,\ca_n;\cb)\)-modules is a subcategory of \(\nOp n_1\), 
whose objects are \((\ca_1,\dots,\ca_n;\cp;\cb)\) (shortly $\cp$) and 
morphisms are \((1_{\ca_1},\dots,1_{\ca_n};h;1_\cb)\).
It has an initial object \(\ci=\cb(0)\) with
\[ \ci(k) =
\begin{cases}
\cb(0), &\text{ if } k=0
\\
0, &\text{ if } k\ne0
\end{cases}
\quad,\quad k\in\NN^n.
\]
Actions are given by
\begin{gather*}
\lambda^i_0 =\id: \cb(0) \to \cb(0),
\\
\rho =\mu^\cb_{0,\dots,0}: (\tens^{\mb m}\cb(0))\tens\cb(m) \to \cb(0).
\end{gather*}

A \emph{system of relations} in an operad
\((\ca_1,\dots,\ca_n;\cb)\)-module $\cp$ means a set (of relations) $R$, 
arity function \(a:R\to\NN^n\), grade function \(g:R\to\ZZ\) and for 
each $r\in R$ elements \(x_r,y_r\in\cp(a(r))^{g(r)}\) supposed to be 
identified by the relation \(x_r=y_r\).
A system of relations gives rise to a free graded $\kk$\n-module 
\(\kk R\in\gr^{\NN^n}\) with 
\[ \kk R(a)^g =\kk\{ r\in R \mid a(r)=a,\ g(r)=g \},
\quad \text{for } a\in\NN^n,\; g\in\ZZ,
\]
and to a map \(\kk R\to\cp\), \(r\mapsto x_r-y_r\).
Denote by $\cn$ the image of this map in abelian category $\gr^{\NN^n}$.
Let $\ck=(0,\dots,0;\ck;0)$ be the graded ideal of \((\ca_1,\dots,\ca_n;\cp;\cb)\) generated by $\cn$.
If \(\cn\partial\subset\ck\), then $\ck$ is a differential graded ideal 
and the quotient $\cp/\ck$ in $\dg^{\NN^n}$ is an operad 
\((\ca_1,\dots,\ca_n;\cb)\)-module.
This is the quotient of an operad module by a system of relations.

\subsubsection{Induced modules.}\label{sec-Induced-modules}
Assume that \(f_i:\ca_i\to\cb_i\) are morphisms of operads for 
$i\in[n]$.
For any $n\wedge1$-operad module \(((\ca_i)_{i\in[n]};\cp)\) there is a 
\((\cb_i)_{i\in[n]}\)-module 
\(\cq=\bigcirc_{i=0}^n\cb_i\odot^i_{\ca_i}\cp\), which is the colimit of 
the diagram (the coequalizer)
\begin{diagram}
\odot_{\ge0}((\cb_i)_{i\in[n]};\odot_{\ge0}((\ca_i)_{i\in[n]};\cp)) 
&\rTTo^{\odot_{\ge0}(\sS{^{[n]}}1;\alpha)}
&\odot_{\ge0}((\cb_i)_{i\in[n]};\cp)
\\
\dTTo<{\odot_{\ge0}(\sS{^{[n]}}1;\odot_{\ge0}((f_i)_{i\in[n]};1))} 
&\ruTTo^\mu_= &\uTTo>{\odot_{\ge0}(\sS{^{[n]}}\mu;1)}
\\
\odot_{\ge0}((\cb_i)_{i\in[n]};\odot_{\ge0}((\cb_i)_{i\in[n]};\cp)) 
&\rTTo^\sim &\odot_{\ge0}((\cb_i\odot\cb_i)_{i\in[n]};\cp)
\end{diagram}
This can be described also using monads 
\(\aA:\cx\mapsto\odot_{\ge0}((\ca_i)_{i\in[n]};\cx)\), 
\(\BB:\cx\mapsto\odot_{\ge0}((\cb_i)_{i\in[n]};\cx)\) in \(\dg^{\NN^n}\) 
and a morphism \(f=\odot_{\ge0}((f_i)_{i\in[n]};1):\aA\to\BB\) between 
them.
Denote \(\alpha^\cp:\aA\cp\to\cp\) the action on the $\aA$\n-module 
$\cp$.
The induced module $\cq$ is the colimit \(\BB_\aA\cp\) of the diagram in 
the category of $\BB$\n-modules (the coequalizer of)
\begin{diagram}[h=1.5em]
\BB\aA\cp &&\rTTo^{\BB\alpha^\cp} &&\BB\cp
\\
&\rdTTo_{\BB f} &&\ruTTo_{\mu^\BB} &
\\
&&\BB\BB\cp
\end{diagram}

If \(\cp=\odot_{\ge0}(\ca_1,\dots,\ca_n;\cx;\ca_0)\) is a free 
$n\wedge1$-module, then the $n\wedge1$-module
	\(\bigcirc_{i=1}^n\cb_i\odot^i_{\ca_i}\cp\odot^0_{\ca_0}\cb_0
	=\odot_{\ge0}(\cb_1,\dots,\cb_n;\cx;\cb_0)\)
is free as well.
It is also generated by $\cx$.
In fact, the coequalizer in question has the structure of a contractible 
coequalizer \cite[Definition~3.3.3]{BarrWells:TopTT} with 
\(d^1=\BB f\cdot\mu^\BB:\BB\aA\aA\cx\to\BB\aA\cx\) and 
\(d=\BB f\cdot\mu^\BB:\BB\aA\cx\to\BB\cx\)
\begin{diagram}
\BB\aA\aA\cx &&\pile{\rTTo^{d^0=\BB\mu^\aA} \\ \lTTo_{t=\BB\aA\eta}} 
&&\BB\aA\cx &\rTTo^{\BB f} &\BB\BB\cx
\\
&\rdTTo_{\BB f} &&\ruTTo_{\mu^\BB} && \luTTo_{s=\BB\eta} 
&\dTTo>{\mu^\BB}
\\
&&\BB\BB\aA\cx &&&&\BB\cx
\end{diagram}

The following is a part of the theory of modules in general categories, 
not necessarily linear.

\begin{proposition}\label{pro-left-adjoint-to-base-change}
So defined functor 
\((\ca_i)_{i\in[n]}\modul\to(\cb_i)_{i\in[n]}\modul\), 
\(\cp\mapsto\bigcirc_{i=0}^n\cb_i\odot^i_{\ca_i}\cp\), is left adjoint 
to the base change functor 
\((\cb_i)_{i\in[n]}\modul\to(\ca_i)_{i\in[n]}\modul\), 
\(\cR\mapsto\sS{_{f_1,\dots,f_n}}\cR_{f_0}\).
Thus, the mapping
\[ (\cb_i)_{i\in[n]}\modul(\bigcirc_{i=0}^n\cb_i\odot^i_{\ca_i}\cp,\cR) 
\rto\sim (\ca_i)_{i\in[n]}\modul(\cp,\sS{_{f_1,\dots,f_n}}\cR_{f_0}), 
\qquad h \mapsto \eta\cdot h,
\]
is a bijection, where
\(\eta:\cp\to\bigcirc_{i=0}^n\cb_i\odot^i_{\ca_i}\cp\) comes from the 
unit of the monad \(\odot_{\ge0}((\cb_i);-)\).
\end{proposition}

\begin{proof}
Since \(\cq=\bigcirc_{i=0}^n\cb_i\odot^i_{\ca_i}\cp\) is a colimit
(a coequalizer), the following set of morphisms is a limit (an equalizer 
of a pair of maps):
\begin{multline*}
(\cb_i)_{i\in[n]}\modul(\cq,\cR)
\\
=\lim\left[ 
\begin{diagram}[inline,nobalance]
(\cb_i)\modul(\odot_{\ge0}((\cb_i);\cp),\cR) 
&\rTTo^{\hspace*{-1em}(\odot_{\ge0}((1);\alpha),1)} 
&(\cb_i)\modul(\odot_{\ge0}((\cb_i);\odot_{\ge0}((\ca_i);\cp)),\cR)
\\
\dTTo~{(\cb_i)\modul(\odot_{\ge0}((\mu_i);1),1)} 
&&\uTTo~{(\cb_i)\modul(\odot_{\ge0}((1);\odot_{\ge0}((f_i);1)),1)}
\\
(\cb_i)\modul(\odot_{\ge0}((\cb_i\odot\cb_i);\cp),\cR) &\rTTo^\sim 
&(\cb_i)\modul(\odot_{\ge0}((\cb_i);\odot_{\ge0}((\cb_i);\cp)),\cR)
\end{diagram}
\right] 
\\
=\lim\left[ 
\begin{diagram}[inline,nobalance]
\dg^{\NN^n}(\cp,\cR) &\rTTo^{\dg^{\NN^n}(\alpha,1)} 
&\dg^{\NN^n}(\odot_{\ge0}((\ca_i);\cp),\cR)
\\
\dTTo<{\odot_{\ge0}((\cb_i);-)}
&&\uTTo~{\dg^{\NN^n}(\odot_{\ge0}((f_i);1),1)}
\\
\dg^{\NN^n}(\odot_{\ge0}((\cb_i);\cp),\odot_{\ge0}((\cb_i);\cR))
&\rTTo^{\dg^{\NN^n}(1,\alpha)}
&\dg^{\NN^n}(\odot_{\ge0}((\cb_i);\cp),\cR)
\end{diagram}
\right].
\end{multline*}
Thus, a morphism \(u:\cp\to\cR\in\dg^{\NN^n}\) belongs to the equalizer 
iff
\begin{multline*}
\bigl[ \odot_{\ge0}((\ca_i)_{i\in[n]};\cp) \rTTo^{\odot_{\ge0}((f_i);1)} 
\odot_{\ge0}((\cb_i)_{i\in[n]};\cp) \rTTo^{\odot_{\ge0}((1);u)} 
\odot_{\ge0}((\cb_i)_{i\in[n]};\cR) \rTTo^\alpha \cR \bigr]
\\
=\bigl[\odot_{\ge0}((\ca_i)_{i\in[n]};\cp) \rTTo^\alpha \cp \rTTo^u 
\cR\bigr].
\end{multline*}
Equivalently, the mapping \(u:\cp\to\sS{_{f_1,\dots,f_n}}\cR_{f_0}\) is 
a homomorphism of \((\ca_i)_{i\in[n]}\)-modules.
\end{proof}

Together with \lemref{lem-((Ai)A(0))-initial} this proposition implies 
the following

\begin{corollary}\label{cor-AA0-AP-BB0-QO}
In assumptions of \secref{sec-Induced-modules} denote 
\(\cq=\bigcirc_{i=0}^n\cb_i\odot^i_{\ca_i}\cp\). 
Then the diagram in the category \(\nOp n_1\)
\begin{diagram}[w=7em]
(\ca_1,\dots,\ca_n;\ca_0(0);\ca_0) &\rTTo^{(1,\dots,1;\rho_\emptyset;1)} 
&(\ca_1,\dots,\ca_n;\cp;\ca_0)
\\
\dTTo<{(f_1,\dots,f_n;f_0(0);f_0)} &&\dTTo>{(f_1,\dots,f_n;\eta;f_0)}
\\
(\cb_1,\dots,\cb_n;\cb_0(0);\cb_0) &\rTTo^{(1,\dots,1;\rho_\emptyset;1)} 
&(\NWpbk\cb_1,\dots,\cb_n;\cq;\cb_0)
\end{diagram}
is a pushout square.
\end{corollary}

\subsection{The operad \texorpdfstring{$T\cn\sqcup\co$}{TNuO}.}\label{sec-operad-TNuO}
As explained in \cite[Proposition~1.8]{Lyu-Ainf-Operad} the operad $T\cn\sqcup\co$ is a direct sum over ordered rooted trees $t$ with inputs whose vertices are coloured with $\cn$ and $\co$ that are terminal the following sense.
Two conditions hold:
\begin{itemize}
\item[$\lightning$)] $t$ contains no edge whose both ends are coloured with $\co$;

\item[$\lightning\lightning$)] insertion of a unary vertex coloured with $\co$ into an arbitrary edge of $t$ breaks condition $\lightning$).
\end{itemize}
By the way, these conditions imply that $t$ contains no edge whose both ends are coloured with $\cn$.
The second condition is related to inserting a unit \(\eta:\kk\to T\cn\sqcup\co\) identified with the unit \(\eta:\kk\to\co\) in an arbitrary summand of $T\cn\sqcup\co$ represented by a coloured tree.
The first condition means that one can not apply binary composition in $\co$ inside an element of
\[ C(\cn,\co;t) \simeq \tens^{v\in(\IV(t),\le)}c(v)(|v|) \subset C(\cn\oplus\co;t) \subset T(\cn\oplus\co)
\]
represented by $t$.
Here $\le$ is an admissible order on $\IV(t)$, the set of internal vertices of an ordered rooted trees $t$ with inputs $\Inp t$, see Section~1.5 of \cite{Lyu-Ainf-Operad}.
The set $\IV(t)$ is coloured by the function \(c:\IV(t)\to\{\cn,\co\}\), which singles out an individual summand of \((\cn\oplus\co)^{\tens\IV(t)}\).

Any sequence of contractions of edges whose ends are coloured with $\co$ and insertions of unary vertices coloured with $\co$ applied to given tree $t'$ may lead to no more than one terminal tree $t$.
The mapping \(t'\mapsto t\) is well defined.
The summand \(C(\cn,\co;t')\) of \(T(\cn\oplus\co)\) corresponding to $t'$ is mapped by binary compositions in $\cb$ and insertions of the unit of $\cb$ to the summand \(C(\cn,\co;t)\).
Associativity and unitality of $\co$ implies that this map is unique.

The requirement of \(\co\to C=T\cn\sqcup\co\) being a morphism of operads (see \corref{cor-biggest-quotient-morphism-T-algebras}) reduces to agreeing with binary compositions and units.
Therefore, in the case of operads equation~\eqref{dia-UFUA-UF(XuUA)} is equivalent to a family of similar diagrams with \(\alpha:T\co\to\co\) replaced with the unit \(1_\co:\kk\to\co\) and with the binary compositions \(\co\odot\co\supset\kk\tdt\kk\tens\co\tens\kk\tdt\kk\tens\co\rto\mu \co\).
That is, the kernel $\cI$ of \(T(\cn\oplus\co)\to T\cn\sqcup\co\) is generated as an ideal by relations coming from binary composition in $\co$ (corresponding to contraction of an $\co$\n-coloured edge) and from inserting a unit of $\co$ (corresponding to insertion of $\co$\n-coloured unary vertex into an edge).
Together with the above this implies that an element of \(C(\cn,\co;t')\) is equivalent to a unique element of \(C(\cn,\co;t)\) modulo $\cI$.
In fact, for elements \(x'\in C(\cn,\co;t')\) and \(x''\in C(\cn,\co;t'')\) related by an elementary relation as above this holds true since there is either a path \((t',t'',\dots,t)\) or a path \((t'',t',\dots,t)\) consisting of contractions or unary insertions.
We conclude that \cite[eq.~(1.2)]{Lyu-Ainf-Operad}
\begin{equation*}
T\cn\sqcup\co =\bigoplus_{t\in\cl} C(\cn,\co;t),
%\label{eq-TNO=C(NOt)}
\end{equation*}
where $\cl$ is the list of terminal trees.

\begin{proposition}
The category $\nOp n_1$ of $n\wedge1$-operad modules with
quasi-isomorphisms as weak equivalences and degreewise surjections as fibrations is a model category.
\end{proposition}

\begin{proof}
Let us apply \thmref{thm-Hinich-model} to the adjunction \(F:\dg^S\leftrightarrows\nOp n_1:U\), where \(S=n\NN\sqcup\NN^n\sqcup\NN\).
By \corref{cor-category-nOp1-complete-cocomplete} the category $\nOp n_1$ is complete and cocomplete.

Let \(\cn\in\Ob\dg^\NN\) or \(\cn\in\Ob\dg^{\NN^n}\) be a complex.
For instance, $\cn=\KK_x[-p]$ for $p\in\ZZ$, $x\in\NN$ or $x\in\NN^n$.
Let \((\ca_1,\dots,\ca_n;\cp;\cb)\) be an $n\wedge1$-operad module.
Denote by $\tilde\cn$ the object \((0,\dots,0;0;\cn)\) or \((0,\dots,0;\cn;0)\) or \((0,\dots,0,\cn,0,\dots,0;0;0)\) ($\cn$ on $i^{\text{th}}$ place) of $\dg^S$.
We shall prove case by case that if $\cn$ is contractible then so is \(U(F\tilde\cn\sqcup(\ca_1,\dots,\ca_n;\cp;\cb))\).

By \corref{cor-biggest-quotient-morphism-T-algebras} the operad module in question \(F\tilde\cn\sqcup(\ca_1,\dots,\ca_n;\cp;\cb)\) is the quotient of \(F(\tilde\cn\oplus(\ca_1,\dots,\ca_n;\cp;\cb))\) by the smallest ideal $\cI$ generated by relations which tell that \((\ca_1,\dots,\ca_n;\cp;\cb)\to F(\tilde\cn\oplus(\ca_1,\dots,\ca_n;\cp;\cb))/\cI\) is a morphism of $n\wedge1$-operad modules.
Notice by the way that
\begin{align*}
F(0,\dots,0;0;\cn)\sqcup(\ca_1,\dots,\ca_n;\cp;\cb) &=(\ca_1,\dots,\ca_n;\cq;T\cn\sqcup\cb),
\\
\cq &=\cp\odot^0_\cb(T\cn\sqcup\cb),
\\
F(0,\dots,0;\cn;0)\sqcup(\ca_1,\dots,\ca_n;\cp;\cb) &=(\ca_1,\dots,\ca_n;\cR;\cb),
\\
\cR &=\cp\oplus\odot_{\ge0}(\ca_1,\dots,\ca_n;\cn;\cb),
\\
F(0,\dots,0,\cn,0,\dots,0;0;0)\sqcup(\ca_1,\dots,\ca_n;\cp;\cb) &=(\ca_1,\dots,T\cn\sqcup\ca_i,\dots,\ca_n;\cs;\cb),
\\
\cs &=(T\cn\sqcup\ca_i)\odot^i_{\ca_i}\cp.
\end{align*}
More important is the presentation of operad polymodules as direct sums over some kind of trees.
This presentation we use for the proof.
Similarly to the operad case considered in \secref{sec-operad-TNuO} we find that
\begin{equation*}
\cq =\oplus_{\tau\in\cl_\cq} \cq(\cp,\cn,\cb;\tau), \qquad \text{where} \qquad \cq(\cp,\cn,\cb;\tau) =\colim_{\le\in\cg_\cq(\tau)} \cq(\cp,\cn,\cb;\tau)(\le).
\end{equation*}
Here objects of the groupoid \(\cg_\cq(\tau)\) are admissible (compatible with the natural partial order $\preccurlyeq$ with the root as the biggest element) total orderings $\le$ of the set of all vertices $\iV(\tau)$.
By definition between any two objects of $\cg(\tau)$ there is precisely one morphism.
Thus \(\cg(\tau)\) is contractible (equivalent to the terminal category with one object and one morphism). 
Therefore the colimit is isomorphic to any of
\[ \cq(\cp,\cn,\cb;\tau)(\le) =\tens^{v\in(\iV(\tau),\le)}c(v)(|v|).
\]
Here \(\tau=(\tau,c,|\cdot|)\) is an ordered rooted tree $\tau$ with inputs $\Inp\tau$, a colouring \(c:\iV(\tau)\to\{\cp,\cn,\cb\}\) such that \(c(\Inp\tau)\subset\{\cp\}\), \(c(\IV(\tau))\subset\{\cn,\cb\}\) and an arbitrary function \(|\cdot|:\Inp(\tau)\to\NN^n\), which complements the valency \(|\cdot|:\IV(\tau)\to\NN\).
The set $\cl_\cq$ of terminal trees consists of $\tau$ such that
\begin{itemize}
\item[$*$)] $\tau$ contains no edge whose both ends are coloured with $\cb$;

\item[$**$)] $\tau$ contains no vertex coloured with $\cb$ whose all entering edges have other ends coloured with $\cp$;

\item[$*{**}$)] insertion of a unary vertex coloured with $\cb$ into an arbitrary edge of $\tau$ breaks condition $*$) or $**$).
\end{itemize}

Respectively,
\begin{equation*}
\cs =\oplus_{\tau\in\cl_\cs} \cs(\cn,\ca_i,\cp;\tau), \qquad \text{where} \qquad \cs(\cn,\ca_i,\cp;\tau) =\colim_{\le\in\cg_\cs(\tau)} \cs(\cn,\ca_i,\cp;\tau)(\le).
\end{equation*}
The colimit over the contractible groupoid \(\cg_\cs(\tau)\) is isomorphic to expression under colimit in any vertex $\le$.
Assuming that $\ell\in\NN^n$, $\ell^i=|\Inp\tau|$, we have
\[ \cs(\cn,\ca_i,\cp;\tau)(\le)(\ell) =\bigl[\tens^{v\in(\IV(\tau)-\{\troot\},\le)}c(v)(|v|)\bigr] \tens\cp(\ell,\ell^i\mapsto q).
\]
Here \(\tau=(\tau,c)\) is an ordered rooted tree $\tau$ with inputs $\Inp\tau$ and a colouring \(c:\IV(\tau)-\{\troot\}\to\{\cn,\ca_i\}\).
The set $\cl_\cs$ of terminal trees consists of $\tau$ such that
\begin{itemize}
\item[$*$)] $\tau$ contains no edge whose both ends are coloured with $\cb$;

\item[$**$)] $\tau$ contains no edge adjacent to the root whose one end is coloured with $\cb$;

\item[$*{**}$)] insertion of a unary vertex coloured with $\cb$ into an arbitrary edge of $\tau$ breaks condition $*$) or $**$).
\end{itemize}

Let us prove existence of contracting homotopy similarly to the case of operads \cite[Proposition~1.8]{Lyu-Ainf-Operad}.
Let \(\cn\in\Ob\dg^\NN\) or \(\cn\in\Ob\dg^{\NN^n}\) be contractible and let \(h:\cn\to\cn\) be a contracting homotopy, \(\deg h=-1\), \(dh+hd=1_\cn\).
Let us show that the operad module morphism \(\alpha=\inj_2:\cm=(\ca_1,\dots,\ca_n;\cp;\cb)\to F\tilde\cn\sqcup(\ca_1,\dots,\ca_n;\cp;\cb)\) is homotopy invertible.
Consider the operad module morphism \(\beta:F\tilde\cn\sqcup(\ca_1,\dots,\ca_n;\cp;\cb)\to(\ca_1,\dots,\ca_n;\cp;\cb)\) which restricts to \(\beta\big|_{(\ca_1,\dots,\ca_n;\cp;\cb)}=\id\) and \(\beta\big|_{F\tilde\cn}\), adjunct to \(0:\tilde\cn\to U(\ca_1,\dots,\ca_n;\cp;\cb)\).
Then \(\alpha\cdot\beta=\id\) and \(g=\beta\cdot\alpha\) is homotopic to \(f=\id_{F\tilde\cn\sqcup(\ca_1,\dots,\ca_n;\cp;\cb)}\) in the $\dg$\n-category $\dg^S$, as we show next.
The homotopy $h$ is extended by $0$ to the map \(h'=h\oplus0:\tilde\cn\oplus U(\ca_1,\dots,\ca_n;\cp;\cb)\to\tilde\cn\oplus U(\ca_1,\dots,\ca_n;\cp;\cb)\), which satisfies \(dh'+h'd=f|-g|:\tilde\cn\oplus U(\ca_1,\dots,\ca_n;\cp;\cb)\to\tilde\cn\oplus U(\ca_1,\dots,\ca_n;\cp;\cb)\).
In all three cases the endomorphisms $f$, $g$ of \(F\tilde\cn\sqcup(\ca_1,\dots,\ca_n;\cp;\cb)\) lift to endomorphisms of \(F(\tilde\cn\oplus(\ca_1,\dots,\ca_n;\cp;\cb))\) obtained by applying \(f\big|_{\tilde\cn}=1:\tilde\cn\to\tilde\cn\), \(f\big|_\cm=1:\cm\to\cm\), \(g\big|_{\tilde\cn}=0:\tilde\cn\to\tilde\cn\), \(g\big|_\cm=1:\cm\to\cm\) to each $\tens$\n-factor corresponding to a vertex of the tree.
For an arbitrary pair of trees $(t,\tau)$ choose admissible orderings \((\le,\le)\).
Then the summands of \(F\tilde\cn\sqcup(\ca_1,\dots,\ca_n;\cp;\cb)\) are preserved by $f$ and $g$ and the restriction to the summand is $f\tdt f$ and $g\tdt g$ respectively.
Define a $\kk$\n-endomorphism \(\hat{h}=\sum_{v\in(\IV(t),\le)}f\tdt f\tens h'\tens g\tdt g\) of degree $-1$, where $h'$ is applied on place indexed by $v$.
Then
\[ d\hat{h} +\hat{h}d =\sum_{v\in(\IV(t),\le)} f\tdt f\tens(f-g)\tens g\tdt g =f\tdt f -g\tdt g =f-g.
\]
Therefore, $f$ and $g$ are homotopic to each other and $\alpha$ is homotopy invertible.
\end{proof}

\propref{pro-left-adjoint-to-base-change} gives a recipe of computing 
colimits in the category $\nOp n_1$ of $n\wedge1$-operad modules in two 
steps.
First of all compute colimits $\cb_i$ on each of $n+1$ operadic places 
$i\in[n]$.
Take induced module over \((\cb_1,\dots,\cb_n;\cb_0)\) on each node of 
the diagram and consider the obtained diagram in the category
\((\cb_1,\dots,\cb_n)\imodi\cb_0\).
Secondly find the colimit of the latter diagram, by finding its colimit 
in \(\dg^{\NN^n}\), then generating by it a free 
\((\cb_1,\dots,\cb_n;\cb_0)\)-module $F$, dividing it precisely by such 
relations that canonical mapping from any module to the quotient of $F$ 
were a morphism of \((\cb_1,\dots,\cb_n;\cb_0)\)-modules.

\subsection{Some colimits of operad modules.}
\begin{lemma}\label{lem-((Ai)A(0))-initial}
Let \(f_i:\ca_i\to\cb_i\) be morphisms of operads for $i\in[n]$.
Let $\cp$ be a \((\cb_i)_{i\in[n]}\)-module. 
Then there is a unique morphism $\phi$ such that 
\[ ((f_i)_{i\in[n]};\phi):
((\ca_i)_{i\in[n]};\ca_0(0))\to((\cb_i)_{i\in[n]};\cp)\in\nOp n_1.
\]
\end{lemma}

\begin{proof}
The morphism $\phi$ is recovered from the equation
\[ \bigl( \ca_0(0) \rTTo^{\rho_\emptyset^{\ca_0(0)}}_1 \ca_0(0) 
\rTTo^\phi \cp(0) \bigr)
=\bigl(\ca_0(0) \rTTo^{f_0(0)} \cb_0(0) \rTTo^{\rho_\emptyset^\cp} 
\cp(0)\bigr)
\]
in the unique possible form \(\phi=f_0(0)\cdot\rho_\emptyset^\cp\).
It is compatible with the action $\rho$ because the diagram
\begin{diagram}
\ca_0(0)^{\tens m}\tens\ca_0(m) &\rTTo^{f_0(0)^{\tens m}\tens f_0(m)} 
&\cb_0(0)^{\tens m}\tens\cb_0(m) &\rTTo^{\rho_\emptyset^{\tens m}\tens1}
&\cp_0(0)^{\tens m}\tens\cb_0(m)
\\
\dTTo<{\mu^{\ca_0}} &= &\dTTo<{\mu^{\cb_0}} &= 
&\dTTo>{\rho^m_{\sS{^m}0}}
\\
\ca_0(0) &\rTTo^{f_0(0)} &\cb_0(0) &\rTTo^{\rho_\emptyset} &\cp_0(0)
\end{diagram}
commutes.
For the second square this follows by associativity of the action.
\end{proof}

\begin{corollary}\label{cor-(Bi)U(Ai)-(BiUAi)}
For arbitrary operads $\cc_i$, $\ca_i$, $i\in[n]$, there is an 
isomorphism in $\nOp n_1$
\[ ((\cc_i)_{i\in[n]};\cc_0(0))\sqcup((\ca_i)_{i\in[n]};\ca_0(0)) \simeq
((\cc_i\sqcup\ca_i)_{i\in[n]};(\cc_0\sqcup\ca_0)(0)).
\]
\end{corollary}

\begin{proposition}\label{pro-A<Mbeta>}
Let \(A=(\ca_1,\dots,\ca_n;\cp;\ca_0)\) be an $n\wedge1$-operad module.
Let \(\beta_i:\cm_i\to\ca_i\in\dg^\NN\) be chain maps for $i\in[n]$.
Denote \(M=(\cm_1,\dots,\cm_n;0;\cm_0)\) and 
 \(\beta=(\beta_1,\dots,\beta_n;0;\beta_0):(\cm_1,\dots,\cm_n;0;\cm_0)
 \to(\ca_1,\dots,\ca_n;\cp;\ca_0)\).
Then \(A\langle M,\beta\rangle\) defined in 
diagram~\eqref{dia-pushout-A<M-alpha>} is isomorphic to
 \(((\ca_i\langle\cM_i,\beta_i\rangle)_{i\in[n]};
 \bigcirc_{i=0}^n\ca_i\langle\cM_i,\beta_i\rangle\odot^i_{\ca_i}\cp)\).
\end{proposition}

\begin{proof}
Denote
	\(\cR=\bigcirc_{i=0}^n\ca_i\langle\cM_i,
	\beta_i\rangle\odot^i_{\ca_i}\cp\) and \(\cc_i=T(\cm_i[1])\).
As we know from \cite[Section~1.10]{Lyu-Ainf-Operad} in
$\CG=\nOp n_1^{\gr}$ the graded module underlying
\(A\langle M,\beta\rangle\) is isomorphic to
\(((\cc_i)_{i\in[n]};\cc_0(0))\sqcup((\ca_i)_{i\in[n]};\cp)\).
Clearly, this coproduct is also a colimit of the following diagram (a 
pushout)
\begin{diagram}[nobalance]
&((\ca_i)_{i\in[n]};\ca_0(0)) &\rTTo^{((1);!)} &((\ca_i)_{i\in[n]};\cp)
\\
&\dTTo<{\inj_2} &&\dTTo>{\inj_2}
\\
((\cc_i\sqcup\ca_i)_{i\in[n]};(\cc_0\sqcup\ca_0)(0)) =\,
&((\cc_i)_{i\in[n]};\cc_0(0))\sqcup((\ca_i)_{i\in[n]};\ca_0(0)) 
&\rTTo^{1\sqcup((1);!)} &A\NWpbk\langle M,\beta\rangle
\end{diagram}
Here the equation is due to \corref{cor-(Bi)U(Ai)-(BiUAi)}.
Denote \(\cb_i=\cc_i\sqcup\ca_i\simeq\ca_i\langle\cM_i,\beta_i\rangle\) 
in $\Op^{\gr}$.
Applying \corref{cor-AA0-AP-BB0-QO} to the canonical embeddings 
\(\inj_2:\ca_i\to\cb_i\) we deduce an isomorphism
	\(\psi_3:A\langle M,\beta\rangle\to
	((\ca_i\langle\cM_i,\beta_i\rangle)_{i\in[n]};
	\bigcirc_{i=0}^n\ca_i\langle\cM_i,
	\beta_i\rangle\odot^i_{\ca_i}\cp)\)
in $\nOp n_1^{\gr}$.

In order to show that the isomorphism is actually in $\nOp n_1^{\dg}$ we 
consider diagram~\vpageref{dia-4x4-wheels}, where 
\(\phi_k=\psi_k^{-1}\), \(1\le k\le3\).
\begin{figure}
\begin{center}
%\boldmath
\resizebox{!}{\texthigh}{\rotatebox{90}{%
\begin{diagram}[height=2.6em,inline,nobalance]
((T\ca_i)_{i\in[n]};\odot_{\ge0}((T\ca_i)_{i\in[n]};\cp)) 
&\rTTo^{((\alpha_i);\odot_{\ge0}((\alpha_i);1))}
&((\ca_i)_{i\in[n]};\odot_{\ge0}((\ca_i)_{i\in[n]};\cp))
&\rTTo^{((1)_{i\in[n]};\alpha)} &((\ca_i)_{i\in[n]};\cp)
\\
\dTTo<{((T\bar\imath_i);\odot_{\ge0}((T\bar\imath_i);1))} &\ovalbox{1} 
&\dTTo<{((\bar\jmath_i)_{i\in[n]};
	\odot_{\ge0}((\bar\jmath_i)_{i\in[n]};1))}
&\ovalbox{2} &\dTTo>{((\bar\jmath_i)_{i\in[n]};\eta)}
\\
((TC_i)_{i\in[n]};\odot_{\ge0}((TC_i)_{i\in[n]};\cp)) 
&\rTTo^{((g_i);\odot_{\ge0}((g_i);1))}
&((\ca_i\langle\cM_i,\beta_i\rangle);
	\odot_{\ge0}((\ca_i\langle\cM_i,\beta_i\rangle);\cp)) 
&\rTTo^{((1)_{i\in[n]};z)} &((\ca_i\langle\cM_i,\beta_i\rangle);\cR)
\\
\dTTo<{\phi_1}>\wr \uTTo>{\psi_1} &\ovalbox{3} 
&\dTTo<{\phi_2}>\wr \uTTo>{\psi_2} &\ovalbox{4} 
&\dTTo<{\phi_3}>\wr \uTTo>{\psi_3}
\\
\begin{array}{l}
((\cc_i)_{i\in[n]};\cc_0(0))\sqcup
\\[2pt]
((T\ca_i)_{i\in[n]};\odot_{\ge0}((T\ca_i)_{i\in[n]};\cp))
\end{array}
&\rTTo^{1\sqcup((\alpha_i);\odot_{\ge0}((\alpha_i);1))} &
\begin{array}{l}
((\cc_i)_{i\in[n]};\cc_0(0))\sqcup
\\[2pt]
((\ca_i)_{i\in[n]};\odot_{\ge0}((\ca_i)_{i\in[n]};\cp))
\end{array}
&\rTTo^{1\sqcup((1);\alpha)} &
\begin{array}{l}
((\cc_i)_{i\in[n]};\cc_0(0))
\\[2pt]
\sqcup((\ca_i)_{i\in[n]};\cp)
\end{array}
\end{diagram}
}}
\end{center}
%\caption{}
\label{dia-4x4-wheels}
\end{figure}
Notice that the isomorphism $\psi_2$ is nothing else but the isomorphism 
$\psi_3$, written for the operad module \(((\ca_i)_{i\in[n]};\cp)\) 
instead of $\cp$, taking into account that
\[ 
\bigcirc_{i=0}^n\ca_i\langle\cM_i,\beta_i\rangle\odot^i_{\ca_i}
\odot_{\ge0}((\ca_i)_{i\in[n]};\cp))
\simeq \odot_{\ge0}((\ca_i\langle\cM_i,\beta_i\rangle)_{i\in[n]};\cp),
\]
see \secref{sec-Induced-modules}.
The isomorphism $\psi_1$ follows from the fact that 
$\bar F:\nOp n_1^{\gr}\to\gr^{n\NN\sqcup\NN^n\sqcup\NN}$ preserves 
colimits. 
The map
\(z:\odot_{\ge0}((\ca_i\langle\cM_i,\beta_i\rangle)_{i\in[n]};\cp)\to
	\bigcirc_{i=0}^n\ca_i\langle\cM_i,\beta_i\rangle\odot^i_{\ca_i}\cp\)
is the canonical projection.

We claim that the diagram commutes.
Its two top squares are in $\nOp n_1^{\dg}$, while the bottom vertical 
isomorphisms are constructed only in $\CG=\nOp n_1^{\gr}$.
Thus, squares \ovalbox{3} and \ovalbox{4} make sense in $\CG$.
First of all, square~\ovalbox{1} commutes due to definition
\eqref{dia-pushout-A<M-alpha>} of \(\ca_i\langle\cM_i,\beta_i\rangle\).
Commutativity of square~\ovalbox{2} follows from the equation
\begin{diagram}[LaTeXeqno]
\odot_{\ge0}((\ca_i)_{i\in[n]};\cp) &\rTTo^\alpha &\cp
\\
\dTTo<{\odot_{\ge0}((\bar\jmath_i)_{i\in[n]};1)} &= &\dTTo>\eta
\\
\odot_{\ge0}((\ca_i\langle\cM_i,\beta_i\rangle)_{i\in[n]};\cp) &\rTTo^z 
&\bigcirc_{i=0}^n\ca_i\langle\cM_i,\beta_i\rangle\odot^i_{\ca_i}\cp
\label{dia-alpha-eta-j1z}
\end{diagram}
proven directly from the definition of the induced operad module.
Namely paths in the following diagram
\begin{diagram}[w=2em,nobalance]
\odot_{\ge0}((\cb_i)_{i\in[n]};\odot_{\ge0}((\ca_i)_{i\in[n]};\cp)) 
&\lTTo^{\odot_{\ge0}((\eta);1)} &\odot_{\ge0}((\ca_i)_{i\in[n]};\cp) 
&\rTTo^\alpha &\cp &\0
\\
\dTTo<{\odot_{\ge0}((1);\odot_{\ge0}((\bar\jmath_i);1))} 
&\rdTTo[hug]_{\odot_{\ge0}((1);\alpha)}
&\dTTo<{\odot_{\ge0}((\bar\jmath_i);1)} 
&\ldTTo[hug]_{\odot_{\ge0}((\eta);1)} &\dTTo>\eta
\\
\odot_{\ge0}((\cb_i)_{i\in[n]};\odot_{\ge0}((\cb_i)_{i\in[n]};\cp)) 
&\rTTo_\mu &\odot_{\ge0}((\cb_i)_{i\in[n]};\cp) &\rTTo_z 
&\bigcirc_{i=0}^n\cb_i\odot^i_{\ca_i}\cp
\end{diagram}
satisfy the relations
\begin{multline*}
\alpha\cdot\eta =\alpha\cdot\odot_{\ge0}((\eta);1)\cdot z 
=\odot_{\ge0}((\eta);1)\cdot\odot_{\ge0}((1);\alpha)\cdot z
\\
=\odot_{\ge0}((\eta);1)\cdot
	\odot_{\ge0}((1);\odot_{\ge0}((\bar\jmath_i);1))\cdot\mu\cdot z 
=\odot_{\ge0}((\bar\jmath_i);1)\cdot z.
\end{multline*}
Commutativity of squares \ovalbox{3} and \ovalbox{4} with the vertical 
arrows $\psi_k$ is proven separately on each of summands of the source 
of the square.
On \(((\cc_i)_{i\in[n]};\cc_0(0))\) commutativity holds due to 
\lemref{lem-((Ai)A(0))-initial}.
On \(((\ca_i)_{i\in[n]};\odot_{\ge0}((\ca_i)_{i\in[n]};\cp))\) 
verification reduces to diagram~\eqref{dia-alpha-eta-j1z}.

One easily finds out that all three columns of 
diagram~\vpageref{dia-4x4-wheels} compose to $\inj_2$:
\[((T\bar\imath_i);\odot_{\ge0}((T\bar\imath_i);1))\cdot\phi_1=\inj_2,\;
((\bar\jmath_i)_{i\in[n]};\odot_{\ge0}
	((\bar\jmath_i)_{i\in[n]};1))\cdot\phi_2 =\inj_2,\;
((\bar\jmath_i)_{i\in[n]};\eta)\cdot\phi_3 =\inj_2.
\]
Therefore, the exterior of this diagram drawn with isomorphisms $\phi_k$ 
is a pushout square.
Hence, the pasting \ovalbox{$1\cup2$} of squares \ovalbox{1} and 
\ovalbox{2} (a diagram in $\nOp n_1^{\dg}$) is a pushout square in
$\nOp n_1^{\gr}$.
However, a cone for a diagram \(\cd\to\nOp n_1^{\dg}\) is a colimiting 
cone iff its image under the forgetful functor
\(\nOp n_1^{\dg}\to\nOp n_1^{\gr}\) is a colimiting cone for the 
composition \(\cd\to\nOp n_1^{\dg}\to\nOp n_1^{\gr}\) (see
\secref{sec-Peculiar-features}).
Thus, the pasting \ovalbox{$1\cup2$} is a pushout square in
$\nOp n_1^{\dg}$, and the colimit \(A\langle M,\beta\rangle\) is 
isomorphic to
	\(((\ca_i\langle\cM_i,\beta_i\rangle)_{i\in[n]};
	\bigcirc_{i=0}^n\ca_i\langle\cM_i,
	\beta_i\rangle\odot^i_{\ca_i}\cp)\)
in $\nOp n_1^{\dg}$.
\end{proof}

\subsection{The lax \texorpdfstring{$\Cat$}{Cat}-multifunctor 
	\texorpdfstring{$\hoM$}{hom}.}
	\label{lax-Cat-multifunctor-HOM}
Starting with an arbitrary braided $\dg$\n-multicategory $\mcC$ we 
construct a lax $\Cat$\n-multifunctor $\hoM:\mcB\to\mcDG$, where the 
$\Cat$\n-multicategory $\mcB$ has \(\Ob\mcB=\Ob\mcC\).
Any symmetric $\dg$\n-multicategory is obviously braided.
For any sequence $(A_i)_{i\in I}$, $B$ of objects of $\mcB$ the category 
\(\mcB((A_i)_{i\in I};B)\) is the terminal category $\1$.
An arbitrary lax $\Cat$\n-multifunctor $\mcB\to\mcDG$ assigns an object 
of \(\dg^{\NN^I}\) to a sequence $(A_i)_{i\in I}$, $B$.
In the case of $\hoM$ this is the object 
\(\hoM((A_i)_{i\in I};B)\in\Ob\dg^{\NN^I}\) given by
\[ \hoM((A_i)_{i\in I};B)((n^i)_{i\in I})
=\mcC\bigl((\sS{^{n^i}}A_i)_{i\in I};B\bigr).
\]

A 2\n-morphism has to be given for each labelled $[I]$\n-tree 
\((t,\ell):[I]\to\co_\sk(\Ob\mcC)\) and each $t$\n-tree 
\(\tau:t\to\co_\sk\):
\begin{multline*}
\comp^I_\tau:
\bigotimes^{h\in I} \bigotimes^{b\in t(h)} \bigotimes^{p\in\tau(h,b)}
\hoM\bigl((A_{h-1}^a)_{a\in t_h^{-1}b};A_h^b\bigr)
\bigl((|\tau_{(h-1,a)\to(h,b)}^{-1}(p)|)_{a\in t_h^{-1}b}\bigr) 
\\
=\bigotimes^{h\in I} \bigotimes^{b\in t(h)} \bigotimes^{p\in\tau(h,b)} 
\mcC\bigl((\sS{^{|\tau_{(h-1,a)\to(h,b)}^{-1}(p)|}}A_{h-1}^a
	)_{a\in t_h^{-1}b};A_h^b\bigr) 
\\
\to \mcC\bigl((\sS{^{|\tau(0,a)|}}A_0^a)_{a\in t(0)};A_l^1\bigr) =
\hoM\bigl((A_0^a)_{a\in t(0)};A_l^1\bigr)
\bigl((|\tau(0,a)|)_{a\in t(0)}\bigr).
\end{multline*}
Here \(0=\min[I]\), \(l=\max[I]\).
As such we take multiplication \(\mu_\mcC^{\tilde\tau}\) in 
$\dg$\n-multicategory $\mcC$ associated with the labelled braided tree 
$\tilde\tau$ with \(\tilde\tau(h)=\sqcup_{b\in t(h)}\tau(h,b)\) and with 
the successor map
\[ S_{\tilde\tau} =\tilde\tau_h: \sqcup_{a\in t(h-1)}\tau(h-1,a)
\to \sqcup_{b\in t(h)}\tau(h,b)
\]
induced by the map \(t_h:t(h-1)\to t(h)\) of indexing sets and by the 
maps \(\tau_{(h-1,a)\to(h,t_h(a))}:\tau(h-1,a)\to\tau(h,t_h(a))\) of 
summands.
Label mappings \(\ell:\tilde\tau(h)\to\Ob\mcC\) associate $A_h^b$ to any 
\(p\in\tau(h,b)\).
Equations~\eqref{lax-sym-Mon-functor-f} follow from 
\cite[equation~(2.25.1)]{BesLyuMan-book} written for the algebra $\mcC$ 
in the lax Monoidal category \(\sS{^{\Ob\mcC}}\pmQuiver_{\dg}\).
The procedure of forming $\tilde\tau$ commutes with restricting trees 
along maps \(\psi:[J]\to[I]\).

An example of such symmetric $\dg$\n-multicategory $\mcC$ comes from 
$\Com$ -- the closed symmetric multicategory of complexes of
$\kk$\n-modules and their chain maps
\cite[Example~3.18]{BesLyuMan-book}.
It is representable by the symmetric Monoidal category $\dg$ of 
complexes and chain maps \cite[Example~3.27]{BesLyuMan-book}.
We take for $\mcC$ the associated enriched symmetric multicategory 
$\uCom$, which is a $\Com$\n-multicategory, or equivalently, a
$\dg$\n-multicategory.
The composition in $\uCom$ has a natural meaning: this is a composition 
(of tensor products) of homogeneous maps, taking into account the Koszul 
rule.

The $\Cat$\n-multicategory $\mcB$ from \secref{lax-Cat-multifunctor-HOM} 
is also a $\Cat$-span multicategory with the discrete category 
$\tar\mcB$, \(\Ob\tar\mcB=\Ob\mcC\).
The lax $\Cat$\n-multifunctor $\hoM:\mcB\to\mcDG$ is simultaneously a 
lax $\Cat$-span multifunctor $\hoM:\mcB\to\mcDG$.
The prism equation from \defref{def-lax-Cat-span-multifunctor} takes in 
notation of \secref{sec-Cat-operad-graded-k-modules} the form
\begin{multline}
\bigl[\circledast_\mcDG(t)
\bigl(\hoM((A_{h-1}^a)_{a\in t_h^{-1}b};A_h^b)\bigr)_{(h,b)\in\IV(t)} 
\\
\rTTo^{\lambda_\mcDG^f} \circledast_\mcDG(t_\psi)\bigl(
\circledast_\mcDG(t^{|c}_{[\psi(g-1),\psi(g)]})
(\hoM((A_{h-1}^a)_{a\in t_h^{-1}b};A_h^b)
	)_{(h,b)\in\IV(t^{|c}_{[\psi(g-1),\psi(g)]})}
\bigr)_{(g,c)\in\IV(t_\psi)}
\\
\rTTo^{\circledast_\mcDG(t_\psi)\comp(t^{|c}_{[\psi(g-1),\psi(g)]})} 
\circledast_\mcDG(t_\psi)\bigl( 
\hoM((A_{\psi(g-1)}^a)_{a\in t_{\psi,g}^{-1}c};A_{\psi g}^c) 
\bigr)_{(g,c)\in\IV(t_\psi)}
\\
\rTTo^{\comp(t_\psi)} \hoM((A_0^a)_{a\in t(0)};A_{\max[I]}^1) \bigr] 
=\comp(t).
\label{eq-hom-comp(t)}
\end{multline}

Let us discuss the relationship between the suspension and the 
composition \(\mu_\uCom^T\) for a symmetric free $T:[l]\to\cS_\sk$.

Let \(g:U\to W\), \(f_i:X_i\to Y_i\), \(1\le i\le k\) be homogeneous 
maps of certain degrees.
For any \(1\le j\le k\) the maps
\[ \uCom\bigl((1)_{i<j},f_j,(1)_{i>j};1\bigr):
\uCom\bigl((X_i)_{i<j},(Y_i)_{i\ge j};U\bigr) \to
\uCom\bigl((X_i)_{i\le j},(Y_i)_{i>j};U\bigr)
\]
are defined as the precomposition with $f_j$,
\(h\mapsto(-1)^{h\cdot f_j}(1^{j-1}\times f_j\times1^{k-j})\cdot h\).
The map
\[ \uCom\bigl((1)_{i=1}^k;g\bigr): 
\uCom\bigl((X_i)_{i=1}^k;U\bigr) \to \uCom\bigl((X_i)_{i=1}^k;W\bigr)
\]
is defined as the postcomposition with $g$, \(h\mapsto h\cdot g\).
By convention, the map
\[ \uCom\bigl(f_1,f_2,\dots,f_k;g\bigr):
\uCom\bigl((Y_i)_{i=1}^k;U\bigr) \to \uCom\bigl((X_i)_{i=1}^k;W\bigr)
\]
is the composition (in this order)
\[ \uCom\bigl(f_1,1,\dots,1;1\bigr) \cdot
\uCom\bigl(1,f_2,\dots,1;1\bigr) \cdot \ldots \cdot 
\uCom\bigl(1,\dots,1,f_k;1\bigr) \cdot \uCom\bigl(1,\dots,1,1;g\bigr).
\]
Factors of this product commute up to the sign depending on parity of 
the product of degrees of factors.

\begin{lemma}
For arbitrary complexes \(A_h^b\in\Ob\uCom\) the following square 
commutes up to the sign \((-1)^{c(T)}\)
\begin{diagram}
\bigotimes^{h\in I} \bigotimes^{b\in T(h)}
\uCom\bigl((sA_{h-1}^a)_{a\in T_h^{-1}b};sA_h^b\bigr) 
&\rTTo^{\mu_\uCom^T}
&\uCom\bigl((sA_0^a)_{a\in T(0)};sA_l^1\bigr)
\\
\dTTo<{\tens^{h\in I}\tens^{b\in T(h)}
	\uCom(\sS{^{T_h^{-1}b}}\sigma;\sigma^{-1})}
&\sss (-1)^{c(T)}
&\dTTo>{\uCom(\sS{^{T(0)}}\sigma;\sigma^{-1})}
\\
\bigotimes_{h\in I} \bigotimes_{b\in T(h)}
\uCom\bigl((A_{h-1}^a)_{a\in T_h^{-1}b};A_h^b\bigr) 
&\rTTo_{\mu_\uCom^T}
&\uCom\bigl((A_0^a)_{a\in T(0)};A_l^1\bigr)
\end{diagram}
where
\[ c(T) =\sum_{h=1}^{l-1} \sum_{b=1}^{|T(h)|} (b-1)(1-|T_h^{-1}b|)
+\sum_{h=1}^{l-1} |\{(x,y)\in T(h-1)^2 \mid x<y,\; T_h(x)>T_h(y) \}|.
\]
\end{lemma}

\begin{proof}
The sign coincides with the sign of permutation of the expression
	\(\tens^{h\in I}\tens^{b\in T(h)}
	(\sigma^{\tens T_h^{-1}b}\tens\sigma^{-1})\),
followed by cancellation of matching \(\sigma^{-1}\) and $\sigma$, 
resulting in \((-1)^{c(T)}\sigma^{\tens T(0)}\tens\sigma^{-1}\).
This permutation can be performed in two steps applied to each floor of 
the tree starting from the root.
At the first step factors of \(\sigma^{\tens T_h^{-1}b}\) (starting from 
the right) are moved to the left towards the matching \(\sigma^{-1}\).
This explains appearance of the first sum in $c(T)$.
At the second step factors \(\tens^{b\in T(h)}\sigma^{\tens T_h^{-1}b}\) 
are permuted to \(\sigma^{\tens T(h-1)}\) accordingly to the map 
\(T_h:T(h-1)\to T(h)\).
This is reflected by the second sum in $c(T)$.
\end{proof}

Let \(g:U\to W\), \(f_i:X_i\to Y_i\), \(1\le i\le k\) be homogeneous 
maps of certain degrees.
Then
\[ \hoM\bigl((f_i)_{i\in I};g\bigr):
\hoM\bigl((Y_i)_{i\in I};U\bigr) \to \hoM\bigl((X_i)_{i\in I};W\bigr)
\]
denotes the collection of homogeneous maps
\[ \hoM\bigl((f_i)_{i\in I};g\bigr)
=\uCom\bigl((\sS{^{n^i}}f_i)_{i\in I};g\bigr):
\uCom\bigl((\sS{^{n^i}}Y_i)_{i\in I};U\bigr) \to 
\uCom\bigl((\sS{^{n^i}}X_i)_{i\in I};W\bigr).
\]

\begin{corollary}\label{cor-square-commutes-up-to-sign}
For each $t$\n-tree $\tau$ as above the following square commutes up to 
the sign \((-1)^{c(\tilde\tau)}\)
\begin{diagram}[nobalance]
\bigotimes^{h\in I} \bigotimes^{b\in t(h)} \bigotimes^{p\in\tau(h,b)}
\hoM\bigl((sA_{h-1}^a)_{a\in t_h^{-1}b};sA_h^b\bigr)
\bigl((|\tau_{(h-1,a)}^{-1}(p)|)_{a\in t_h^{-1}b}\bigr) \hspace*{-8em}
&&\phantom{\hoM\bigl((sA_0^a)_{a\in t(0)};sA_l^1\bigr)
	\bigl((|\tau(0,a)|)_{a\in t(0)}\bigr)}
\\
&\rdTTo_{\comp^I_\tau} &
\\
\dTTo~{\tens^{h\in I}\tens^{b\in t(h)}\tens^{p\in\tau(h,b)}
\hoM(\sS{^{t_h^{-1}b}}\sigma;\sigma^{-1})
((|\tau_{(h-1,a)}^{-1}(p)|)_{a\in t_h^{-1}b})\hspace*{-1.5em}}
&&\hoM\bigl((sA_0^a)_{a\in t(0)};sA_l^1\bigr)
	\bigl((|\tau(0,a)|)_{a\in t(0)}\bigr)
\\
&\sss (-1)^{c(\tilde\tau)}
&\dTTo~{\hoM(\sS{^{t(0)}}\sigma;\sigma^{-1})((|\tau(0,a)|)_{a\in t(0)})}
\\
\bigotimes^{h\in I} \bigotimes^{b\in t(h)} \bigotimes^{p\in\tau(h,b)}
\hoM\bigl((A_{h-1}^a)_{a\in t_h^{-1}b};A_h^b\bigr)
\bigl((|\tau_{(h-1,a)}^{-1}(p)|)_{a\in t_h^{-1}b}\bigr) \hspace*{-8em}
\\
&\rdTTo_{\comp^I_\tau} &
\\
&&\hoM\bigl((A_0^a)_{a\in t(0)};A_l^1\bigr)
	\bigl((|\tau(0,a)|)_{a\in t(0)}\bigr)
\end{diagram}
where $t:[l]\to\cO_\sk$, $\tau:t\to\cO_\sk$,
\begin{multline}
c(\tilde\tau) =\sum_{h=1}^{l-1} \sum_{x\in\tilde\tau(h)}
(1-|\tilde\tau_h^{-1}x|) |\{y\in\tilde\tau(h)\mid y<x\}|
\\
+\sum_{h=1}^{l-1} \sum_{\substack{a,b\in t(h-1)\\a<b,\;t_ha=t_hb}}
|\{(u,v)\in\tau(h-1,a)\times\tau(h-1,b) \mid 
\tau_{(h-1,a)}u>\tau_{(h-1,b)}v \}|.
\label{eq-c(tilde-tau)SS-SS}
\end{multline}
\end{corollary}

\begin{proof}
Let \(a,b\in t(h-1)\), \(u\in\tau(h-1,a)\), \(v\in\tau(h-1,b)\).
Then \(x=(a,u)\), \(y=(b,v)\in\tilde\tau(h)\).
Assume that $x<y$.
Several cases occur.
If $a=b$, $u<v$, then \(\tilde\tau_h(x)\le\tilde\tau_h(y)\).
If $a<b$, $t_ha<t_hb$, then \(\tilde\tau_h(x)<\tilde\tau_h(y)\).
If $a<b$, $t_ha=t_hb$, then \(\tilde\tau_h(x)>\tilde\tau_h(y)\) 
$\Longleftrightarrow$ \(\tau_{(h-1,a)}u>\tau_{(h-1,b)}v\).
\end{proof}

\subsection{Multicategory of operad modules.}
Note that $n\wedge1$-operad modules form a $\Cat$-span multiquiver 
$\mcM$ with $\tar\mcM=\Op$, the category of operads, with the functors
\[ (\ca_i)_{i\in I} \lMapsTo^\src ((\ca_i)_{i\in I};\cp;\cb) 
\rMapsTo^\tgt \cb.
\]

The $\Cat$-span multiquiver $\mcM$ becomes a weak $\Cat$-span 
multicategory when equipped with the tensor product \(\circledast_\mcM\) 
defined as follows.
Consider a tree \(t:[p]\to\co_\sk\) like in \eqref{eq-t-t(0)-t(l)-1}.
If $p=0$ the morphism
\(\circledast^\emptyset:\boxdot^\emptyset\mcM=\tar\mcM\to\mcM\) takes an 
operad $\ca$ to the regular $\ca$\n-bimodule \((\ca;\ca;\ca)\).
If $p=1$ the morphism \(\circledast^{\mb1}:\boxdot^{\mb1}\mcM\to\mcM\) 
is the natural isomorphism.
If $p>1$, $g\in\NN$, $0<g<p$, introduce a new tree by doubling the 
$g$\n-th level of $t$:
\begin{equation*}
t_+^g =\bigl(t(0) \rto{t_1} \dots\to t(g-1) \rto{t_g} t(g) \rto\id t(g) 
\rto{t_{g+1}} t(g+1) \to\dots \rto{t_p} t(p)=\mb1\bigr).
\end{equation*}
Consider a diagram of the form
\begin{diagram}[LaTeXeqno,h=1.6em,w=1.6em]
&&\bullet
\\
&\cdots &\dTTo\dTTo &\cdots 
\\
\bullet &\pile{\rTTo \\ \rTTo} &\sss\blacklozenge &\pile{\lTTo \\ \lTTo} &\bullet
\label{dia-double-T}
\end{diagram}
where $p-1$ pairs of parallel arrows are
\begin{equation}
\bullet =\circledast_\mcDG(t_+^g)((\cp_h^b)_{(h,b)\in\IV(t)},
(\ca_g^c)_{c\in t(g)}) \pile{\rTTo \\ \rTTo}
\circledast_\mcDG(t)(\cp_h^b)_{(h,b)\in\IV(t)} =\;\sss\blacklozenge
\label{eq-O(t+g)toto-O(t)}
\end{equation}
for $0<g<p$.
The first arrow corresponds to simultaneous right actions $\rho$ of 
$\ca_g^c$ on $\cp_g^c$ and the second arrow consists of simultaneous 
left actions $\lambda$ of \((\ca_g^c)_{c\in t_{g+1}^{-1}b}\) on 
$\cp_{g+1}^b$, \(b\in t(g+1)\).
In detail: we take the two embeddings \(\psi:[p]\to[p+1]\) missing $g$ 
or $g+1$, then \((t_+^g)_\psi=t\).
We use isomorphism~\eqref{eq-circledast(tpsi)circledast}
\begin{multline*}
\circledast_\mcDG(t_+^g)
\bigl((\cp_h^b)_{(h,b)\in\IV(t)},(\ca_g^c)_{c\in t(g)}\bigr) \rTTo^\sim
\circledast_\mcDG(t)\bigl((\cp_h^b)_{h<g}^{b\in t(h)},
(X_b)_{b\in t(g)},(\cp_h^b)_{h>g}^{b\in t(h)}\bigr)
\\
\rTTo \circledast_\mcDG(t)(\cp_h^b)_{(h,b)\in\IV(t)}
\end{multline*}
where $X_b$ stands for \(\circledast_\mcDG(t_g^{-1}b\to\{b\}\to\{b\})
(\cp_g^b,\ca_g^b)\) (resp. for
	\(\circledast_\mcDG(t_g^{-1}b\rto\id t_g^{-1}b\to\{b\})
	((\ca_g^c)_{c\in t_{g+1}^{-1}b},\cp_{g+1}^b)\)),
mapped by $\rho$ (resp. $\lambda$) to $\cp_g^b$ (resp. $\cp_{g+1}^b$).
The summand of the source of \eqref{eq-O(t+g)toto-O(t)} indexed by 
\(\tau_+^g:t_+^g\to\co_\sk\) is mapped by the first (resp. the second) 
arrow to the summand indexed by $t$\n-tree \(\tau_\rho^g\) (resp. 
\(\tau_\lambda^g\)), given by
\(\bigl(t\rTTo^{\psi-} t_+^g\rTTo^{\tau_+^g} \co_\sk\bigr)\), where the 
functor \(\psi-:t\to t_+^g\) corresponds to the embedding $\psi$ of 
vertices.
Thus,
\begin{equation}
\pi:\circledast_\mcDG(t)(\cp_h^b)_{(h,b)\in\IV(t)} \to 
\circledast_\mcM(t)(\cp_h^b)_{(h,b)\in\IV(t)}
\label{eq-circledastG-circledastM}
\end{equation}
determines the colimiting cone of \eqref{dia-double-T}.

Another definition is based on diagram~\eqref{dia-double-T} with pairs 
of arrows
\begin{equation*}
\circledast_\mcDG(t_+^g)((\cp_h^b)_{(h,b)\in\IV(t)},\ca_g^c)
\pile{\rTTo \\ \rTTo} \circledast_\mcDG(t)(\cp_h^b)_{(h,b)\in\IV(t)}
\end{equation*}
indexed by elements $(g,c)$ of \(\IV'(t)=\IV(t)-\{\text{root}\}\) -- the 
set of internal edges of $t$ (edges from $t_1$ are excluded).
This pair is obtained from pair~\eqref{eq-O(t+g)toto-O(t)} by 
precomposing with units \(\eta:\kk\to\ca_g^b\) for
\(b\in t(g)-\{c\}\) ($\kk$ being the unit operad).
The operad $\ca_g^c$ still acts via $\rho$ on $\cp_g^c$ on the right and 
via $\lambda^c$ on $\cp_{g+1}^{t_{g+1}c}$ on the left.
The colimit of \eqref{dia-double-T} is
\(\circledast_\mcM(t)(\cp_h^b)_{(h,b)\in\IV(t)}\).
Equivalence of the two definitions follows from the possibility to take 
elements of the operads equal to the unity for all operads $\ca_g^b$ but 
one.

Notice that for $p>0$ the collection
	\(\bigl((\ca_0^b)_{b\in t(0)};
	\circledast_\mcDG(t)(\cp_h^b)_{(h,b)\in\IV(t)};\ca_p^1)\)
has an operad polymodule structure.
The right action is constructed with the help of two increasing 
injections: \(\psi_1:[2]\to[p+1]\), \(0\mapsto0\), \(1\mapsto p\), 
\(2\mapsto p+1\), and \(\psi_2:[p]\to[p+1]\) missing $p$.
Namely,
\begin{multline}
\rho =\bigl[\bigl( 
\circledast_\mcDG(t)(\cp_h^b)_{(h,b)\in\IV(t)} \bigr)\odot_0\ca_p^1 
\rTTo^{\lambda_{\psi_1}^{-1}}_\sim
\circledast_\mcDG(t_+^p)((\cp_h^b)_{(h,b)\in\IV(t)},\ca_p^1)
\\
\rTTo^{\lambda_{\psi_2}}_\sim \circledast_\mcDG(t)
((\cp_h^b)_{1\le h<p}^{b\in t(h)},\cp_p^1\odot_0\ca_p^1) 
\rTTo^{\circledast_\mcDG(t)((1),\rho)}
\circledast_\mcDG(t)(\cp_h^b)_{(h,b)\in\IV(t)} \bigr],
\label{eq-rho-t+p}
\end{multline}
where
	\(\cp_p^1\odot_0\ca_p^1
	=\circledast_\mcDG(t(p-1)\to\mb1\to\mb1)(\cp_p^1,\ca_p^1)\).
The left action uses the following increasing injections:
\(\psi_1:[2]\to[p+1]\), \(0\mapsto0\), \(1\mapsto1\), \(2\mapsto p+1\), 
and \(\psi_2:[p]\to[p+1]\) missing $1$.
Namely,
\begin{multline}
\lambda =\bigl[ \odot_{>0} \bigl( (\ca_0^b)_{b\in t(0)}; 
\circledast_\mcDG(t)(\cp_h^b)_{(h,b)\in\IV(t)} \bigr) 
\rTTo^{\lambda_{\psi_1}^{-1}}_\sim
\circledast_\mcDG(t_+^0)( (\ca_0^b)_{b\in t(0)};
(\cp_h^b)_{(h,b)\in\IV(t)} ) \rTTo^{\lambda_{\psi_2}}_\sim
\\
\circledast_\mcDG(t)\bigl[ 
\bigl(\odot_{>0}((\ca_0^c)_{c\in t_1^{-1}b};\cp_1^b) \bigr)_{b\in t(1)},
(\cp_h^b)_{1<h\le p}^{b\in t(h)} \bigr]
\rTTo^{\circledast_\mcDG(t)((\lambda),1)}
\circledast_\mcDG(t)(\cp_h^b)_{(h,b)\in\IV(t)} \bigr],
\label{eq-lambda-t+0}
\end{multline}
where
	\(\odot_{>0}((\ca_0^c)_{c\in t_1^{-1}b};\cp_1^b)
	=\circledast_\mcDG(t_1^{-1}b\rTTo^\id t_1^{-1}b\to\{b\})
	((\ca_0^c)_{c\in t_1^{-1}b};\cp_1^b)\).
These (outer) actions commute with inner actions of $\ca_g^c$, thus they 
project to the quotient
\(\circledast_\mcM(t)(\cp_h^b)_{(h,b)\in\IV(t)}\), making
	\(\bigl((\ca_0^b)_{b\in t(0)};
	\circledast_\mcM(t)(\cp_h^b)_{(h,b)\in\IV(t)};\ca_p^1\bigr)\)
into an operad polymodule.

There is a lax $\Cat$-span multifunctor \(\mcM\to\mcDG\), 
\(((\ca_i)_{i\in I};\cp;\cb)\mapsto\cp\).
The component
	\(\pi:\circledast_\mcDG(t)(\cp_h^b)_{(h,b)\in\IV(t)}\to 
	\circledast_\mcM(t)(\cp_h^b)_{(h,b)\in\IV(t)}\)
is that of definition~\eqref{eq-circledastG-circledastM}.
We lift $\hoM$ to a lax $\Cat$-span multifunctor $\HOM:\mcB\to\mcM$ so 
that \(\bigl(\mcB\rTTo^\HOM \mcM\rTTo \mcDG\bigr)=\hoM\).
Namely, \(\tar\HOM:\tar\mcB\to\tar\mcM\), \(B\mapsto\END B\), and
$\HOM:\mcB\to\mcM$,
\[ ((A_i)_{i\in I};B) \mapsto
((\END A_i)_{i\in I};\hoM((A_i)_{i\in I};B);\END B)
=\HOM((A_i)_{i\in I};B).
\]
Prism equation~\eqref{eq-hom-comp(t)} projects to
\begin{multline}
\bigl[\circledast_\mcM(t)
\bigl(\HOM((A_{h-1}^a)_{a\in t_h^{-1}b};A_h^b)\bigr)_{(h,b)\in\IV(t)} 
\\
\rTTo^{\lambda_\mcM^f}
\circledast_\mcM(t_\psi)\bigl(
\circledast_\mcM(t^{|c}_{[\psi(g-1),\psi(g)]})
(\HOM((A_{h-1}^a)_{a\in t_h^{-1}b};A_h^b))
	_{(h,b)\in\IV(t^{|c}_{[\psi(g-1),\psi(g)]})} 
	\bigr)_{(g,c)\in\IV(t_\psi)}
\\
\rTTo^{\circledast_\mcM(t_\psi)\comp(t^{|c}_{[\psi(g-1),\psi(g)]})} 
\circledast_\mcM(t_\psi)\bigl( 
\HOM((A_{\psi(g-1)}^a)_{a\in t_{\psi,g}^{-1}c};A_{\psi g}^c) 
\bigr)_{(g,c)\in\IV(t_\psi)}
\\
\rTTo^{\comp(t_\psi)} \HOM((A_0^a)_{a\in t(0)};A_{\max[I]}^1) \bigr] 
=\comp(t).
\label{eq-comp(t)-decomposed}
\end{multline}

\begin{example}
Let us consider in detail the particular case of $p=2$.
Then $g=1$, \(t=\{\sqcup_{c=1}^n\mb{l}^c\to\bn\to\mb1\}\), 
\(t_+^1=\{\sqcup_{c=1}^n\mb{l}^c\to\bn\rto\id \bn\to\mb1\}\).
An arbitrary tree \(\tau_+^1:t_+^1\to\co_\sk\) has the form
\begin{equation}
\tau_+^1 =\quad
\begin{diagram}[h=0.8em,w=3em,nobalance,inline]
\sqcup_{p=1}^{u^n}\mb r_p^{n,l^n} \\
\cdots &\rdTTo \\
\sqcup_{p=1}^{u^n}\mb r_p^{n,1} &\rTTo
&\mb u^n=\sqcup_{q=1}^{k^n}\mb j_q^n &\rTTo &\mb k^n \\
\cdots &\cdots &&&&\rdTTo(2,4) \\
\sqcup_{p=1}^{u^c}\mb r_p^{c,l^c} &&\cdots &\rTTo &\cdots \\
&\rdTTo &&&&\rdTTo \\
\cdots &&\mb u^c=\sqcup_{q=1}^{k^c}\mb j_q^c &\rTTo &\mb k^c &\rTTo 
&\ \ \mb1\ \ . \\
&\ruTTo &&&&\ruTTo \ruTTo(2,4) \\
\sqcup_{p=1}^{u^c}\mb r_p^{c,1} &&\cdots &\rTTo &\cdots \\
\cdots &\cdots \\
\sqcup_{p=1}^{u^1}\mb r_p^{1,l^1} &\rTTo
&\mb u^1=\sqcup_{q=1}^{k^1}\mb j_q^1 &\rTTo &\mb k^1 \\
\cdots &\ruTTo \\
\sqcup_{p=1}^{u^1}\mb r_p^{1,1}
\end{diagram}
\label{dia-rocket-uk}
\end{equation}
The resulting $t$\n-trees \(\tau_\lambda^1\) and \(\tau_\rho^1\) are 
presented below:
\begin{equation}
\!\!\! \tau_\lambda^1 =
\begin{diagram}[inline,h=0.8em,w=3em,nobalance]
\sqcup_{p=1}^{u^n}\mb r_p^{n,l^n} \\
\cdots &\rdTTo \\
\sqcup_{p=1}^{u^n}\mb r_p^{n,1} &\rTTo &\mb u^n \\
\cdots &\cdots &&\rdTTo(2,4) \\
\sqcup_{p=1}^{u^c}\mb r_p^{c,l^c} &&\cdots \\
&\rdTTo &&\rdTTo \\
\cdots &&\mb u^c &\rTTo &\ \ \mb1\ , \\
&\ruTTo &&\ruTTo \ruTTo(2,4) \\
\sqcup_{p=1}^{u^c}\mb r_p^{c,1} &&\cdots \\
\cdots &\cdots \\
\sqcup_{p=1}^{u^1}\mb r_p^{1,l^1} &\rTTo &\mb u^1 \\
\cdots &\ruTTo \\
\sqcup_{p=1}^{u^1}\mb r_p^{1,1}
\end{diagram}
\qquad \tau_\rho^1 =
\begin{diagram}[inline,h=0.8em,w=3em,nobalance]
\sqcup_{q=1}^{k^n}\mb s_q^{n,l^n} \\
\cdots &\rdTTo \\
\sqcup_{q=1}^{k^n}\mb s_q^{n,1} &\rTTo &\mb k^n \\
\cdots &\cdots &&\rdTTo(2,4) \\
\sqcup_{q=1}^{k^c}\mb s_q^{c,l^c} &&\cdots \\
&\rdTTo &&\rdTTo \\
\cdots &&\mb k^c &\rTTo &\ \ \mb1\ . \\
&\ruTTo &&\ruTTo \ruTTo(2,4) \\
\sqcup_{q=1}^{k^c}\mb s_q^{c,1} &&\cdots \\
\cdots &\cdots \\
\sqcup_{q=1}^{k^1}\mb s_q^{1,l^1} &\rTTo &\mb k^1 \\
\cdots &\ruTTo \\
\sqcup_{q=1}^{k^1}\mb s_q^{1,1}
\end{diagram}
\label{eq-tau-lambda-tau-rho}
\end{equation}
In the latter case we use the notation \(r_{q,v}^{b,a}=r_p^{b,a}\) for 
\(q\in\mb k^b\), \(v\in\mb j_q^b\), where 
\(p=j_1^b+\dots+j_{q-1}^b+v\in\mb u^b\).
Furthermore,
\[ s_q^{b,a} =\sum_{v=1}^{j_q^b} r_{q,v}^{b,a} 
=\sum_{p=j_1^b+\dots+j_{q-1}^b+1}^{j_1^b+\dots+j_{q-1}^b+j_q^b}r_p^{b,a}.
\]
Direct summands corresponding to these trees are related by the two 
maps.
The second arrow is
\begin{multline}
\circledast_\mcG(t_+^1)((\cp^b)_{b\in\bn},(\cb^b)_{b\in\bn},
\cq(\tau_+^1) 
=\bigl[\bigotimes^{b\in\bn}
\bigotimes^{p\in\mb u^b} \cp^b((r_p^{b,a})_{a\in\mb l^b})\bigr]\tens
\bigl(\bigotimes^{b\in\bn} \bigotimes^{q\in\mb k^b}
\cb^b(j_q^b)\bigr) \tens \cq((k^b)_{b\in\bn})
\\
\rTTo^{1\tens\lambda} \bigl[\bigotimes^{b\in\bn}
\bigotimes^{p\in\mb u^b} \cp^b((r_p^{b,a})_{a\in\mb l^b})\bigr]
\tens \cq((u^b)_{b\in\bn})
=\circledast_\mcDG(t)((\cp^b)_{b\in\bn},\cq)(\tau_\lambda^1).
\label{eq-125--}
\end{multline}
The first arrow includes
isomorphism~\eqref{eq-circledast(tpsi)circledast}:
\begin{multline}
\bigl[\bigotimes^{b\in\bn} \bigotimes^{p\in\mb u^b} 
\cp^b((r_p^{b,a})_{a\in\mb l^b})\bigr]\tens \bigl(\bigotimes^{b\in\bn} 
\bigotimes^{q\in\mb k^b} \cb^b(j_q^b)\bigr) \tens \cq((k^b)_{b\in\bn}) 
\rTTo^\sim
\\
\bigl\{ \bigotimes^{b\in\bn} \bigotimes^{q\in\mb k^b} \bigl[ \bigl(
\bigotimes^{v\in\mb j_q^b} 
\cp^b((r_{q,v}^{b,a})_{a\in\mb l^b})\tens \cb^b(j_q^b)\bigr]\bigr\} 
\tens \cq((k^b)_{b\in\bn})
\rTTo^{(\tens^{b\in\bn}\tens^{q\in\mb k^b}\rho)\tens1}
\\
\bigl[ \bigotimes^{b\in\bn} \bigotimes^{q\in\mb k^b} 
\cp^b((s_q^{b,a})_{a\in\mb l^b})\bigr] \tens \cq((k^b)_{b\in\bn})
=\circledast_\mcDG(t)((\cp^b)_{b\in\bn},\cq)(\tau_\rho^1).
\label{eq-125-}
\end{multline}
In the special case with distinguished \(c\in\bn\) we take $j_q^b=1$ for 
$b\ne c$, $q\in\mb k^b$.
The above maps are precomposed with the units \(\eta:\kk\to\cb^b(1)\) 
for $b\ne c$.
For $b=c$ the action $\rho$ of $\cb^c$ on $\cp^c$ is the essential part 
of the first arrow, and the action $\lambda^c$ of $\cb^c$ on $\cq$ is 
left in the second arrow.

For \(t=\{\sqcup_{c=1}^n\mb{l}^c\to\bn\to\mb1\}\) there is 
\(t_+^2=\{\sqcup_{c=1}^n\mb{l}^c\to\bn\to\mb1\to\mb1\}\).
A \(t_+^2\)\n-tree has the form
\begin{equation*}
\tau_+^2 =\quad
\begin{diagram}[inline,h=0.8em,w=3em,nobalance]
\sqcup_{p=1}^{u^n}\mb r_p^{n,l^n} \\
\cdots &\rdTTo \\
\sqcup_{p=1}^{u^n}\mb r_p^{n,1} &\rTTo &\mb u^n=\sqcup_{v=1}^w\mb u_v^n 
\\
\cdots &\cdots &&\rdTTo(2,4) \\
\sqcup_{p=1}^{u^c}\mb r_p^{c,l^c} &&\cdots \\
&\rdTTo &&\rdTTo \\
\cdots &&\mb u^c=\sqcup_{v=1}^w\mb u_v^c &\rTTo &\ \ \mb w &\rTTo
&\mb1\ \ . \\
&\ruTTo &&\ruTTo \ruTTo(2,4) \\
\sqcup_{p=1}^{u^c}\mb r_p^{c,1} &&\cdots \\
\cdots &\cdots \\
\sqcup_{p=1}^{u^1}\mb r_p^{1,l^1} &\rTTo
&\mb u^1=\sqcup_{v=1}^w\mb u_v^1 \\
\cdots &\ruTTo \\
\sqcup_{p=1}^{u^1}\mb r_p^{1,1}
\end{diagram}
\end{equation*}
Respectively, right action \eqref{eq-rho-t+p} of $\cc$ on 
\(\circledast_\mcG(t)((\cp^c)_{c\in\bn},\cq)\) is given on the direct 
summand corresponding to $\tau_+^2$ by
\begin{multline*}
\!\!\! \bigl[\bigl( \circledast_\mcG(t)((\cp^c)_{c\in\bn},\cq) \bigr)
\odot_0\cc \bigr](\tau_+^2)
=\biggl[\bigotimes^{v\in\mb w} \Bigl( \bigl[\bigotimes^{c\in\bn} 
\bigotimes_{p=u_1^c+\dots+u_{v-1}^c+1}^{u_1^c+\dots+u_{v-1}^c+u_v^c} \!\!\!
\cp^c((r_p^{c,g})_{g\in\mb l^c})\bigr]\tens\cq(u_v)\Bigr) 
\biggr]\tens\cc(w)
\\
\rTTo^\sim
\Bigl[\bigotimes^{c\in\bn} \bigotimes^{p\in\mb u^c} 
\cp^c((r_p^{c,g})_{g\in\mb l^c})\Bigr]\tens
\Bigl[\bigotimes^{v\in\mb w}\cq(u_v) \Bigl]\tens\cc(w)
\rTTo^{1\tens\rho}
\Bigl[\bigotimes^{c\in\bn} \bigotimes^{p\in\mb u^c} 
\cp^c((r_p^{c,g})_{g\in\mb l^c})\Bigr]\tens \cq(u).
\end{multline*}

There is also
	\(t_+^0=\{\sqcup_{c=1}^n\mb{l}^c\rTTo^\id 
	\sqcup_{c=1}^n\mb{l}^c\to\bn\to\mb1\}\).
A \(t_+^0\)\n-tree has the form
\begin{equation*}
\tau_+^0 =\quad
\begin{diagram}[inline,h=0.8em,w=3em,nobalance]
\sqcup_{q=1}^{k^n}\sqcup_{\yi=1}^{s_q^{n,l^n}}
\boldsymbol\nu_{q,\yi}^{n,l^n} 
&\rTTo &\mb j^{n,l^n}=\sqcup_{q=1}^{k^n}\mb s_q^{n,l^n} \\
&\cdots &\cdots &\rdTTo \\
\sqcup_{q=1}^{k^n}\sqcup_{\yi=1}^{s_q^{n,1}}\boldsymbol\nu_{q,\yi}^{n,1} 
&\rTTo &\mb j^{n,1}=\sqcup_{q=1}^{k^n}\mb s_q^{n,1} &\rTTo &\mb k^n \\
&\cdots &\cdots &\cdots &&\rdTTo(2,4) \\
\sqcup_{q=1}^{k^c}\sqcup_{\yi=1}^{s_q^{c,l^c}}
\boldsymbol\nu_{q,\yi}^{c,l^c}
&\rTTo &\mb j^{c,l^c}=\sqcup_{q=1}^{k^c}\mb s_q^{c,l^c} &&\cdots \\
&&&\rdTTo &&\rdTTo \\
&\cdots &\cdots &&\mb k^c &\rTTo &\ \ \mb1\ . \\
&&&\ruTTo &&\ruTTo \ruTTo(2,4) \\
\sqcup_{q=1}^{k^c}\sqcup_{\yi=1}^{s_q^{c,1}}\boldsymbol\nu_{q,\yi}^{c,1} 
&\rTTo &\mb j^{c,1}=\sqcup_{q=1}^{k^c}\mb s_q^{c,1} &&\cdots \\
&\cdots &\cdots &\cdots \\
\sqcup_{q=1}^{k^1}\sqcup_{\yi=1}^{s_q^{1,l^1}}
\boldsymbol\nu_{q,\yi}^{1,l^1} &\rTTo
&\mb j^{1,l^1}=\sqcup_{q=1}^{k^1}\mb s_q^{1,l^1} &\rTTo &\mb k^1 \\
&\cdots &\cdots &\ruTTo \\
\sqcup_{q=1}^{k^1}\sqcup_{\yi=1}^{s_q^{1,1}}\boldsymbol\nu_{q,\yi}^{1,1} 
&\rTTo &\mb j^{1,1}=\sqcup_{q=1}^{k^1}\mb s_q^{1,1}
\end{diagram}
\end{equation*}
Respectively, left action \eqref{eq-lambda-t+0} of operads $\ca_g^c$ on 
\(\circledast_\mcG(t)((\cp^c)_{c\in\bn},\cq)\) is given on the direct 
summand corresponding to $\tau_+^0$ by
\begin{multline*}
\odot_{>0}\bigl((\ca_g^c)^{c\in\bn}_{g\in\mb{l}^c};
\circledast_\mcG(t)((\cp^c)_{c\in\bn},\cq) \bigr)(\tau_+^0)
\\
=\Bigl[\bigotimes^{c\in\bn} \bigotimes^{g\in\mb l^c} 
\bigotimes^{q\in\mb k^c} \bigotimes^{\yi\in\mb s_q^{c,g}} 
\ca_g^c(\boldsymbol\nu_{q,\yi}^{c,g}) \Bigr] \tens 
\Bigl(\bigotimes^{c\in\bn} \bigotimes^{q\in\mb k^c} 
\cp^c((s_q^{c,g})_{g\in\mb l^c})\Bigr) \tens \cq(k)
\\
\rTTo^{\lambda_{\psi_2}}_\sim
\Bigl\{\bigotimes^{c\in\bn} \bigotimes^{q\in\mb k^c} \Bigl[
\Bigl( \bigotimes^{g\in\mb l^c} \bigotimes^{\yi\in\mb s_q^{c,g}} 
\ca_g^c(\boldsymbol\nu_{q,\yi}^{c,g}) \Bigr) \tens 
\cp^c((s_q^{c,g})_{g\in\mb l^c}) \Bigr] \Bigr\} \tens \cq(k)
\\
\rTTo^{(\tens\tens\lambda)\tens1}
\Bigl[\bigotimes^{c\in\bn} \bigotimes^{q\in\mb k^c} \cp^c\Bigl( 
\Bigl(\sum_{\yi=1}^{s_q^{c,g}}\nu_{q,\yi}^{c,g}\Bigr)_{g\in\mb l^c} 
\Bigr) \Bigr] \tens \cq(k).
\end{multline*}
\end{example}

\section{Morphisms with several entries}
	\label{sec-Morphisms-several-entries}
Here we give support to the observation that morphisms with $n$ entries 
of algebras over operads form an $n\wedge1$-operad module.
In particular, we find this module for \ainf-algebras.

\subsection{Main source of \texorpdfstring{$n\wedge1$}{n1}-operad 
	modules.}
	\label{sec-source-operad-modules}
Starting with an arbitrary symmetric $\dg$\n-multicategory $\mcC$ we get 
a $\dg$\n-operad $\ce(X)=\END X$ for any object $X$ and an
$n\wedge1$-module $\HOM=(\END A_1,\dots,\END A_n;\cH;\END B)$ for any 
family $A_1$, \dots, $A_n$, $B$ in $\Ob\mcC$
\begin{gather*}
(\END X)(v) =\mcC(\sS{^v}X;X),
\\
\cH(j^1,\dots,j^n) =\HOM(A_1,\dots,A_n;B)(j^1,\dots,j^n)
=\mcC\bigl((\sS{^{j^i}}A_i)_{i=1}^n;B\bigr).
\end{gather*}
The right action
\begin{multline*}
\rho_{(j_p^i)}: \Bigl[\bigotimes_{p\in\mb k}
\cH\bigl((j_p^i)_{i\in\bn}\bigr)\Bigr] \tens(\END B)(k)
=\Bigl[\bigotimes_{p\in\mb k}\mcC\bigl((\sS{^{j_p^i}}A_i)_{i=1}^n;B
\bigr)\Bigr] \tens\mcC(\sS{^k}B;B)
\\
\to\mcC\bigl((\sS{^{\ell^i}}A_i)_{i=1}^n;B\bigr)
=\cH\bigl((\ell^i)_{i=1}^n\bigr),
\end{multline*}
where \(\ell^i=\sum_{p=1}^kj_p^i\), equals to the multicategory 
composition $\mu_\phi$, which corresponds to the map
\(\phi:\mb l^1\sqcup\dots\sqcup\mb l^n\to\mb k\), whose restriction to 
$\mb l^i$ is isotonic and sends exactly $j_p^i$ elements to
\(p\in\mb k\).

The left action
\[ \lambda_{(j_p^i)}:
\Bigl[\bigotimes_{i\in\bn} \bigotimes_{p=1}^{k^i}(\END A_i)(j_p^i)\Bigr]
\tens\cH\bigl((k^i)_{i=1}^n\bigr)
\to \cH\biggl(\Bigl(\sum_{p=1}^{k^i}j_p^i\Bigr)_{i=1}^n\biggr),
\]
that is,
\begin{equation*}
\lambda_{(j_p^i)}: \Bigl[\bigotimes_{i\in\bn} \bigotimes_{p=1}^{k^i} 
\mcC(\sS{^{j_p^i}}A_i;A_i) \Bigr] 
\tens\mcC\bigl((\sS{^{k^i}}A_i)_{i=1}^n;B\bigr)
\to \mcC\bigl((\sS{^{\ell^i}}A_i)_{i=1}^n;B\bigr)
\end{equation*}
with \(\ell^i=\sum_{p=1}^{k^i}j_p^i\), equals to the multicategory 
composition $\mu_\psi$, corresponding to the isotonic map
 \(\psi=\sqcup\sqcup\con:\sqcup_{i=1}^n\sqcup_{p=1}^{k^i}\mb j_p^i
 \to\sqcup_{i=1}^n\sqcup_{p=1}^{k^i}\mb1\),
which sends exactly $j_p^i$ elements to the element of the target 
indexed by $(i,p)$. 
Notice that $\rho_\emptyset:(\END Y)(0)=\mcC(;Y)=\cH(0,\dots,0)$ is the 
identity map.

\begin{example}
In particular, reasonings of \secref{sec-source-operad-modules} apply to 
the symmetric $\dg$\n-multicategory $\mcC=\uCom$ and for any
$(n+1)$-tuple \((X_1,\dots,X_n;Y)\) of complexes give an
$n\wedge1$-operad module
\begin{gather*}
(\END X_1,\dots,\END X_n;\HOM(X_1,\dots,X_n;Y);\END Y),
\\
\cH(j^1,\dots,j^n) =\HOM(X_1,\dots,X_n;Y)(j^1,\dots,j^n)
=\uCom\bigl((\sS{^{j^i}}X_i)_{i=1}^n;Y\bigr).
\end{gather*}
The case of $n=0$ gives $\cH=Y$.
The left action
\begin{equation*}
\lambda_{(j_p^i)}: \Bigl[\bigotimes_{i\in\bn} \bigotimes_{p=1}^{k^i} 
\uCom(\sS{^{j_p^i}}A_i;A_i) \Bigr] 
\tens\uCom\bigl((\sS{^{k^i}}A_i)_{i=1}^n;B\bigr)
\to \uCom\bigl((\sS{^{\ell^i}}A_i)_{i=1}^n;B\bigr),
\quad (\tens_i\tens_p g_p^i)\tens f \mapsto h
\end{equation*}
with \(\ell^i=\sum_{p=1}^{k^i}j_p^i\) is found as
\[ h =\bigl[ \boxt^{i\in\bn}T^{\ell^i}\ca_i
\rTTo^{\boxt^{i\in\bn}\lambda^{\gamma_i}}
\boxt^{i\in\bn}\tens^{p\in\mb k^i}T^{j_p^i}\ca_i
\rTTo^{\boxt^{i\in\bn}\tens^{p\in\mb k^i}g_p^i} 
\boxt^{i\in\bn}T^{k^i}\ca_i \rTTo^f \cb \bigr]
\]
For the scope of this article there is no distinction between $\boxt$ 
and $\tens$, see \secref{sec-Motivation}.
The right action
\begin{multline*}
\rho_{(j_p^i)}: \Bigl[\bigotimes_{p\in\mb k}
\uCom\bigl((\sS{^{j_p^i}}A_i)_{i=1}^n;B\bigr)\Bigr]\tens\uCom(\sS{^k}B;B)
\to \uCom\bigl((\sS{^{\ell^i}}A_i)_{i=1}^n;B\bigr),
\quad (\tens_pf^p)\tens g \mapsto h,
\end{multline*}
where \(\ell^i=\sum_{p=1}^kj_p^i\), is found as (see
\cite[Eq.~(6.1.1)]{BesLyuMan-book} for isomorphism 
$\overline{\varkappa}$)
\[ h =\bigl[ \boxt^{i\in\bn}T^{\ell^i}A_i
\rTTo^{\boxt^{i\in\bn}\lambda^{\gamma_i}}
\boxt^{i\in\bn}\tens^{p\in\mb k}T^{j_p^i}A_i 
\rTTo^{\overline{\varkappa}^{-1}}
\tens^{p\in\mb k}\boxt^{i\in\bn}T^{j_p^i}A_i
\rTTo^{\tens^{p\in\mb k}f^p} \tens^{p\in\mb k}\cb \rTTo^g \cb \bigr].
\]
\end{example}

Given an operad $\co$ and an $n\wedge1$-operad $\co$\n-module $\cF_n$ 
for each $n\ge0$ we define a morphism of $\co$\n-algebras with $n$ 
arguments \(X_1,\dots,X_n\to Y\) as a morphism of $\nOp n_1$
\[ (\co,\dots,\co;\cF_n;\co) \to
(\END X_1,\dots,\END X_n;\hoM(X_1,\dots,X_n;Y);\END Y).
\]

\begin{example}\label{exa-Com-FAsn}
Produce an $n\wedge1$-operad $\As$\n-module $\FAs_n$ from the symmetric
$\dg$\n-multicategory $\COM$ with one object $*$ -- the symmetric 
$\dg$\n-operad of associative non--unital commutative algebras. 
It has \(\COM(k)=\kk\) for $k>0$, and \(\COM(0)=0\). 
The compositions are given by multiplication in $\kk$. 
Hence, \(\END_{\COM}*=\As\) and \(\FAs_n=\hoM_{\COM}(\sS{^n}*;*)\) has 
\(\FAs_n(j^1,\dots,j^n)=\kk=\kk u_j\) for all non--vanishing
\((j^1,\dots,j^n)\in\NN^n\), while \(\FAs_n(0,\dots,0)=0\).
In particular, \(\FAs_0=0\).
The actions for $\FAs_n$ are given by multiplication in $\kk$. 
A morphism of $n\wedge1$-operad modules
\[ (\As,\dots,\As;\FAs_n;\As) \to
(\END A_1,\dots,\END A_n;\hoM(A_1,\dots,A_n;B);\END B)
\]
amounts to a family of morphisms \(f_i:A_i\to B\) of associative 
differential graded $\kk$\n-algebras without units, \(i\in\bn\), such 
that the following diagrams commute for all \(1\le i<j\le n\):
\begin{diagram}[LaTeXeqno]
A_i\tens A_j &\rTTo^c_\sim &A_j\tens A_i &\rTTo^{f_j\tens f_i} &B\tens B
\\
\dTTo<{f_i\tens f_j} &&&&\dTTo>{m_B}
\\
B\tens B &&\rTTo^{m_B} &&B
\label{dia-fifj-cfjfi}
\end{diagram}

In fact, morphisms $f_i=f_{(e_i)}$ are images of the unit under the 
action map
\[ \dot{f}_{(e_i)}: \kk =\FAs_n(e_i) \to \uCom(A_i;B),
\quad 1\mapsto f_{(e_i)},
\]
where \(e_i=(0,\dots,0,1,0,\dots,0)\in\NN^n\) has 1 on $i$\n-th place. 
The equations hold for all \(1\le i<j\le n\):
\begin{diagram}
\FAs_n(e_i)\tens\FAs_n(e_j)\tens\As(2) &\rTTo^{\text{mult}} 
&\FAs_n(e_i+e_j)
\\
\dTTo &= &\dTTo
\\
\cH(e_i)\tens\cH(e_j)\tens(\END B)(2) &\rTTo^{\mu_{\id}} &\cH(e_i+e_j)
\end{diagram}
\begin{diagram}
\FAs_n(e_j)\tens\FAs_n(e_i)\tens\As(2) &\rTTo^{\text{mult}} 
&\FAs_n(e_i+e_j)
\\
\dTTo &= &\dTTo
\\
\cH(e_j)\tens\cH(e_i)\tens(\END B)(2) &\rTTo^{\mu_{(12)}} &\cH(e_i+e_j)
\end{diagram}
Here the compositions $\mu_{\id}$ and $\mu_{(12)}$ in $\uCom$ correspond 
to the two maps
\begin{align*}
\id:  \mb0\sqcup\dots\sqcup\mb0\sqcup\mb1\sqcup\mb0\sqcup\dots\sqcup 
\mb0\sqcup\mb1\sqcup\mb0\sqcup\dots\sqcup\mb0  &=\mb2 \to \mb2,
\\
(12): \mb0\sqcup\dots\sqcup\mb0\sqcup\mb1\sqcup\mb0\sqcup\dots
\sqcup\mb0\sqcup\mb1\sqcup\mb0\sqcup\dots\sqcup\mb0  &=\mb2 \to \mb2.
\end{align*}
The equations are more explicit in the form
\begin{diagram}
\kk\tens\kk\tens\kk &\rTTo^{\text{mult}} &\kk
\\
\dTTo<{\dot{f}_{(e_i)}\tens\dot{f}_{(e_j)}\tens\dot m_B} &=
&\dTTo>{\dot{f}_{(e_i+e_j)}}
\\
\uCom(A_i;B)\tens\uCom(A_j;B)\tens\uCom(B,B;B) &\rTTo^{\mu_{\id}} 
&\uCom(A_i,A_j;B)
\end{diagram}
\begin{diagram}
\kk\tens\kk\tens\kk &\rTTo^{\text{mult}} &\kk
\\
\dTTo<{\dot{f}_{(e_j)}\tens\dot{f}_{(e_i)}\tens\dot m_B} &=
&\dTTo>{\dot{f}_{(e_i+e_j)}}
\\
\uCom(A_j;B)\tens\uCom(A_i;B)\tens\uCom(B,B;B) &\rTTo^{\mu_{(12)}} 
&\uCom(A_i,A_j;B)
\end{diagram}
The same equations can be written as
\[(f_{(e_i)}\tens f_{(e_j)})m_B =f_{(e_i+e_j)} 
=c(f_{(e_j)}\tens f_{(e_i)})m_B: A_i\tens A_j \to B,
\]
which coincides with condition~\eqref{dia-fifj-cfjfi}.

Notice that if $A_i$, $i\in\bn$, $B$ are unital $\dg$\n-algebras, a 
collection of unital morphisms \(f_i:A_i\to B\) that satisfy 
equation~\eqref{dia-fifj-cfjfi} is the same as a single unital morphism 
\(f:A_1\tdt A_n\to B\). 
In fact, such $f$ determines \(f_i(x)=f(1\tdt1\tens x\tens1\tdt1)\) and 
can be recovered from the whole collection of $f_i$'s.

So defined multimorphisms turn the class of associative non--unital 
$\dg$\n-algebras into a multicategory $\AS$. 
In fact, any totally ordered finite set $I$ can be used instead of 
$\bn$. 
Corresponding to an isotonic map $\phi:I\to J$, composition of
multimorphisms
	\((f^j:(A_i)_{i\in\phi^{-1}j}\to B_j)
	=(f^j_i:A_i\to B_j)_{i\in\phi^{-1}j}\),
$j\in J$, \((g:(B_j)_{j\in J}\to C)=(g_j:B_j\to C)_{j\in J}\) is given 
by
\[ [(f^j:
(A_i)_{i\in\phi^{-1}j}\to B_j)_{j\in J}\cdot(g:(B_j)_{j\in J}\to C)]_k
=f^{\phi(k)}_k\cdot g_{\phi(k)}: A_k \to C.
\]
One can check condition~\eqref{dia-fifj-cfjfi} for this collection of 
morphisms.
\end{example}

\begin{definition}
An $n\wedge1$-operad module homomorphism
\[ (f_1,\dots,f_n;h;g): (\ca_1,\dots,\ca_n;\cp;\cb) \to
(\cc_1,\dots,\cc_n;\cq;\cd),
\]
of degree \(r\in\ZZ\) is a family of $\dg$\n-operad homomorphisms
\(g:\cb\to\cd\), \(f_i:\ca_i\to\cc_i\), \(0\le i\le n\), of degree $r$ 
and a collection of homogeneous $\kk$\n-linear maps
\(h(j):\cp(j)\to\cq(j)\), $j\in\NN^n$, of degree \(r(1-|j|)\) such that
\begin{itemize}
\item for all \(l\in\NN\), \((k_q\in\NN^n\mid1\le q\le l)\), the 
following square commutes up to the sign
\begin{diagram}[h=2.8em,LaTeXeqno]
\Bigl(\bigotimes_{q=1}^l\cp(k_q)\Bigr)\tens\cB(l) &\rTTo^\rho 
&\cP\biggl(\sum_{q=1}^lk_q\biggr)
\\
\dTTo<{(\tens_{q=1}^lh(k_q))\tens g(l)} &\sss (-1)^{c_\rho} 
&\dTTo>{h(\sum_{q=1}^lk_q)}
\\
\Bigl(\bigotimes_{q=1}^l\cq(k_q)\Bigr)\tens\cd(l) &\rTTo^\rho
&\cq\biggl(\sum_{q=1}^lk_q\biggr)
\label{dia-PPPB-rho-P}
\end{diagram}
\begin{multline}
c_\rho =r\sum_{q=1}^l(q-1)(1-|k_q|) 
+r\sum_{1\le b<a\le l}^{1\le c<d\le n}k_a^ck_b^d
\\
+\frac{r(r-1)}2\biggl\{ (1-l)\sum_{q=1}^l (1-|k_q|)
+\sum_{1\le q<s\le l} (1-|k_q|)(1-|k_s|) \biggl\};
\label{eq-c-rho-oper-mod-hom}
\end{multline}

\item for all \(k\in\NN^n\),
\((j_p^i\in\NN\mid1\le i\le n,\;0\le p\le k^i)\), the following square 
commutes up to the sign
\begin{diagram}[h=3em,LaTeXeqno]
\biggl[\bigotimes_{i=1}^n \bigotimes_{p=1}^{k^i} \cA_i(j_p^i)\biggr]
\tens\cP\bigl((k^i)_{i=1}^n\bigr) &\rTTo^\lambda
&\cP\biggl(\Bigl(\sum_{p=1}^{k^i}j_p^i\Bigr)_{i=1}^n\biggr)
\\
\dTTo<{[\tens_{i=1}^n\tens_{p=1}^{k^i}f_i(j_p^i)]\tens h(k)}
&\sss (-1)^{c_\lambda} &\dTTo>{h((\sum_{p=1}^{k^i}j_p^i)_{i=1}^n)}
\\
\biggl[\bigotimes_{i=1}^n \bigotimes_{p=1}^{k^i} \cC_i(j_p^i)\biggr]
\tens\cq\bigl((k^i)_{i=1}^n\bigr) &\rTTo^\lambda
&\cq\biggl(\Bigl(\sum_{p=1}^{k^i}j_p^i\Bigr)_{i=1}^n\biggr)
\label{dia-AAAAP-lambda-P}
\end{diagram}
\begin{multline*}
c_\lambda =r\sum_{i=1}^n \sum_{p=1}^{k^i}
(1-j_p^i)\Bigl(p-1+\sum_{q=1}^{i-1} k^q\Bigr)
+\frac{r(r-1)}2\biggl\{ (1-|k|)\sum_{i=1}^n \sum_{p=1}^{k^i} (1-j_p^i)
\\
+\sum_{1\le i<l\le n} \Bigl[\sum_{p=1}^{k^i}(1-j_p^i)\Bigr]
\Bigl[\sum_{q=1}^{k^l}(1-j_q^l)\Bigr] 
+\sum_{i=1}^n \sum_{1\le p<q\le k^i} (1-j_p^i)(1-j_q^i) \biggl\};
\end{multline*}

\item for all $j\in\NN^n$
\[ d\cdot h(j) =(-1)^{r(1-|j|)}h(j)\cdot d: \cp(j) \to \cq(j).
\]
\end{itemize}
\end{definition}

The second (shuffle) part of $c_\rho$, $c_\lambda$ proportional to 
\(r(r-1)/2\) makes sure that the composition of morphisms of degrees $r$ 
and $r'$ be a morphism of degree $r+r'$.
The first part coincides with \(rc(\tilde{\tau}_\rho)\), 
\(rc(\tilde{\tau}_\lambda)\).

The last condition using $\lambda$ can be replaced with $n$ conditions 
using $\lambda^i$, \(1\le i\le n\):
\begin{diagram}[h=3em]
\biggl[\bigotimes_{p=1}^{k^i} \cA_i(j_p)\biggr]\tens\cP(k) 
&\rTTo^{\lambda^i} &\cP\biggl(k,k^i\mapsto\sum_{p=1}^{k^i}j_p\biggr)
\\
\dTTo<{[\tens_{p=1}^{k^i}f_i(j_p)]\tens h(k)} &\sss (-1)^{c_{\lambda^i}} 
&\dTTo>{h(k,k^i\mapsto\sum_{p=1}^{k^i}j_p)}
\\
\biggl[\bigotimes_{p=1}^{k^i} \cC_i(j_p)\biggr]\tens\cq(k) 
&\rTTo^{\lambda^i} &\cq\biggl(k,k^i\mapsto\sum_{p=1}^{k^i}j_p\biggr)
\end{diagram}
\begin{equation}
c_{\lambda^i}
=r\sum_{p=1}^{k^i} (1-j_p)\Bigl(p-1+\sum_{q=1}^{i-1} k^q\Bigr)
+\frac{r(r-1)}2\biggl\{ (1-|k|)\sum_{p=1}^{k^i} (1-j_p)
+\sum_{1\le p<q\le k^i} (1-j_p)(1-j_q) \biggl\}.
\label{eq-c-lambda-i-oper-mod-hom}
\end{equation}

\begin{example}
For all complexes $A_1$, \dots, $A_n$, $B$ the collection
\begin{multline*}
\varSigma
=(\sS{^n}\hoM(\sigma;\sigma^{-1});\hoM(\sS{^n}\sigma;\sigma^{-1});
\hoM(\sigma;\sigma^{-1})):
\\
\HOM(sA_1,\dots,sA_n;sB) \to \HOM(A_1,\dots,A_n;B)
\end{multline*}
is an $n\wedge1$-operad morphism of degree 1.
In fact, equations for $\varSigma$ involving $\lambda$ and $\rho$ are 
particular cases of \corref{cor-square-commutes-up-to-sign}.
\end{example}

\subsection{\texorpdfstring{$A_\infty$}{A8}-morphisms with several 
	entries.}
Produce an $n\wedge1$-operad $\Ass$\n-module $\FAss_n$ from the 
symmetric $\dg$\n-multicategory $\COMM$ with one object $*$ -- the 
symmetric $\dg$\n-operad of associative unital commutative algebras. 
It has \(\COMM(p)=\kk\) for $p\ge0$. 
The compositions are given by multiplication in $\kk$. 
Hence, \(\END_{\COMM}*=\Ass\) and \(\FAss_n=\hoM_{\COMM}(\sS{^n}*;*)\) 
has \(\FAss_n(j^1,\dots,j^n)=\kk\) for all \((j^1,\dots,j^n)\in\NN^n\).
In particular, \(\FAss_0=\kk\).
The actions for $\FAss_n$ are given by multiplication in $\kk$. 
A morphism of $n\wedge1$-operad modules
\[ (\Ass,\dots,\Ass;\FAss_n;\Ass) \to
(\END A_1,\dots,\END A_n;\hoM(A_1,\dots,A_n;B);\END B)
\]
amounts to a family of unital morphisms \(f_i:A_i\to B\) of associative 
unital differential graded $\kk$\n-algebras, \(i\in\bn\), such that
diagrams~\eqref{dia-fifj-cfjfi} commute for all \(1\le i<j\le n\). 
These data are in bijection with unital homomorphisms
\(f:A_1\tdt A_n\to B\), where $A_1$, \dots, $A_n$, $B$ are unital 
associative $\dg$\n-algebras.

In fact, each complex $A_1$, \dots, $A_n$, $B$ acquires a unital 
associative $\dg$\n-algebra structure through morphisms
\(\Ass\to\END A_i\), \(\Ass\to\END B\). 
Particular cases of actions
\begin{align*}
\lambda_{e_i}: \ca_i(0)\tens\cp(e_i) &\to \cp(0),
\\
\rho_\emptyset: \cb(0) =\kk\tens\cb(0) &\to \cp(0),
\end{align*}
for the module \((\Ass,\dots,\Ass;\FAss_n;\Ass)\) take unity to unity:
\begin{align*}
\lambda_{e_i}:
\Ass(0)\tens\FAss_n(e_i) \ni 1\tens1 &\mapsto 1 \in \FAss_n(0),
\\
\rho_\emptyset: \Ass(0) \ni 1 &\mapsto 1 \in \FAss_n(0).
\end{align*}
Commutative diagram~\cite[(2.2)]{Lyu-Ainf-Operad} %{dia-As1-As1-EndB-Hom} 
with \(\HOM(A_1,\dots,A_n;B)(0)\) in place of \(\HOM(A;B)(0)\) shows 
that \(1\in\FAss_n(0)\) is represented by $1_B$.
Since the representation agrees with \(\lambda_{e_i}\) the equation
\(1_{A_i}.f_i=1_B\) holds, thus, $f_i$ is unital.

\begin{proposition}\label{pro-Fn-FAsn-Fn}
There is the $n\wedge1$-operad $A_\infty$\n-module
	\(F_n=\boxdot_{\ge0}(\sS{^n}A_\infty;
	\kk\{f_j\mid j\in\NN^n-0\};A_\infty)\)
freely generated as a graded module by elements
$f_{j^1,\dots,j^n}\in F_n(j^1,\dots,j^n)$,
\((j^1,\dots,j^n)\in\NN^n-0\), of degree $0$.
The differential for it is given by
\begin{equation}
f_\ell\partial =\sum_{q=1}^n\sum_{r+x+t=\ell^q}^{x>1} 
\lambda^q_{(\sS{^r}1,x,\sS{^t}1)}
	(\sS{^r}1,b_x,\sS{^t}1;f_{\ell-(x-1)e_q})
-\sum_{\substack{j_1,\dots,j_k\in\NN^n-0\\j_1+\dots+j_k=\ell}}^{k>1} 
\rho_{(j_p^i)}((f_{j_p})_{p=1}^k;b_k).
\label{eq-f-ell-partial-bs}
\end{equation}
The first arguments of $\lambda$ are all $1\in A_\infty(1)$ except $b_x$ 
on the only place $p=r+1$.
$F_n$-maps are \ainf-algebra morphisms \(A_1,\dots,A_n\to B\) (for 
algebras written with operations $b_n$).
\end{proposition}

\begin{proof}
Notice that \(F_0=A_\infty(0)=0\) by \lemref{lem-((Ai)A(0))-initial}.

The following lemma is verified straightforwardly.

\begin{lemma}\label{lem-dg-category-claim}
For $\dg$\n-operads $\ca_1$, \dots, $\ca_n$ there is a $\dg$\n-category
\(\ca_1\text-\cdots\text-\ca_n\modul\), whose objects are left
$n$\n-operad \(\ca_1\text-\cdots\text-\ca_n\)-modules and degree $t$ 
morphisms \(f:\cP\to\cq\) are collections of $\kk$\n-linear maps
\(f(k^1,\dots,k^n):\cP(k^1,\dots,k^n)\to\cq(k^1,\dots,k^n)\) of degree 
$t$ such that
\begin{diagram}[height=3em]
\biggl[\bigotimes_{i\in\bn} \bigotimes_{p=1}^{k^i} \cA_i(j_p^i)\biggr]
\tens\cP\bigl((k^i)_{i=1}^n\bigr)
&\rTTo^{\lambda_{(j_p^i)}}
&\cP\biggl(\Bigl(\sum_{p=1}^{k^i}j_p^i\Bigr)_{i=1}^n\biggr)
\\
\dTTo<{[\tens\tens1]\tens f} &= &\dTTo>f
\\
\biggl[\bigotimes_{i\in\bn} \bigotimes_{p=1}^{k^i} \cA_i(j_p^i)\biggr]
\tens\cq\bigl((k^i)_{i=1}^n\bigr)
&\rTTo^{\lambda_{(j_p^i)}}
&\cq\biggl(\Bigl(\sum_{p=1}^{k^i}j_p^i\Bigr)_{i=1}^n\biggr)
\end{diagram}
The differential is \(f\mapsto[f,\partial]=f\partial-(-1)^f\partial f\).
\end{lemma}

A \emph{connection} on a graded $n\wedge1$-operad module $\cp$ over 
$\dg$\n-operads $\ca_1$, \dots, $\ca_n$, $\cb$ is a collection of 
$\kk$\n-linear maps \(\partial:\cp(j)\to\cp(j)\) of degree 1, 
\(j\in\NN^n\), which can be viewed as functors \(\ZZ\to\kk\modul\), 
\(p\mapsto\cp(j)^p\), where the category $\ZZ$ comes from the ordered 
set $\ZZ$.
All action maps \(\lambda^i\), $\rho$ from \eqref{eq-rho-rhokr} are 
required to be natural (chain) transformations with respect to the sum 
of maps \(1^{\tens a}\tens\partial\tens1^{\tens b}\) in the source, 
where $\partial$ denotes the connection on the module or the 
differential in an operad. 
Equivalently, action maps \(\lambda\), $\rho$ are chain transformations, 
or, equivalently, action maps $\alpha$ from \eqref{eq-alpha-full-action} 
are chain transformations.
A connection on a freely generated module $\cp$ is unambiguously fixed 
by its value on generators.

The square $\partial^2$ of a connection $\partial$ is also a connection 
(of degree 2).
It makes all actions into chain transformations with respect to the sum 
of maps \(1^{\tens a}\tens\partial^2\tens1^{\tens b}\) in the source 
(where $\partial^2$ vanishes if applied to an operad). 
In particular, \(\partial^2:\cp\to\cp\) is a morphism of graded left 
$n$\n-operad \(\ca_1\text-\cdots\text-\ca_n\)-modules of degree 2 as 
defined in \lemref{lem-dg-category-claim}.
If $\partial^2$ vanishes, $(\cp,\partial)$ becomes an $n\wedge1$-operad 
$\dg$\n-module.

\begin{lemma}
$F_n$ is an $n\wedge1$-operad $\dg$\n-module.
\end{lemma}

\begin{proof}
Recall that the differential in the operad $A_\infty$ is given by
\[ b_n.\partial 
=-\sum_{a+p+c=n}^{p>1,\;a+c>0} \mu(\sS{^a}1,b_p,\sS{^c}1;b_{a+1+c}).
\]
Let us prove that \(\partial^2=0\) for connection $\partial$ given by 
\eqref{eq-f-ell-partial-bs}.
Let us verify this on generators:
\begin{align}
&f_\ell\partial^2 
=\sum_{q=1}^n \sum_{k+y+m=\ell^q}^{y>1}
\lambda^q(\sS{^k}1,b_y,\sS{^m}1;f_{\ell-(y-1)e_q}.\partial)
+\sum_{q=1}^n \sum_{u+c+v=\ell^q}^{c>1}
\lambda^q(\sS{^u}1,b_c.\partial,\sS{^v}1;f_{\ell-(c-1)e_q}) \notag
\\
&-\sum_{\substack{j_1,\dots,j_k\in\NN^n-0\\j_1+\dots+j_k=\ell}}^{k>1}
\rho((f_{j_p})_{p=1}^k;b_k.\partial)
+\sum_{\substack{j_1,\dots,j_k\in\NN^n-0\\j_1+\dots+j_k=\ell}}^{k>1} 
\sum_{p=1}^k
\rho(f_{j_1},\dots,f_{j_{p-1}},f_{j_p}.\partial,
f_{j_{p+1}},\dots,f_{j_k};b_k)
\notag
\\
&=\sum_{q=1}^n \sum_{k+y+m=\ell^q}^{y>1} \sum_{p=1}^n 
\sum_{u+h+v=\ell^q-y+1}^{y>1}
\lambda^q(\sS{^k}1,b_y,\sS{^m}1;
\lambda^p(\sS{^u}1,b_h,\sS{^v}1;f_{\ell-(y-1)e_q-(h-1)e_p}))
\label{eq-bbf1}
\\
&-\sum_{q=1}^n \sum_{c>1} \sum_{r+c+t=\ell^q}
\sum_{
 \substack{j_1,\dots,j_k\in\NN^n-0\\j_1+\dots+j_k=\ell-(c-1)e_q}}^{k>1}
\lambda^q(\sS{^r}1,b_c,\sS{^t}1;\rho((f_{j_p})_p;b_k))
\label{eq-bfb1}
\\
&-\sum_{q=1}^n \sum_{u+c+v=\ell^q}^{c>1} \sum_{x+y+z=c}^{y>1}
\lambda^q(\sS{^{u+x}}1,b_y,\sS{^{z+v}}1;
\lambda^q(\sS{^u}1,b_{x+1+z},\sS{^v}1;f_{\ell-(c-1)e_q})) 
\label{eq-bbf2}
\\
&-\sum_{\substack{j_1,\dots,j_s\in\NN^n-0\\j_1+\dots+j_s=\ell}}^{s>1}
\sum_{x+m+z=s}^{m>1}
\rho((f_{j_p})_{p=1}^s;\mu(\sS{^x}1,b_m,\sS{^z}1;b_{x+1+z}))
\label{eq-fbb1}
\\
&+\sum_{q=1}^n
\sum_{\substack{j_1,\dots,j_k\in\NN^n-0\\j_1+\dots+j_k=\ell}}^{k>1} 
\sum_{p=1}^k \sum_{x+c+z=j_p^q}^{c>1} 
\rho(f_{j_1},\dots,f_{j_{p-1}},\lambda^q(\sS{^x}1,b_c,
\sS{^z}1;f_{j_p-(c-1)e_q}),f_{j_{p+1}},\dots,f_{j_k};b_k)
\label{eq-bfb2}
\\
&-\sum_{\substack{y_1,\dots,y_k\in\NN^n-0\\y_1+\dots+y_k=\ell}}^{k>1} 
\sum_{p=1}^k
\sum_{\substack{t_1,\dots,t_m\in\NN^n-0\\t_1+\dots+t_m=y_p}}^{m>1} 
\rho(f_{y_1},\dots,f_{y_{p-1}},\rho(f_{t_1},\dots,f_{t_m};b_m),
f_{y_{p+1}},\dots,f_{y_k};b_k).
\label{eq-fbb2}
\end{align}
Summands of \eqref{eq-bbf1} pairwise cancel each other if $p\ne q$.
Also summands of \eqref{eq-bbf1} with $p=q$ pairwise cancel each other 
if output of $b_y$ does not become an input of $b_h$.
The remainder of \eqref{eq-bbf1} cancels with sum~\eqref{eq-bbf2}.
Sums \eqref{eq-bfb1} and \eqref{eq-bfb2} cancel each other.
Identifying in sums \eqref{eq-fbb1} and \eqref{eq-fbb2} the index $s$ 
with \(k+m-1\) and the sequence \((j_1,\dots,j_s)\) with the sequence 
\((y_1,\dots,y_{p-1},t_1,\dots,t_m,y_{p+1},\dots,y_k)\) we see that they 
cancel.
Therefore, $\partial^2$ vanishes.
\end{proof}

The image of \(f_\ell\partial\) in \(\hoM((sA_i)_i;sB)\) is
\begin{multline*}
\sum_{q=1}^n \sum_{r+x+t=\ell^q}^{x>1}
\bigl[ \boxt^{i\in\bn}T^{\ell^i}sA_i
\rTTo^{1^{\boxt(q-1)}\boxt(1^{\tens r}\tens b_x
\tens1^{\tens t})\boxt1^{\boxt(n-q)}}
\\
T^{\ell^1}sA_1\boxt\dots\boxt T^{\ell^{q-1}}sA_{q-1}\boxt T^{r+1+t}sA_q 
\boxt T^{\ell^{q+1}}sA_{q+1}\boxt\dots\boxt T^{\ell^n}sA_n
\rTTo^{f_{\ell-(x-1)e_q}} sB \bigr]
\\
-\sum_{\substack{j_1,\dots,j_k\in\NN^n-0\\j_1+\dots+j_k=\ell}}^{k>1}
\bigl[\boxt^{i\in\bn}T^{\ell^i}sA_i 
\rTTo^{\boxt^{i\in\bn}\lambda^{\gamma_i}}
\boxt^{i\in\bn}\tens^{p\in\mb k}T^{j_p^i}sA_i 
\rTTo^{\overline{\varkappa}^{-1}}
\tens^{p\in\mb k}\boxt^{i\in\bn}T^{j_p^i}sA_i
\\
\rTTo^{\tens^{p\in\mb k}f_{j_p}} \tens^{p\in\mb k}sB \rTTo^{b_k} 
sB\bigr].
\end{multline*}
Isomorphisms \(\lambda^{\gamma_i}\) and \(\overline{\varkappa}\) are the 
obvious ones, see \cite{BesLyuMan-book} for details.

An $F_n$\n-algebra map is specified by \ainf-algebras $A_1$, \dots, 
$A_n$, $B$, and a collection of $\kk$\n-linear degree 0 maps
\(f_j:\boxt^{i\in\bn}T^{j^i}sA_i\to sB\) assigned to generators 
\((f_j)_{j\in\NN^n-0}\).
It suffices to satisfy on generators the only requirement that 
\(F_n\to\hoM((A_i[1])_{i=1}^n;B[1])\) be a chain map.
The latter means that the equation holds for all $\ell\in\NN^n-0$:
\[ f_\ell b_1 -\biggl[ \sum_{q=1}^n \sum_{r+1+t=\ell^q}
1^{\boxt(q-1)}\boxt(1^{\tens r}\tens
b_1\tens1^{\tens t})\boxt1^{\boxt(n-q)} \biggr]f_\ell =f_\ell\partial.
\]
Explicitly this equation says
\begin{multline}
\sum_{q=1}^n \sum_{r+x+t=\ell^q}^{x>0}
\bigl[ \boxt^{i\in\bn}T^{\ell^i}sA_i
\rTTo^{1^{\boxt(q-1)}\boxt(1^{\tens r}\tens b_x
\tens1^{\tens t})\boxt1^{\boxt(n-q)}}
\\
T^{\ell^1}sA_1\boxt\dots\boxt T^{\ell^{q-1}}sA_{q-1}\boxt T^{r+1+t}sA_q
\boxt T^{\ell^{q+1}}sA_{q+1}\boxt\dots\boxt T^{\ell^n}sA_n
\rTTo^{f_{\ell-(x-1)e_q}} sB \bigr]
\\
\hskip\multlinegap
=\sum_{\substack{j_1,\dots,j_k\in\NN^n-0\\j_1+\dots+j_k=\ell}}^{k>0}
\bigl[ \boxt^{i\in\bn}T^{\ell^i}sA_i 
\rTTo^{\boxt^{i\in\bn}\lambda^{\gamma_i}}
\boxt^{i\in\bn}\tens^{p\in\mb k}T^{j_p^i}sA_i \hfill
\\
\rTTo^{\overline{\varkappa}^{-1}} 
\tens^{p\in\mb k}\boxt^{i\in\bn}T^{j_p^i}sA_i
\rTTo^{\tens^{p\in\mb k} f_{j_p}}\tens^{p\in\mb k}sB \rTTo^{b_k} 
sB\bigr].
\label{eq-A8-f-shifted}
\end{multline}
Collections \((f_j)_{j\in\NN^n-0}\) are in bijection with augmented 
coalgebra morphisms \(f:\boxt^{i\in\bn}TsA_i\to TsB\). 
Coherence with augmentation means that
\[ \bigl( \kk \rEq \boxt^{i\in\bn}T^0sA_i \rTTo^{f|} TsB \bigr)
=\bigl( \kk \rEq T^0sB \rMono TsB \bigr).
\]
Tensor quivers $TsB$ of \ainf-algebras $B$ are $\dg$\n-coalgebras, whose
differential \(b:TsB\to TsB\) has the components \(b_k:T^ksB\to sB\).
Equation~\eqref{eq-A8-f-shifted} can be rewritten as
\begin{equation*}
\bigl( \boxt^{i\in\bn}TsA_i \xrightarrow{f} TsB \xrightarrow b sB \bigr)
= \bigl( \boxt^{i\in\bn}TsA_i
\rTTo^{\sum_{i=1}^n1^{\boxt(i-1)}\boxt b\boxt1^{\boxt(n-i)}}
\boxt^{i\in\bn}TsA_i \xrightarrow{f} sB \bigr).
\end{equation*}
In other terms, $f$ is an augmented $\dg$\n-coalgebra morphism.
These are \ainf-morphisms $A_1,\dots,A_n\to B$ by definition, see 
\cite{BesLyuMan-book}.
\end{proof}

\begin{proposition}\label{pro-Fn-FAsn-rmFn}
There is an $n\wedge1$-operad module \((\mainf,\rmF_n)\) freely 
generated as graded module by elements
$\sff_{j^1,\dots,j^n}\in\rmF_n(j^1,\dots,j^n)$,
\((j^1,\dots,j^n)\in\NN^n-0\), of degree $1-j^1-\dots-j^n=1-|j|$.
The differential for it is given by 
\begin{multline}
\sff_\ell\partial =\sum_{q=1}^n \sum_{r+x+t=\ell^q}^{x>1}
(-1)^{(1-x)(\ell^1+\dots+\ell^{q-1}+r)+1-|\ell|}
\lambda^q_{(\sS{^r}1,x,\sS{^t}1)}
	(\sS{^r}1,m_x,\sS{^t}1;\sff_{\ell-(x-1)e_q})
\\
+\sum_{\substack{j_1,\dots,j_k\in\NN^n-0\\j_1+\dots+j_k=\ell}}^{k>1}
(-1)^{k+\sum_{1\le b<a\le k}^{1\le c<d\le n}j_a^cj_b^d
	+\sum_{p=1}^k(p-1)(|j_p|-1)}
\rho_{(j_p^i)}((\sff_{j_p})_{p=1}^k;m_k).
\label{eq-fldel-prop}
\end{multline}
There is an invertible morphism of degree 1 between these
$n\wedge1$-operad modules
\begin{equation}
(\varSigma,\varSigma): (A_\infty,F_n) \to (\mainf,\rmF_n), 
\qquad b_i \mapsto m_i, \quad f_j \mapsto \sff_j.
\label{eq-(Sigma-Sigma)-(AF)-(AF)}
\end{equation}
$\rmF_n$-maps are \ainf-algebra morphisms \(A_1,\dots,A_n\to B\)
(for algebras written with operations $m_n$).
The two notions of \ainf-morphisms agree in the sense that the square of 
$n\wedge1$-operad module maps
\begin{diagram}[nobalance,LaTeXeqno,bottom]
(\sS{^n}A_\infty;F_n;A_\infty) &\rTTo
&((\END A_i[1])_{i=1}^n;\hoM((A_i[1])_{i=1}^n;B[1]);\END B[1])
\\
\dTTo<{(\sS{^n}\varSigma;\varSigma;\varSigma)}
&&\dTTo>{(\sS{^n}\hoM(\sigma;\sigma^{-1});
\hoM(\sS{^n}\sigma;\sigma^{-1});\hoM(\sigma;\sigma^{-1}))}
\\
(\sS{^n}\mainf;\rmF_n;\mainf) &\rTTo
&((\END A_i)_{i=1}^n;\hoM((A_i)_{i=1}^n;B);\END B)
\label{dia-nA8-Fn-A8}
\end{diagram}
commutes.
\end{proposition}

\begin{proof}
The existence of $F_n$ implies the existence of $\rmF_n$ as the 
following lemma shows.

\begin{lemma}
Let \((\ca_1,\dots,\ca_n;\cp;\cb)\) be a $\dg$\n-$n\wedge1$-operad 
module, \((\cc_1,\dots,\cc_n;\cq;\cd)\) be a graded $n\wedge1$-operad 
module and
\[ (f_1,\dots,f_n;h;g): (\ca_1,\dots,\ca_n;\cp;\cb) \to
(\cc_1,\dots,\cc_n;\cq;\cd),
\]
be an invertible graded $n\wedge1$-operad module homomorphism of degree 
$r$ (equations \eqref{dia-PPPB-rho-P} and \eqref{dia-AAAAP-lambda-P} 
hold).
Then $\cc_1$, \dots, $\cc_n$, $\cd$ are $\dg$\n-operads (see
\remref{rem-homogeneous-operad-homomorphisms}) and $\cp$ has a unique 
differential $d$ which turns it into a $\dg$\n-$n\wedge1$-operad module 
and \((f_1,\dots,f_n;h;g)\) into a $\dg$\n-$n\wedge1$-operad module 
isomorphism of degree $r$.
\end{lemma}

\begin{proof}
The differential is given by a unique expression
\[ d =\bigl(\cq(j) \rTTo^{h(j)^{-1}} \cp(j) \rTTo^{(-1)^{r(1-|j|)}d} 
\cp(j) \rTTo^{h(j)} \cq(j)\bigr).
\]
Clearly, $\deg d=1$ and $d^2=0$.
Verification that $\rho$ and $\lambda$ for $\cq$ are chain maps is 
straightforward.
\end{proof}

Let us compute the value of the differential on generators $\sff_\ell$:
\begin{align*}
\sff_\ell\partial &=(f_\ell.\varSigma(\ell))\partial 
=(-1)^{1-|\ell|}(f_\ell.\partial)\varSigma(\ell)
\\
&=(-1)^{1-|\ell|}\sum_{q=1}^n \sum_{r+x+t=\ell^q}^{x>1}
\lambda^q_{(\sS{^r}1,x,\sS{^t}1)}
	(\sS{^r}1,b_x,\sS{^t}1;f_{\ell-(x-1)e_q}).\varSigma(\ell)
\\
&\qquad +(-1)^{|\ell|}
\sum_{\substack{j_1,\dots,j_k\in\NN^n-0\\j_1+\dots+j_k=\ell}}^{k>1}
\rho_{(j_p^i)}((f_{j_p})_{p=1}^k;b_k).\varSigma(\ell)
\\
&=\sum_{q=1}^n \sum_{r+x+t=\ell^q}^{x>1} 
(-1)^{c(\tilde\tau_{\lambda^q})+1-|\ell|}
\lambda^q_{(\sS{^r}1,x,\sS{^t}1)}
	(\sS{^r}1,m_x,\sS{^t}1;\sff_{\ell-(x-1)e_q}) \\
&\qquad
+\sum_{\substack{j_1,\dots,j_k\in\NN^n-0\\j_1+\dots+j_k=\ell}}^{k>1} 
(-1)^{k+c(\tilde\tau_\rho)} \rho_{(j_p^i)}((\sff_{j_p})_{p=1}^k;m_k),
\end{align*}
which coincides with \eqref{eq-fldel-prop}, if one plugs in expressions 
\(c(\tilde\tau_{\lambda^q})=c_{\lambda^q}\) from
\eqref{eq-c-lambda-i-oper-mod-hom} and \(c(\tilde\tau_\rho)=c_\rho\) 
from \eqref{eq-c-rho-oper-mod-hom} for $r=1$.

The image of \(\sff_\ell\partial\) in \(\hoM((A_i)_i;B)\) is
\begin{multline*}
\sum_{q=1}^n \sum_{r+x+t=\ell^q}^{x>1}
(-1)^{(1-x)(\ell^1+\dots+\ell^{q-1}+r)+1-|\ell|}
\bigl[ \boxt^{i\in\bn}T^{\ell^i}A_i
\rTTo^{1^{\boxt(q-1)}\boxt(1^{\tens r}\tens m_x \tens1^{\tens t})
	\boxt1^{\boxt(n-q)}}
\\
T^{\ell^1}A_1\boxt\dots\boxt T^{\ell^{q-1}}A_{q-1}\boxt 
T^{r+1+t}A_q\boxt T^{\ell^{q+1}}A_{q+1}\boxt\dots\boxt T^{\ell^n}A_n
\rTTo^{\sff_{\ell-(x-1)e_q}} B \bigr]
\\
+\sum_{\substack{j_1,\dots,j_k\in\NN^n-0\\j_1+\dots+j_k=\ell}}^{k>1}
(-1)^{k+\sum_{1\le b<a\le k}^{1\le c<d\le n}j_a^cj_b^d
	+\sum_{p=1}^k(p-1)(|j_p|-1)}
\bigl[ \boxt^{i\in\bn}T^{\ell^i}A_i \rTTo^{\boxt^{i\in\bn}\lambda^{\gamma_i}}
\boxt^{i\in\bn}\tens^{p\in\mb k}T^{j_p^i}A_i
\\
\rTTo^{\overline{\varkappa}^{-1}}
\tens^{p\in\mb k}\boxt^{i\in\bn}T^{j_p^i}A_i
\rTTo^{\tens^{p\in\mb k}\sff_{j_p}} \tens^{p\in\mb k}B \rTTo^{m_k} B 
\bigr].
\end{multline*}
$\rmF_n$\n-algebra maps consist of \ainfm-algebras $A_1$, \dots, $A_n$, 
$B$, and a collection \((\sff_j)_{j\in\NN^n-0}\) that satisfies the 
following equation for all $\ell\in\NN^n-0$:
\[ \sff_\ell m_1
+(-1)^{|\ell|} \biggl[ \sum_{q=1}^n \sum_{r+1+t=\ell^q}
1^{\boxt(q-1)}\boxt(1^{\tens r}\tens m_1\tens1^{\tens t})
\boxt1^{\boxt(n-q)} \biggr]\sff_\ell =\sff_\ell\partial.
\]
In expanded form the equation says:
\begin{multline}
\sum_{q=1}^n \sum_{r+x+t=\ell^q}^{x>0}
(-1)^{(1-x)(\ell^1+\dots+\ell^{q-1}+r)-|\ell|}
\bigl[ \boxt^{i\in\bn}T^{\ell^i}A_i
\rTTo^{1^{\boxt(q-1)}\boxt(1^{\tens r}\tens m_x
\tens1^{\tens t})\boxt1^{\boxt(n-q)}}
\\
T^{\ell^1}A_1\boxt\dots\boxt T^{\ell^{q-1}}A_{q-1}\boxt T^{r+1+t}A_q
\boxt T^{\ell^{q+1}}A_{q+1}\boxt\dots\boxt T^{\ell^n}A_n
\rTTo^{\sff_{\ell-(x-1)e_q}} B \bigr]
\\
=\sum_{\substack{j_1,\dots,j_k\in\NN^n-0\\j_1+\dots+j_k=\ell}}^{k>0}
(-1)^{k+\sum_{1\le b<a\le k}^{1\le c<d\le n}j_a^cj_b^d
	+\sum_{p=1}^k(p-1)(|j_p|-1)}
\bigl[ \boxt^{i\in\bn}T^{\ell^i}A_i 
\rTTo^{\boxt^{i\in\bn}\lambda^{\gamma_i}}
\boxt^{i\in\bn}\tens^{p\in\mb k}T^{j_p^i}A_i
\\
\rTTo^{\overline{\varkappa}^{-1}} 
\tens^{p\in\mb k}\boxt^{i\in\bn}T^{j_p^i}A_i
\rTTo^{\tens^{p\in\mb k}\sff_{j_p}} \tens^{p\in\mb k}B \rTTo^{m_k} 
B\bigr].
\label{eq-A8-morphism-signes}
\end{multline}
This is actually the definition of an \ainfm-algebra morphism
\(A_1,\dots,A_n\to B\) for algebras written with operations $m_n$, 
adopted in the current article.
%Previous version used different generators $\sff'_j$ represented by
%\((s^{\tens j^1}\boxt\dots\boxt s^{\tens j^n})\cdot f_j\cdot s^{-1}\) in \(\hoM((A_i)_{i=1}^n;B)(j)\) instead of $\sff_j$ represented by \((\sigma^{\tens j^1}\boxt\dots\boxt\sigma^{\tens j^n})\cdot f_j\cdot\sigma^{-1}\).
%Equations for $\sff'_j$ coincide with \eqref{eq-A8-morphism-signes} if one plugs in \(\sff'_j=(-1)^{1-|j|}\sff_j\).

Relationship between $f_j$ and $\sff_j$ in \(\hoM((A_i)_{i=1}^n;B)(j)\),
\begin{diagram}
T^{j^1}A_1\boxt\dots\boxt T^{j^n}A_n &\rTTo^{\sff_j} &B
\\
\dTTo<{\sigma^{\tens j^1}\boxt\dots\boxt\sigma^{\tens j^n}} 
&&\dTTo>\sigma
\\
T^{j^1}sA_1\boxt\dots\boxt T^{j^n}sA_n &\rTTo^{f_j} &sB
\end{diagram}
shows that diagram~\eqref{dia-nA8-Fn-A8} commutes on generators.
Therefore, it is commutative.
\end{proof}

Reducing the data used in \secref{sec-n1-operad-modules} or
\defref{def-M-operad-polymodules} we call an \emph{$n$\n-dimensional 
right operad module} the pair \((\cp;\cb)\) consisting of a
$\dg$\n-operad $\cB$ and an object \(\cp\in\dg^{\NN^n}\), equipped with 
a unital associative action
\[ \rho: \Bigl(\bigotimes_{q=1}^l\cp(k_q)\Bigr)\tens\cB(l)
\rTTo \cP\biggl(\sum_{q=1}^lk_q\biggr) \in\dg.
\]

\begin{definition}
An \emph{$n$\n-dimensional right operad module homomorphism}
\((h;g):(\cp;\cb)\to(\cq;\cd)\) of degree $(p;0)$, \(p\in\ZZ^n\), is a 
$\dg$\n-operad homomorphism \(g:\cb\to\cd\) of degree $0$ and a 
collection of homogeneous $\kk$\n-linear maps \(h(j):\cp(j)\to\cq(j)\), 
$j\in\NN^n$, of degree \((p|j)=\sum_{i=1}^np^ij^i\) such that
\begin{itemize}
\item for all \(l\in\NN\), \((k_q\in\NN^n\mid1\le q\le l)\), the 
following square commutes up to the sign
\begin{diagram}[h=2.8em,nobalance,LaTeXeqno]
\Bigl(\bigotimes_{q=1}^l\cp(k_q)\Bigr)\tens\cB(l) &\rTTo^\rho 
&\cP\biggl(\sum_{q=1}^lk_q\biggr)
\\
\dTTo<{(\tens_{q=1}^lh(k_q))\tens g(l)} &\sss (-1)^{c(k_1,\dots,k_l)} 
&\dTTo>{h(\sum_{q=1}^lk_q)}
\\
\Bigl(\bigotimes_{q=1}^l\cq(k_q)\Bigr)\tens\cd(l) &\rTTo^\rho
&\cq\biggl(\sum_{q=1}^lk_q\biggr)
\label{dia-PPPB-rho-P-c}
\end{diagram}
\begin{equation}
c(k_1,\dots,k_l) =\sum_{1\le t<q\le l} \chi(k_t,k_q),
\label{eq-c(kkkkk)-chi(kk)}
\end{equation}
where \(\chi:\NN^n\times\NN^n\to\ZZ/2\) is an arbitrary bilinear form 
(it is specified by a matrix \(\chi\in\Mat(n,\ZZ/2)\));

\item for all $j\in\NN^n$
\begin{equation}
d\cdot h(j) =(-1)^{(p|j)}h(j)\cdot d: \cp(j) \to \cq(j).
\label{eq-dh-(-1)hd}
\end{equation}
\end{itemize}
\end{definition}

\begin{lemma}\label{lem-PB=hg-QD}
Let \((\cp;\cb)\) be an $n$\n-dimensional right $\dg$\n-operad module, 
let \(g:\cb\to\cd\) be a $\dg$\n-operad isomorphism of degree $0$.
Let \(h(j):\cp(j)\to\cq(j)\), $j\in\NN^n$, be a collection of invertible 
homogeneous $\kk$\n-linear maps of degree \((p|j)\) for some 
\(p\in\ZZ^n\).
Let \(\chi:\NN^n\times\NN^n\to\ZZ/2\) be a bilinear form.
Then $\cq$ admits a unique structure of an $n$\n-dimensional right 
$\cd$\n-module such that \((h;g):(\cp;\cb)\to(\cq;\cd)\) is a 
homomorphism of degree $(p;0)$ with respect to $\chi$.
\end{lemma}

\begin{proof}
The value of the differential in $\cq$ is fixed by \eqref{eq-dh-(-1)hd}.
The unique candidate $\rho$ for action of $\cd$ on $\cq$ is found from 
diagram~\eqref{dia-PPPB-rho-P-c}.
This $\rho$ is a chain map, as follows from a cubical diagram consisting 
of two faces~\eqref{dia-PPPB-rho-P-c} joined by differentials.
Opposite faces of the cube commute up to the same sign, since
\((p|\sum_{q=1}^lk_q)=\sum_{q=1}^l(p|k_q)\).
Therefore, the both squares expressing commutation of $\rho$ with the 
differential commute simultaneously.

Associativity of the action of $\cd$ on $\cq$ is expressed by the 
pentagon
\begin{diagram}[h=2.8em,nobalance]
\bigotimes_{q=1}^l\Bigl(\bigotimes_{t=1}^{n_q}\cq(\sS{_t}k_q)\tens
\cd(n_q)\Bigr)\tens\cd(l)
&\rTTo^{\hspace*{-0.9em}\tens_{q=1}^l\rho\tens1} 
&\Bigl(\bigotimes_{q=1}^l\cq\Bigl(\sum_{t=1}^{n_q}\sS{_t}k_q\Bigr)
\Bigr)\tens\cd(l)
&\rTTo^\rho &\cq\Bigl(\sum_{q=1}^l\sum_{t=1}^{n_q}\sS{_t}k_q\Bigr)
\\
\dTTo<\wr &&&\ruTTo<\rho
\\
\Bigl(\bigotimes_{q=1}^l\bigotimes_{t=1}^{n_q}\cq(\sS{_t}k_q)\Bigr)
\tens\Bigl(\bigotimes_{q=1}^l\cd(n_q)\Bigr)\tens\cd(l)
&\rTTo^{1\tens\mu_\cd}
&\Bigl(\bigotimes_{q=1}^l\bigotimes_{t=1}^{n_q}\cq(\sS{_t}k_q)\Bigr)
\tens\cd\Bigl(\sum_{q=1}^ln_q\Bigr) \hspace*{-2.6em}
\end{diagram}
lying at the bottom of a rectangular prism, whose top face is the 
pentagon, expressing associativity of the action of $\cb$ on $\cp$.
Vertical maps are tensor products of $h$ and $g$.
The walls commute up to sign.
The product of these signs is $+1$, since
\[ c\Bigl(\bigl((\sS{_t}k_q)_{t=1}^{n_q}\bigr)_{q=1}^l\Bigr) 
=c\biggl(\Bigl(\sum_{t=1}^{n_q}\sS{_t}k_q\Bigr)_{q=1}^l\biggr)
+\sum_{q=1}^lc\bigl((\sS{_t}k_q)_{t=1}^{n_q}\bigr)
\]
due to definition~\eqref{eq-c(kkkkk)-chi(kk)} of $c$ and bilinearity of 
$\chi$.

Unitality of the action of $\cd$ on $\cq$ follows from that for $\cb$ 
and $\cp$, since \(c(k)=0\), \(k\in\NN^n\).
\end{proof}

Cofibrant replacement of an $n\wedge1$-operad module
$(\co,\cp)\overset{\text{def}}=(\co,\dots,\co;\cp;\co)$ is a trivial 
fibration \((\ca,\cF)\to(\co,\cp)\) such that the only map from the 
initial $n\wedge1$-operad module \((\1,0)\to(\ca,\cF)\) is a cofibration 
in $\nOp n_1$.

\begin{theorem}\label{thm-Fn-FAsn-cofibrant-replacement}
The $n\wedge1$-operad module \((\mainf,\rmF_n)\) is a cofibrant 
replacement of $(\As,\FAs_n)$.
Moreover, \((\mainf,\rmF_n)\to(\As,\FAs_n)\) is a homotopy isomorphism 
in \(\dg^{\NN\sqcup\NN^n}\).
\end{theorem}

\begin{proof}
Generate a free $n\wedge1$-operad $\As$\n-module $\overline{\rmF}_n$ by 
elements $\sff_{j^1,\dots,j^n}\in\overline{\rmF}_n(j^1,\dots,j^n)$,
\((j^1,\dots,j^n)\in\NN^n-0\), of degree $1-j^1-\dots-j^n$. 
Actually, \((\As,\overline{\rmF}_n)\) is the coequalizer in $\nOp n_1$ 
of the pair of morphisms of collections
\[ 0,\inj: (\kk\{(m_2\tens1)m_2-(1\tens m_2)m_2,m_n\mid n\ge3\},0)
\rightrightarrows (\mainf,\rmF_n),
\]
the second arrow is just the embedding. 
Therefore, the differential in $\overline{\rmF}_n$ reduces to
\begin{equation*}
\sff_\ell\partial =\sum_{q=1}^n \sum_{r+2+t=\ell^q}
(-1)^{1+t+\ell^{q+1}+\dots+\ell^n}
\lambda^q(\sS{^r}1,m,\sS{^t}1;\sff_{\ell-e_q})
-\sum_{q+r=\ell}^{q,r\in\NN^n-0}
(-1)^{|r|+\sum_{c>d}q^cr^d}\rho(\sff_q,\sff_r;m),
\end{equation*}
$m=m_2$, and the equation $\partial^2=0$ follows. 
Notice that the quadratic part of the differential
\begin{equation}
\sff_\ell\bar\partial =\sum_{q+r=\ell}^{q,r\in\NN^n-0}
(-1)^{1-|r|+\sum_{c>d}q^cr^d}\rho(\sff_q,\sff_r;m)
\label{eq-fl-bar-partial}
\end{equation}
is a differential itself, $\bar\partial^2=0$.

The $n\wedge1$-operad module morphism in question decomposes as
\[ (\mainf,\rmF_n) \rFib^{htis} (\As,\overline{\rmF}_n) \rFib^{(1,p)} 
(\As,\FAs_n).
\]
The first epimorphism is a homotopy isomorphism, since $\mainf\to\As$ 
is. 
Let us describe the second epimorphism and prove for it the same 
property, that is, the $n\wedge1$-operad $\As$\n-module epimorphism
$p:\overline{\rmF}_n\to\FAs_n$ is a homotopy isomorphism. 
We prove more: the zero degree cycle $p:\overline{\rmF}_n\to\FAs_n$ is 
homotopy invertible in the $\dg$\n-category \(n\text-\!\As\modul\).

Any left $n$\n-operad $\As$\n-module $\cp$ decomposes into a direct sum 
of submodules.
Any subset $S\subset\bn$ with the induced total ordering is viewed as 
the isomorphic ordinal with $|S|$ elements.
For any \(k\in\NN^n\) denote by \(\supp k=\{i\in\bn \mid k^i\ne0\}\) its 
support.
Consider the $n$\n-operad $\As$\n-submodule
\[ \cp^S(k) =
\begin{cases}
\cp(k) &\quad \text{if } \supp k =S,
\\
0 &\quad \text{otherwise}.
\end{cases}
\]
Then \(\cp=\oplus_{S\subset\bn}\cp^S\). 
Since \(\As(0)=0\), the $\ZZ^n$\n-graded collection $\cp^S$ is a left 
$n$\n-operad $\As$\n-module.
This structure boils down to a $\ZZ^S$\n-graded collection $\cp^S$, 
which is a left $S$\n-operad $\As$\n-module (that is, a $|S|$\n-operad 
$\As$\n-module).
A left $n$\n-operad $\As$\n-module $\cp$ is freely generated iff left 
$S$\n-operad $\As$\n-modules $\cp^S$ are freely generated for all 
\(S\subset\bn\).

Let \(e_S\in\NN^n\) have the coordinates
\(e_S^i=\chi(i\in S)\in\{0,1\}\), \(e_i\overset{\text{def}}=e_{\{i\}}\). 
For \(j\in\NN^n\), $j\ne0$, consider the basic element 
\(u_j=1\in\kk=\FAs_n(j)\).
For $S\ne\emptyset$ the element \(u_{e_S}=1\in\kk=\FAs_n(e_S)\) freely 
generates the left $S$\n-operad $\As$\n-module $\FAs_n^S$, while 
$\FAs_n^\emptyset=0$. 
Namely, for any \(j\in\NN^n\), $j\ne0$, with support \(S=\supp j\) we 
have \(u_j=\lambda((m^{(j^i)})_{i\in S};u_{e_S})\).

The left $n$\n-operad $\As$\n-module $\overline{\rmF}_n$ is also freely 
generated.
Its basis is given by elements
\(\rho(\sff_{j_1},\dots,\sff_{j_k};m^{(k)})\), where $k>0$ and 
\(j_t\in\NN^n-0\) for all $t$.

The $n\wedge1$-operad $\As$\n-module map $p$ is specified on the 
generators as follows:
\[ \sff_j.p =
\begin{cases}
u_j, &\quad \text{if } |j| =1,
\\
0, &\quad \text{otherwise}.
\end{cases}
\]
On the basis of the left $n$\n-operad $\As$\n-module $\overline{\rmF}_n$ 
the map $p$ is computed as
\[ \rho(\sff_{j_1},\dots,\sff_{j_k};m^{(k)}).p =
\begin{cases}
u_j, &\quad \text{if}\quad |j_1| =\dots =|j_k| =1,\quad j=\sum_r j_r,
\\
0, &\quad \text{otherwise}.
\end{cases}
\]
In order to prove that $p$ is a chain map it suffices to prove that 
\(\sff_{2e_a}.\partial p=0\), \(1\le a\le n\), and 
\(\sff_{e_a+e_b}.\partial p=0\) for all \(1\le a<b\le n\).
These equations are verified straightforwardly:
\begin{align*}
\sff_{2e_a}.\partial =-\lambda(m;\sff_{e_a}) +\rho(\sff_{e_a},
\sff_{e_a};m) &\rMapsTo^p -\lambda(m;u_{e_a}) +u_{2e_a} =0,
\\
\sff_{e_a+e_b}.\partial =\rho(\sff_{e_a},\sff_{e_b};m) 
-\rho(\sff_{e_b},\sff_{e_a};m) &\rMapsTo^p u_{e_a+e_b} -u_{e_b+e_a} =0.
\end{align*}

A zero degree cycle \(\beta:\FAs_n\to\overline{\rmF}_n\) in 
\(n\text-\!\As\modul\) is given on generators $u_j$ of free
$\kk$\n-modules $\FAs_n(j)$ by the formula
\[ u_j.\beta =\rho\bigl((\lambda(m^{(j^i)};\sff_{e_i}))_{i\in\supp j};
m^{(|\supp j|)}\bigr).
\]
The composition
\begin{diagram}
\FAs_n &\rTTo^\beta &\overline{\rmF}_n &\rTTo^p &\FAs_n
\\
u_{e_S} &\rMapsTo &\rho((\sff_{e_i})_{i\in S};m^{(|S|)})
=(\tens_{i\in S}\sff_{e_i})m^{(S)} &\rMapsTo
&(\tens_{i\in S}u_{e_i})m^{(S)} =u_{e_S}
\end{diagram}
is the identity map.
Let us prove that \(p\beta\) is homotopy invertible. 
These two statements would imply that $p$ is homotopy invertible in 
\(n\text-\!\As\modul\) and $\beta$ is its homotopy inverse.

Let $\overline{\rmF}_n^{(q)}$ be a $n\text-\!\As$\n-submodule generated 
by \(\rho(\sff_{j_1},\dots,\sff_{j_k};m^{(k)})\), $k\le q$, 
$\overline{\rmF}_n^{(0)}=0$.
This filtration induces the graded $n\text-\!\As$\n-module with the 
components 
	\(\overline{\rmF}_n^{\{k\}}
	=\overline{\rmF}_n^{(k)}/\overline{\rmF}_n^{(k-1)}\).
Since the differential in $\As$ vanishes, the differential
\(\partial:\overline{\rmF}_n^{(q)}\to\overline{\rmF}_n^{(q+1)}\) is a 
left $n$\n-operad $\As$\n-module map.
We look for a left $n$\n-operad $\As$\n-module map 
$h:\overline{\rmF}_n\to\overline{\rmF}_n$ of degree $-1$ such that 
\(\overline{\rmF}_n^{(q)}.h\subset\overline{\rmF}_n^{(q-1)}\).
Consider the zero degree cycle 
\[ N=1-p\beta+h\partial+\partial h\overline{\rmF}_n\to\overline{\rmF}_n.
\]
It satisfies 
\(\overline{\rmF}_n^{(q)}.N\subset\overline{\rmF}_n^{(q)}\). 
We are going to choose $h$ in such a way that $N$ be locally nilpotent. 
Thus, $1-N$ is invertible with the (well--defined) inverse 
\(\sum_{a=0}^\infty N^a\). Therefore,
\[ p\beta =1 -N +h\partial +\partial h:
\overline{\rmF}_n \to \overline{\rmF}_n
\]
is homotopy invertible.

Since
	\(\rho\bigl(\overline{\rmF}_n^{(q_1)}(j_1)\tdt
	\overline{\rmF}_n^{(q_k)}(j_k)\tens\As(k)\bigr) 
	\subset\overline{\rmF}_n^{(q_1+\dots+q_k)}(j_1+\dots+j_k)\)
there is an induced map between quotients:
\[ \overline\rho: \overline{\rmF}_n^{\{q_1\}}
(j_1)\tdt\overline{\rmF}_n^{\{q_k\}}(j_k)\tens\As(k) 
\to \overline{\rmF}_n^{\{q_1+\dots+q_k\}}(j_1+\dots+j_k)
\]
The actions $\overline\rho$ assemble to an action of $\As$ on the sum 
	\(\overline{\rmF}_n^{\{\}}
	=\oplus_{q=0}^\infty\overline{\rmF}_n^{\{q\}}\).
The quadratic differential
	\(\bar{\partial}:\overline{\rmF}_n^{\{q\}}(j)
	\to\overline{\rmF}_n^{\{q+1\}}(j)\)
from \eqref{eq-fl-bar-partial} induces a differential \(\bar{\partial}\) 
in \(\overline{\rmF}_n^{\{\}}\), thereby making it into a differential 
graded $n\wedge1$-$\As$-module.
As a left $n\text-\!\As$\n-module it is generated by its
$n$\n-dimensional right $\As$-$\dg$-submodule 
\(\overline{\sff}_n^{\{\}}\):
\[ \overline{\sff}_n^{\{\}}
=\oplus_{k=0}^\infty \overline{\sff}_n^{\{k\}}, \qquad
\overline{\sff}_n^{\{k\}}(j)
=\kk\bigl\{ \rho(\sff_{j_1},\dots,\sff_{j_k};m^{(k)}) \mid
j_1+\dots+j_k =j,\ \forall\;q\le k\; j_q\in\NN^n-0 \bigr\}.
\]
The matrix coefficients of
	\(\bar{\partial}:\overline{\sff}_n^{\{k\}}(j)
	\to\overline{\sff}_n^{\{k+1\}}(j)\)
are integers and we shall find
\(h:\overline{\sff}_n^{\{k\}}(j)\to\overline{\sff}_n^{\{k-1\}}(j)\) with 
the same property.
Thus, instead of working over a general ring $\kk$ we can assume that 
$\kk=\ZZ$, and we do it till the end of the proof.
Any such map $h$ extends to a morphism of left $n\text-\!\As$\n-modules 
in a unique way.

The operator induced by $N$ in the graded $n\text-\!\As$\n-module 
$\overline{\rmF}_n^{\{\}}$ is denoted
\(\bar{N}:\overline{\rmF}_n^{\{\}}\to\overline{\rmF}_n^{\{\}}\).
It can be described via a simplified formula
\begin{equation}
\bar{N} =1 -\overline{p\beta} +h\bar\partial +\bar\partial h: 
\overline{\rmF}_n^{\{k\}} \to \overline{\rmF}_n^{\{k\}},
\label{eq-N1-pd-hd-dh}
\end{equation} 
where \(\bar\partial\) is given by \eqref{eq-fl-bar-partial}, 
\(h:\overline{\rmF}_n^{\{p\}}\to\overline{\rmF}_n^{\{p-1\}}\), and 
\(\rho(\sff_{j_1},\dots,\sff_{j_k};m^{(k)}).\overline{p\beta}\) vanishes 
unless \(|j_1|=\dots=|j_k|=1\) and \(\supp j_q\) are all distinct for 
\(1\le q\le k\).
When \((j_1,\dots,j_k)\) is a permutation of \((e_{a_1},\dots,e_{a_k})\) 
with \(1\le a_1<\dots<a_k\le n\), then 
\[ \rho(\sff_{j_1},\dots,\sff_{j_k};m^{(k)}).\overline{p\beta}
=\rho(\sff_{e_{a_1}},\dots,\sff_{e_{a_k}};m^{(k)}).
\]
Otherwise,
\(\rho(\sff_{j_1},\dots,\sff_{j_k};m^{(k)}).\overline{p\beta}\) 
vanishes.
The operator $N$ is locally nilpotent iff $\bar{N}$ is.
We shall achieve $\bar{N}=0$.

Let us define a family of graded abelian groups \(\tilde{f}_n(j)\), 
\(j\in\NN^n\),
\[ \tilde{f}_n(j)^k =\ZZ\bigl\{ x(j_1,\dots,j_k) \mid
j_q\in\NN^n-0,\; j_1+\dots+j_k =j \bigr\}.
\]
The family \(\tilde{f}_n\) has an obvious structure of a graded
$n$\n-dimensional right $\As$\n-module, namely,
\[ \tilde\rho\bigl( x\bigl((\sS{_t}j_1)_{t=1}^{n_1}\bigr), \dots, 
x\bigl((\sS{_t}j_k)_{t=1}^{n_k}\bigr);m^{(k)} \bigr)
=x\bigl((\sS{_t}j_1)_{t=1}^{n_1}\bigr), \dots, 
(\sS{_t}j_k)_{t=1}^{n_k}\bigr).
\]
This structure is completely fixed by the requirement
\begin{equation}
\tilde\rho\bigl(x(j_1),\dots,x(j_k);m^{(k)}\bigr)=x(j_1,\dots,j_k)
\label{eq-rho(xxm)=x(jj)}
\end{equation}
Consider the bilinear form \(\chi:\NN^n\times\NN^n\to\ZZ\), 
\(\chi(t,p)=\sum_{c\le d}t^cp^d\), and define the corresponding $c$ by 
\eqref{eq-c(kkkkk)-chi(kk)}.
There are invertible mappings
\(\psi(j):\overline{\sff}_n^{\{\}}(j)\to\tilde{f}_n(j)\) of degree
\(|j|=((1,1,\dots,1)|j)\) such that
\begin{itemize}
\item \((\sff_j).\psi(j)=x(j)\);

\item the right $As$\n-module structure obtained from
\((\psi,\id_{\As})\) and the bilinear form $\chi$ as in
\lemref{lem-PB=hg-QD} satisfies condition~\eqref{eq-rho(xxm)=x(jj)}.
\end{itemize}
Existence and uniqueness of $\psi$ is shown in the following computation 
in square~\eqref{dia-PPPB-rho-P-c}:
\begin{diagram}[nobalance,h=2.2em]
\sff_{j_1}\tdt\sff_{j_k}\tens m^{(k)} &\rMapsTo^{\overline\rho} 
&\rho(\sff_{j_1},\dots,\sff_{j_k};m^{(k)}) 
&\rMapsTo^{\psi(j_1+\dots+j_k)}
&\rho(\sff_{j_1},\dots,\sff_{j_k};m^{(k)}).\psi
\\
\dMapsTo<{\psi(j_1)\tdt\psi(j_k)\tens1} 
&&&&\dMapsTo>{(-1)^{\sum_{q<r}^{c\le d}j_q^cj_r^d}}
\\
(-1)^{\sum_{q<r}|j_q|(1-|j_r|)}x(j_1)\tdt x(j_k)\tens m^{(k)} 
\hspace*{-5em} &&\overset{\tilde\rho}\longmapsto \hspace*{-1em} 
&&\hspace*{-7em}
(-1)^{\sum_{q=1}^k(k-q)|j_q|-\sum_{q<r}|j_q|\cdot|j_r|}x(j_1,\dots,j_k)
\end{diagram}
wherefore
\[ \rho(\sff_{j_1},\dots,\sff_{j_k};m^{(k)}).\psi =
(-1)^{\sum_{q=1}^k(k-q)|j_q|-\sum_{q<r}^{c>d}j_q^cj_r^d}x(j_1,\dots,j_k)
\]
and \(\deg\psi(j)=|j|\) as claimed.
We conclude that for this $\psi$ and $\chi$ the induced (by
\lemref{lem-PB=hg-QD}) right action of $\As$ on $\tilde{f}_n$ is the 
natural one.

Let us compute the differential \(\tilde{\partial}\) in $\tilde{f}_n$.
For \(\ell=j_1+\dots+j_k\) the expression
\begin{multline*}
(-1)^{|\ell|}\rho(\sff_{j_1},\dots,\sff_{j_k};m^{(k)}).
\bar{\partial}\psi
\\
=\sum_{q=1}^k \sum_{t+p=j_q}^{t,p\ne0}
(-1)^{\sum_{r=q+1}^k(1-|j_r|)+\sum_{c>d}t^cp^d+|p|+1+|\ell|}
\rho(\sff_{j_1},\dots,\sff_{j_{q-1}},\sff_t,\sff_p,\sff_{j_{q+1}},
\dots,\sff_{j_k};m^{(k+1)}).\psi
\\
\hskip\multlinegap =\sum_{q=1}^k \sum_{t+p=j_q}^{t,p\ne0}
(-1)^{\sum_{r=q+1}^k(1-|j_r|) +\sum_{c>d}t^cp^d+|p|+1+|\ell| 
+\sum_{r=1}^{q-1}(k+1-r)|j_r| +(k+1-q)|t| +(k-q)|p|} \hfill
\\
\times(-1)^{\sum_{r=q+1}^k(k-r)|j_r| -\sum_{u<r}^{c>d}j_u^cj_r^d 
+\sum_{c>d}t^cp^d} x(j_1,\dots,j_{q-1},t,p,j_{q+1},\dots,j_k)
\end{multline*}
has to coincide with
\[ \rho(\sff_{j_1},\dots,\sff_{j_k};m^{(k)}).\psi \tilde{\partial}
=(-1)^{\sum_{q=1}^k(k-q)|j_q|-\sum_{q<r}^{c>d}j_q^cj_r^d}
x(j_1,\dots,j_k).\tilde{\partial}.
\]
This gives the differential \(\tilde{\partial}\):
\begin{equation}
x(j_1,\dots,j_k).\tilde\partial =\sum_{q=1}^k (-1)^{k+1-q} 
\sum_{t+p=j_q}^{t,p\ne0} x(j_1,\dots,j_{q-1},t,p,j_{q+1},\dots,j_k).
\label{eq-y(j1jk)-sum}
\end{equation}

Note that the differential
	\(\bar{\partial}:\overline{\sff}_n^{\{k\}}(\ell)
	\to\overline{\sff}_n^{\{k+1\}}(\ell)\)
makes
\begin{equation}
0 \rTTo \overline{\sff}_n^{\{1\}}(\ell) \rTTo^{\bar\partial} \dots 
\rTTo^{\bar\partial} \overline{\sff}_n^{\{k\}}(\ell) 
\rTTo^{\bar\partial} \dots
\rTTo^{\bar\partial} \overline{\sff}_n^{\{|\ell|\}}(\ell) \rTTo 0
\label{eq-0F1ddFkddFl0}
\end{equation}
into a bounded complex of abelian groups. 
The term \(\overline{\sff}_n^{\{k\}}(\ell)\) is placed in degree
\(k-|\ell|\).

Consider the operad morphism $\As\to\ZZ$, \(m^{(k)}\mapsto0\) for 
$k\ge2$, where $\ZZ$ is the unit operad, \(\ZZ(1)=\ZZ\), \(\ZZ(n)=0\) 
for $n\ne1$.
We may view \(\oplus_{\ell\in\NN^n-0}\overline{\sff}_n^{\{k\}}(\ell)\) 
as a left $n$\n-operad $\ZZ$\n-module, quotient of 
\(\overline{\sff}_n^{\{k\}}\) by the submodule spanned by images of all 
left actions of elements \(m^{(k)}\) for $k\ge2$.
Applying the same quotient procedure to $\FAs_n$ we get
\[ \overline\FAs_n(\ell) =
\begin{cases}
\ZZ =\ZZ u(\ell) =\ZZ u(e_{\supp\ell}),
&\quad \text{ if } |\ell| =|\supp\ell|,
\\
0, &\quad \text{ if } |\ell| >|\supp\ell|.
\end{cases}
\]
We are going to prove that complex~\eqref{eq-0F1ddFkddFl0} is homotopy 
isomorphic via $\bar{p}$ and $\bar{\beta}$ to its cohomology 
\(\overline\FAs_n(\ell)\). 
If \(\ell=e_S\) for some \(S\subset\bn\), then the cohomology is 
concentrated in degree~0 and equals
$\overline\FAs_n(e_S)=\ZZ=\ZZ u(e_S)$.
If \(|\ell|>|\supp\ell|\), then the cohomology vanishes.

We construct mappings of abelian groups 
	\(h:\overline{\sff}_n^{\{p\}}(\ell)
	\to\overline{\sff}_n^{\{p-1\}}(\ell)\)
such that
	\(\bar{N}:\overline{\sff}_n^{\{k\}}(\ell)
	\to\overline{\sff}_n^{\{k\}}(\ell)\)
given by \eqref{eq-N1-pd-hd-dh} vanishes. 
These $h$ induce the left $n\text-\!\As$\n-module morphism
$h:\overline{\rmF}_n\to\overline{\rmF}_n$ compatibly with the 
generator--to--generator mapping 
\(\overline{\sff}_n^{\{k\}}\to\overline{\rmF}_n^{(k)}\). 
Thus, vanishing of
	\(\bar{N}:\overline{\sff}_n^{\{k\}}(\ell)
	\to\overline{\sff}_n^{\{k\}}(\ell)\)
implies vanishing of 
\(\bar{N}:\overline{\rmF}_n^{\{\}}\to\overline{\rmF}_n^{\{\}}\) and 
local nilpotency of \(N:\overline{\rmF}_n\to\overline{\rmF}_n\). 

We have reduced the proposition to proving that the chain maps $\bar p$, 
$\bar\beta$ in
\begin{diagram}[LaTeXeqno,w=2em]
0 &\rTTo &\overline{\sff}_n^{\{1\}}(\ell) &\rTTo^{\bar\partial} &\dots 
&\rTTo^{\bar\partial} &\overline{\sff}_n^{\{k\}}(\ell) 
&\rTTo^{\bar\partial} 
&\dots &\rTTo^{\bar\partial} &\overline{\sff}_n^{\{|\ell|-1\}}(\ell) 
&\rTTo^{\bar\partial} &\overline{\sff}_n^{\{|\ell|\}}(\ell) &\rTTo &0
\\
&&&&&& &&&&&& \uTTo<{\bar\beta} \dTTo>{\bar p}
\\
&&&&& &&&&& 0 &\rTTo &\overline{\FAs}_n(\ell) &\rTTo &0
\label{dia-Fn1l-Fnkl-FAsnl}
\end{diagram}
are homotopy inverse to each other for any \(\ell\in\NN^n-0\).
We add formally the case of $\ell=0$ by defining the top and the bottom 
rows as complexes \(\overline{\sff}_n^{\{0\}}(0)=\ZZ\) and 
\(\overline{\FAs}_n(0)=\ZZ\) concentrated in degree~0.
Here $\bar p$, $\bar\beta$ are defined as the identity maps.

Chain maps $\bar p$, $\bar\beta$ give rise to other chain maps
$\tilde p$, $\tilde\beta$ in the commutative diagram
\begin{diagram}[h=1.4em,w=2em]
0\to &\overline{\sff}_n^{\{1\}}(\ell) &\rTTo^{\bar\partial} &\dots 
&\rTTo^{\bar\partial} &\overline{\sff}_n^{\{k\}}(\ell) 
&\rTTo^{\bar\partial} &\dots &\rTTo^{\bar\partial} 
&\overline{\sff}_n^{\{|\ell|\}}(\ell) &\to0
\\
&\dTTo<\psi &&&&\dTTo<\psi &&&&\dTTo<\psi
&\rdFromTo^{\bar p}_{\bar\beta}
\\
&&&&& &&&&&0\to &\overline{\FAs}_n(\ell) &\to0
\\
0\to &\tilde{f}_n(\ell)^1 &\rTTo^{\tilde\partial} &\dots 
&\rTTo^{\tilde\partial} &\tilde{f}_n(\ell)^k &\rTTo^{\tilde\partial} 
&\dots &\rTTo^{\tilde\partial} &\tilde{f}_n(\ell)^{|\ell|} &\to0 &\dEq
\\
&&&&& &&&&&\rdFromTo^{\tilde p}_{\tilde\beta}
\\
&&&&& &&&&&0\to &\overline{\FAs}_n(\ell) &\to0
\end{diagram}
In fact, the maps $\tilde p$, $\tilde\beta$ have to be defined if 
\(l^i\in\{0,1\}\) for all $1\le i\le n$.
For \(k=|\ell|\) we find
\begin{equation}
x(e_{a_1},\dots,e_{a_k}).\tilde p= \sign(a_1,\dots,a_k) u(\ell)
\overset{\text{def}}= (-1)^{\sum_{q<r}\chi(a_q>a_r)} u(\ell),
\label{eq-x(ea1eak)p}
\end{equation}
where \(\chi(b>c)\) is 1 or 0 depending on the case whether the 
inequality holds or not.
The exponent is the number of inversions in the sequence
$(a_1,\dots,a_k)$.
If \(\ell=e_{\bn}=(1,1,\dots,1)\), then $k=n$ and the sign is just the 
sign of the permutation $(a_1,\dots,a_k)$.
The map
	$\tilde\beta=\bar\beta\psi:\overline{\FAs}_n(\ell)
	\to\tilde{f}_n(\ell)^{|\ell|}$
satisfies
\begin{equation}
u(\ell).\tilde\beta =x(e_{c_1},\dots,e_{c_k}),
\label{eq-u(l)beta}
\end{equation}
where \(\{c_1<c_2<\dots<c_k\}=\supp\ell\). 
In particular, if \(\ell=e_{\bn}\), then $k=n$ and 
\(u(e_{\bn})\tilde\beta=x(e_1,\dots,e_n)\).

Consider the augmented coalgebra \(C_n=\ZZ\{\NN^n\}\) with the 
comultiplication \(j.\Delta=\sum_{q+r=j}q\tens r\) where 
\(j,q,r\in\NN^n\).
Generators $j\in\NN^n$ of the free abelian group $C_n$ are denoted also 
$x(j)$.
The augmentation is \(\eta:\ZZ\to C_n\), $1\mapsto x(0)$.
The counit is $\eps:C_n\to\ZZ$, $x(j)\mapsto\delta_{j0}$.
The reduced comultiplication is defined as
\[ \bar\Delta =\Delta -\eta\tens\id -\id\tens\eta +\eps\eta\tens\eta, 
\qquad \vec0.\bar\Delta =0, \quad 
j.\bar\Delta=\sum_{q+r=j}^{q,r\ne0}q\tens r \quad \text{for } j\ne0.
\]
The abelian subgroup $\bar C_n=\Ker\eps=\ZZ\{\NN^n-0\}\subset C_n$ 
equipped with the comultiplication $\bar\Delta$ is a coassociative 
coalgebra, which is not counital.
The complex \(\tilde{f}_n\) is nothing else but the cohomology complex 
$K'(\bar C_n)$ of the coassociative coalgebra $\bar C_n$, which is the 
upper row of the diagram
\begin{diagram}[w=2em,LaTeXeqno]
0\to\!\!\!\!\! &\ZZ &\rTTo^0 &\bar C_n &\rTTo^{\bar\Delta}
&\bar C_n^{\tens2} &\rTTo^{\tilde\partial} &\dots 
&\rTTo^{\tilde\partial} &\bar C_n^{\tens k} &\rTTo^{\tilde\partial} 
&\dots &\rTTo^{\tilde\partial} &\bar C_n^{\tens n} 
&\rTTo^{\tilde\partial} &\dots
\\
&\dEq &&\uTTo<{\tilde\beta}\dTTo>{\tilde p} 
&&\uTTo<{\tilde\beta}\dTTo>{\tilde p} 
&&&&\uTTo<{\tilde\beta}\dTTo>{\tilde p} 
&&&&\uTTo<{\tilde\beta}\dTTo>{\tilde p}
\\
0\to\!\!\!\!\! &\ZZ &\rTTo^0 &\wedge_\ZZ^1(\ZZ^n) &\rTTo^0 
&\wedge_\ZZ^2(\ZZ^n) &\rTTo^0 &\dots &\rTTo^0 &\wedge_\ZZ^k(\ZZ^n) 
&\rTTo^0 &\dots &\rTTo^0 &\wedge_\ZZ^n(\ZZ^n) &\rTTo^0 &0
\label{dia-ZCnCn2-ZL1ZnL2Zn}
\end{diagram}
We identify \(x(j_1,\dots,j_k)\in\tilde{f}_n\) with
\(j_1\tdt j_k\in\bar C_n^{\tens k}\).
The exterior algebra
$\wedge_\ZZ(\ZZ^n)=T_\ZZ(\ZZ^n)/(x\tens x\mid x\in\ZZ^n)$ has the basis 
	$(e_{\{c_1<c_2<\dots<c_k\}}
	=e_{c_1}\wedge e_{c_2}\wedge\dots\wedge e_{c_k})$,
where \(1\le c_1<c_2<\dots<c_k\le n\).
The mappings in this diagram are
\begin{gather*}
(j_1\tdt j_k).\tilde\partial =\sum_{q=1}^k (-1)^{k-q+1}
j_1\tdt j_{q-1}\tens j_q.\bar\Delta\tens j_{q+1}\tdt j_k,
\\
x(j_1,\dots,j_k).\tilde p =0 \text{ \ unless \ } 
\Big|\sum_{q=1}^kj_q\Big| =k= \Big|\supp\sum_{q=1}^kj_q\Big|,
\quad j_q\in\NN^n-0,
\\
x(e_{a_1},\dots,e_{a_k}).\tilde p
=(-1)^{\sum_{q<r}\chi(a_q>a_r)} e_{\{a_1,\dots,a_k\}},
\quad a_1,\dots,a_k \text{ -- distinct},
\\
e_{\{c_1<c_2<\dots<c_k\}}.\tilde\beta =x(e_{c_1},\dots,e_{c_k}),
\end{gather*}
which coincides with \eqref{eq-y(j1jk)-sum}, \eqref{eq-x(ea1eak)p} and 
\eqref{eq-u(l)beta}.

It remains to prove that the maps $\tilde p$, $\tilde\beta$ are homotopy 
inverse to each other.
Clearly, $\tilde\beta\tilde p=1$.

\begin{lemma}\label{lem-p-beta-homotopy-inverse}
For $n=1$ the maps $\tilde p$, $\tilde\beta$ are homotopy inverse to 
each other.
\end{lemma}

\begin{proof}
For $n=1$ the chain maps in question become
\begin{diagram}[w=2em]
0\to\!\!\!\!\! &\ZZ &\rTTo^0 &\bar C_1 &\rTTo^{\bar\Delta}
&\bar C_1^{\tens2} &\rTTo^{\tilde\partial} &\dots 
&\rTTo^{\tilde\partial} &\bar C_1^{\tens k} &\rTTo^{\tilde\partial} 
&\dots 
\\
&\dEq &&\uTTo<{\tilde\beta}\dTTo>{\tilde p} &&\uTTo \dTTo
&&&&\uTTo \dTTo
\\
0\to\!\!\!\!\! &\ZZ &\rTTo^0 &\ZZ &\rTTo &0 &\rTTo &\dots &\rTTo &0 
&\rTTo &\dots 
\end{diagram}
\[ x(j).\tilde p =\delta_{j1}, \qquad 1.\tilde\beta =x(1).
\]

Define a map of graded abelian groups \(h:K'(\bar C_1)\to K'(\bar C_1)\) 
of degree $-1$ by the formula
\[ x(j_1,\dots,j_{k-1},j_k).h =
\begin{cases}
x(j_1,\dots,j_{k-2},j_{k-1}+1), &\quad \text{if } k>1,\, j_k=1,
\\
0, &\quad \text{otherwise}.
\end{cases}
\]
We claim that the chain map
	\(E=\tilde p\tilde\beta-h\tilde\partial-\tilde\partial h:
	K'(\bar C_1)\to K'(\bar C_1)\)
is the identity map.
In fact, \(x(1).E=x(1)\), and for $j\ge2$ we have
\[ x(j).E =\sum_{q+r=j}^{q,r>0} x(q,r).h =x(j-1+1) =x(j).
\]
For $k>1$ and $j_k\ge2$ we find
\[ x(j_1,\dots,j_k).E =\sum_{q+r=j_k} x(j_1,\dots,j_{k-1},q,r).h 
=x(j_1,\dots,j_{k-1},j_k).
\]
It remains to consider for $k>1$ the value
\begin{multline*}
x(j_1,\dots,j_{k-1},1).E =-x(j_1,\dots,j_{k-2},j_{k-1}+1).\tilde\partial 
-x(j_1,\dots,j_{k-1},1).\tilde\partial h
\\
=\sum_{q+r=j_{k-1}+1} x(j_1,\dots,j_{k-1},q,r)
-\sum_{q+t=j_{k-1}} x(j_1,\dots,j_{k-2},q,t,1).h
=x(j_1,\dots,j_{k-2},j_{k-1},1).
\end{multline*}
Hence, $E=\id$, and \(\tilde p\tilde\beta\) is homotopic to the identity 
map.
\end{proof}

\subsubsection{Homology of augmented algebras.}
Let $\cc$ be a symmetric monoidal category with the tensor product 
$\tens$ and the unit object $\1$.
Assume that \(A=(A,\mu:A\tens A\to A,\eta:\1\to A,\eps:A\to\1)\) is an 
augmented unital associative algebra in $\cc$.
There is an associated simplicial object $S(A)$:
\[ \cdots A\tens A\tens A 
\begin{diagram}[w=3.2em,inline]
&\pile{\rTTo^{d_3=1\tens1\tens\eps} \\ \lTTo~{s_2=1\tens1\tens\eta} 
\\ \rTTo~{d_2=1\tens\mu} \\ \lTTo~{s_1=1\tens\eta\tens1} 
\\ \rTTo~{d_1=\mu\tens1} \\ \lTTo~{s_0=\eta\tens1\tens1} 
\\ \rTTo_{d_0=\eps\tens1\tens1}} &
\end{diagram}
A\tens A
\begin{diagram}[w=2.6em,inline]
&\pile{\rTTo^{d_2=1\tens\eps} \\ \lTTo~{s_1=1\tens\eta}
\\ \rTTo~{d_1=\mu}
\\ \lTTo~{s_0=\eta\tens1} \\ \rTTo_{d_0=\eps\tens1}} &
\end{diagram}
A 
\begin{diagram}[w=2em,inline]
&\pile{\rTTo^{d_1=\eps} \\ \lTTo~{s_0=\eta} \\ \rTTo_{d_0=\eps}} &
\end{diagram}
\1,
\]
where $d_i$ and $s_i$ are face maps and degeneracy maps respectively.
When $B$ is another augmented algebra in $\cc$, the Cartesian product 
\(S(A)\times S(B)\) of simplicial objects 
\cite[Section~VIII.8]{MacLane-Homology} is naturally isomorphic to the 
simplicial object $S(A\tens B)$.

Assume also that $\cc$ is abelian and the tensor product $\tens$ is 
bilinear.
A complex \(K(A)\overset{\text{def}}=K(S(A))\) is associated with the 
simplicial object $S(A)$. 
It has the differential 
\(\partial=\sum_{i=0}^q(-1)^id_i:A^{\tens q}\to A^{\tens q-1}\).
Homology of \(K(A)\) gives the torsion objects \(\Tor^A_\bull(\1,\1)\), 
where the left and the right $A$\n-module $\1$ obtains its structure via 
$\eps:A\to\1$.
Given two augmented algebras $A$ and $B$ in $\cc$ we can form a 
bisimplicial object in $\cc$, whose terms are
\(A^{\tens p}\tens B^{\tens q}\).
By Eilenberg--Zilber theorem \cite[Theorem~8.5.1]{Weibel} the complexes 
\(K(A\tens B)=K(S(A\tens B))\simeq K(S(A)\times S(B))\) and
\(K(A)\tens K(B)\) are quasi-isomorphic.

Consider associative algebras 
\(\bar A=(\bar A,\mu:\bar A\tens\bar A\to\bar A)\) in $\cc$, which are 
not required to have a unit.
Such an algebra gives rise to a unital one \(A=\1\oplus\bar A\) for 
which \(\eta=\inj_\1:\1\to A\) is the unit and \(\eps=\pr_\1:A\to\1\) is 
an augmentation.
Introduce another monoidal product in $\cc$ (not bilinear) via the 
formula
\[ \bar A\circledast\bar B =\bar A\oplus\bar B\oplus(\bar A\tens\bar B).
\]
There is an obvious isomorphism
\[ \1\oplus(\bar A\circledast\bar B) 
= (\1\oplus\bar A)\tens(\1\oplus\bar B).
\]
If $\bar A$, $\bar B$ are associative algebras in $(\cc,\tens)$, then 
$\bar A\circledast\bar B$ obtains an associative algebra structure in 
$(\cc,\tens)$ via this isomorphism, namely, 
\(\1\oplus(\bar A\circledast\bar B)=A\tens B\).

There is a normalised chain complex $K_N(\bar A)=K_N(S(A))$ of the 
simplicial complex $S(A)$:
\[ \dots \rTTo \bar A^{\tens3} \rTTo^{-\mu\tens1+1\tens\mu}
\bar A^{\tens2} \rTTo^{-\mu} \bar A \rTTo^0 \1 \to 0,
\]
with the differential 
 \(\partial=\sum_{i=0}^{q-2}(-1)^{i+1}1^{\tens i}\tens\mu
 \tens1^{\tens q-i-2}:\bar A^{\tens q}\to\bar A^{\tens q-1}\),
where $\1$ is placed in degree 0. 
By a generalization of normalization theorem of Eilenberg and Mac Lane 
\cite[Theorem~VIII.6.1]{MacLane-Homology} the natural projection 
\(K(A)\to K_N(\bar A)\) is a homotopy isomorphism.
As a corollary, we get the following

\begin{proposition}\label{pro-K(AB)-K(A)K(B)}
For associative algebras $\bar A$, $\bar B$ in $(\cc,\tens)$ there is a 
natural quasi-isomorphism 
\[ K_N(\bar A\circledast\bar B) \rightleftarrows
K_N(\bar A)\tens K_N(\bar B).
\]
\end{proposition}

\subsubsection{Conclusion of the proof of
	\thmref{thm-Fn-FAsn-cofibrant-replacement}.}
Let us take for $(\cc,\tens)$ the category $(\Ab^\op,\tens_\ZZ^\op)$ 
opposite to the category of abelian groups with the opposite tensor 
product.
Clearly, the category $\Ab^\op$ is abelian.
An associative algebra in this monoidal category is a coassociative 
coalgebra over $\ZZ$ in the ordinary sense.
In particular, such is $\bar C_n=\ZZ\{\NN^n-0\}$.
Adding a unit to it in $\Ab^\op$ gives \(C_n=\ZZ\{\NN^n\}\).
Since \(C_n\tens_\ZZ C_m\simeq C_{n+m}\), we conclude that 
\(\bar C_n\circledast\bar C_m\simeq\bar C_{n+m}\).
The homological complex \(K_N(\bar C_n)\) in $\Ab^\op$ and the top line 
of \eqref{dia-ZCnCn2-ZL1ZnL2Zn}, the cohomological complex
$K'(\bar C_n)$ in $\Ab$ are identified: the $n$\n-th abelian groups and 
the differentials between them coincide.
Thus, the results of the previous section apply to $K'(\bar C_n)$.

We claim that the cohomology of $K'(\bar C_n)$ is isomorphic to 
\(\wedge_\ZZ(\ZZ^n)\).
In fact, using induction we deduce from
\lemref{lem-p-beta-homotopy-inverse} and \propref{pro-K(AB)-K(A)K(B)} 
the quasi-isomorphism
\[ K'(\bar C_{n+1}) \rTTo_{qis} K'(\bar C_n)\tens K'(\bar C_1) 
\rTTo^{q\tens\tilde p} \wedge_\ZZ(\ZZ^n)\tens\wedge_\ZZ(\ZZ) \simeq 
\wedge_\ZZ(\ZZ^{n+1}).
\]
Here $q$, $\tilde p$ are quasi-isomorphisms. 
So is their tensor product \(q\tens\tilde p\), since complexes
$K'(\bar C_n)$, \(\wedge_\ZZ(\ZZ^n)\) consist of free abelian groups, 
\(\wedge_\ZZ(\ZZ^n)\) are bounded and $K'(\bar C_n)$ are direct sums of 
bounded complexes.

Both rows of diagram~\eqref{dia-ZCnCn2-ZL1ZnL2Zn} have the same 
homology, which coincides with the bottom row and consists of finitely 
generated free abelian groups.
Since $H(\tilde\beta)H(\tilde p)=1$, the matrices of $H(\tilde\beta)$ 
and $H(\tilde p)$ are invertible.
Thus, $\tilde\beta$ and $\tilde p$ induce isomorphisms in homology.
They are quasi-isomorphisms of complexes consisting of free abelian 
groups.
Therefore, their cones are acyclic complexes consisting of free abelian 
groups.
They split into short exact sequences whose terms are also free abelian 
groups (as subgroups of such).
Hence, these short exact sequences split and the cones are contractible.
Thus, $\tilde p$ and $\tilde\beta$ are homotopy isomorphisms.
Clearly, they are homotopy inverse to each other.
This implies the same conclusion for $\bar p$ and $\bar\beta$ and for 
$p$ and $\beta$.
\end{proof}

\begin{corollary}[to \propref{pro-Fn-FAsn-rmFn},
\thmref{thm-Fn-FAsn-cofibrant-replacement}]
The polymodule $F_n$ is homotopy isomorphic to its cohomology and 
\(H^\bull(F_n(j))=\kk[1-|j|]\) for \(j\in\NN^n-0\).
\end{corollary}

This is due to existence of a degree 1 isomorphism 
\(\varSigma:H^\bull(F_n)\to\FAs_n\).

\subsection{Homotopy unital \texorpdfstring{$A_\infty$}{A8}-morphisms.}
Consider the free $n\wedge1$-$\mainf^\su$-module
\[ \tilde{\rmF}_n=\bigcirc_{i=1}^n\mainf^\su\odot^i_{\mainf}\rmF_n 
\odot^0_{\mainf}\mainf^\su
=\odot_{\ge0}(\sS{^n}\mainf^\su;\kk\{\sff_j\mid j\in\NN^n-0\};
\mainf^\su).
\]
In particular, \(\tilde{\rmF}_0=\mainf^\su(0)=\kk\one^\su\) by 
\lemref{lem-((Ai)A(0))-initial}.
The graded ideal generated by the following system of relations in it
\[ \rho_\emptyset(\one^\su) =\lambda^i_{e_i}(\one^\su;\sff_{e_i}),
\ \forall\,i, \quad
\lambda^i_\ell(\sS{^a}1,\one^\su,\sS{^b}1;\sff_\ell) =0 
\text{ \ if \ } a+1+b=\ell^i, \; |\ell|>1,
\]
is stable under the differential, as one easily verifies.
Therefore the quotient \(\rmF_n^\su\) of \(\tilde{\rmF}_n\) by these 
relations is an $n\wedge1$-$\mainf^\su$-module.
We still have \(\rmF_0^\su=\mainf^\su(0)=\kk\one^\su\).
Note that \(\rmF_n^\su\)\n-algebra maps coincide with \emph{strictly 
unital \ainfm-algebra morphisms}, which are by
\cite[Definition~9.2]{BesLyuMan-book} \ainfm-morphisms
$\sff:(A_1,\dots,A_n)\to B$ between strictly unital \ainfm-algebras such 
that all components of $\sff$ vanish if any of its entries is 
$\one^\su_{A_i}$, except $\one^\su_{A_i}\sff_{e_i}=\one^\su_B$.

The rows of the following diagram in \(\dg^{\NN\sqcup\NN^n}\)
\begin{diagram}[LaTeXeqno]
0 &\rTTo &(\mainf,\rmF_n) &\rTTo &(\mainf^\su,\rmF_n^\su) &\rTTo 
&(\kk\one^\su,\kk\one^\su\rho_\emptyset) &\rTTo &0
\\
&& \dFib<{htis} &&\dFib<{htis}>{p'} &&\dEq
\\
0 &\rTTo &(\As,\FAs_n) &\rTTo &(\Ass,\FAss_n) &\rTTo
&(\kk\one^\su,\kk\one^\su\rho_\emptyset) &\rTTo &0
\label{dia-(AF)-(AsuFsu)-(k1k1)}
\end{diagram}
are exact sequences, split in the obvious way. 
Therefore, the middle vertical arrow $p'$ is a homotopy isomorphism.

Consider the embedding of free graded operads 
\(\mainf^\su\to\mainf^\su\langle\sfi,\sfj\rangle\), where $\sfi$, $\sfj$ 
are two nullary operations, $\deg \sfi=0$, $\deg \sfj=-1$.
Assuming \(\sfi\partial=0\), \(\sfj\partial=\one^\su-\sfi\), we make the 
second operad differential graded and the embedding becomes a chain map.
It is proven in \cite{Lyu-Ainf-Operad} (end of proof of Proposition~1.8) 
that this embedding is a homotopy isomorphism.
Or, the reader can simplify the lines of the proof given below and adopt 
it to the case of \(\mainf^\su\to\mainf^\su\langle\sfi,\sfj\rangle\).

\begin{proposition}
The embedding
	\(\iota:(\mainf^\su,\rmF_n^\su)
	\to(\mainf^\su,\rmF_n^\su)\langle\sfi,\sfj\rangle\)
is a homotopy isomorphism.
\end{proposition}

\begin{proof}
An arbitrary chain $n\wedge1$-module map
\(\phi:(\mainf^\su,\rmF_n^\su)\langle\sfi,\sfj\rangle\to(\ca,\cp)\) is 
fixed by specifying a chain $n\wedge1$-module map
\((\mainf^\su,\rmF_n^\su)\to(\ca,\cp)\) and the image 
\(\phi(\sfj)\in\ca(0)^{-1}\).
In particular, there is a unique chain $n\wedge1$-module map
\[ \pi: (\mainf^\su,\rmF_n^\su)\langle\sfi,\sfj\rangle \to
(\mainf^\su,\rmF_n^\su), \qquad \sfi \mapsto \one^\su, 
\qquad \sfj \mapsto 0,
\]
whose restriction to \((\mainf^\su,\rmF_n^\su)\) is identity.
Let us prove that $\pi$ is homotopy inverse to $\iota$.

The restrictions of the above chain maps
\(\iota':\kk\one^\su\hookrightarrow\kk\{\one^\su,\sfi,\sfj\}\) and
\(\pi':\kk\{\one^\su,\sfi,\sfj\}\to\kk\one^\su\), 
\(\one^\su\mapsto\one^\su\), \(\sfi\mapsto\one^\su\), \(\sfj\mapsto0\), 
are homotopy isomorphisms: the homotopy
\(h:\kk\{\one^\su,\sfi,\sfj\}\to\kk\{\one^\su,\sfi,\sfj\}\), 
\(\one^\su.h=0\), \(\sfi.h=\sfj\), \(\sfj.h=0\), satisfies
\(\partial h+h\partial=\pi'\iota'-1\).
We know from \propref{pro-A<Mbeta>} that the $n\wedge1$-module 
\((\mainf^\su,\rmF_n^\su)\langle\sfi,\sfj\rangle\) coincides with 
	\((\mainf^\su\langle\sfi,\sfj\rangle,
	\cp=\bigcirc_{k=0}^n\mainf^\su\langle\sfi,
	\sfj\rangle\odot^k_{\mainf^\su}\rmF_n^\su)\).
The component $\cp(l)$ is obtained from components \(\rmF_n^\su(r)\) 
with \(r^q\ge l^q\) for all \(q\in\bn\) by plugging the unused
\(r^q-l^q\) entries with $\sfi$ and $\sfj$ in all possible ways 
determined by injections \(\psi^q:\mb l^q\hookrightarrow\mb r^q\):
\[ \cp(l) =\bigoplus_{(\psi^q:\mb l^q\hookrightarrow\mb r^q)_{q=1}^n} 
\Bigl( \bigotimes_{q\in\bn} \bigotimes_{\mb r^q-\im\psi^q}
\kk\{\sfi,\sfj\} \Bigr)\tens \rmF_n^\su(r).
\]

There is a split surjection
\[ \lambda: \cq(l)
=\bigoplus_{(\psi^q:\mb l^q\hookrightarrow\mb r^q)_{q=1}^n} \Bigl(
\bigotimes_{q\in\bn}
\bigotimes_{\mb r^q-\im\psi^q} (\kk\{\sfi,\sfj\}\oplus\kk\one^\su) 
\Bigr)\tens\rmF_n^\su(r) \rTTo \cp(l),
\]
obtained by acting with all \(\one^\su\) on \(\rmF_n^\su\) on the left 
via $\lambda$.
This reduces the quantities $r^q$ by the number of factors \(\one^\su\).

Denote
	\(f=\pi'\iota',g=\id:\kk\{\one^\su,\sfi,\sfj\}
	\to\kk\{\one^\su,\sfi,\sfj\}\).
Equip the set \(S=\sqcup_{q\in\bn}(\mb r^q-\im\psi^q)\) with the 
lexicographic order, $q<y\in\bn$ implies \((q,c)<(y,z)\).
The maps $\partial$ and $\pi\iota$ satisfy
\begin{diagram}
\cq(l)
&\rTTo^{\oplus_{(\psi^q)_q}(1\tens\partial+\sum_{(y,z)\in S}
	(\tens^{(q,c)<(y,z)}1)\tens\partial\tens(\tens^{(q,c)>(y,z)}1)
	\tens1)}_{\wh\partial}
&\cq(l)
\\
\dEpi<\lambda &= &\dEpi>\lambda
\\
\cp(l) &\rTTo^\partial &\cp(l)
\end{diagram}
\begin{diagram}
\cq(l)
&\rTTo^{\oplus_{(\psi^q)_q}(\tens^{(q,c)\in S}f)\tens1}_{\wh{\pi\iota}} 
&\cq(l)
\\
\dEpi<\lambda &= &\dEpi>\lambda
\\
\cp(l) &\rTTo^{\pi\iota} &\cp(l)
\end{diagram}
Since \(\one^\su.f=\one^\su.g=\one^\su\), $\one^\su.h=0$, there is a 
unique map \(H:\cp(l)\to\cp(l)\) of degree $-1$ such that
\begin{diagram}
\cq(l)
&\rTTo^{\oplus_{(\psi^q)_q}\sum_{(y,z)\in S}(\tens^{(q,c)<(y,z)}f)\tens 
	h\tens(\tens^{(q,c)>(y,z)}g)\tens1}_{\wh H}
&\cq(l)
\\
\dEpi<\lambda &= &\dEpi>\lambda
\\
\cp(l) &\rTTo^H &\cp(l)
\end{diagram}
In order to find the commutator \(\partial H+H\partial\) we can compute
\begin{multline*}
\wh\partial\wh H+\wh H\wh\partial
=\bigoplus_{(\psi^q:\mb l^q\hookrightarrow\mb r^q)_{q=1}^n}
\sum_{(y,z)\in S}\bigl(\bigotimes_{(q,c)<(y,z)}f\bigr)\tens
	(f-g)\tens\bigl(\bigotimes_{(q,c)>(y,z)}g\bigr)\tens1
\\
=\bigoplus_{(\psi^q:\mb l^q\hookrightarrow\mb r^q)_{q=1}^n}\bigl
	(\bigotimes_{(q,c)\in S}f-\bigotimes_{(q,c)\in S}g\bigr)\tens1
=\wh{\pi\iota} -1.
\end{multline*}
Therefore, \(\partial H+H\partial=\pi\iota-1\).
\end{proof}

The projection $p'$ decomposes into a standard trivial cofibration and 
an epimorphism $p''$
\[ p' =\bigl( (\mainf^\su,\rmF_n^\su) \rCof~\Sim^{htis}
(\mainf^\su,\rmF_n^\su)\langle \sfi,\sfj\rangle \rTTo^{p''}
(\Ass,\FAss_n) \bigr),
\]
where \(p''(\one^\su)=\one^\su\), \(p''(\sfi)=\one^\su\), 
\(p''(m_2)=m_2\), \(p''(\sff_{e_i})=1\in\FAss_n(e_i)\), and other 
generators go to 0. 
As a corollary \(p''(\one^\su\rho_\emptyset)=\one^\su\rho_\emptyset\). 
Hence, the projection $p''$ is a homotopy isomorphism as well.

Generators $\sff_\ell$ of the $n\wedge1$-operad module $\rmF_n$ are 
interpreted as maps \(\sff_\ell:\boxt^{k\in\bn}T^{\ell^k}A_k\to B\) of 
degree \(\deg\sff_\ell=1-|\ell|\).
A cofibrant replacement \((\mainf^{hu},\rmF_n^{hu})\to(\Ass,\FAss_n)\) 
is constructed as a $\gr$-submodule of 
\((\mainf^\su,\rmF_n^\su)\langle\sfi,\sfj\rangle\) generated in operadic 
part by $\sfi$ and $g$\n-ary operations of degree $4-g-2k$
\[ m_{g_1;g_2;\dots;g_k} =(1^{\tens g_1}\tens\sfj\tens1^{\tens g_2}\tens
\sfj\tdt1^{\tens g_{k-1}}\tens \sfj\tens1^{\tens g_k})m_{g+k-1},
\]
where \(g=\sum_{q=1}^kg_q\), $k\ge1$, $g_q\ge0$, \(g+k\ge3\) and in 
module part by the nullary elements 
	\(\sfv_k=\lambda^k_{e_k}(\sfj;
	\sff_{e_k})-\sfj\rho_\emptyset=\sfj\sff_{e_k}-\sfj\rho_\emptyset\),
$k\in\bn$, \(\deg \sfv_k=-1\), and by elements 
\begin{multline}
\sff_{(\ell^k_1;\ell^k_2;\dots;\ell^k_{t^k})_{k\in\bn}} 
=\lambda_{\hat\ell}\bigl((\sS{^{\ell^k_1}}1,\sfj,\sS{^{\ell^k_2}}1,\sfj,
\dots,\sS{^{\ell^k_{t^k-1}}}1,\sfj,\sS{^{\ell^k_{t^k}}}1)_{k\in\bn}; 
\sff_{\hat\ell}\bigr)
\\
=\bigl[\boxt^{k\in\bn}T^{\ell^k}A_k
\rTTo^{\boxt^{k\in\bn}(1^{\tens\ell^k_1}\tens\sfj\tens1^{\tens\ell^k_2} 
	\tens\sfj\tdt1^{\tens\ell^k_{t^k-1}}\tens\sfj
	\tens1^{\tens\ell^k_{t^k}})}
\boxt^{k\in\bn}T^{\hat\ell^k}A_k \rTTo^{\sff_{\hat\ell}} B\bigr]
\label{eq-f(llll)=lambda-(1j1j1j1)f}
\end{multline}
of arity \(\ell=\bigl(\sum_{p=1}^{t^k}\ell^k_p\bigr)_{k\in\bn}\), where 
the intermediate arity is
	\(\hat\ell=\bigl(t_k-1+\sum_{p=1}^{t^k}\ell^k_p\bigr)_{k\in\bn}
	=\bigl(-1+\sum_{p=1}^{t^k}(\ell^k_p+1)\bigr)_{k\in\bn}\),
and of degree
	\(\deg\sff_{(\ell^k_1;\dots;\ell^k_{t^k})_{k\in\bn}}
	=1+2n-\sum_{k=1}^n\sum_{p=1}^{t^k}(\ell^k_p+2)\).
We assume that $t_k\ge1$ for all $k\in\bn$ and either 
	\(|\hat\ell|
	=\sum_{k=1}^n(t^k-1)+\sum_{k=1}^n\sum_{p=1}^{t^k}\ell^k_p\ge2\),
or all $t_k=1$ and there is $m\in\bn$ such that \(\ell^k_1=\delta^k_m\).
The last condition eliminates from the list the summands 
\(\sff_{0,\dots,0,(0;0),0,\dots,0}=\sfj\sff_{e_k}\) of $\sfv_k$.
Setting \(\sfi\partial=0\), \(\sfj\partial=\one^\su-\sfi\), we turn 
\((\mainf^\su,\rmF_n^\su)\langle \sfi,\sfj\rangle\) into a $\dg$-module 
and \((\mainf^{hu},\rmF_n^{hu})\) into its $\dg$-submodule.
Note that \(\sfv_k\partial=\sfi\rho_\emptyset-\sfi\sff_{e_k}\).

Let us prove that the graded $n\wedge1$-module
$(\mainf^{hu},\rmF_n^{hu})$ is free over $(\kk,0)$. 
The graded $n\wedge1$-module \((\mainf,\rmF_n)\langle \sfj\rangle\) can 
be presented as
\begin{multline}
\bigl(\mainf, \odot_{\ge0}(\sS{^{[n]}}\mainf;
\kk\{\sff_\ell \mid \ell\in\NN^n-0\})\bigr)\langle \sfj\rangle
\\
\simeq \bigl( \mainf\langle \sfj\rangle, \odot_{\ge0}(
\sS{^n}\kk\langle m_{n_1;\dots;n_k}
\mid {k+\ttt\sum_{q=1}^k} n_q\ge3\rangle; 
\kk\{\sff_{(\ell^k_1;\dots;\ell^k_{t^k})_{k\in\bn}}\mid|\hat\ell|\ge1\}; 
\mainf\langle \sfj\rangle) \bigr).
\label{eq-(AAokfoA)j}
\end{multline}
The free graded $n\wedge1$-operad module generated by
\(m_{n_1;\dots;n_k}\) and
\(\sff_{(\ell^k_1;\dots;\ell^k_{t^k})_{k\in\bn}}\) has the form
\begin{multline*}
K= F\bigl(\kk\{m_{n_1;\dots;n_k} \mid {k+\ttt\sum_{q=1}^k} n_q\ge3\},
\kk\{\sff_{(\ell^k_1;\dots;\ell^k_{t^k})_{k\in\bn}}
\mid|\hat\ell|\ge1\}\bigr)
\\
\hskip\multlinegap
=\bigl(\kk\langle m_{n_1;\dots;n_k} 
\mid {k+\ttt\sum_{q=1}^k} n_q\ge3\rangle,
\hfill
\\
\odot_{\ge0}(
\sS{^{[n]}}\kk\langle m_{n_1;\dots;n_k}
\mid{k+\ttt\sum_{q=1}^k} n_q\ge3\rangle; 
\kk\{\sff_{(\ell^k_1;\dots;\ell^k_{t^k})_{k\in\bn}}
\mid |\hat\ell|\ge1\}) \bigr).
\end{multline*}
It is a direct summand of \eqref{eq-(AAokfoA)j}, so we have a split 
exact sequence in \(\gr^{\NN\sqcup\NN^n}\)
\[ 0 \rTTo K \pile{\rTTo^\alpha \\ \lTTo_\pi}
(\mainf,\rmF_n)\langle \sfj\rangle
\pile{\rTTo^\varkappa \\ \lTTo_\omega} (\kk \sfj,\kk \sfj\rho_\emptyset) 
\rTTo 0,
\]
where $\omega$ takes $\sfj\rho_\emptyset$ to the nullary generator 
\(\sfj\rho_\emptyset\).
Consider also the graded $n\wedge1$-module 
\[ L= F\bigl(\kk\{m_{n_1;\dots;n_k} \mid {k+\ttt\sum_{q=1}^k} n_q\ge3\},
\kk\{\sfv_k, \sff_{e_k}, \sff_{(\ell^k_1;\dots;\ell^k_{t^k})_{k\in\bn}} 
\mid |\hat\ell|\ge2 \} \bigr).
\]
Notice that the map $L\to K$, \(\sfv_k\mapsto \sfj\sff_{e_k}\), which 
maps other generators identically, identifies the $n\wedge1$-modules $L$ 
and $K$.

Consider the graded module morphism
\begin{align*}
\beta: L 
=\bigl(\kk\langle m_{n_1;\dots;n_k}
\mid {k+\ttt\sum_{q=1}^k} n_q\ge3\rangle,\overline{L}\bigr) 
&\to (\mainf,\rmF_n)\langle \sfj\rangle
=\bigl(\mainf\langle \sfj\rangle,
\overline{(\mainf,\rmF_n)\langle \sfj\rangle}\bigr), 
\\
\sfv_k &\mapsto \sfj\sff_{e_k} -\sfj\rho_\emptyset,
\end{align*}
which maps other generators identically.
The morphism $\beta$ extends to basic elements so that each factor 
$\sfj\rho_\emptyset$ arising from a vertex of type $\sfv_k$ gives its 
$\sfj$ to subsequent $m_{;;;}$ adding another semicolon to its indexing 
sequence.
This follows by associativity of $\rho$.
The basic elements $\sfv_k$ are mapped by \(\beta\varkappa\) to 
$-\sfj\rho_\emptyset$.
For any other basic element \(b(t)\in L\) we have
\(b(t).\beta\varkappa=0\).

The map \(\beta-\beta\varkappa\omega\in\gr^{\NN\sqcup\NN^n}\) factors 
through $\alpha$ as the following diagram shows:
\begin{diagram}
0 &\rTTo &K &\pile{\rTTo^\alpha \\ \lTTo_\pi}
&(\mainf,\rmF_n)\langle \sfj\rangle
&\pile{\rTTo^\varkappa \\ \lTTo_\omega}
&(\kk \sfj,\kk \sfj\rho_\emptyset) &\rTTo &0
\\
&&\uTTo<{\exists!\gamma} &\ruTTo>{\beta-\beta\varkappa\omega} &
\\
&&L
\end{diagram}
The unique map
 \(\gamma=(\beta-\beta\varkappa\omega)\pi=\beta\pi:
 L\to K\in\gr^{\NN\sqcup\NN^n}\), 
such that \(\beta-\beta\varkappa\omega=\gamma\alpha\), has a triangular 
matrix.
In fact, $L$ and $K$ have an $\NN$\n-grading, \(L^q=\oplus\kk b(t)\), 
\(K^q=\oplus\kk b(t)\), where the summation is over forests $t$ with $q$ 
vertices labelled by one of $\sfv_k$ (resp. one of \(\sfj\sff_{e_k}\)).
The map $\gamma$ takes the filtration \(L_q=L^0\oplus\dots\oplus L^q\) 
to the filtration \(K_q=K^0\oplus\dots\oplus K^q\). 
The diagonal entries \(\gamma^{qq}:L^q\to K^q\) are identity maps.
Thus, the matrix of \(\gamma\) equals \(1-N\), where $N$ is locally 
nilpotent, and $\gamma$ is invertible.
We obtained a split exact sequence
\begin{equation}
0 \rTTo L \rTTo^{\beta-\beta\varkappa\omega}
(\mainf,F_1)\langle \sfj\rangle
\pile{\rTTo^\varkappa \\ \lTTo_\omega} (\kk \sfj,\kk \sfj\rho_\emptyset) 
\rTTo 0.
\label{eq-0L-beta-beta-kappa-omega-(AF)-kk0}
\end{equation}
Let us decompose the first two terms into direct sums 
\begin{gather*}
\overline{L}=\kk\{\sfv_s\mid s\in\bn\}\oplus
\bigl(\overline{L}\ominus\kk\{\sfv_s\mid s\in\bn\}\bigr),
\\
\overline{(\mainf,\rmF_n)\langle \sfj\rangle}
=\kk\{\sfj\rho_\emptyset,\sfj\sff_{e_s} \mid s\in\bn\}\oplus
\bigl(\overline{(\mainf,\rmF_n)\langle \sfj\rangle}\ominus
\kk\{\sfj\rho_\emptyset,\sfj\sff_{e_s} \mid s\in\bn\}\bigr),
\end{gather*}
where the complements are spanned by all basic elements except those 
listed in the first summands.
The maps $\beta$ and \(\beta-\beta\varkappa\omega\) preserve this 
decomposition.
Their restriction to the second summand coincide and this is an 
isomorphism due to exactness of
\eqref{eq-0L-beta-beta-kappa-omega-(AF)-kk0}.
If we drop out complements, this split exact sequence takes the form
\[ 0 \rTTo \kk\{\sfv_s\mid s\in\bn\} \rTTo^{\beta-\beta\varkappa\omega} 
\kk\{\sfj\rho_\emptyset,\sfj\sff_{e_s} \mid s\in\bn\}
\pile{\rTTo^\varkappa \\ \lTTo_\omega} \kk \sfj\rho_\emptyset \rTTo 0,
\]
where \(\sfv_s.(\beta-\beta\varkappa\omega)=\sfj\sff_{e_s}\).
Let us replace it with another split exact sequence
\[ 0 \rTTo \kk\{\sfv_s\mid s\in\bn\} \pile{\rTTo^\beta \\ \lTTo_\tau}
\kk\{\sfj\rho_\emptyset,\sfj\sff_{e_s} \mid s\in\bn\}
\pile{\rTTo^\theta \\ \lTTo_\omega} \kk \sfj\rho_\emptyset \rTTo 0,
\]
where 
\(\sfj\sff_{e_s}.\theta=\sfj\rho_\emptyset.\theta=\sfj\rho_\emptyset\) 
and \(\sfj\rho_\emptyset.\tau=0\), \(\sfj\sff_{e_s}.\tau=\sfv_s\).
Restoring back the dropped isomorphism of second summands we obtain from 
the above the split exact sequence
\[ 0 \rTTo L \rTTo^\beta (\mainf,F_1)\langle \sfj\rangle
\pile{\rTTo^\theta \\ \lTTo_\omega}
(\kk \sfj,\kk \sfj\rho_\emptyset) \rTTo 0,
\]
such that $\theta$ vanishes on the complement.

Adding freely $\sfi$ we deduce the split exact sequence in 
\(\gr^{\NN\sqcup\NN^n}\)
\begin{equation}
0 \to L\langle \sfi\rangle \to (\mainf,\rmF_n)\langle \sfj,\sfi\rangle
\to (\kk \sfj,\kk \sfj\rho_\emptyset) \to 0.
\label{eq-Li-(AF)ji-(kjkjr)}
\end{equation}
The image of the embedding is precisely $(\mainf^{hu},\rmF_n^{hu})$, 
thus the latter graded $n\wedge1$-module is free.
In the particular case of $n=0$ the module part is generated by the 
empty set of generators.
Therefore, \(\rmF_0^{hu}=\mainf^{hu}(0)\) by
\lemref{lem-((Ai)A(0))-initial}.

Furthermore, from the top row of 
diagram~\eqref{dia-(AF)-(AsuFsu)-(k1k1)} we deduce a splittable exact 
sequence in \(\gr^{\NN\sqcup\NN^n}\)
\[ 0 \to (\mainf,\rmF_n)\langle \sfi\rangle \to 
(\mainf^\su,\rmF_n^\su)\langle \sfi\rangle \to
(\kk\one^\su,\kk\one^\su\rho_\emptyset) \to 0.
\]
We may choose the splitting of this exact sequence as indicated below:
\[ 0 \to (\mainf,\rmF_n)\langle \sfi\rangle \to 
(\mainf^\su,\rmF_n^\su)\langle \sfi\rangle \to
(\kk\{\one^\su-\sfi\},\kk\{\one^\su\rho_\emptyset-\sfi\rho_\emptyset\}) 
\to 0.
\]
Adding freely $\sfj$ we get the split exact sequence
\begin{equation}
0 \to (\mainf,\rmF_n)\langle \sfi,\sfj\rangle \to 
(\mainf^\su,\rmF_n^\su)\langle \sfi,\sfj\rangle \to
(\kk\{\one^\su-\sfi\},\kk\{(\one^\su-\sfi)\rho_\emptyset\}) \to 0.
\label{eq-(AF)ij-(AsuFsu)ij-(k1ik1ir)}
\end{equation}
Combining \eqref{eq-Li-(AF)ji-(kjkjr)} with 
\eqref{eq-(AF)ij-(AsuFsu)ij-(k1ik1ir)} we get a split exact sequence 
\begin{equation}
0 \to (\mainf^{hu},\rmF_n^{hu}) \xrightarrow{i'}
(\mainf^\su,\rmF_n^\su)\langle \sfi,\sfj\rangle \to
(\kk\{\one^\su-\sfi,\sfj\},\kk\{(\one^\su-\sfi)\rho_\emptyset,
\sfj\rho_\emptyset\}) \to 0.
\label{eq-0-(Ahu8F1hu)-ij-0}
\end{equation}

The differential in $(\mainf^{hu},\rmF_n^{hu})$ is computed through that 
of \((\mainf^\su,\rmF_n^\su)\langle \sfi,\sfj\rangle\). 
Actually, \eqref{eq-0-(Ahu8F1hu)-ij-0} is a split exact sequence in 
\(\dg^{\NN\sqcup\NN^n}\), where the third term obtains the differential 
\(\sfj.\partial=\one^\su-\sfi\), 
	\(\sfj\rho_\emptyset.\partial
	=\one^\su\rho_\emptyset-\sfi\rho_\emptyset\). 
The third term is contractible, which shows that the inclusion $i'$ is a 
homotopy isomorphism in \(\dg^{\NN\sqcup\NN^n}\). 
Hence, the epimorphism 
\(p=i'\cdot p'':(\mainf^{hu},\rmF_n^{hu})\to(\Ass,\FAss_n)\) is a 
homotopy isomorphism as well.

In order to prove that $(\1,0)\to(\mainf^{hu},\rmF_n^{hu})$ is a 
standard cofibration we present it as a colimit of sequence of 
elementary cofibrations
\[ (\1,0)\to\cd_0=F(\kk\{\sfi,m_2\},\kk\{\sff_{e_s} \mid s\in\bn\})
\to\cd_1\to\cd_2\to\dots,
\]
where for $r>0$
\[ \cd_r =F\bigl(
\kk\{\sfi,m_{n_1;\dots;n_k} \mid \deg m_{n_1;\dots;n_k} \ge-r\},
\kk\{\sfv_s,\sff_{(\ell^k_1;\dots;\ell^k_{t^k})_{k\in\bn}} \mid s\in\bn,
\; \deg\sff_{(\ell^k_1;\dots;\ell^k_{t^k})} \ge-r\} \bigr).
\]

Algebra maps over $(\mainf^{hu},\rmF_n^{hu})$ are identified with 
homotopy unital \ainf-morphisms, which we define in the spirit of 
Fukaya's approach:

\begin{definition}\label{def-homotopy-unital-structure-A8-morphism}
A \emph{homotopy unital structure} of an \ainfm-morphism
$\sff:A_1,\dots,A_n\to B$ is an \ainfm-morphism 
	\(\sff^+:(A_k^+)_{k\in\bn}
	=(A_k\oplus\kk\one^\su_{A_k}\oplus\kk\sfj^{A_k})_{k\in\bn}
	\to B\oplus\kk\one^\su_B\oplus\kk \sfj^B=B^+\)
between given homotopy unital \ainfm-algebra structures such that:
\begin{myitemize}
\item{(1)} $\sff^+$ is a strictly unital: for all \(1\le k\le n\)
\[ \one^\su_{A_k}\sff_{e_k}^+=\one^\su_B, \quad \bigl[
1^{\boxt(k-1)}\boxt(1^{\tens a}\tens\one^\su_{A_k}\tens1^{\tens b})
\boxt1^{\boxt(n-k)}\bigr]\sff_\ell^+=0
\text{ \ if \ } a+1+b=\ell^k,\;|\ell|>1.
\]

\item{(2)} the element \(\sfv_k^B=\sfj^{A_k}\sff_{e_k}^+-\sfj^B\) is 
contained in $B$;

\item{(3)} the restriction of $\sff^+$ to $A_1,\dots,A_n$ gives $\sff$;

\item{(4)}
	\(\bigl[\boxt^{k\in\bn}(A_k\oplus\kk \sfj^{A_k})^{\tens\ell^k}\bigr]
	\sff_\ell^+\subset B\),
for each \(\ell\in\NN^n\), $|\ell|>1$.
\end{myitemize}
\end{definition}

Homotopy unital structure of an \ainfm-morphism $f$ means a 
\emph{choice} of such $f^+$.
There is another notion of unitality which is a \emph{property} of an 
\ainfm-morphism:

\begin{definition}[See {\cite[Proposition~9.13]{BesLyuMan-book}}]
An \ainfm-morphism $\sff:A_1,\dots,A_n\to B$ between unital
\ainfm-algebras is \emph{unital} if the cycles \(\sfi^{A_k}\sff_{e_k}\) 
and \(\sfi^B\) differ by a boundary for all \(1\le k\le n\).
\end{definition}

For a homotopy unital \ainfm-morphism $\sff:A_1,\dots,A_n\to B$ the 
equation holds 
	\(\sfv_k^Bm_1=\sfv_k\partial=\sfi\rho_\emptyset-\sfi\sff_{e_k}
	=\sfi^B-\sfi^{A_k}\sff_{e_k}\). 
Thus an \ainfm-morphism with a homotopy unital structure is unital.

\begin{conjecture}
Unitality of an \ainfm-morphism is equivalent to homotopy unitality: any 
unital \ainfm-morphism admits a homotopy unital structure.
\end{conjecture}

All reasoning of this section can be applied to $F_n$ in place of 
$\rmF_n$.
A nullary degree $-1$ cycle $\bone^\su$ subject to relations 
\eqref{eq-(1bone)b2-1-(bone1)b2-1} is added to $A_\infty$.
The resulting operad is denoted $A_\infty^\su$.
We consider the $A_\infty^\su$\n-module
\[ \tilde{F}_n=\bigcirc_{i=1}^nA_\infty^\su\odot^i_{A_\infty}F_n 
\odot^0_{A_\infty}A_\infty^\su
=\odot_{\ge0}(\sS{^n}A_\infty^\su;\kk\{f_j\mid j\in\NN^n-0\};
A_\infty^\su).
\]
It is divided by the graded ideal generated by the following system of 
relations
\[ \rho_\emptyset(\bone^\su) =\lambda^i_{e_i}(\bone^\su;f_{e_i}),
\ \forall\,i, \quad
\lambda^i_\ell(\sS{^a}1,\bone^\su,\sS{^b}1;f_\ell) =0 
\text{ \ if \ } a+1+b=\ell^i, \; |\ell|>1.
\]
The quotient is denoted $F_n^\su$.
Similarly to the above we add two nullary operations $\bi$, $\bj$ to 
$A_\infty^\su$ with $\deg\bi=-1$, $\deg\bj=-2$, \(\bi\partial=0\), 
\(\bj\partial=\bi-\bone^\su\).
The obtained \(A_\infty^\su\langle\bi,\bj\rangle\)-module 
\(F_n^\su\langle\bi,\bj\rangle\) contains an $A_\infty^{hu}$\n-submodule 
$F_n^{hu}$ spanned by the nullary elements
	\(\bv_k=\lambda^k_{e_k}(\bj;f_{e_k})-\bj\rho_\emptyset
	=\bj f_{e_k}-\bj\rho_\emptyset\),
$k\in\bn$, \(\deg\bv_k=-2\), and by elements
\(f_{(\ell^k_1;\ell^k_2;\dots;\ell^k_{t^k})_{k\in\bn}}\) similar to 
\eqref{eq-f(llll)=lambda-(1j1j1j1)f}.
There are invertible operad module homomorphisms $\varSigma$ of degree 1 
sending $f_j\mapsto\sff_j$, $\bone^\su\mapsto\one^\su$, 
$\bi\mapsto\sfi$, $\bj\mapsto\sfj$, $\bv_k\mapsto\sfv_k$:
\begin{gather*}
\varSigma: (A_\infty^\su,F_n^\su) \to (\mainf^\su,\rmF_n^\su), \qquad
\varSigma: (A_\infty^\su,F_n^\su)\langle\bi,\bj\rangle \to
(\mainf^\su,\rmF_n^\su)\langle\sfi,\sfj\rangle,
\\
\varSigma: (A_\infty^{hu},F_n^{hu}) \to (\mainf^{hu},\rmF_n^{hu}).
\end{gather*}

\section{Composition of morphisms with several arguments}
The mechanism which provides an associative composition of morphisms 
with several arguments is that of convolution product in the module of 
linear maps from a coalgebra to an algebra.
The part of a coalgebra is played by a colax $\Cat$-span multifunctor.
A lax $\Cat$-span multifunctor $\HOM$ comes in place of an algebra.
The convolution product of these multifunctors gives a multicategory 
structure to the collection of \ainf-algebras and \ainf-morphisms with 
several arguments.

\subsection{Colax \texorpdfstring{$\Cat$}{Cat}-span multifunctors.}
\begin{definition}
A \emph{colax $\Cat$-span multifunctor}
	\((F,\psi^I):(\mcL,\circledast_\mcL^I,\lambda_\mcL^f,\rho_\mcL)
	\to(\mcM,\circledast_\mcM^I,\lambda_\mcM^f,\rho_\mcM)\) 
between lax $\Cat$-span multicategories is
\begin{enumerate}
\renewcommand{\labelenumi}{\roman{enumi})}
\item a 1\n-morphism
	\(F=(\tar F,F,\tar F):\mcL=(\tar\mcL,\mcL,\tar\mcL)
	\to(\tar\mcM,\mcM,\tar\mcM)=\mcM\)
in \(\smQuiver\);

\item a 2\n-morphism for each set $I\in\Ob\co_\sk$
\begin{diagram}
\boxdot^I\mcL &\rTTo^{\boxdot^IF} &\boxdot^I\mcM
\\
\dTTo<{\circledast_\mcL^I} &\ruTwoar^{\psi^I} 
&\dTTo>{\circledast_\mcM^I}
\\
\mcL &\rTTo^F &\mcM
\end{diagram}
\end{enumerate}
such that
\begin{equation*}
\begin{diagram}[inline,w=4em]
\boxdot^{\mb1}\mcL &\rTTo^{\boxdot^{\mb1}F} &\boxdot^{\mb1}\mcM
\\
\dTTo<{\circledast_\mcL^{\mb1}} &\ruTwoar^{\psi^{\mb1}} 
&\dTTo<{\circledast_\mcM^{\mb1}} \overset{\rho_\mcM}\Rightarrow 
\dTTo>{\Rho}
\\
\mcL &\rTTo^F &\mcM
\end{diagram}
\;=\;
\begin{diagram}[inline,w=4em]
\boxdot^{\mb1}\mcL &\rTTo^{\boxdot^{\mb1}F} &\boxdot^{\mb1}\mcM
\\
\dTTo<{\circledast_\mcL^{\mb1}} \overset{\rho_\mcL}\Rightarrow 
\dTTo>{\Rho} &= &\dTTo>{\Rho}
\\
\mcL &\rTTo^F &\mcM
\end{diagram}
\end{equation*}
and for every map $f:I\to J$ of $\co_\sk$ the following equation holds:
\begin{equation*}
\begin{diagram}[inline,nobalance,tight,width=2.6em]
\boxdot^I\mcL &&\rTTo^{\boxdot^IF} &&\boxdot^I\mcM
\\
&\rdTTo^{\circledast_\mcL^f} &&\ruTwoar^{\boxdot^{j\in J}\psi^{f^{-1}j}} 
&&\rdTTo^{\circledast_\mcM^f} &
\\
\dTTo<{\circledast_\mcL^I} &\rTwoar^{\lambda_\mcL^f} &\boxdot^J\mcL 
&&\rTTo^{\boxdot^JF} &&\boxdot^J\mcM
\\
&\ldTTo^{\circledast_\mcL^J} &&&&\ruTwoar(6,2)^{\psi^J} 
\ldTTo_{\circledast_\mcM^J} & 
\\
\mcL &&&\rTTo^F &\mcM
\end{diagram}
\qquad=\quad
\begin{diagram}[inline,nobalance,tight,width=2.6em]
\boxdot^I\mcL &\rTTo^{\boxdot^IF} &\boxdot^I\mcM
\\
&\ruTwoar(2,4)<{\psi^I} &&\rdTTo^{\circledast_\mcM^f} 
\\
\dTTo<{\circledast_\mcL^I} &&\dTTo<{\circledast_\mcM^I}  
&\rTwoar^{\lambda_\mcM^f} &\boxdot^J\mcM 
\\
&&&\ldTTo>{\circledast_\mcM^J} 
\\
\mcL &\rTTo^F &\mcM
\end{diagram}
\end{equation*}
\end{definition}

Here 2\n-morphism \(\boxdot^{j\in J}\psi^{f^{-1}j}\) means the pasting
\begin{diagram}
\boxdot^I\mcL &\rTTo^{\boxdot^IF} &\boxdot^I\mcM
\\
\dTTo<{\Lambda^f} &= &\dTTo>{\Lambda^f}
\\
\boxdot^{j\in J}\boxdot^{f^{-1}j}\mcL
&\rTTo^{\boxdot^{j\in J}\boxdot^{f^{-1}j}F}
&\boxdot^{j\in J}\boxdot^{f^{-1}j}\mcM
\\
\dTTo<{\boxdot^{j\in J}\circledast_\mcL^{f^{-1}j}} 
&\ruTwoar^{\boxdot^{j\in J}\psi^{f^{-1}j}} 
&\dTTo>{\boxdot^{j\in J}\circledast_\mcM^{f^{-1}j}}
\\
\boxdot^J\mcL &\rTTo^{\boxdot^JF} &\boxdot^J\mcM
\end{diagram}

We shall show that \ainf-modules $F_n$ form a polymodule cooperad, that 
is a colax $\Cat$\n-multifunctor \(F:\mcF\to\mcM\), where the category 
$\mcM$ of operad polymodules is described in
\defref{def-M-operad-polymodules}. 
Here (strict) $\Cat$\n-operad $\mcF$ has (1\n-element set of objects), 
\(\mcF(I)=\1\) is the terminal category for any $I\in\co_\sk$,
1\n-morphism \(\circledast^I:\boxdot^I\mcF\to\mcF\) is the unique one, 
2\n-morphisms $\lambda^f$ and $\rho$ are identity morphisms.
$\mcF$ is also a $\Cat$-span operad.

A general \emph{polymodule cooperad}, that is, a colax
$\Cat$\n-multifunctor \(F:\mcF\to\mcM\) (equivalently, a colax
$\Cat$-span multifunctor) amounts to the following data: an operad 
\(\ca=F(*)\), for each $I\in\co_\sk$ an $I\wedge1$-$\ca$\n-module $F_I$, 
for each tree \(t:[I]\to\co_\sk\) a morphism
	\(\Delta(t):F_{t(0)}\to
	\circledast_\mcM(t)(F_{t_h^{-1}b})_{(h,b)\in\IV(t)}\)
of $t(0)\wedge1$-$\ca$\n-modules such that for all corollas
\(t:[1]\to\co_\sk\) 
\begin{equation}
\bigl(F_{t(0)} \rTTo^{\Delta(t)} \circledast_\mcM(t)(F_{t(0)}) 
\rTTo^\sim F_{t(0)}\bigr) =1,
\label{eq-Ft(0)-Delta(t)-1}
\end{equation}
for all non-decreasing maps \(f:I\to J\), the induced
\(\psi=[f]:[J]\to[I]\) as in \eqref{eq-x<fy-fx<y}, and for all trees
\(t:[I]\to\co_\sk\)
\begin{diagram}[nobalance,LaTeXeqno]
F_{t(0)} &\rTTo^{\Delta(t_\psi)}
&\circledast_\mcM(t_\psi)(F_{t_{\psi,g}^{-1}c})_{(g,c)\in\IV(t_\psi)} 
\hspace*{3em}
\\
\dTTo<{\Delta(t)} &=
&\dTTo>{\circledast_\mcM(t_\psi)(\Delta(t^{|c}_{[\psi(g-1),\psi(g)]}))}
\\
\hspace*{-3em} \circledast_\mcM(t)(F_{t_h^{-1}b})_{(h,b)\in\IV(t)} &\rTTo^{\lambda^f}
&\circledast_\mcM(t_\psi)\bigl(
\circledast_\mcM(t^{|c}_{[\psi(g-1),\psi(g)]})
(F_{t_h^{-1}b})_{(h,b)\in\IV(t^{|c}_{[\psi(g-1),\psi(g)]})} 
\bigr)_{(g,c)\in\IV(t_\psi)}
\label{dia-Delta-Delta-Delta-lambda}
\end{diagram}
We are interested in the cases of \(\ca=A_\infty\), $\mainf$, 
\(A_\infty^{hu}\) and \(\mainf^{hu}\).

\begin{exercise}
Write down explicitly equation~\eqref{dia-Delta-Delta-Delta-lambda} for 
the both non-decreasing surjections \(\psi:[2]\to[1]\) and a tree
\(t:[1]\to\co_\sk\).
Conclude that for the tree \(r:[0]\to\co_\sk\) the operad
$\ca$\n-bimodule map \(\Delta(r):F_1\to\ca\) plays the part of a counit 
for $\Delta$.
\end{exercise}

Viewing (system of $\ca$\n-modules) \(F:\mcF\to\mcM\) as a coalgebra and 
$\HOM:\mcB\to\mcM$ (coming from a symmetric multicategory $\mcC$) as an 
algebra we consider homomorphisms between them (in the sense of $\mcM$) 
and they have to form an algebra as well.
So we define a multiquiver $\mcH$ whose objects are $\ca$\n-algebras 
\((B,\alpha_B:\ca\to\END B)\) with
\begin{multline*}
\mcH((A_i,\alpha_{A_i})_{i\in I};(B,\alpha_B)) 
\\
=\{ ((\alpha_{A_i})_{i\in I};\phi;\alpha_B) \in 
\mcM\bigl((\sS{^I}\ca;F_I;\ca),
((\END A_i)_{i\in I};\hoM((A_i)_{i\in I};B);\END B)\bigr) \}.
\end{multline*}
Let us define a multicategory composition for it.
For any tree $t$ and any collection of $\ca$\n-algebras 
\(\alpha_h^b:\ca\to\END A_h^b\), \((h,b)\in\iV(t)\), assume given 
\(t_h^{-1}b\wedge1\)\n-operad module morphisms for \((h,b)\in\IV(t)\):
\begin{multline*}
\bigl((\alpha_{h-1}^a)_{a\in t_h^{-1}b};g_h^b;\alpha_h^b\bigr): 
\bigl(\sS{^{t_h^{-1}b}}\ca;F_{t_h^{-1}b};\ca\bigr) 
\\
\to \bigl((\END A_{h-1}^a)_{a\in t_h^{-1}b}; 
\hoM((A_{h-1}^a)_{a\in t_h^{-1}b};A_h^b); \END A_h^b\bigr).
\end{multline*}
Then their composition is defined as 
\(\bigl((\alpha_0^a)_{a\in t(0)};\comp(t)(g_h^b);
\alpha_{\max[I]}^1\bigr)\), where
\begin{multline}
\comp(t)(g_h^b) =\bigl[F_{t(0)} \rTTo^{\Delta(t)} 
\circledast_\mcM(t)(F_{t_h^{-1}b})_{(h,b)\in\IV(t)} 
\rTTo^{\circledast_\mcM(t)(g_h^b)}
\\
\circledast_\mcM(t)
\bigl(\HOM((A_{h-1}^a)_{a\in t_h^{-1}b};A_h^b)\bigr)_{(h,b)\in\IV(t)} 
\rTTo^{\comp(t)} \HOM((A_0^a)_{a\in t(0)};A_{\max[I]}^1) \bigr].
\label{eq-comp(t)(ghb)-defined}
\end{multline}

\begin{proposition}\label{pro-composition-H-strictly-associative}
The composition in $\mcH$ is strictly associative.
\end{proposition}

\begin{proof}
For any $f:I\to J$ and \((y,c)\in\IV(t_\psi)\) in notation from 
\secref{sec-Cat-operad-graded-k-modules} denote 
\begin{multline*}
h_y^c =\bigl[F_{t_{\psi,g}^{-1}c}
\rTTo^{\Delta(t^{|c}_{[\psi(y-1),\psi y]})}
\circledast_\mcM(t^{|c}_{[\psi(y-1),\psi y]})
(F_{t_x^{-1}b})_{(x,b)\in\IV(t^{|c}_{[\psi(y-1),\psi y]})} 
\\
\rTTo^{\circledast_\mcM(t^{|c}_{[\psi(y-1),\psi y]})(g_x^b)}
\circledast_\mcM(t^{|c}_{[\psi(y-1),\psi y]})
\bigl(\HOM((A_{x-1}^a)_{a\in t_x^{-1}b};A_x^b)
\bigr)_{(x,b)\in\IV(t^{|c}_{[\psi(y-1),\psi y]})} 
\\
\rTTo^{\comp(t^{|c}_{[\psi(y-1),\psi y]})} 
\HOM((A_{\psi(y-1)}^a)_{a\in t_{\psi,g}^{-1}c};A_{\psi y}^c) \bigr].
\end{multline*}
We plug in this expression into
\begin{multline*}
\bigl[F_{t(0)} \rTTo^{\Delta(t_\psi)} 
\circledast_\mcM(t_\psi)(F_{t_{\psi,g}^{-1}c})_{(y,c)\in\IV(t_\psi)} 
\rTTo^{\circledast_\mcM(t_\psi)(h_y^c)}
\\
\hfill \circledast_\mcM(t_\psi) \bigl(
\HOM((A_{\psi(y-1)}^a)_{a\in t_{\psi,g}^{-1}c};
	A_{\psi y}^c)\bigr)_{(y,c)\in\IV(t_\psi)} 
\rTTo^{\comp(t_\psi)} \HOM((A_0^a)_{a\in t(0)};A_{\max[I]}^1) 
\bigr]\quad
\\
\hskip\multlinegap =\bigl[F_{t(0)} \rTTo^{\Delta(t_\psi)} 
\circledast_\mcM(t_\psi)(F_{t_{\psi,g}^{-1}c})_{(y,c)\in\IV(t_\psi)} 
\rTTo^{\circledast_\mcM(t_\psi)(\Delta(t^{|c}_{[\psi(y-1),\psi y]}))} 
\hfill
\\
\circledast_\mcM(t_\psi) \bigl(
\circledast_\mcM(t^{|c}_{[\psi(y-1),\psi y]})
(F_{t_x^{-1}b})_{(x,b)\in\IV(t^{|c}_{[\psi(y-1),\psi y]})} 
\bigr)_{(y,c)\in\IV(t_\psi)}
\rTTo^{\circledast_\mcM(t_\psi)
\circledast_\mcM(t^{|c}_{[\psi(y-1),\psi y]})(g_x^b)}
\\
\circledast_\mcM(t_\psi) \bigl(
\circledast_\mcM(t^{|c}_{[\psi(y-1),\psi y]})
\bigl(\HOM((A_{x-1}^a)_{a\in t_x^{-1}b};A_x^b)
\bigr)_{(x,b)\in\IV(t^{|c}_{[\psi(y-1),\psi y]})} 
\bigr)_{(y,c)\in\IV(t_\psi)}
\\
\rTTo^{\circledast_\mcM(t_\psi) \comp(t^{|c}_{[\psi(y-1),\psi y]})} 
\circledast_\mcM(t_\psi) \bigl(
\HOM((A_{\psi(y-1)}^a)_{a\in t_{\psi,g}^{-1}c};
	A_{\psi y}^c)\bigr)_{(y,c)\in\IV(t_\psi)} 
\\
\hfill \rTTo^{\comp(t_\psi)}
\HOM((A_0^a)_{a\in t(0)};A_{\max[I]}^1)\bigr]\quad
\\
\hskip\multlinegap =\bigl[F_{t(0)} \rTTo^{\Delta(t)} 
\circledast_\mcM(t)(F_{t_x^{-1}b})_{(x,b)\in\IV(t)} 
\rTTo^{\lambda_\mcM^f} \hfill
\\
\circledast_\mcM(t_\psi) \bigl(
\circledast_\mcM(t^{|c}_{[\psi(y-1),\psi y]})
(F_{t_x^{-1}b})_{(x,b)\in\IV(t^{|c}_{[\psi(y-1),\psi y]})} 
\bigr)_{(y,c)\in\IV(t_\psi)}
\rTTo^{\circledast_\mcM(t_\psi)
\circledast_\mcM(t^{|c}_{[\psi(y-1),\psi y]})(g_x^b)}
\\
\circledast_\mcM(t_\psi) \bigl(
\circledast_\mcM(t^{|c}_{[\psi(y-1),\psi y]})
\bigl(\HOM((A_{x-1}^a)_{a\in t_x^{-1}b};A_x^b)
\bigr)_{(x,b)\in\IV(t^{|c}_{[\psi(y-1),\psi y]})} 
\bigr)_{(y,c)\in\IV(t_\psi)}
\\
\rTTo^{\circledast_\mcM(t_\psi) \comp(t^{|c}_{[\psi(y-1),\psi y]})} 
\circledast_\mcM(t_\psi) \bigl(
\HOM((A_{\psi(y-1)}^a)_{a\in t_{\psi,g}^{-1}c};
	A_{\psi y}^c)\bigr)_{(y,c)\in\IV(t_\psi)} 
\\
\hfill \rTTo^{\comp(t_\psi)}
\HOM((A_0^a)_{a\in t(0)};A_{\max[I]}^1)\bigr] \quad
\\
\hskip\multlinegap =\bigl[F_{t(0)} \rTTo^{\Delta(t)} 
\circledast_\mcM(t)(F_{t_x^{-1}b})_{(x,b)\in\IV(t)} 
\rTTo^{\circledast_\mcM(t)(g_x^b)}
\circledast_\mcM(t) \bigl( \HOM((A_{x-1}^a)_{a\in t_x^{-1}b};A_x^b) 
\bigr)_{(x,b)\in\IV(t)} \hfill
\\
\rTTo^{\lambda_\mcM^f} 
\circledast_\mcM(t_\psi) \bigl(
\circledast_\mcM(t^{|c}_{[\psi(y-1),\psi y]})
\bigl(\HOM((A_{x-1}^a)_{a\in t_x^{-1}b};A_x^b)
\bigr)_{(x,b)\in\IV(t^{|c}_{[\psi(y-1),\psi y]})} 
\bigr)_{(y,c)\in\IV(t_\psi)}
\\
\rTTo^{\circledast_\mcM(t_\psi) \comp(t^{|c}_{[\psi(y-1),\psi y]})} 
\circledast_\mcM(t_\psi) \bigl(
\HOM((A_{\psi(y-1)}^a)_{a\in t_{\psi,g}^{-1}c};
A_{\psi y}^c)\bigr)_{(y,c)\in\IV(t_\psi)} 
\\
\hfill \rTTo^{\comp(t_\psi)}
\HOM((A_0^a)_{a\in t(0)};A_{\max[I]}^1)\bigr] \quad
\\
\hskip\multlinegap =\bigl[F_{t(0)} \rTTo^{\Delta(t)} 
\circledast_\mcM(t)(F_{t_x^{-1}b})_{(x,b)\in\IV(t)} 
\rTTo^{\circledast_\mcM(t)(g_x^b)}
\circledast_\mcM(t) \bigl( \HOM((A_{x-1}^a)_{a\in t_x^{-1}b};A_x^b) 
\bigr)_{(x,b)\in\IV(t)} \hfill
\\
\rTTo^{\comp(t)} \HOM((A_0^a)_{a\in t(0)};A_{\max[I]}^1) \bigr] 
=\comp(t)(g_h^b).
\end{multline*}
Here we have used \eqref{dia-Delta-Delta-Delta-lambda}, naturality of 
$\lambda_\mcM$, equation~\eqref{eq-comp(t)-decomposed} and 
definition~\eqref{eq-comp(t)(ghb)-defined}.
\end{proof}

Below we shall prove that the convolution $\mcH$ of \ainfm-module 
cooperad $\rmF$ and the lax $\Cat$\n-multifunctor $\HOM$ built from 
$\uCom$ gives a multicategory of \ainfm-algebras and \ainfm-morphisms.
Its objects are \ainfm-algebras and morphisms
\((A_i)_{i\in I}\to B\in\mcH\) are morphisms of $n\wedge1$-operad 
modules
\[ (\sS{^I}\mainf;\rmF_n;\mainf) \to
((\END A_i)_{i\in I};\HOM((A_i)_{i\in I};B);\END B),
\]
which are precisely \ainfm-morphisms with several arguments.
Their composition is the composition in $\mcH$.

\subsection{Comultiplication under \texorpdfstring{$A_\infty$}{A8}.}
Taking tensor coalgebra of a graded $\kk$-module gives morphism of 
multiquivers to multicategory of differential graded augmented counital 
coassociative coalgebras \(Ts:\mcH\rMono \dgac\).
We wish to define a colax $\Cat$-span multifunctor \(F:\mcF\to\mcM\) 
such that $Ts$ becomes a multifunctor.
The following statements follow from results of \cite{BesLyuMan-book}.

\begin{proposition}[See Proposition~6.8 of \cite{BesLyuMan-book}]
	\label{pro-Tge1-coalgebra}
A $T^{\ge1}$\n-coalgebra $C$ in $\dg$ is a coassociative coalgebra 
$(C,\overline\Delta:C\to C\tens C)$ in $\dg$ such that 
\begin{equation*}
C =\colim_{k\to\infty}
\Ker\bigl(\overline\Delta^{(k)}:C\to C^{\tens k}\bigr).
%\label{eq-C(XY)-colim-Ker(Delta)}
\end{equation*}
\end{proposition}

\begin{corollary}[See Corollary~6.11 of \cite{BesLyuMan-book}]
	\label{cor-Tge1-coalgebra}
Let $C$ be a $T^{\ge1}$\n-coalgebra in $\dg$, and let $B\in\Ob\dg$ be a 
complex. 
Then there is a natural bijection
\begin{gather*}
\dg_{T^{\ge1}}(C,T^{\ge1}B) \to \dg(C,B), \qquad
(f:C\to T^{\ge1}B) \mapsto
\bigl(C \rTTo^f T^{\ge1}B \rTTo^{\pr_1} B\bigr),
\end{gather*}
where \(\dg_{T^{\ge1}}\) is the category of $T^{\ge1}$\n-coalgebras in 
$\dg$.
\end{corollary}

\begin{proposition}[See Corollary~6.18 and Proposition~6.19 of 
\cite{BesLyuMan-book}]
The full and faithful functor
\begin{multline*}
T^{\le1}:\dg_{T^{\ge1}}\to\dgac, \qquad (C,\overline\Delta) \mapsto
\\
\bigl(\kk\oplus C,
\Delta_0 =\pr_1\cdot\overline\Delta\cdot(\inj_1\tens\inj_1) 
+\id\tens\inj_0 +\inj_0\tens\id -\pr_0\cdot(\inj_0\tens\inj_0), 
\eps=\pr_0, \eta=\inj_0\bigr)
\end{multline*}
makes \(\dg_{T^{\ge1}}\) into a symmetric Monoidal subcategory of 
$\dgac$.
\end{proposition}

\begin{corollary}
An arbitrary augmented $\dg$\n-coalgebra \(A=\tens^{i\in I}TA_i\) comes 
from a $T^{\ge1}$\n-coalgebra \(A\ominus\kk\) and there is a natural 
bijection
\begin{align*}
\dgac(\tens^{i\in I}TA_i,TB) &\to \dg((\tens^{i\in I}TA_i)\ominus\kk,B),
\\
(f:\tens^{i\in I}TA_i\to TB) &\mapsto
\bigl((\tens^{i\in I}TA_i)\ominus\kk \rTTo^{f|} T^{\ge1}B
\rTTo^{\pr_1} B\bigr).
\end{align*}
\end{corollary}

Let us denote by $\alinf$ the multiquiver of \ainfm-algebras and their 
morphisms.
It admits a full and faithful embedding \(Ts:\alinf\rMono \dgac\) into 
the multiquiver of differential graded augmented counital coassociative 
coalgebras over $\kk$.
Actually the latter is a multicategory and $\alinf$ is isomorphic to its 
submulticategory.
In this way $\rmF$ obtains a unique colax $\Cat$-span multifunctor 
structure \(\Delta(t)\).
Let us describe the details.

Objects of $\alinf=\mcH$, that is, \ainfm-algebras
\((B,\alpha_B:\mainf\to\END B)\) are taken by $Ts$ to the tensor 
coalgebra \((TsB,\Delta_0,\eps,\eta)\), where
\(TsB=\oplus_{n=0}^\infty T^nsB=\oplus_{n=0}^\infty(B[1])^{\tens n}\), 
$\Delta_0$ is the cut comultiplication
\[\Delta_0(x_1\tdt x_n)=\sum_{i=0}^n(x_1\tdt x_i)\tens(x_{i+1}\tdt x_n),
\]
\(\eps=\pr_0:TsB\rEpi T^0sB=\kk\) is the counit and
\(\eta=\inj_0:\kk=T^0sB\rMono TsB\) is the augmentation.

The differential \(b:TsB\to TsB\), \(\deg b=1\), has matrix entries 
\(b^{n,k}:T^nsB\to T^ksB\),
\[ b^{n,k}=\sum_{a+p+c=n}^{a+1+c=k}1^{\tens a}\tens b_p\tens1^{\tens c},
\quad\text{where}\quad 
b_p=(-1)^n(\sigma^{\tens p})^{-1}\cdot m_p\cdot\sigma: T^psB\to sB
\]
for $p\ge1$, $b_0=0$, and \(s:B\to sB=B[1]\), $x\mapsto x$, is the shift 
map (the suspension), \(\deg s=-1\), \(\deg b_p=1\).
Here \(m_p:T^pB\to B\), \(\deg m_p=2-p\), are linear maps representing 
generators \(m_p\in\mainf(p)\) for $p\ge2$ and $m_1:B\to B$ is the 
differential in the complex $B$.

Morphisms of $\alinf=\mcH$, \ainfm-algebra morphisms
\(\sff:(A_i)_{i\in\bn}\to B\), are taken by $Ts$ to the augmented 
coalgebra chain homomorphisms \(f:\tens^{i\in\bn}TsA_i\to TsB\), whose 
compositions with the projections \(\pr_l:TsB\to T^lsB\) are given by
\begin{multline}
f\cdot\pr_l =\bigl[ TsA_1\tdt TsA_n 
\rTTo^{\Delta_0^{(l)}\tdt\Delta_0^{(l)}}
(TsA_1)^{\tens l}\tdt(TsA_n)^{\tens l} 
\\
\rTTo^{\overline{\varkappa}_{n,l}} (TsA_1\tdt TsA_n)^{\tens l} 
\rTTo^{\check{f}^{\tens l}} (sB)^{\tens l} \bigr],
\label{eq-fprl-Delta0l-sigma-fl}
\end{multline}
where the restriction of $\check{f}$ to \(T^{j^1}sA_1\tdt T^{j^n}sA_n\) 
is given by the component
\begin{equation}
f_j 
=(\sigma^{\tens j^1}\tdt\sigma^{\tens j^n})^{-1}\cdot\sff_j\cdot\sigma: 
T^{j^1}sA_1\tdt T^{j^n}sA_n \to sB,
\label{eq-fj-(ss)fjs}
\end{equation}
$\sff_j$ being linear maps that represent the generators 
\(\sff_j\in\rmF_n(j^1,\dots,j^n)\).
The symmetry $\overline{\varkappa}_{n,l}=c_{s_{n,l}}$ corresponds to the 
permutation $s_{n,l}$ of the set $\{1,2,\dots,nl\}$,
\begin{equation*}
s_{n,l}(1+t+kl)=1+k+tn \text{ \ for \ } 0\le t<l,\; 0\le k<n.
%\label{eq-s-permute}
\end{equation*}
Detailed description of map~\eqref{eq-fprl-Delta0l-sigma-fl} on direct 
summands is
\begin{multline}
\Bigl[ \bigotimes_{a\in\bn}T^{j^a}sA_a \rTTo^{\tens^{\bn}\Delta_0^{(l)}}
\bigotimes_{a\in\bn} \bigoplus_{\sum_{q=1}^lr_q^a=j^a} 
\bigotimes_{p\in\mb l} T^{r_p^a}sA_a \rto\sim
\bigoplus_{\sum_{q=1}^lr_q^a=j^a} \bigotimes_{a\in\bn} 
\bigotimes_{p\in\mb l} T^{r_p^a}sA_a 
\\
\rTTo^{\oplus\overline{\varkappa}_{n,l}}
\bigoplus_{\sum_{q=1}^lr_q^a=j^a} \bigotimes_{p\in\mb l} 
\bigotimes_{a\in\bn} T^{r_p^a}sA_a
\rTTo^{\sum\tens^{p\in\mb l}f_{(r_p^a)_{a\in\bn}}}
\bigotimes_{p\in\mb l}sB = T^lsB \Bigr].
\label{eq-map-on-direct-summands}
\end{multline}

Due to reasoning after equation~(8.20.2) in \cite{BesLyuMan-book} the 
map
\[ Ts: \alinf((A_i)_{i\in\bn};B) \to \dgac(\tens^{i\in\bn}TsA_i,TsB), 
\quad \sff \mapsto f,
\]
is a bijection.
Thus, for an arbitrary tree \(t:[I]\to\co_\sk\) and arbitrary vertex 
\((h,b)\in\IV(t)\)
\[ Ts: \alinf((A_{h-1}^a)_{a\in t_h^{-1}b};A_h^b) \to 
\dgac(\tens^{a\in t_h^{-1}b}TsA_{h-1}^a,TsA_h^b), \quad 
g_h^b =(g_{h,j}^b)_{j\in\NN^n} \mapsto \hat{g}_h^b,
\]
is a bijection.
Define composition in $\alinf$ as
\begin{multline*}
\bigl[ \prod_{(h,b)\in\IV(t)} \alinf((A_{h-1}^a)_{a\in t_h^{-1}b};A_h^b) 
\rTTo^{\prod Ts}_\sim 
\prod_{(h,b)\in\IV(t)} \dgac(\tens^{a\in t_h^{-1}b}TsA_{h-1}^a,TsA_h^b) 
\\
\rTTo^{\comp_\dgac(t)} \dgac(\tens^{a\in t(0)}TsA_0^a,TsA_{\max[I]}^1) 
\rTTo_\sim \alinf((A_0^a)_{a\in t(0)};A_{\max[I]}^1) \bigr],
\end{multline*}
\[ (g_h^b)_{(h,b)\in\IV(t)} \mapsto \comp(t)(g_h^b), \quad \text{where}
\]
\begin{multline}
\comp(t)(g_h^b)_j =\bigl[ \tens^{a\in t(0)}T^{j^a}A_0^a 
\rTTo^{\tens^{a\in t(0)}\sigma^{\tens j^a}}
\tens^{a\in t(0)}T^{j^a}sA_0^a
\rTTo^{\tens^{b_1\in t(1)}\wh{g_1^{b_1}}} 
\tens^{b_1\in t(1)}T^{j_1^{b_1}}sA_1^{b_1}
\\
\rTTo^{\tens^{b_2\in t(2)}\wh{g_2^{b_2}}}
\tens^{b_2\in t(2)}T^{j_2^{b_2}}sA_2^{b_2} \to \dots 
\rTTo^{g_{\max[I]}^1} 
sA_{\max[I]}^1 \rTTo^{\sigma^{-1}} A_{\max[I]}^1 \bigr],
\label{eq-comp(t)(ghb)j}
\end{multline}
and $g_h^b$ is given via its components
\[ \sfg_{h,j}^b: \rmF_{t_h^{-1}b}(j) \to 
\und\dg(\tens^{a\in t_h^{-1}b}T^{j^a}A_{h-1}^a,A_h^b), 
\quad \sff_j \mapsto
\bigl( \sfg_{h,j}^b: \tens^{a\in t_h^{-1}b}T^{j^a}A_{h-1}^a 
\to A_h^b \bigr).
\]
Here $j$ belongs to \(\NN^{t_h^{-1}b}\).

\begin{figure}
\begin{equation}
%\boldmath
\resizebox{!}{.5\texthigh}{\rotatebox[origin=c]{90}{%
\begin{diagram}[height=4em,inline,nobalance]
F_{t(0)}(j) &\rTTo^{\!\!\!\!\varDelta(t)}
&\circledast_\mcG(t)(F_{t_h^{-1}b})_{(h,b)\in\IV(t)}(j)
&\rTTo^{\circledast_\mcG(t)(g_h^b)}
&\circledast_\mcG(t)
	(\hoM((sA_{h-1}^a)_{a\in t_h^{-1}b};sA_h^b))_{(h,b)\in\IV(t)}(j)
&\rTTo^{\comp(t)} &\hoM((sA_0^a)_{a\in t(0)};sA_{\max[I]}^1)(j)
\\
\dTTo<{\varSigma(j)} &\sss ((-1)^{c(\tilde\tau)})_\tau \hspace*{-1.8em} 
&\dTTo>{\circledast_\mcG(t)(\varSigma)(j)} &=&
\dTTo~{\vphantom{1_{\big|}}\circledast_\mcG(t)
	(\hoM((\sigma)_{a\in t_h^{-1}b};\sigma^{-1}))_{(h,b)\in\IV(t)}(j)}
&\!\! \sss ((-1)^{c(\tilde\tau)})_\tau \!\!
&\dTTo~{\hoM((\sigma)_{a\in t(0)};\sigma^{-1})(j)}
\\
\rmF_{t(0)}(j) &\rTTo^{\!\!\!\!\Delta(t)}
&\circledast_\mcG(t)(\rmF_{t_h^{-1}b})_{(h,b)\in\IV(t)}(j)
&\rTTo^{\circledast_\mcG(t)(\sfg_h^b)}
&\circledast_\mcG(t)
	(\hoM((A_{h-1}^a)_{a\in t_h^{-1}b};A_h^b))_{(h,b)\in\IV(t)}(j)
&\rTTo^{\comp(t)} &\hoM((A_0^a)_{a\in t(0)};A_{\max[I]}^1)(j)
\end{diagram}
}}
\label{eq-Ft(0)-DeltaG(t)}
\end{equation}
%\caption{}
%\label{dia-4x4-wheels}
\end{figure}

We are going to show that expression~\eqref{eq-comp(t)(ghb)j} is the 
image of \(f_j\) under the left-bottom path
\begin{multline*}
\bigl[ F_{t(0)}(j) \rTTo^{\varSigma(j)} \rmF_{t(0)}(j) 
\rTTo^{\Delta^\mcG(t)}
\circledast_\mcG(t)(\rmF_{t_h^{-1}b})_{(h,b)\in\IV(t)}(j)
\rTTo^{\circledast_\mcG(t)(\sfg_h^b)}
\\
\circledast_\mcG(t)
	(\hoM((A_{h-1}^a)_{a\in t_h^{-1}b};A_h^b))_{(h,b)\in\IV(t)}(j)
\rTTo^{\comp(t)} \hoM((A_0^a)_{a\in t(0)};A_{\max[I]}^1)(j) \bigr]
\end{multline*}
of \eqref{eq-Ft(0)-DeltaG(t)}.
This diagram uses the invertible homomorphism \(\varSigma:F_n\to\rmF_n\) 
of degree 1 defined by \eqref{eq-(Sigma-Sigma)-(AF)-(AF)}.

Let \(j\in\NN^{t(0)}\) and let $\tau$ denote a $t$\n-tree 
\(t\to\co_\sk\) such that $|\tau(0,a)|=j^a$ for all $a\in t(0)$.
Shorten \(\tau_{(h-1,a)\to(h,b)}:\tau(h-1,a)\to\tau(h,b)\) to 
\(\tau_{(h-1,a)}\).
By definition
\begin{gather}
\circledast_\mcG(t)(F_{t_h^{-1}b})_{(h,b)\in\IV(t)}(j)
=\bigoplus^{t-\text{tree }\tau}_{\forall a\in t(0)\,|\tau(0,a)|=j^a}
\bigotimes^{h\in I} \bigotimes^{b\in t(h)} \bigotimes^{p\in\tau(h,b)}
F_{t_h^{-1}b}\Bigl(|\tau_{(h-1,a)}^{-1}(p)|_{a\in t_h^{-1}b}\Bigr),
\label{eq-oxxxF()}
\\
\begin{split}
\circledast_\mcG(t)
	(\hoM((A_{h-1}^a &)_{a\in t_h^{-1}b};A_h^b))_{(h,b)\in\IV(t)}(j)
\\
&=\bigoplus^{t-\text{tree }\tau}_{\forall a\in t(0)\,|\tau(0,a)|=j^a}
\bigotimes^{h\in I} \bigotimes^{b\in t(h)} \bigotimes^{p\in\tau(h,b)}
\und\dg(\tens^{a\in t_h^{-1}b}
T^{|\tau_{(h-1,a)}^{-1}(p)|}A_{h-1}^a,A_h^b).
\end{split}
\notag
\end{gather}
Diagram~\eqref{dia-nA8-Fn-A8} implies commutativity of 
\begin{diagram}
F_{t_h^{-1}b} &\rTTo^{g_h^b} &\hoM((sA_{h-1}^a)_{a\in t_h^{-1}b};sA_h^b)
\\
\dTTo<\varSigma &&\dTTo>{\hoM((\sigma)_{a\in t_h^{-1}b};\sigma^{-1})}
\\
\rmF_{t_h^{-1}b} &\rTTo^{\sfg_h^b}
&\hoM((A_{h-1}^a)_{a\in t_h^{-1}b};A_h^b)
\end{diagram}
Here
	\(\hoM((\sigma)_{a\in t_h^{-1}b};\sigma^{-1})
	=\hoM((\sigma)_{a\in t_h^{-1}b};1)\cdot
	\hoM((1)_{a\in t_h^{-1}b};\sigma^{-1})\)
is the product of right operators.
Hence, the middle square of \eqref{eq-Ft(0)-DeltaG(t)} is commutative.

The second and the third term of both rows of
diagram~\eqref{eq-Ft(0)-DeltaG(t)} are direct sums over $\tau$.
The diagram splits into a direct sum over $\tau$ of 3\n-squares-diagrams 
in which the second square is commutative, while the third square 
commutes up to the sign \((-1)^{c(\tilde\tau)}\).
The mapping \(\Delta^\mcG(t)\) is defined so that the first square of 
the diagram commutes up to the same sign \((-1)^{c(\tilde\tau)}\).
Thus, the exterior of diagram~\eqref{eq-Ft(0)-DeltaG(t)} is commutative.

\begin{proposition}\label{pro-comultiplication-varDelta}
Define for a tree $t:[I]\to\co_\sk$ the degree 0 graded 
$t(0)\wedge1$-\ainf-module homomorphism
	$\varDelta^\mcG(t)(j):F_{t(0)}(j)\to
	\circledast_\mcG(t)(F_{t_h^{-1}b})_{(h,b)\in\IV(t)}(j)$
on generators as
\begin{equation}
\varDelta^\mcG(t)(f_j)
=\sum^{t-\text{tree }\tau}_{\forall a\in t(0)\,|\tau(0,a)|=j^a}
\tens^{h\in I} \tens^{b\in t(h)} \tens^{p\in\tau(h,b)}
f_{|\tau_{(h-1,a)}^{-1}(p)|_{a\in t_h^{-1}b}}.
\label{eq-varDeltaG(t)(fj)}
\end{equation}
In particular, for the only tree $t$ with empty $I$ the \ainf-bimodule 
homomorphism $\varDelta^\mcG(t)(j):F_1(j)\to\kk(j)=\delta_{j1}\kk$ is 
determined by \(\varDelta^\mcG(t)(f_j)=\delta_{j1}\), $j\ge1$.
Then this comultiplication is coassociative and
expression~\eqref{eq-comp(t)(ghb)j} is the image of \(f_j\) under the
left-bottom path of \eqref{eq-Ft(0)-DeltaG(t)}.
Thus, the canonical actions of the unit operad $\kk$ on $F_\bull$ make 
\((\kk,F_\bull,\varDelta^\mcG)\) into a graded polymodule cooperad.
\end{proposition}

A tree $r:[I]\to\co_\sk$ is \emph{surjective} if mappings 
\(r_h:r(h-1)\to r(h)\) are surjective for all $h\in I$.
The summation in expression~\eqref{eq-oxxxF()} and 
formula~\eqref{eq-varDeltaG(t)(fj)} extends precisely over $t$\n-trees 
$\tau$ such that all trees \(\tilde\tau\) are surjective.
In fact, \(F_\emptyset=0\) and the summand of \eqref{eq-oxxxF()} 
corresponding to $\tau$ does not vanish iff
\[ \forall\, h\in I\; \forall\, b\in t(h)\; \forall\, p\in\tau(h,b) \;
\exists\, a\in t_h^{-1}b \quad \tau_{(h-1,a)}^{-1}(p) \ne \emptyset.
\]
Equivalently, for all \(h\in I\)
\[ \forall\, (b,p)\in\bigsqcup_{b\in t(h)}\tau(h,b)=\tilde\tau(h) \;
\exists\, (a,q)\in \!\!\bigsqcup_{a\in t(h-1)}\!\! 
\tau(h-1,a)=\tilde\tau(h-1) \,\; (t_ha,\tau_{(h-1,a)}q) =(b,p).
\]
The last equation means that \(\tilde\tau_h(a,q)=(b,p)\), that is, 
\(\tilde\tau\) is surjective.

The number of surjective trees $\tilde\tau$ is finite.
Hence, the number of tree mappings \(\tilde\tau\to t\) is finite, and 
the number of $t$\n-trees $\tau$ with this surjectivity property is 
finite as well.
Thus the sum is finite.

Deduced comultiplication
	$\Delta^\mcG(t)(j):\rmF_{t(0)}(j)\to
	\circledast_\mcG(t)(\rmF_{t_h^{-1}b})_{(h,b)\in\IV(t)}(j)$
has degree 0:
\begin{equation}
\Delta^\mcG(t)(\sff_j)
=\sum^{t-\text{tree }\tau}_{\forall a\in t(0)\,|\tau(0,a)|=j^a}
(-1)^{c(\tilde\tau)} \tens^{h\in I} \tens^{b\in t(h)} 
\tens^{p\in\tau(h,b)} \sff_{|\tau_{(h-1,a)}^{-1}(p)|_{a\in t_h^{-1}b}}.
\label{eq-Delta-DG(t)(f)}
\end{equation}

\begin{proof}[Proof of \propref{pro-comultiplication-varDelta}]
Using \eqref{eq-map-on-direct-summands} and identifying $j_h^b$ with 
\(|\tau(h,b)|\) we write down
\begin{multline*}
\bigl[ \wh{g_h^b}\cdot\pr_{j_h^b}: 
\tens^{a\in t_h^{-1}b}T^{j_{h-1}^a}sA_{h-1}^a \to T^{j_h^b}sA_h^b \bigr]
\\
\hskip\multlinegap =\bigl[ 
\tens^{a\in t_h^{-1}b}T^{|\tau(h-1,a)|}sA_{h-1}^a 
\rTTo^{\tens^{a\in t_h^{-1}b}\Delta_0^{(|\tau(h,b)|)}} \hfill
\\
\bigotimes_{a\in t_h^{-1}b} \;
\bigoplus_{\sum_{p\in\tau(h,b)}|\tau_{(h-1,a)}^{-1}(p)|=|\tau(h-1,a)|} 
\; \bigotimes_{p\in\tau(h,b)} T^{|\tau_{(h-1,a)}^{-1}(p)|}sA_{h-1}^a 
\\
\rto\sim
\bigoplus_{\sum_{p\in\tau(h,b)}|\tau_{(h-1,a)}^{-1}(p)|=|\tau(h-1,a)|} 
\; \bigotimes_{a\in t_h^{-1}b} \; \bigotimes_{p\in\tau(h,b)} 
T^{|\tau_{(h-1,a)}^{-1}(p)|}sA_{h-1}^a 
\\
\rTTo^{\oplus\overline{\varkappa}_{t_h^{-1}b,\tau(h,b)}}
\bigoplus_{\sum_{p\in\tau(h,b)}|\tau_{(h-1,a)}^{-1}(p)|=|\tau(h-1,a)|} 
\; \bigotimes_{p\in\tau(h,b)} \; \bigotimes_{a\in t_h^{-1}b}
T^{|\tau_{(h-1,a)}^{-1}(p)|}sA_{h-1}^a 
\\
\rTTo^{\sum\tens^{p\in\tau(h,b)}
	g^b_{h,|\tau_{(h-1,a)}^{-1}(p)|_{a\in t_h^{-1}b}}}
\tens^{p\in\tau(h,b)} sA_h^b = T^{|\tau(h,b)|}sA_h^b \bigr].
\end{multline*}

In the particular case of \(j_h^b=|\tau(h,b)|=1\) we get $g_h^b$.
This is realized for the root \((h,b)=(\max[I],1)=\troot_t\) of $t$, 
since \(\tau(\troot_t)=\mb1\).

Tensor product of these expressions gives a factor of
\eqref{eq-comp(t)(ghb)j} for each $h\in I$:
\begin{multline*}
\tens^{b\in t(h)}\wh{g_h^b} =\bigl[ 
\tens^{a\in t(h-1)}T^{|\tau(h-1,a)|}sA_{h-1}^a \rTTo^{\lambda^{t_h}}
\tens^{b\in t(h)}\tens^{a\in t_h^{-1}b}T^{|\tau(h-1,a)|}sA_{h-1}^a 
\\
\rTTo^{\tens^{b\in t(h)}\tens^{a\in t_h^{-1}b}\Delta_0^{(|\tau(h,b)|)}} 
\bigotimes_{b\in t(h)} \; \bigotimes_{a\in t_h^{-1}b} \;
\bigoplus_{\sum_{p\in\tau(h,b)}|\tau_{(h-1,a)}^{-1}(p)|=|\tau(h-1,a)|} 
\; \bigotimes_{p\in\tau(h,b)} T^{|\tau_{(h-1,a)}^{-1}(p)|}sA_{h-1}^a 
\\
\rto\sim \bigotimes_{b\in t(h)} \; 
\bigoplus_{\sum_{p\in\tau(h,b)}|\tau_{(h-1,a)}^{-1}(p)|=|\tau(h-1,a)|} 
\; \bigotimes_{a\in t_h^{-1}b} \; \bigotimes_{p\in\tau(h,b)} 
T^{|\tau_{(h-1,a)}^{-1}(p)|}sA_{h-1}^a 
\\
\rTTo^{\tens^{b\in t(h)}\oplus
	\overline{\varkappa}_{t_h^{-1}b,\tau(h,b)}}
\bigotimes_{b\in t(h)} \; 
\bigoplus_{\sum_{p\in\tau(h,b)}|\tau_{(h-1,a)}^{-1}(p)|=|\tau(h-1,a)|} 
\; \bigotimes_{p\in\tau(h,b)} \; \bigotimes_{a\in t_h^{-1}b}
T^{|\tau_{(h-1,a)}^{-1}(p)|}sA_{h-1}^a 
\\
\rTTo^{\tens^{b\in t(h)}\sum
 \tens^{p\in\tau(h,b)}g^b_{h,|\tau_{(h-1,a)}^{-1}(p)|_{a\in t_h^{-1}b}}}
\tens^{b\in t(h)}\tens^{p\in\tau(h,b)} sA_h^b 
=\tens^{b\in t(h)} T^{|\tau(h,b)|}sA_h^b \bigr].
\end{multline*}

Plugging these expressions into \eqref{eq-comp(t)(ghb)j} we obtain an 
explicit form of the latter.
On the other hand, we compute the image of \(\sff_j\) under the bottom 
row of \eqref{eq-Ft(0)-DeltaG(t)} via the top exterior path of this 
diagram
\begin{multline*}
\sff_j =(\tens^{a\in t(0)}\sigma^{\tens j^a}) \cdot f_j\cdot \sigma^{-1} 
\\
\mapsto (\tens^{a\in t(0)}\sigma^{\tens j^a}) \cdot 
\sum^{t-\text{tree }\tau}_{\forall a\in t(0)\,|\tau(0,a)|=j^a}
\tens^{h\in I} \tens^{b\in t(h)} \tens^{p\in\tau(h,b)}
f_{|\tau_{(h-1,a)}^{-1}(p)|_{a\in t_h^{-1}b}} \cdot \sigma^{-1}
\\
\mapsto (\tens^{a\in t(0)}\sigma^{\tens j^a}) \cdot 
\sum^{t-\text{tree }\tau}_{\forall a\in t(0)\,|\tau(0,a)|=j^a}
\tens^{h\in I} \tens^{b\in t(h)} \tens^{p\in\tau(h,b)}
g^b_{h,|\tau_{(h-1,a)}^{-1}(p)|_{a\in t_h^{-1}b}} \cdot \sigma^{-1}
\\
\mapsto (\tens^{a\in t(0)}\sigma^{\tens j^a}) \cdot \comp_{\dg}(t)\Bigl[
\sum^{t-\text{tree }\tau}_{\forall a\in t(0)\,|\tau(0,a)|=j^a}
\tens^{h\in I} \tens^{b\in t(h)} \tens^{p\in\tau(h,b)}
g^b_{h,|\tau_{(h-1,a)}^{-1}(p)|_{a\in t_h^{-1}b}} 
\Bigr]\cdot\sigma^{-1}.
\end{multline*}
The last expression coincides with \(\comp(t)(g_h^b)_j\).
In fact, multicategory composition in $\und\dg$ restricted to degree 0 
cycles $g^b_h$ coincides with the multicategory composition
$\comp_{\dg}(t)$ in $\dg$.
The latter in our case is the map
\begin{multline*}
\tens^{h\in I} \tens^{b\in t(h)} \tens^{p\in\tau(h,b)}
\dg(\tens^{a\in t_h^{-1}b}T^{|\tau_{(h-1,a)}^{-1}(p)|}sA_{h-1}^a,sA_h^b)
\\
\rTTo^{\tens^{h\in I}\tens^{b\in t(h)}\tens^{\tau(h,b)}}
\tens^{h\in I} \tens^{b\in t(h)} 
\dg(\tens^{p\in\tau(h,b)}\tens^{a\in t_h^{-1}b}
	T^{|\tau_{(h-1,a)}^{-1}(p)|}sA_{h-1}^a,\tens^{p\in\tau(h,b)}sA_h^b)
\\
\rTTo^{\tens^{h\in I}\tens^{b\in t(h)}\dg((\tens^{a\in t_h^{-1}b}
	\lambda^{\tau_{(h-1,a)}})\cdot
	\overline{\varkappa}_{t_h^{-1}b,\tau(h,b)},1)}
\tens^{h\in I} \tens^{b\in t(h)}
\dg(\tens^{a\in t_h^{-1}b}T^{|\tau_{(h-1,a)}|}sA_{h-1}^a,
	T^{|\tau(h,b)|}sA_h^b)
\\
\rTTo^{\tens^{h\in I}\tens^{t(h)}}
\tens^{h\in I}
\dg(\tens^{b\in t(h)} \tens^{a\in t_h^{-1}b}
	T^{|\tau_{(h-1,a)}|}sA_{h-1}^a,
	\tens^{b\in t(h)} T^{|\tau(h,b)|}sA_h^b)
\\
\rTTo^{\tens^{h\in I}\dg(\lambda^{t_h},1)}
\tens^{h\in I}
\dg(\tens^{a\in t(h-1)} T^{|\tau_{(h-1,a)}|}sA_{h-1}^a,
	\tens^{b\in t(h)} T^{|\tau(h,b)|}sA_h^b)
\\
\rTTo^{\mu_{\dg}^I} 
\dg(\tens^{a\in t(0)} T^{|\tau_{(0,a)}|}sA_0^a,sA_{\max[I]}^1).
\end{multline*}
Thus, the considered expression coincides with \eqref{eq-comp(t)(ghb)j}.

Recall the bijection
$\tau\mapsto(\sS{_\psi}\tau,(g,c,q)\mapsto\tR gcq\tau)$ from 
\eqref{eq-(tau-tau)} with the inverse mapping given by
\eqref{eq-T(hb)=bigsqcup}.
It is used in formula~\eqref{eq-circledast(tpsi)circledast} for 
$\lambda^f$.
Summation in \eqref{eq-varDeltaG(t)(fj)} extends precisely over
$t$\n-trees $\tau$ with surjective $\tilde\tau$.
Each summand implements an obvious isomorphism \(\kk\to\kk^{\tens-}\).
We have
\[ \wt{\sS{_\psi}\tau}=\tilde\tau_\psi: [J] \to \co_\sk
\quad\text{ and }\quad
\wt{\tR gcq\tau} =\tilde\tau_{[\psi(g-1),\psi(g)]}^{|(c,g,q)}:
[\psi(g)-\psi(g-1)] \to \co_\sk.
\]
Therefore, $\tilde\tau$ is surjective iff $\sS{_\psi}\tau$ and all
$\tR gcq\tau$ are surjective trees.
This implies that images of $f_j\in F_{t(0)}$ under both paths in the 
following diagram
\begin{diagram}[nobalance,LaTeXeqno]
F_{t(0)} &\rTTo^{\varDelta(t_\psi)}
&\circledast_\mcG(t_\psi)(F_{t_{\psi,g}^{-1}c})_{(g,c)\in\IV(t_\psi)} 
\hspace*{3em}
\\
\dTTo<{\varDelta(t)} &=
&\dTTo>{\circledast_\mcG(t_\psi)
	(\varDelta(t^{|c}_{[\psi(g-1),\psi(g)]}))}
\\
\hspace*{-3em} \circledast_\mcG(t)(F_{t_h^{-1}b})_{(h,b)\in\IV(t)} 
&\rTTo^{\lambda^f}
&\circledast_\mcG(t_\psi)\bigl(
\circledast_\mcG(t^{|c}_{[\psi(g-1),\psi(g)]})
(F_{t_h^{-1}b})_{(h,b)\in\IV(t^{|c}_{[\psi(g-1),\psi(g)]})} 
\bigr)_{(g,c)\in\IV(t_\psi)}
\label{dia-Ft(0)-OG-3-10}
\end{diagram}
are sums indexed by the same set of indices with equal summands.
Thus this diagram is commutative and $\varDelta$ is coassociative.
Notice also that the tensor product $\varDelta^\mcM$ over the unit 
operad $\kk$ coincides with $\varDelta^\mcG$.
\end{proof}

\begin{proposition}
Define comultiplication $\varDelta^\mcM(t)$ for the \ainf-polymodule 
$F_\bull$ and a tree \(t:[p]\to\co_\sk\), $p>0$, via composition with 
$\pi$ from \eqref{eq-circledastG-circledastM}
\[ \varDelta^\mcM(t)(j) =\bigl[ F_{t(0)}(j) \rTTo^{\varDelta^\mcG(t)} 
\circledast_\mcG(t)(F_{t_h^{-1}b})_{(h,b)\in\IV(t)}(j)
\rTTo^\pi \circledast_\mcM(t)(F_{t_h^{-1}b})_{(h,b)\in\IV(t)}(j) \bigr]
\]
(on generators it is still given by \eqref{eq-varDeltaG(t)(fj)}).
For the only tree \(t:[0]\to\co_\sk\) define
$\varDelta^\mcM(t)(j):F_1(j)\to A_\infty(j)$, \(f_j\mapsto\delta_{j1}\), 
$j\ge1$.
Then all \(\varDelta^\mcM(t)(j)\) are chain maps, thus, 
\((A_\infty,F_\bull,\varDelta^\mcM)\) is a $\dg$\n-polymodule cooperad.
\end{proposition}

\begin{proof}
First of all, comultiplication $\varDelta^\mcM$ is coassociative --- it 
satisfies \eqref{dia-Delta-Delta-Delta-lambda}, since $\varDelta^\mcG$ 
satisfies \eqref{dia-Ft(0)-OG-3-10}.
For trees \(t:[1]\to\co_\sk\) the morphism 
$\varDelta^\mcM(t)=\varDelta^\mcG(t)$ satisfies
\eqref{eq-Ft(0)-Delta(t)-1}.

It suffices to prove the commutation of $\varDelta^\mcM(t)$ and 
$\partial$ on generators:
\begin{equation}
f_j.\varDelta^\mcM(t)\partial=f_j.\partial\varDelta^\mcM(t)
\text{ for all trees \(t:[p]\to\co_\sk\) and all indices } 
j\in\NN^{t(0)}-0.
\label{eq-fj-Delta-partial-partial-Delta}
\end{equation}
In fact, any element of $F_n$, $n=|t(0)|$, is a sum of elements of the 
form
	\(\alpha\bigl((\otimes_{i=1}^n\otimes_{p=1}^{k_1^i+\dots+k_m^i}
	\omega_p^i)\tens(\otimes_{r=1}^mf_{k_r})\tens\omega\bigr)\), 
where \(\omega_p^i\in A_\infty(j_p^i)\), \(\omega\in A_\infty(m)\), see 
the first row of the following diagram
\begin{diagram}[nobalance]
\Bigl(\bigotimes_{i=1}^n \bigotimes_{p=1}^{k_1^i+\dots+k_m^i} 
A_\infty(j_p^i)\Bigr)\tens
\Bigl(\bigotimes_{r=1}^m F_{t(0)}(k_r)\Bigr) \tens
A_\infty(m) \hspace*{2.5em}
&&\hphantom{F_{t(0)}\biggl(
\Bigl(\sum_{p=1}^{k_1^i+\dots+k_m^i}j_p^i\Bigr)_{i=1}^n\biggr)}
\\
&\rdTTo^\alpha
\\
\dTTo<{1\tens(\otimes_{r=1}^m\varDelta^\mcM(t))\tens1}
&&F_{t(0)}\biggl(
\Bigl(\sum_{p=1}^{k_1^i+\dots+k_m^i}j_p^i\Bigr)_{i=1}^n\biggr)
\\
\Bigl(\bigotimes_{i=1}^n \bigotimes_{p=1}^{k_1^i+\dots+k_m^i} 
A_\infty(j_p^i)\Bigr)\tens
\Bigl(\bigotimes_{r=1}^m \circledast_\mcM(t)
(F_{t_h^{-1}b})_{(h,b)\in\IV(t)}(k_r)\Bigr)
\tens A_\infty(m) \hspace*{-4em}
&&\dTTo>{\varDelta^\mcM(t)}
\\
&\rdTTo_\alpha
\\
&&\hspace*{-6.5em} \circledast_\mcM(t)(F_{t_h^{-1}b})_{(h,b)\in\IV(t)} 
\biggl(\Bigl(\sum_{p=1}^{k_1^i+\dots+k_m^i}j_p^i\Bigr)_{i=1}^n\biggr)
\end{diagram}
This square commutes since $\varDelta^\mcM(t)$ is a
$t(0)\wedge1$-\ainf-module homomorphism of degree 0.
Use this square as the top and the bottom faces of a cubical diagram 
whose vertical edges are given by the differential $\partial$.
We know that $\alpha$ is a chain map.
Apply all 3-arrow-paths in this cube to the element
	\(x=(\otimes_{i=1}^n\otimes_{p=1}^{k_1^i+\dots+k_m^i}\omega_p^i)
	\tens(\otimes_{r=1}^mf_{k_r})\tens\omega\) 
from the top vertex.
The equation
 \(x.\bigl(1\tens(\otimes_{r=1}^m\varDelta^\mcM(t))\tens1\bigr)\partial
 =x.\partial(1\tens(\otimes_{r=1}^m\varDelta^\mcM(t))\tens1\bigr)\)
holds by assumption, hence, the equation
\(x.\alpha\varDelta^\mcM(t)\partial=x.\alpha\partial\varDelta^\mcM(t)\) 
holds as well.
	
For the only tree of height 0, \(t:[0]\to\co_\sk\),
\(0\mapsto t(0)=\mb1\), and and a positive integer $j$ we have
\(f_j.\varDelta^\mcM(t)\partial=\delta_{j1}.\partial=0\).
On the other hand,
\begin{multline*}
f_j.\partial\varDelta^\mcM(t)
=\sum_{r+n+t=j}^{n>1} (1^{\tens r}\tens b_n\tens1^{\tens t})
f_{r+1+t}.\varDelta^\mcM(t)
-\sum^{l>1}_{i_1+\dots+i_l=j}
(f_{i_1}\tens f_{i_2}\tdt f_{i_l}) b_l.\varDelta^\mcM(t)
\\
=b_j -b_j =0
\end{multline*}
by \eqref{eq-fk-partial-1b1f-fffb}.
For $p=1$ the map $\varDelta^\mcM(t)$ is the identity morphism.
The case of $p>2$ can be reduced to trees of height smaller than $p$.
In fact, for \(\psi:[2]\to[p]\), $0\mapsto0$, $1\mapsto2$, $2\mapsto p$, 
coassociativity equation~\eqref{dia-Delta-Delta-Delta-lambda} represents 
$\varDelta^\mcM(t)$ as composition of $\varDelta^\mcM(t_\psi)$ and 
tensor product of $\varDelta^\mcM(t')$ for trees $t'$ of height 2 or 
$p-2$.
Therefore, it suffices to prove that $\varDelta^\mcM(t)$ is a chain map 
for height $p=2$ of $t$.

Let us prove equation~\eqref{eq-fj-Delta-partial-partial-Delta} for 
\(t=\{\sqcup_{c=1}^n\mb{l}^c\to\bn\to\mb1\}\), \(|t(0)|=l^1+\dots+l^n\).
Elements of \(\NN^{t(0)}\) are written as
\(j=(j^{c,g}\mid c\in\bn,\; g\in\mb{l}^c)\).
Summands of \eqref{eq-varDeltaG(t)(fj)} are indexed by $t$\n-trees 
\(\tau=\tau_\lambda^1\) from \eqref{eq-tau-lambda-tau-rho}.
The tree $\tau$ occurs in expansion of \(\varDelta^\mcM(t)(f_j)\) if 
\(\sum_{p=1}^{u^c}r_p^{c,g}=j^{c,g}\).
Denote \(r_p^c=(r_p^{c,g})_{g\in\mb l^c}\).
Thus,
\[ f_j.\varDelta^\mcM(t) =\sum_{u\in\NN^n} \sum_{\forall c\in\bn\;
\forall g\in\mb l^c\;\sum_{p=1}^{u^c}r_p^{c,g}=j^{c,g}}
\Bigl( \bigotimes_{c\in\bn} \bigotimes_{p\in\mb u^c} 
f_{r_p^c}\Bigr)\tens f_u.
\]
Equally well we could use \(\tau=\tau_\rho^1\) from
\eqref{eq-tau-lambda-tau-rho}, which results in replacing $u$ with $k$, 
$p$ with $q$, and $r$ with $s$ in the above formula.
Below we use both presentations.
Using the lexicographic order on the set \(\sqcup_{c=1}^n\mb u^c\) 
($c<h$ implies that \((c,p)<(h,y)\)), we find
\begin{align}
\begin{split}
&f_j.\varDelta^\mcM(t)\partial 
=\sum_{k\in\NN^n} \sum_{h\in\bn} \sum_{y\in\mb k^h}
\sum_{\forall c\in\bn\;\forall g\in\mb l^c\;
\sum_{q=1}^{k^c}s_q^{c,g}=j^{c,g}} \sum_{z\in\mb l^h} 
\sum_{a+x+m=s_y^{h,z}}
\\
&\hspace*{3em} \Bigl( \bigotimes^{(c,q)<(h,y)\in\sqcup_b\mb k^b} 
f_{s_q^c}\Bigr)\tens 
\lambda^z(\sS{^a}1,b_x,\sS{^m}1;f_{s_y^h-(x-1)e_z})\tens\Bigl(
	\bigotimes^{(c,q)>(h,y)\in\sqcup_b\mb k^b}f_{s_q^c}\Bigr)\tens f_k
\end{split}
\label{eq-fj-Delta-M(t)partial-first}
\\
\begin{split}
&-\sum_{k\in\NN^n} \sum_{h\in\bn} \sum_{y\in\mb k^h}
\sum_{\forall c\in\bn\;\forall g\in\mb l^c\;
\sum_{q=1}^{k^c}s_q^{c,g}=j^{c,g}}
\sum^{i_1,\dots,i_w\in\NN^{l^h}-0}_{i_1+\dots+i_w=s_y^h}
\\
&\hspace*{3em} \Bigl( \bigotimes^{(c,q)<(h,y)\in\sqcup_b\mb k^b} 
f_{s_q^c}\Bigr)\tens \rho((f_{i_s})_{s=1}^w;b_w)
\tens\Bigl( \bigotimes^{(c,q)>(h,y)\in\sqcup_b\mb k^b} 
f_{s_q^c}\Bigr)\tens f_k
\end{split}
\label{eq-fj-Delta-M(t)partial-second}
\\
&+\sum_{u\in\NN^n} \sum_{\forall c\in\bn\;\forall g\in\mb l^c\;
\sum_{p=1}^{u^c}r_p^{c,g}=j^{c,g}} \sum_{h\in\bn}
\sum_{a+w+m=u^h} \Bigl( \bigotimes_{c\in\bn} \bigotimes_{p\in\mb u^c} 
f_{r_p^c}\Bigr)\tens
\lambda^h(\sS{^a}1,b_w,\sS{^m}1;f_{u-(w-1)e_h})
\label{eq-fj-Delta-M(t)partial-third}
\\
&-\sum_{u\in\NN^n} \sum_{\forall c\in\bn\;\forall g\in\mb l^c\;
\sum_{p=1}^{u^c}r_p^{c,g}=j^{c,g}}
\sum^{u_1,\dots,u_w\in\NN^n-0}_{u_1+\dots+u_w=u}
\Bigl( \bigotimes_{c\in\bn} \bigotimes_{p\in\mb u^c} 
f_{r_p^c}\Bigr)\tens \rho((f_{u_v})_{v=1}^w;b_w).
\label{eq-fj-Delta-M(t)partial-fourth}
\end{align}
Sums \eqref{eq-fj-Delta-M(t)partial-second} and
\eqref{eq-fj-Delta-M(t)partial-third} cancel each other.
In fact, consider the tree \(\tau_+^1:t_+^1\to\co_\sk\) from
\eqref{dia-rocket-uk} with $j_q^c=1$ for $c\ne h$, $j_q^h=1$ for
$q\ne y$, and $j_y^h=w$.
This implies relations $u^c=k^c$ for $c\ne h$, $u^h=k^h+w-1$, 
\(s_q^{c,g}=r_q^{c,g}\) for $c\ne h$, \(s_q^{h,g}=r_q^{h,g}\) if $q<y$, 
\(s_y^{h,g}=\sum_{p=y}^{y+w-1}r_p^{h,g}\), and 
\(s_q^{h,g}=r_{q+w-1}^{h,g}\) if $q>y$.
The images of the element
\begin{multline*}
\Bigl( \bigotimes^{(c,q)<(h,y)} f_{s_q^c}\Bigr)\tens \Bigl( 
\bigotimes^{v\in\mb w} f_{i_v}\Bigr)
\tens\Bigl( \bigotimes^{(c,q)>(h,y)} f_{s_q^c}\Bigr) \tens
\Bigl( \bigotimes^{(c,q)<(h,y)} 1\Bigr)\tens b_w\tens \Bigl( 
\bigotimes^{(c,q)>(h,y)} 1\Bigr) \tens f_k
\\
=\Bigl( \bigotimes^{c\in\bn} \bigotimes^{p\in\mb u^c} 
f_{r_p^c}\Bigr)\tens
\Bigl( \bigotimes^{(c,q)<(h,y)\in\sqcup_b\mb k^b} 1\Bigr)\tens b_w\tens 
\Bigl( \bigotimes^{(c,q)>(h,y)\in\sqcup_b\mb k^b} 1\Bigr) \tens f_k,
\end{multline*}
where \(i_v=r_{y,v}^{h,g}=r_p^{h,g}\) for \(p=y-1+v\), $v\in\mb w$, 
under maps \eqref{eq-125-} and \eqref{eq-125--} are, respectively,
\begin{gather*}
\Bigl( \bigotimes^{(c,q)<(h,y)\in\sqcup_b\mb k^b} f_{s_q^c}\Bigr)\tens 
\rho((f_{i_v})_{v=1}^w;b_w)
\tens\Bigl( \bigotimes^{(c,q)>(h,y)\in\sqcup_b\mb k^b} 
f_{s_q^c}\Bigr)\tens f_k,
\\
\Bigl( \bigotimes^{c\in\bn} \bigotimes^{p\in\mb u^c} 
f_{r_p^c}\Bigr)\tens \lambda^h(\sS{^{y-1}}1,b_w,\sS{^{k^h-y}}1;f_k).
\end{gather*}
These are identified with summands of sums
\eqref{eq-fj-Delta-M(t)partial-second} and
\eqref{eq-fj-Delta-M(t)partial-third}.
Namely, we identify $a$ with $y-1$, $m$ with $k^h-y$ and notice that 
\(u-(w-1)e_h=k\).
Therefore, summands of these sums pairwise cancel each other.

We claim that the difference of sums
\eqref{eq-fj-Delta-M(t)partial-first} and
\eqref{eq-fj-Delta-M(t)partial-fourth} equals to
\(f_j.\partial\varDelta^\mcM(t)\).
In fact, by \eqref{eq-f-ell-partial-bs}
\begin{align}
&f_j.\partial\varDelta^\mcM(t) =\sum_{h=1}^n\sum_{z=1}^{l^h} 
\sum_{\gamma+x+m=j^{h,z}}^{x>1}
\lambda^{h,z}(\sS{^\gamma}1,b_x,\sS{^m}1;f_{j-(x-1)e_{h,z}}.
\varDelta^\mcM(t)) \notag
\\
&\hspace*{5em} -\sum^{j_1,\dots,j_w\in\NN^{t(0)}-0}_{j_1+\dots+j_w=j} 
\rho((f_{j_v}.\varDelta^\mcM(t))_{v=1}^w;b_w) \notag
\\
&=\sum_{h=1}^n\sum_{z=1}^{l^h} \sum_{\gamma+x+m=j^{h,z}}^{x>1} 
\sum_{k\in\NN^n} \sum_{\forall c\in\bn\;\forall g\in\mb l^c\;
\sum_{q=1}^{k^c}s_q^{c,g}=j^{c,g}-(x-1)\delta_h^c\delta_z^g} 
\hspace*{-1em}
\lambda^{h,z}\Bigl(\sS{^\gamma}1,b_x,\sS{^m}1; \Bigl( 
\bigotimes_{c\in\bn} \bigotimes_{q\in\mb k^c} f_{s_q^c}\Bigr)\tens 
f_k\Bigr)
\label{eq-fj-partial-Delta-M(t)first}
\\
&-\hspace*{-0.6em}
\sum^{j_1,\dots,j_w\in\NN^{t(0)}-0}_{j_1+\dots+j_w=j} 
\rho\Biggl(\biggl(\sum_{u_v\in\NN^n}
\sum_{\forall c\in\bn\;\forall g\in\mb l^c\;
\sum_{p=u_1^c+\dots+u_{v-1}^c+1}^{u_1^c
	+\dots+u_{v-1}^c+u_v^c}r_p^{c,g}=j_v^{c,g}}
\hspace*{-1em}
\Bigl( \bigotimes^{c\in\bn} 
\bigotimes_{p=1+\sum_{\alpha=1}^{v-1}u_\alpha^c}^{\sum_{\alpha=1}^v
	u_\alpha^c} \hspace*{-0.5em} f_{r_p^c}\Bigr) 
\tens f_{u_v}\biggr)_{v=1}^w \hspace*{-0.2em} ;b_w\Biggr).
\label{eq-fj-partial-Delta-M(t)second}
\end{align}
Sum \eqref{eq-fj-partial-Delta-M(t)second} coincides with
\eqref{eq-fj-Delta-M(t)partial-fourth}.
Sums \eqref{eq-fj-Delta-M(t)partial-first} and
\eqref{eq-fj-partial-Delta-M(t)first} are equal after the following 
identification.
The sum in \eqref{eq-fj-partial-Delta-M(t)first} over $\gamma\in\NN$ 
such that \(1\le\gamma+1\le j^{h,z}-x+1=\sum_{q=1}^{k^h}s_q^{h,z}\) is 
equivalent to summing over pairs \((y,\yi)\), \(y\in\mb k^h\), 
\(\yi\in\mb s_y^{h,z}\), namely, 
\(\gamma+1=\yi+\sum_{q=1}^{y-1}s_q^{h,z}\).
The expression $\yi-1$ is denoted $a$ in
\eqref{eq-fj-Delta-M(t)partial-first}.
Variables $s_q^{c,g}$ coincide in both expressions, except $s_y^{h,z}$ 
whose values differ by $x-1$.
This finishes the proof.
\end{proof}

\begin{definition}
A polymodule cooperad morphism \((h,p):(\ca,F)\to(\cb,G)\) of degree 1 
is a family of $n\wedge1$-operad module homomorphisms of degree 1
\[ (\sS{^n}h;p_n;h): (\sS{^n}\ca;F_n;\ca) \to (\sS{^n}\cb;G_n;\cb)
\]
such that for all non-decreasing maps \(\phi:I\to J\), the induced 
\(\psi=[\phi]:[J]\to[I]\) as in \eqref{eq-x<fy-fx<y}, and for all trees 
\(t:[I]\to\co_\sk\) each wall (vertical face) of cube at
\figref{fig-dia-Ft(0)-Delta(t)-cube} commutes up to prescribed sign, 
varying with the summand of $\circledast_\mcM$.
Here the sign $(-1)^{s(f)}$ comes from the permutation of factors $p(k)$ 
of degree $1-|k|$ according to
rule~\eqref{eq-circledast(tpsi)circledast}.
The floor and the ceiling of cube at
\figref{fig-dia-Ft(0)-Delta(t)-cube} commute.
\end{definition}

\begin{figure}
\begin{equation*}
%\boldmath
\resizebox{!}{.48\texthigh}{\rotatebox[origin=c]{90}{%
\begin{diagram}[nobalance,inline,h=1.7em]
F_{t(0)} &&\rTTo^{\varDelta(t_\psi)} 
&&\circledast_\mcM(t_\psi)(F_{t_{\psi,g}^{-1}c})_{(g,c)\in\IV(t_\psi)} 
\hspace*{-8em} &&\hspace*{24em}
\\
&\rdTTo(2,4)>{\varDelta(t)} &&&\dLine 
&\rdTTo(2,4)>{\circledast_\mcM(t_\psi)
	(\varDelta(t^{|c}_{[\psi(g-1),\psi(g)]}))}
\\
&&&\makebox[0mm][c]{$\sss (-1)^{c(\wt{\sS{_\psi}\tau})}$}
\\
&&&&&\makebox[0mm][l]{$\sss (-1)^{\sum_{g,c,q}c(\wt{\tR gcq\tau})}$}
\\
\dTTo<p &&\hspace*{-3em}
\circledast_\mcM(t)(F_{t_h^{-1}b})_{(h,b)\in\IV(t)} 
&\rTTo^{\lambda^f}_\sim &\HonV
&&\circledast_\mcM(t_\psi)\bigl(
\circledast_\mcM(t^{|c}_{[\psi(g-1),\psi(g)]})
(F_{t_h^{-1}b})_{(h,b)\in\IV(t^{|c}_{[\psi(g-1),\psi(g)]})} 
\bigr)_{(g,c)\in\IV(t_\psi)}
\\
\\
&&\dTTo<{\circledast_\mcM(t_\psi)(p)}
&&\dTTo>{\circledast_\mcM(t_\psi)(p)}
\\
\\
G_{t(0)} &\rLine &\VonH &\rTTo^{\varDelta(t_\psi)}
&\circledast_\mcM(t_\psi)(G_{t_{\psi,g}^{-1}c})_{(g,c)\in\IV(t_\psi)} 
\hspace*{-8em}
&&\dTTo>{\circledast_\mcM(t_\psi)
	(\circledast_\mcM(t^{|c}_{[\psi(g-1),\psi(g)]})(p))}
\\
&\rdTTo(2,4)<{\varDelta(t)} \makebox[0mm][l]{$\sss(-1)^{c(\wt\tau)}$}
&&&&\rdTTo(2,4)<{\circledast_\mcM(t_\psi)
	(\varDelta(t^{|c}_{[\psi(g-1),\psi(g)]}))}
\\
&&&\makebox[0mm][c]{$\sss (-1)^{s(f)}$}
\\
\\
&&\hspace*{-3em} \circledast_\mcM(t)(G_{t_h^{-1}b})_{(h,b)\in\IV(t)} 
&&\rTTo^{\lambda^f}_\sim
&&\circledast_\mcM(t_\psi)\bigl(
\circledast_\mcM(t^{|c}_{[\psi(g-1),\psi(g)]})
(G_{t_h^{-1}b})_{(h,b)\in\IV(t^{|c}_{[\psi(g-1),\psi(g)]})} 
\bigr)_{(g,c)\in\IV(t_\psi)}
\end{diagram}
}}
%\label{eq-Ft(0)-Delta(t)-cube}
\end{equation*}
\caption{Cube equations for a polymodule cooperad morphism of degree 1}
\label{fig-dia-Ft(0)-Delta(t)-cube}
\end{figure}

\begin{remark}\label{rem-signs-on-walls}
For each summand of the target of cube at
\figref{fig-dia-Ft(0)-Delta(t)-cube} (the right lower vertex) the 
product of all signs on walls of this diagram is $+1$.
Indeed, the same signs occur in diagram
\begin{diagram}[nobalance,h=1.5em,w=1.7em,LaTeXeqno]
H'(\dd{t}) &&\rTTo^{\comp(t_\psi)\hspace*{-2em}} 
&&\circledast_\mcG(t_\psi)H'(t_\psi) \hspace*{-4em} &&\hspace*{21em}
\\
&\rdTTo(2,4)>{\comp(t)} &&&\dLine
&\rdTTo(2,4)>{\circledast_\mcG(t_\psi)
	(\comp(t^{|c}_{[\psi(g-1),\psi(g)]}))}
\\
&&&\makebox[0mm][c]{$\sss (-1)^{c(\wt{\sS{_\psi}\tau})}$}
\\
&&&&&\makebox[0mm][l]{$\sss (-1)^{\sum_{g,c,q}c(\wt{\tR gcq\tau})}$}
\\
\dTTo<{h(\dd{t})} &&\hspace*{-3em} \circledast_\mcG(t)H'(t) 
&\rTTo^{\lambda^f}_\sim &\HonV
&&\circledast_\mcG(t_\psi)\bigl(
\circledast_\mcG(t^{|c}_{[\psi(g-1),\psi(g)]})
H'(t^{|c}_{[\psi(g-1),\psi(g)]}) \bigr)_{(g,c)\in\IV(t_\psi)}
\\
\\
&&\dTTo<{\circledast_\mcG(t_\psi)h(t)} 
&&\dTTo>{\circledast_\mcG(t_\psi)h(t_\psi)}
\\
\\
H(\dd{t}) &\rLine &\VonH 
&\rTTo^{\hspace*{-1em}\comp(t_\psi)\hspace*{-1em}} 
&\circledast_\mcG(t_\psi)H(t_\psi) \hspace*{-4em}
&&\dTTo~{\circledast_\mcG(t_\psi)(\circledast_\mcG
	(t^{|c}_{[\psi(g-1),\psi(g)]})h(t^{|c}_{[\psi(g-1),\psi(g)]}))}
\\
&\rdTTo(2,4)<{\comp(t)} \makebox[0mm][l]{$\sss(-1)^{c(\wt\tau)}$}
&&&&\rdTTo(2,4)<{\circledast_\mcG(t_\psi)
	(\comp(t^{|c}_{[\psi(g-1),\psi(g)]}))}
\\
&&&\makebox[0mm][c]{$\sss (-1)^{s(f)}$}
\\
\\
&&\hspace*{-3em} \circledast_\mcG(t)H(t) &\rTTo^{\lambda^f}_\sim
&&&\circledast_\mcG(t_\psi)\bigl(
	\circledast_\mcG(t^{|c}_{[\psi(g-1),\psi(g)]})
H(t^{|c}_{[\psi(g-1),\psi(g)]}) \bigr)_{(g,c)\in\IV(t_\psi)}
\label{dia-hom-shift-cube}
\end{diagram}
relating composition $\comp$ for functors $\hoM$ and the shifts 
$\sigma$.
The following abbreviations for families of operad polymodules and their 
degree 1 homomorphisms are used in diagram~\eqref{dia-hom-shift-cube}:
\begin{align*}
H'(t) &= 
\bigl(\hoM((sA_{h-1}^a)_{a\in t_h^{-1}b};sA_h^b)\bigr)_{(h,b)\in\IV(t)},
\\
h(t) &= \bigl(
\hoM((\sigma)_{a\in t_h^{-1}b};\sigma^{-1})\bigr)_{(h,b)\in\IV(t)},
\\
H(t) &= 
\bigl(\hoM((A_{h-1}^a)_{a\in t_h^{-1}b};A_h^b)\bigr)_{(h,b)\in\IV(t)}.
\end{align*}
For \(t:[I]\to\co_\sk\) the tree $\dd{t}$ means the corolla 
\(t(0)\to\mb1=t_{\max[I]}\) labelled with $(A_0^a)_{a\in t(0)}$ and 
$A_{\max[I]}^1$.
Thus \(H(\dd{t})=\hoM((A_0^a)_{a\in t(0)};A_{\max[I]}^1)\) etc.
The commuting floor of diagram~\eqref{dia-hom-shift-cube} is precisely 
equation~\eqref{eq-hom-comp(t)}.
The commuting ceiling is similar, with $sA_h^b$ instead of $A_h^b$.
Commutativity of the floor and the ceiling implies the claim on signs.
\end{remark}

\begin{lemma}\label{lem-degree-1-isomorphism-cooperads}
Let \((\sS{^n}h;p_n;h):(\sS{^n}\ca;F_n;\ca)\to(\sS{^n}\cb;G_n;\cb)\), 
$n\ge0$, be a family of invertible $n\wedge1$-operad module 
homomorphisms of degree 1.
If the family \((\ca,F_\bull)\) or \((\cb,G_\bull)\) has a structure of 
a polymodule cooperad, then the other family has a unique cooperad 
structure such that \((h,p_\bull)\) is a degree 1 isomorphism of 
cooperads.
\end{lemma}

\begin{proof}
Follows from diagram at \figref{fig-dia-Ft(0)-Delta(t)-cube} and the 
equation between signs on walls proven in \remref{rem-signs-on-walls}.
\end{proof}

\begin{corollary}
The family of $n\wedge1$-operad modules \((\mainf,\rmF_n)\), $n\ge0$, 
equipped with comultiplication~\eqref{eq-Delta-DG(t)(f)} is a polymodule 
cooperad.
\end{corollary}

\subsection{Comultiplication for homotopy unital case.}
Let multiquiver \(\alinf^{hu}=\mcH\) be convolution of
\(\rmF^{hu}:\mcF\to\mcM\) and $\HOM:\mcB\to\mcM$ coming from $\uCom$.
Objects of \(\alinf^{hu}\) are homotopy unital \ainfm-algebras and 
morphisms are homotopy unital \ainfm-morphisms.

There is a multiquiver map \(\text-^+:\alinf^{hu}\to\alinf\), 
	\((A,\sfi,m_1,m_{n_1;n_2;\dots;n_k}\mid k+\sum_{q=1}^kn_q\ge3)
	\to(A^+,m^+_n\mid n\ge1)\),
where \(A^+=A\oplus\kk\one^\su\oplus\kk \sfj\) is strictly unital with 
the strict unit $\one^\su$, \(m^+_n\big|_{A^{\tens n}}=m_n\),
\(\sfj m^+_1=\one^\su-\sfi\) and
\[ (1^{\tens n_1}\tens \sfj\tens1^{\tens n_2}\tens
\sfj\tdt1^{\tens n_{k-1}}\tens \sfj\tens1^{\tens n_k})m^+_{n+k-1} 
=m_{n_1;n_2;\dots;n_k}: A^{\tens n+k-1} \to A
\]
for $k\ge1$, $n_q\ge0$, \(n=\sum_{q=1}^kn_q\), \(n+k\ge3\).
On morphisms with $n$ arguments we have
\[ \sff
=(\sfv_k, \sff_{(\ell^k_1;\ell^k_2;\dots;\ell^k_{t^k})_{k\in\bn}})
\mapsto \sff^+ =(\sff^+_j \mid j\in\NN^n-0),
\]
where \(\sfj\sff^+_{e_k}=\sfv_k+\sfj\rho_\emptyset\) and for
$|\hat\ell|\ge2$
\begin{multline*}
\bigl[\tens^{k\in\bn}T^{\ell^k}A_k
\rTTo^{\tens^{k\in\bn}(1^{\tens\ell^k_1}\tens\sfj\tens1^{\tens\ell^k_2} 
	\tens\sfj\tdt1^{\tens\ell^k_{t^k-1}}\tens\sfj
	\tens1^{\tens\ell^k_{t^k}})}
\tens^{k\in\bn}T^{\hat\ell^k}A^+_k \rTTo^{\sff^+_{\hat\ell}} B\bigr]
\\
=\lambda_{\hat\ell}\bigl((\sS{^{\ell^k_1}}1,\sfj,\sS{^{\ell^k_2}}1,\sfj,
\dots,\sS{^{\ell^k_{t^k-1}}}1,\sfj,\sS{^{\ell^k_{t^k}}}1)_{k\in\bn};
\sff^+_{\hat\ell}\bigr)
=\sff_{(\ell^k_1;\ell^k_2;\dots;\ell^k_{t^k})_{k\in\bn}}.
\end{multline*}
This multiquiver map is injective on morphisms and the conditions of 
\defref{def-homotopy-unital-structure-A8-morphism} describe its image.
The image is closed under composition in $\alinf$, hence, it is a 
submulticategory.
In this way \(\alinf^{hu}\) becomes a multicategory and 
\(\text-^+:\alinf^{hu}\to\alinf\) becomes a multifunctor.
Composing it with the multifunctor \(Ts:\alinf\to\dgac\) we get again a 
full and faithful embedding \(Ts(\text-)^+:\alinf^{hu}\to\dgac\).
Its image is described by conditions parallel to that of 
\defref{def-homotopy-unital-structure-A8-morphism}:
\begin{enumerate}
\renewcommand{\labelenumi}{(\arabic{enumi})}
\item $\sff^+$ is a strictly unital;

\item
\(\wh{\sff^+}(1\tdt1\tens(A_k+\sfj^{A_k})\tens1\tdt1)\subset B+\sfj^B\);

\item $\wh{\sff^+}(\tens^{k\in\bn}TA_k)\subset TB$;

\item
	\(\wh{\sff^+}(\tens^{k\in\bn}T^{\ell^k}(A_k\oplus\kk \sfj^{A_k}))
	\subset B\oplus T^{>1}(B\oplus\kk \sfj^B)\)
for each \(\ell\in\NN^n\), $|\ell|>1$.
\end{enumerate}
One checks directly that the set of such coalgebra morphisms is closed 
under composition.

\subsubsection{Actions of operads in the tensor product of operad 
	modules.}
Assume given a tree \(t:[l]\to\co_\sk\) and a family 
\((\cp_h^b)_{(h,b)\in\IV(t)}\) of \(t_h^{-1}b\wedge1\)-operad modules.
Their tensor product \(\cp=\circledast(t)(\cp_h^b)_{(h,b)\in\IV(t)}\) is 
a \(t(0)\wedge1\)-operad module.
Operads acting on the left are \((\ca_a)_{a\in t(0)}\), where 
\((\ca_a)_{a\in t_1^{-1}b}\) act on $\cp_1^b$ for $b\in t(1)$.
The summand indexed by $t$\n-tree $\tau$, \(\tau(0,a)=k^a\) is mapped by 
$\lambda$ from \eqref{eq-rho-rhokr} to the $t$\n-tree $\bar\tau$ 
obtained as follows.
There is a tree \(t':[-1,l]\cap\ZZ\to\co_\sk\), \(t'\big|_{[l]}=t\), 
\(t'_0=\id_{t(0)}\), and a $t'$\n-tree $\tau'$ such that 
\(\tau'\big|_{[l]}=\tau\), and the maps \(\tau'_{(-1,a)}\) correspond to 
partition \(\sum_{p=1}^{k^a}j_p^a\) into $k^a$ summands for
\(a\in t(0)\).
The $t$\n-tree $\bar\tau$ is obtained from $\tau'$ by dropping the 
intermediate level $0$, so that level $-1$ becomes level $0$.
The mapping in question is the tensor product over \(b\in t(1)\) of 
actions $\lambda$ of \((\ca_a)_{a\in t_1^{-1}b}\) on $\cp_1^b$.

The operad acting on the right of $\cp$ is the operad $\cb$ acting on 
the right of $\cp_l^1$.
The action map $\rho$ given by \eqref{eq-rho-rhokr} sends the summand 
\(\bigl(\otimes_{r=1}^m\cp(\tau_r)\bigr)\tens\cB(m)\to\cP(\tau')\), 
where $\tau'$ and $\tau_r$, \(1\le r\le m\), are $t$\n-trees, 
\(\cp=\oplus_{t-\text{tree }\tau}\cp(\tau)\), and $\tau'$ is constructed 
as \(\tau_1\sqcup\tau_2\sqcup\dots\sqcup\tau_m\) with identified roots.
The action map $\rho$ for $\cp$ is the identity map tensored with the 
action $\rho$ for $\cp_l^1$.

In particular, the map \(\rho_\emptyset^\cp:\cb(0)\to\cp(0)\) sends 
$\cb(0)$ to the summand \(\cp(\tau_0)\simeq\cp_l^1(0)\) via 
\(\rho_\emptyset^{\cp_l^1}:\cb(0)\to\cp_l^1(0)\), where 
\(\tau_0(h,b)=\emptyset\) for \((h,b)\in\IV(t)-\{\troot\}\), while 
\(\tau_0(\troot)=\tau_0(l,1)=\mb1\).

Comultiplication~\eqref{eq-Delta-DG(t)(f)} extends in a unique way to 
\((\mainf^\su,\rmF_n^\su)\), which differs from \((\mainf,\rmF_n)\) by a 
direct summand \((\kk\one^\su,\kk\one^\su\rho_\emptyset)\), see 
\eqref{dia-(AF)-(AsuFsu)-(k1k1)}.
In fact, for a tree \(t:[I]\to\co_\sk\) the equation
\begin{equation}
\rho_\emptyset =\bigl[\mainf^\su(0) \rTTo^{\rho_\emptyset}_\sim 
\rmF_{t(0)}^\su(0) \rTTo^{\Delta(t)}
\circledast_\mcM(t)(\rmF_{t_h^{-1}b}^\su)_{(h,b)\in\IV(t)}(0) \bigr]
\label{eq-rho-empty-set-rho-empty-set}
\end{equation}
is one of those saying that \(\Delta(t)\) agree with $\rho$ (see 
\eqref{dia-PPPB-rho-P} with $l=0$).
So we set 
	\(\Delta(t)\bigl(\rho_\emptyset(\one^\su)\bigr)
	=\rho_\emptyset(\one^\su)\).
For non-empty $I$ the image of \eqref{eq-rho-empty-set-rho-empty-set} is 
contained in the image of the summand
\(\rmF_{t(|I|-1)}^\su(0)\) of
\(\circledast_\mcG(t)(\rmF_{t_h^{-1}b}^\su)_{(h,b)\in\IV(t)}(0)\) 
indexed by the $t$\n-tree $\tau$ with \(\tau(h,b)=\emptyset\) for all 
$(h,b)$ such that \(0\le h<|I|\), while \(\tau(|I|,1)=\mb1\).
For the tree \(t:[0]\to\co_\sk\) \eqref{eq-rho-empty-set-rho-empty-set} 
is the right action in the regular bimodule $\mainf^\su$:
\[ \id =\rho_\emptyset =\bigl[ \mainf^\su(0) \rTTo^{\rho_\emptyset}_\sim 
\rmF_1^\su(0) \rTTo^{\Delta(t)} \mainf^\su(0) \bigr].
\]
We can be more precise in this case: 
	\(\Delta(t)\bigl(\rho_\emptyset(\one^\su)\bigr)
	=\rho_\emptyset(\one^\su)=\one^\su\).

So extended comultiplication obviously agrees with the left action 
$\lambda$ (see \eqref{dia-AAAAP-lambda-P} with $k=0$).
It agrees also with the right action $\rho$, see \eqref{dia-PPPB-rho-P} 
for $l>0$ with \(J=\{q\in\mb l\mid k_q=0\}\).
We may take elements \(\one^\su\rho_\emptyset\) in each place 
\(\cp(0)=\rmF_n^\su(0)\) for $q\in J$.
Then \(\one^\su\rho_\emptyset\) will appear also in 
\(\cq(0)=\circledast_\mcM(t)(\rmF_{t_h^{-1}b}^\su)_{(h,b)\in\IV(t)}(0)\) 
for the same $q$.
Using associativity of $\rho$ we can absorb those \(\one^\su\) into an 
element of \(\mainf^\su\) and get rid of \(\one^\su\)'s completely.
The equation is reduced to the case of \((\mainf,\rmF_n)\), which is 
already verified.
Coassociativity of extended comultiplication is obvious.

Let us extend comultiplication further to
\[ (\mainf^\su,\rmF_n^\su)\langle \sfi,\sfj\rangle \simeq
\bigl(\mainf^\su\langle \sfi,\sfj\rangle,
\bigcirc_{k=0}^n\mainf^\su\langle \sfi,
\sfj\rangle\odot^k_{\mainf^\su}\rmF_n^\su \bigr)
\]
using \propref{pro-A<Mbeta>}.
Let \(n=|t(0)|\).
The comultiplication is the lower diagonal in
\[ \hspace*{-3em}
\begin{diagram}[inline,nobalance]
\circledast_\mcM(t)(\mainf^\su,\rmF^\su_{t_h^{-1}b})_{(h,b)\in\IV(t)} 
&\rMono
&\circledast_\mcM(t)(\mainf^\su,\rmF^\su_{t_h^{-1}b})_{(h,b)\in\IV(t)}
\langle\sfi,\sfj\rangle
\\
\uTTo<{\Delta^\mcM(t)} &&\dTTo>\wr
\\
(\mainf^\su,\rmF_n^\su) &&\bigl(\mainf^\su\langle \sfi,\sfj\rangle,
\bigcirc_{k=0}^n\mainf^\su\langle \sfi,\sfj\rangle\odot^k_{\mainf^\su}
\circledast_\mcM(t)
(\mainf^\su,\rmF^\su_{t_h^{-1}b})_{(h,b)\in\IV(t)}\bigr)
\\
\dMono &\ruTTo[hug]^{\Delta^\mcM(t)\langle \sfi,\sfj\rangle} &\dMono
\\
(\mainf^\su,\rmF_n^\su)\langle \sfi,\sfj\rangle 
&&\circledast_\mcM(t)(\mainf^\su\langle \sfi,\sfj\rangle,
\bigcirc_{k\in[t_h^{-1}b]}
\mainf^\su\langle \sfi,\sfj\rangle\odot^k_{\mainf^\su}
\rmF^\su_{t_h^{-1}b})_{(h,b)\in\IV(t)}
\\
&\rdTTo_{\Delta^\mcM(t)} &\dTTo>\wr
\\
&&\circledast_\mcM(t)(\mainf^\su,\rmF^\su_{t_h^{-1}b})
\langle \sfi,\sfj\rangle_{(h,b)\in\IV(t)}
\end{diagram}
\hspace*{3em}
\]

Proof of coassociativity is contained in the following diagram.
The operad module \((\mainf^\su,\rmF^\su_n)\) is short-handed to 
$\rmF^\su_n$.
Similarly $\rmF^\su_n\langle \sfi,\sfj\rangle$ stands for
\((\mainf^\su,\rmF^\su_n)\langle \sfi,\sfj\rangle\).
Being diagram~\eqref{dia-Delta-Delta-Delta-lambda} the top square 
commutes.
The middle square parallel to the top face is obtained by adding freely 
operations $\sfi$ and $\sfj$.
Hence it also commutes.
The vertical faces commute as well, therefore, the bottom quadrangle is 
commutative.
\[ \hspace*{0.1em}
\begin{diagram}[nobalance,inline]
&&\hspace*{-0.2em} \circledast_\mcM^{(h,b)}(t)\rmF^\su_{t_h^{-1}b}
&&\rTTo^{\lambda^f} 
&&\circledast_\mcM^{(g,c)}(t_\psi)\bigl( \circledast_\mcM^{(h,b)}
(t^{|c}_{[\psi(g-1),\psi(g)]})\rmF^\su_{t_h^{-1}b} \bigr) 
\hspace*{4em}
\\
&\ruTTo^{\Delta(t)} &\dGyrlyga &&&\ruTTo>{\circledast_\mcM^{(g,c)}
(t_\psi)(\Delta(t^{|c}_{[\psi(g-1),\psi(g)]}))}
\\
\rmF^\su_{t(0)} &&\HonV &\rTTo^{\hspace*{-2.2em}\Delta(t_\psi)} 
&\circledast_\mcM^{(g,c)}(t_\psi)\rmF_{t_{\psi,g}^{-1}c} &&\dMono
\\
&&\dTTo 
\\
\dMono &&\hspace*{-2.4em} 
[\circledast_\mcM^{(h,b)}(t)\rmF^\su_{t_h^{-1}b}]\langle\sfi,\sfj\rangle 
&\rLine^{\lambda^f\langle \sfi,\sfj\rangle\hspace*{-2.5em}} &\dMono 
&\rTTo 
&\bigl[\circledast_\mcM^{(g,c)}(t_\psi)\bigl( \circledast_\mcM^{(h,b)}
(t^{|c}_{[\psi(g-1),\psi(g)]})\rmF^\su_{t_h^{-1}b} \bigr)\bigr]
\langle \sfi,\sfj\rangle \hspace*{7em}
\\
&\ruTTo[hug]^{\Delta(t)\langle \sfi,\sfj\rangle} &\dGyrlyga 
&&&\ruTTo>{[\circledast_\mcM^{(g,c)}(t_\psi)
	(\Delta(t^{|c}_{[\psi(g-1),\psi(g)]}))]\langle \sfi,\sfj\rangle}
\\
\rmF^\su_{t(0)}\langle \sfi,\sfj\rangle &&\HonV 
&\rTTo^{\hspace*{-3.3em}\Delta(t_\psi)\langle \sfi,\sfj\rangle}
&[\circledast_\mcM^{(g,c)}(t_\psi)\rmF_{t_{\psi,g}^{-1}c}]
	\langle \sfi,\sfj\rangle
\hspace*{-2em} &&\dMono
\\
&\rdTTo(2,5)<{\Delta(t)} \rdTTo(4,3)[hug]^{\qquad\Delta(t_\psi)} 
&\rule[1em]{0.4pt}{1.5em} &&\dMono
&&\hspace*{-1em} \circledast_\mcM^{(g,c)}(t_\psi)\bigl[\bigl(
\circledast_\mcM^{(h,b)}(t^{|c}_{[\psi(g-1),\psi(g)]})\rmF^\su_{t_h^{-1}b}
	\bigr)\langle \sfi,\sfj\rangle\bigr] \hspace*{6em}
\\
&&&&&\ruTTo>{\circledast_\mcM^{(g,c)}(t_\psi)
[(\Delta(t^{|c}_{[\psi(g-1),\psi(g)]}))\langle\sfi,\sfj\rangle]}
\\
&&\dTTo &&\circledast_\mcM^{(g,c)}(t_\psi)
[\rmF_{t_{\psi,g}^{-1}c}\langle\sfi,\sfj\rangle]
\hspace*{-2em} &&\dMono
\\
&&&&&\rdTTo>{\circledast_\mcM^{(g,c)}(t_\psi)
	(\Delta(t^{|c}_{[\psi(g-1),\psi(g)]}))}
\\
&&\hspace*{-2.4em}
\circledast_\mcM^{(h,b)}(t)[\rmF^\su_{t_h^{-1}b}\langle\sfi,\sfj\rangle] 
&&\rTTo^{\lambda^f} 
&&\circledast_\mcM^{(g,c)}(t_\psi)\bigl(
	\circledast_\mcM^{(h,b)}(t^{|c}_{[\psi(g-1),\psi(g)]})
\bigl[\rmF^\su_{t_h^{-1}b}\langle \sfi,\sfj\rangle \bigr]\bigr) 
\hspace*{7em}
\end{diagram}
\hspace*{-0.1em}
\]

Thus, a polymodule cooperad
\((\mainf^\su,\rmF^\su)\langle\sfi,\sfj\rangle\) is constructed.
By \lemref{lem-degree-1-isomorphism-cooperads} there is a polymodule 
cooperad \((A_\infty^\su,F^\su)\langle\bi,\bj\rangle\) isomorphic to it 
via a degree 1 isomorphism.

\begin{proposition}
The collection of operad submodules 
	\((A_\infty^{hu},F^{hu})\subset
	(A_\infty^\su,F^\su)\langle\bi,\bj\rangle\)
is a subcooperad.
\end{proposition}

\begin{proof}
Assume that \(k\in t(0)\) for a $l$\n-tree $t$.
Let us compute \(\varDelta^\mcM(t)(f_{e_k})\).
Notice that there exists the only $t$\n-tree $\tau$ such that 
$\tilde\tau$ is surjective and \(|\tau(0,a)|=\delta_k^a\) for all
\(a\in t(0)\).
In fact, \(\tilde\tau(0)=\mb1\), hence, \(\tilde\tau(h)=\mb1\) for all 
\(h\in[I]\).
The tree $\tau$ is given by the formula
\[ \tau(h,b)=
\begin{cases}
\mb1, \quad &\text{if } b=t_h\ldots t_2t_1(k),
\\
\emptyset, &\text{otherwise}.
\end{cases}
\]
Denoting \(e_p^S\in\NN^S\) a basis vector for \(p\in S\) we find
\[ \varDelta^\mcM(t)(f_{e_k^{t(0)}})
=f_{e_k^{t_1^{-1}t_1k}}\tens f_{e_{t_1k}^{t_2^{-1}t_2t_1k}}\tens 
f_{e_{t_2t_1k}^{t_3^{-1}t_3t_2t_1k}} \tdt 
f_{e_{t_{l-2}\ldots t_1k}^{t_{l-1}^{-1}t_{l-1}t_{l-2}\ldots t_1k}} \tens 
f_{e_{t_{l-1}\ldots t_1k}^{t(l-1)}}.
\]
Now we compute
\begin{align*}
\varDelta^\mcM(t)(\bv_k^{t(0)})
=\varDelta^\mcM(t)(\lambda^k_{e_k}(\bj;f_{e_k^{t(0)}}) 
-\bj\rho_\emptyset)
=\lambda^k_{e_k}(\bj;\varDelta^\mcM(t)(f_{e_k^{t(0)}}))
-\bj\rho_\emptyset &
\\
=\bj f_{e_k^{t_1^{-1}t_1k}} \tens f_{e_{t_1k}^{t_2^{-1}t_2t_1k}} \tens 
f_{e_{t_2t_1k}^{t_3^{-1}t_3t_2t_1k}} \tdt 
f_{e_{t_{l-2}\ldots t_1k}^{t_{l-1}^{-1}t_{l-1}t_{l-2}\ldots t_1k}} \tens 
f_{e_{t_{l-1}\ldots t_1k}^{t(l-1)}} &
\\
-\bj\rho_\emptyset \tens f_{e_{t_1k}^{t_2^{-1}t_2t_1k}} \tens 
f_{e_{t_2t_1k}^{t_3^{-1}t_3t_2t_1k}} \tdt 
f_{e_{t_{l-2}\ldots t_1k}^{t_{l-1}^{-1}t_{l-1}t_{l-2}\ldots t_1k}} \tens 
f_{e_{t_{l-1}\ldots t_1k}^{t(l-1)}} &
\\
+\bj f_{e_{t_1k}^{t_2^{-1}t_2t_1k}}\tens 
f_{e_{t_2t_1k}^{t_3^{-1}t_3t_2t_1k}} \tdt
f_{e_{t_{l-2}\ldots t_1k}^{t_{l-1}^{-1}t_{l-1}t_{l-2}\ldots t_1k}} \tens 
f_{e_{t_{l-1}\ldots t_1k}^{t(l-1)}} &
\\
-\bj\rho_\emptyset \tens f_{e_{t_2t_1k}^{t_3^{-1}t_3t_2t_1k}} \tdt 
f_{e_{t_{l-2}\ldots t_1k}^{t_{l-1}^{-1}t_{l-1}t_{l-2}\ldots t_1k}} \tens 
f_{e_{t_{l-1}\ldots t_1k}^{t(l-1)}} &
\\
+\bj f_{e_{t_2t_1k}^{t_3^{-1}t_3t_2t_1k}} \tdt 
f_{e_{t_{l-2}\ldots t_1k}^{t_{l-1}^{-1}t_{l-1}t_{l-2}\ldots t_1k}} \tens 
f_{e_{t_{l-1}\ldots t_1k}^{t(l-1)}} &
\\
\cdots \hspace*{8em} &
\\
-\bj\rho_\emptyset \tens f_{e_{t_{l-1}\ldots t_1k}^{t(l-1)}} &
\\
+\bj f_{e_{t_{l-1}\ldots t_1k}^{t(l-1)}} &
\\
-\bj\rho_\emptyset &
\\
=\bv_k^{t_1^{-1}t_1k} \tens f_{e_{t_1k}^{t_2^{-1}t_2t_1k}} \tens 
f_{e_{t_2t_1k}^{t_3^{-1}t_3t_2t_1k}} \tdt 
f_{e_{t_{l-2}\ldots t_1k}^{t_{l-1}^{-1}t_{l-1}t_{l-2}\ldots t_1k}} \tens 
f_{e_{t_{l-1}\ldots t_1k}^{t(l-1)}} &
\\
+\bv_{t_1k}^{t_2^{-1}t_2t_1k} \tens f_{e_{t_2t_1k}^{t_3^{-1}t_3t_2t_1k}} 
\tdt f_{e_{t_{l-2}\ldots t_1k}^{t_{l-1}^{-1}t_{l-1}t_{l-2}\ldots t_1k}} 
\tens f_{e_{t_{l-1}\ldots t_1k}^{t(l-1)}} &
\\
+\bv_{t_2t_1k}^{t_3^{-1}t_3t_2t_1k} \tdt 
f_{e_{t_{l-2}\ldots t_1k}^{t_{l-1}^{-1}t_{l-1}t_{l-2}\ldots t_1k}} \tens 
f_{e_{t_{l-1}\ldots t_1k}^{t(l-1)}} &
\\
\cdots \hspace*{8em} &
\\
+\bv_{t_{l-2}\ldots t_1k}^{t_{l-1}^{-1}t_{l-1}t_{l-2}\ldots t_1k} \tens 
f_{e_{t_{l-1}\ldots t_1k}^{t(l-1)}} &
\\
+\bv_{t_{l-1}\ldots t_1k}^{t(l-1)} &.
\end{align*}

On other generators we transform
\[ 
\varDelta^\mcM(t)(f_{(\ell^k_1;\ell^k_2;\dots;\ell^k_{t^k})_{k\in\bn}})
=\lambda_{\hat\ell}\bigl((\sS{^{\ell^k_1}}1,\bj,\sS{^{\ell^k_2}}1,\bj,
\dots, \sS{^{\ell^k_{t^k-1}}}1,\bj,\sS{^{\ell^k_{t^k}}}1)_{k\in\bn};
\varDelta^\mcM(t)(f^+_{\hat\ell})\bigr)
\]
as follows.
First factors $f_{i_q}$ are replaced with generators
\(f_{a_1;a_2;\dots;a_p}\) accordingly with the set of $\bj$'s appearing 
among the arguments of $f_{i_q}$. 
The only exception is the case of $\bj f_{e_k}$ which is replaced with 
$\bv_k+\bj\rho_\emptyset$. 
In obtained summands all instances of $\bj\rho_\emptyset$ are moved to 
the right as arguments $\bj$ of next $f_p$ due to defining the tensor 
product as a colimit, and this procedure goes on until no 
$\bj\rho_\emptyset$ are left.
Notice that the separate term $\bj\rho_\emptyset$ can not appear 
elsewhere but in the expression
\(\varDelta^\mcM(t)(\bj f_{e_k^{t(0)}})\), which is not considered by 
itself, but only as a summand of \(\bv_k^{t(0)}\).
\end{proof}

\begin{corollary}
The collection of operad submodules
	\((\mainf^{hu},\rmF^{hu})\subset
	(\mainf^\su,\rmF^\su)\langle\sfi,\sfj\rangle\)
is a polymodule subcooperad.
\end{corollary}

\appendix
\section{Colimits of algebras over monads}
Let \(\top:\cc\to\cc\) be a monad, and let \(F:\cc\leftrightarrows\cc^\top:U\) be the associated adjunction.
Assume that $\cc$ is cocomplete and $\cc^\top$ has coequalizers.
The latter condition is satisfied in each of the two following cases:
\begin{itemize}
\item $\cc$ is a complete, regular, regularly co-well-powered category with coequalizers, and $\top$ is a monad which preserves regular epimorphisms \cite[Proposition~9.3.8]{BarrWells:TopTT}.

\item $\cc$ has finite colimits and equalizers of arbitrary sets of maps (with the same source and target), and $\top$ is a monad in $\cc$ which preserves colimits along countable chains  \cite[Theorem~9.3.9]{BarrWells:TopTT}.
\end{itemize}

When $\cc^\top$ has coequalizers, the category $\cc^\top$ is cocomplete by a result of Barr and Wells \cite[Corollary~9.3.3]{BarrWells:TopTT}.
Our goal in this section is to reprove this result expressing the colimit in $\cc^\top$ through the colimit in $\cc$ via sufficiently explicit recipe.

\begin{proposition}
Assume that $\cc$ is cocomplete and $\cc^\top$ has coequalizers.
Then the category $\cc^\top$ is cocomplete.
\end{proposition}

\begin{proof}
Let $I$ be a small category and let \(I\ni i\mapsto P_i\in\cc^\top\) be a diagram in $\cc^\top$.
Denote by $C$ the colimit (coequalizer) of the following diagram in $\cc^\top$
\begin{diagram}[LaTeXeqno,w=4em]
F\colim_i UFUP_i
\\
\dTTo<{F\colim\alpha_i} &\rdTTo>{F\can}
\\
F\colim_i UP_i &&FUF\colim_i UP_i &\rTTo^\eps &F\colim_i UP_i
\\
\dTTo<{F\colim\eta} &\ruTTo>{F\can}
\\
F\colim_i UFUP_i
\label{dia-FcolimiUFUPi}
\end{diagram}
where `can' means any canonical map.
Equip \(C=(C,\can:F\colim_iUP_i\to C)\) with maps in $\cc$ going through the rightmost vertex
\begin{align}
\In_i &=\bigl( UP_i \rTTo^{\inj_i} \colim_i UP_i \rto\eta UF\colim_i UP_i \rTTo^{U\can} UC \bigr) \notag
\\
&=\bigl( UP_i \rto\eta UFUP_i \rTTo^{UF\inj_i} UF\colim_i UP_i \rTTo^{U\can} UC \bigr).
\label{eq-Ini-UPi-UC}
\end{align}
We claim that \(\In_i\in\cc^\top\) and \((C,\In_i:P_i\to C\mid i\in I)\) is the colimiting cocone of the diagram \(i\mapsto P_i\) in $\cc^\top$.

Let us verify that $\In_i$ are morphisms of $\top$\n-algebras.
The exterior of the following diagram commutes
\begin{diagram}[LaTeXeqno]
UFUP_i &\rTTo^{UF\eta} &UFUFUP_i &\rTTo^{UFUF\inj_i} &UFUF\colim_i UP_i &\rTTo^{UFU\can} &UFUC
\\
\dTTo<{\alpha_i} &\rdTTo(4,2)>{UF\inj_i} &&= &\dTTo>{U\eps} &= &\dTTo>{\alpha_C}
\\
UP_i &\rTTo^\eta &UFUP_i &\rTTo^{UF\inj_i} &UF\colim_i UP_i &\rTTo^{U\can} &UC
\label{dia-UFUPi-UC}
\end{diagram}
if and only if
\begin{diagram}[LaTeXeqno]
UFUP_i &\rTTo^{UF\inj_i} &UF\colim_i UP_i &\rTTo^{U\can} &UC
\\
\dTTo<{\alpha_i} &&= &&\uTTo>{U\can}
\\
UP_i &\rTTo^\eta &UFUP_i &\rTTo^{UF\inj_i} &UF\colim_i UP_i
\label{dia-UFUPi-UFcolimiUPi}
\end{diagram}
Schematically this is the equation \(f=\bigr(A\rto gA\rto fC\bigr)\), where \(f=UF\inj_i\cdot U\can:A\to C\in\cc^\top\) but \(g=\alpha_i\cdot\eta\in\cc\).
By the freeness of $\top$\n-algebra $\top A$ (see the proof of \cite[Theorem~3.2.1]{BarrWells:TopTT}) this is equivalent to equation
\[ \bigr(\top A \rto{\alpha_A} A \rto f C\bigr) =\bigr(\top A \rto{\top g} \top A \rto{\alpha_A} A \rto f C\bigr).
\]
In detail it is the equation
\begin{diagram}[LaTeXeqno]
UFUFUP_i &\rTTo^{U\eps} &UFUP_i &\rTTo^{UF\inj_i} &UF\colim_i UP_i &\rTTo^{U\can} &UC
\\
\dTTo<{UF\alpha_i} &&&= &&&\uTTo>{U\can}
\\
UFUP_i &\rTTo^{UF\eta} &UFUFUP_i &\rTTo^{U\eps} &UFUP_i &\rTTo^{UF\inj_i} &UF\colim_i UP_i
\label{dia-UFUFUPi-UFcolimiUPi}
\end{diagram}
Removing the unnecessary $U$ we write it as an equation in $\cc^\top$:
\begin{diagram}[LaTeXeqno]
FUFUP_i &\rTTo^{FUF\inj_i} &FUF\colim_i UP_i &\rTTo^\eps &F\colim_i UP_i &\rTTo^\can &C
\\
\dTTo<{F\alpha_i} &&&= &&&\uTTo>\can
\\
FUP_i &\rTTo^{F\eta} &FUFUP_i &\rTTo^{FUF\inj_i} &FUF\colim_i UP_i &\rTTo^\eps &F\colim_i UP_i
\label{dia-FUFUPi-FcolimiUPi}
\end{diagram}
which holds due to \((C,\can)\) being coequalizer of \eqref{dia-FcolimiUFUPi}.

Clearly, \(\In_i:P_i\to C\) is a cocone from the diagram \(i\mapsto P_i\).
Let us prove that it is an initial one.
Let \(\phi_i:P_i\to Q\in\cc^\top\) be an arbitrary cocone from the diagram \(i\mapsto P_i\).
There is a unique map \(\beta:\colim_iUP_i\to UQ\in\cc\) such that \(U\phi_i=\bigl(UP_i\rTTo^{\inj_i} \colim_iUP_i\rto\beta UQ\bigr)\).
It has an adjunct \(\gamma=\sS{^t}\beta=\bigl(F\colim_iUP_i\rto{F\beta} FUQ\rto\eps Q\bigr)\in\cc^\top\), so that \(\sS{^t}(U\phi_i)=\bigl(FUP_i\rTTo^{F\inj_i} F\colim_iUP_i\rto\gamma Q\bigr)\).
Consequently,
\[ U\phi_i =\bigl(UP_i \rto\eta UFUP_i \rTTo^{UF\inj_i} UF\colim_iUP_i \rto{U\gamma} UQ\bigr).
\]
Since \(\phi_i\in\cc^\top\) the exterior of diagram~\eqref{dia-UFUPi-UC} commutes, where $\can$ and $C$ are replaced with $\gamma$ and $Q$.
Therefore, equation~\eqref{dia-UFUPi-UFcolimiUPi} with the same replacement holds.
As explained above this implies equations \eqref{dia-UFUFUPi-UFcolimiUPi} and \eqref{dia-FUFUPi-FcolimiUPi} with the same modification.
Therefore, both paths in diagram~\eqref{dia-FcolimiUFUPi} postcomposed with \(\gamma:F\colim_iUP_i\to Q\) from the top vertex \(F\colim_iUFUP_i\) to $Q$ are equal to each other.
Hence, $\gamma$ factorizes as \(F\colim_iUP_i\rTTo^\can C\rto\psi Q\) for a unique \(\psi\in\cc^\top\).
\end{proof}

\begin{remark}
It is shown in the proof of the above proposition that $\colim_iP_i$ is the biggest quotient of $F\colim_iUP_i$ via a regular epimorphism \(\can:F\colim_iUP_i\to C=\colim_iP_i\) such that morphisms \(\In_i:UP_i\to UC\in\cc\) from \eqref{eq-Ini-UPi-UC} are morphisms of $\top$\n-algebras.
\end{remark}

\begin{proposition}
Assume that $\cc$ is cocomplete and $\cc^\top$ has coequalizers.
Let \(X\in\Ob\cc\) and \(A=(UA,\alpha:UFUA\to UA)\in\Ob\cc^\top\).
Then the colimit \(C=(C,\can:F(X\sqcup UA)\to C)\) of the diagram in $\cc^\top$
\begin{diagram}[LaTeXeqno,w=4em]
FUFUA
\\
\dTTo<{F\alpha} &\rdTTo>{FUF\inj_2}
\\
FUA &&FUF(X\sqcup UA) &\rTTo^\eps &F(X\sqcup UA)
\\
\dTTo<{F\eta} &\ruTTo>{FUF\inj_2}
\\
FUFUA
\label{dia-FUFUA-F(XuUA)}
\end{diagram}
equipped with the morphisms of $\top$\n-algebras
\begin{align}
\In_1 &=\bigl( FX \rTTo^{F\inj_1} F(X\sqcup UA) \rTTo^\can C \bigr), \notag
\\
\In_2 &=\bigl( UA \rTTo^{\inj_2} X\sqcup UA \rto\eta UF(X\sqcup UA) \rTTo^{U\can} UC \bigr) \notag
\\
&=\bigl( UA \rto\eta UFUA \rTTo^{UF\inj_2} UF(X\sqcup UA) \rTTo^{U\can} UC \bigr)
\label{eq-In1-In2-In2}
\end{align}
is the coproduct $FX\sqcup A$ in $\cc^\top$.
\end{proposition}

\begin{proof}
Let us verify that $\In_2$ is a morphism of $\top$\n-algebras.
This is equivalent to commutativity of the exterior of the following diagram
\begin{diagram}[LaTeXeqno]
UFUA &\rTTo^{UF\eta} &UFUFUA &\rTTo^{UFUF\inj_2} &UFUF(X\sqcup UA) &\rTTo^{UFU\can} &UFUC
\\
\dTTo<\alpha &\rdTTo(4,2)>{UF\inj_2} &&= &\dTTo>{U\eps} &= &\dTTo>{\alpha_C}
\\
UA &\rTTo^\eta &UFUA &\rTTo^{UF\inj_2} &UF(X\sqcup UA) &\rTTo^{U\can} &UC
\label{dia-UFUA-UC}
\end{diagram}
which holds if and only if
\begin{diagram}[LaTeXeqno]
UFUA &\rTTo^{UF\inj_2} &UF(X\sqcup UA) &\rTTo^{U\can} &UC
\\
\dTTo<\alpha &&= &&\uTTo>{U\can}
\\
UA &\rTTo^\eta &UFUA &\rTTo^{UF\inj_2} &UF(X\sqcup UA)
\label{dia-UFUA-UF(XuUA)}
\end{diagram}
Schematically this is the equation \(f=\bigr(B\rto gB\rto fC\bigr)\), where \(f=UF\inj_2\cdot U\can:B\to C\in\cc^\top\) but \(g=\alpha\cdot\eta\in\cc\).
By the freeness of $\top$\n-algebra $\top B$ (see the proof of \cite[Theorem~3.2.1]{BarrWells:TopTT}) this is equivalent to equation
\[ \bigr(\top B \rto{\alpha_B} B \rto f C\bigr) =\bigr(\top B \rto{\top g} \top B \rto{\alpha_B} B \rto f C\bigr).
\]
In detail it is the equation
\begin{diagram}[LaTeXeqno]
UFUFUA &\rTTo^{U\eps} &UFUA &\rTTo^{UF\inj_2} &UF(X\sqcup UA) &\rTTo^{U\can} &UC
\\
\dTTo<{UF\alpha} &&&= &&&\uTTo>{U\can}
\\
UFUA &\rTTo^{UF\eta} &UFUFUA &\rTTo^{U\eps} &UFUA &\rTTo^{UF\inj_2} &UF(X\sqcup UA)
\label{dia-UFUFUA-UF(XuUA)}
\end{diagram}
Removing the unnecessary $U$ we write it as an equation in $\cc^\top$:
\begin{diagram}[LaTeXeqno]
FUFUA &\rTTo^{FUF\inj_2} &FUF(X\sqcup UA) &\rTTo^\eps &F(X\sqcup UA) &\rTTo^\can &C
\\
\dTTo<{F\alpha} &&&= &&&\uTTo>\can
\\
FUA &\rTTo^{F\eta} &FUFUA &\rTTo^{FUF\inj_2} &FUF(X\sqcup UA) &\rTTo^\eps &F(X\sqcup UA)
\label{dia-FUFUA-F(XuUA)-C}
\end{diagram}
which holds due to \((C,\can)\) being coequalizer of \eqref{dia-FUFUA-F(XuUA)}.

Let us prove that \((C,\In_1:FX\to C,\In_2:A\to C)\) is the coproduct $FX\sqcup A$ in $\cc^\top$.
Let \(\phi_1:FX\to Q\in\cc^\top\) and \(\phi_2:A\to Q\in\cc^\top\).
The maps \(\delta=\phi_1^t=\bigl(X\rto\eta UFX\rTTo^{U\phi_1} UQ)\) and \(U\phi_2:UA\to UQ\) determine a unique map \(\beta:X\sqcup UA\to UQ\) in $\cc$.
It has an adjunct \(\gamma=\sS{^t}\beta=\bigl(F(X\sqcup UA)\rto{F\beta} FUQ\rto\eps Q\bigr)\in\cc^\top\), so that \(\phi_1=\sS{^t}\delta=\bigl(FX\rTTo^{F\inj_1} F(X\sqcup UA)\rto\gamma Q\bigr)\), \(\sS{^t}(U\phi_2)=\bigl(FUA\rTTo^{F\inj_2} F(X\sqcup UA)\rto\gamma Q\bigr)\).
Consequently,
\[ U\phi_2 =\bigl(UA \rto\eta UFUA \rTTo^{UF\inj_2} UF(X\sqcup UA) \rto{U\gamma} UQ\bigr).
\]
Since \(\phi_2\in\cc^\top\) the exterior of diagram~\eqref{dia-UFUA-UC} commutes, where $\can$ and $C$ are replaced with $\gamma$ and $Q$.
Therefore, equation~\eqref{dia-UFUA-UF(XuUA)} with the same replacement holds.
As explained above this implies equations \eqref{dia-UFUFUA-UF(XuUA)} and \eqref{dia-FUFUA-F(XuUA)-C} with the same modification.
Therefore, both paths in diagram~\eqref{dia-FUFUA-F(XuUA)} postcomposed with \(\gamma:F(X\sqcup UA)\to Q\) from the top vertex \(FUFUA\) to $Q$ are equal to each other.
Hence, $\gamma$ factorizes as \(F(X\sqcup UA)\rTTo^\can C\rto\psi Q\) for a unique \(\psi\in\cc^\top\).
We get
\begin{align*}
\phi_1 &=\bigl( FX \rTTo^{F\inj_1} F(X\sqcup UA) \rTTo^\can C \rto\psi Q \bigr) =\bigl( FX \rTTo^{\In_1} C \rto\psi Q \bigr),
\\
U\phi_2 &=\bigl( UA \rto\eta UFUA \rTTo^{UF\inj_2} UF(X\sqcup UA) \rTTo^{U\can} UC \rto{U\psi} UQ \bigr),
\end{align*}
hence, \(\phi_2=\bigl(A\rTTo^{\In_2} C \rto\psi Q \bigr)\).
This shows that \((C,\In_1,\In_2)\) is the coproduct $FX\sqcup A$ in $\cc^\top$.
\end{proof}

\begin{corollary}\label{cor-biggest-quotient-morphism-T-algebras}
$FX\sqcup A$ is the biggest quotient of $F(X\sqcup UA)$ via a regular epimorphism \(\can:F(X\sqcup UA)\to C=FX\sqcup A\) such that the morphism \(\In_2:UA\to UC\in\cc\) from \eqref{eq-In1-In2-In2} is a morphism of $\top$\n-algebras.
\end{corollary}

%\bibliography{yuri}

\begin{thebibliography}{FOOO09}

\bibitem[BLM08]{BesLyuMan-book}
Yu.~Bespalov, V.~Lyubashenko, and O.~Manzyuk, \emph{Pretriangulated
  ${A}_\infty$-categories}, Proc. of the Inst. of Math. of NAS of 
  Ukraine. Mathematics and its Applications, vol.~76, 
  Institute of Mathematics of NAS of Ukraine, Kyiv, 2008,
  \url{http://www.math.ksu.edu/~lub/papers.html}.

\bibitem[BW05]{BarrWells:TopTT}
M.~Barr and C.~Wells, \emph{Toposes, triples and theories}, Theory
  Appl. Categ. \textbf{278} (1985, 2005), no.~12, x+288 pp. 
  (electronic), Corrected reprint of the 1985 original, 
  Grundlehren der mathematischen {W}issenschaften,
  \url{http://www.tac.mta.ca/tac/reprints/articles/12/tr12abs.html}.

\bibitem[FOOO09]{FukayaOhOhtaOno:Anomaly}
K.~Fukaya, Yo.-G.~Oh, H.~Ohta, and K.~Ono, \emph{Lagrangian
  intersection {F}loer theory: Anomaly and obstruction}, 
  AMS/IP Studies in Advanced Mathematics Series, 
  vol.~46, Amer. Math. Soc., 2009.

\bibitem[Hin97]{Hinich:q-alg/9702015}
V.~Hinich, \emph{Homological algebra of homotopy algebras}, 
  Comm. Algebra \textbf{25} (1997), no.~10, 3291--3323,
  \href{http://arXiv.org/abs/q-alg/9702015}{{\tt
  arXiv:\linebreak[1]q-alg/9702015}}.

\bibitem[Kel06]{math.RT/0510508}
B.~Keller, \emph{{$A$}-infinity algebras, modules and functor
  categories}, Trends in representation theory of algebras and related 
  topics, Contemp. Math., vol. 406, Amer. Math. Soc., 
  Providence, RI, 2006, pp.~67--93,
  \href{http://arXiv.org/abs/math.RT/0510508}{{\tt 
  arXiv:\linebreak[1]math.RT/0510508}}.

\bibitem[Lei03]{math.CT/0305049}
T.~Leinster, \emph{Higher operads, higher categories}, London Math.
  Soc. Lecture Notes Series, Cambridge University Press, Boston, Basel,
  Berlin, 2003, \href{http://arXiv.org/abs/math.CT/0305049}{{\tt
  arXiv:\linebreak[1]math.CT/0305049}}.

\bibitem[LH03]{Lefevre-Ainfty-these}
K.~Lef\`evre-Hasegawa, \emph{Sur les ${A}_\infty$-cat\'egories}, Ph.D.
  thesis, Universit\'e Paris 7, U.F.R. de Math\'ematiques, 2003,
  \href{http://arXiv.org/abs/math.CT/0310337}{{\tt
  arXiv:\linebreak[1]math.CT/0310337}}.

\bibitem[Lyu11]{Lyu-Ainf-Operad}
V.~Lyubashenko, \emph{Homotopy unital ${A}_\infty$-algebras},
  J. Algebra \textbf{329} (2011), no.~1, 190--212, Special Issue 
  Celebrating the 60th Birthday of Corrado De Concini,
  \url{http://dx.doi.org/10.1016/j.jalgebra.2010.02.009}.

\bibitem[Lyu12]{Lyu-Hinich-Thm}
V.~Lyubashenko, \emph{A model structure on categories related to 
  categories of complexes}, 2012.

\bibitem[M{ac}63]{MacLane-Homology}
S.~M{ac~Lane}, \emph{Homology}, Die Grundlehren der mathematischen
  Wissenschaften, no. 114, Springer-Verlag, Berlin, Heidelberg, 1963.

\bibitem[M{ac}88]{MacLane}
S.~M{ac~Lane}, \emph{Categories for the working mathematician},
  GTM, vol.~5, Springer-Verlag, New York, 1971, 1988.

\bibitem[Mak96]{MakkaiM:avoacc}
M.~Makkai, \emph{Avoiding the axiom of choice in category theory},
  J. Pure Appl. Algebra \textbf{108} (1996), 109--173.

\bibitem[Mar96]{Markl:ModOp}
M.~Markl, \emph{Models for operads}, Commun. in Algebra
  \textbf{24} (1996), no.~4, 1471--1500, 
  \href{http://arXiv.org/abs/hep-th/9411208}{{\tt
  arXiv:\linebreak[1]hep-th/9411208}}.

\bibitem[Mur11]{MR2821434}
F.~Muro, \emph{Homotopy theory of nonsymmetric operads}, Algebr. Geom.
  Topol. \textbf{11} (2011), no.~3, 1541--1599,
  \href{http://arXiv.org/abs/1101.1634}{{\tt arXiv:\linebreak[1]1101.1634}}
  \url{http://dx.doi.org/10.2140/agt.2011.11.1541}.

\bibitem[Sei08]{SeidelP-book-Fukaya}
P.~Seidel, \emph{Fukaya categories and {P}icard--{L}efschetz theory},
  Zurich Lectures in Advanced Math., European Math. Soc. (EMS), 
  Z\"urich, 2008.

\bibitem[Spi01]{math/0101102}
M.~Spitzweck, \emph{Operads, algebras and modules in general model
  categories}, jan 2001, \href{http://arXiv.org/abs/math/0101102}{{\tt
  arXiv:\linebreak[1]math/0101102}}.

\bibitem[Sta63]{Stasheff:HomAssoc}
J.~D. Stasheff, \emph{Homotopy associativity of {H}-spaces I $\&$ {II}},
  Trans. Amer. Math. Soc. \textbf{108} (1963), 275--292, 293--312.

\bibitem[Wei94]{Weibel}
C.~A. Weibel, \emph{An introduction to homological algebra}, Cambridge
  Studies in Adv. Math., vol.~38, Cambridge Univ. Press, Cambridge, 
  New York, Melbourne, 1994.

\end{thebibliography}
\tableofcontents

\ifx\chooseClass1
	\else
\textit{\\
Institute of Mathematics,
National Academy of Sciences of Ukraine, \\
3 Tereshchenkivska st.,
Kyiv-4, 01601 MSP, Ukraine \\
}

Email: \textsf{lub@imath.kiev.ua}
\fi

\end{document}